\documentclass{article}
\usepackage{amssymb,latexsym}
\def\restypeone{{\rm (RT)}\ \ }
\def\restypetwo{{\rm (RT)}\ \ }
\def\remRT{\refer{r: technical reasons residue type}}
\def\thmendRT{\refer{t: inversion for Schwartz}}
\def\eqnreset{\setcounter{equation}{0}}

\newtheorem{Thm}{Theorem}[section]
\newtheorem{Defi}[Thm]{Definition}
\newtheorem{Cor}[Thm]{Corollary}
\newtheorem{Lemma}[Thm]{Lemma}
\newtheorem{Prop}[Thm]{Proposition}
\newtheorem{Rem}[Thm]{Remark}
\newtheorem{Conj}[Thm]{Conjecture}
\newtheorem{Prelim}[Thm]{Preliminary}

\newenvironment{thm}[0]{\begin{Thm}\noindent}%
{\end{Thm}}
\newenvironment{defi}[0]{\begin{Defi}\noindent\rm}%
{\end{Defi}}
\newenvironment{cor}[0]{\begin{Cor}\noindent}%
{\end{Cor}}
\newenvironment{lemma}[0]{\begin{Lemma}\noindent}%
{\end{Lemma}}
\newenvironment{prop}[0]{\begin{Prop}\noindent}%
{\end{Prop}}
\newenvironment{rem}[0]{\begin{Rem}\noindent\rm}%
{\end{Rem}}
{\end{Conj}}
{\end{Prelim}}

\def\medno{\medbreak\noindent}

\def\qed{~\hfill$\square$\medbreak}
\def\proof{\par\noindent{\bf Proof:}{\ }{\ }}

\def\naam{\label}
\def\refer{\ref}

\def\bib#1{\cite{#1}}

\def\text#1{\;\;\;\;{\rm \hbox{#1}}\;\;\;\;}
\def\qquad{\quad\quad}

\def\itema{\item[{\rm (a)}]}
\def\itemb{\item[{\rm (b)}]}
\def\itemc{\item[{\rm (c)}]}
\def\itemd{\item[{\rm (d)}]}
\def\iteme{\item[{\rm (e)}]}

\def\msy#1{{\mathbb #1}}
\def\C{{\msy C}}
\def\N{{\msy N}}

\def\R{{\msy R}}
\def\D{{\msy D}}

\def\ga{\alpha}
\def\gb{\beta}
\def\gd{\delta}
\def\geps{\varepsilon}
\def\gf{\varphi}
\def\gg{\gamma}

\def\gl{\lambda}

\def\gs{\sigma}

\def\gD{\Delta}

\def\gS{\Sigma}

\def\got#1{\mathfrak #1}
\def\fa{{\got a}}
\def\fb{{\got b}}
\def\fc{{\got c}}
\def\fg{{\got g}}
\def\fh{{\got h}}
\def\fj{{\got j}}
\def\fk{{\got k}}
\def\fl{{\got l}}
\def\fm{{\got m}}
\def\fn{{\got n}}

\def\fp{{\got p}}
\def\fq{{\got q}}

\def\ft{{\got t}}

\def\implies{\Rightarrow}
\def\to{\rightarrow}
\def\Re{{\rm Re}\,}

\def\inp#1#2{\langle#1\,,\,#2\rangle}
\def\hinp#1#2{\langle#1\,|\,#2\rangle}
\def\Ad{{\rm Ad}}
\def\End{{\rm End}}
\def\Hom{{\rm Hom}}

\def\ad{{\rm ad}}
\def\after{\,{\scriptstyle\circ}\,}

\def\pr{{\rm pr}}

\def\tr{{\rm tr}\,}

\def\iM{0}
\def\iiM{0}

\def\ik{{\rm k}}
\def\ih{{\rm h}}
\def\iq{{\rm q}}
\def\ip{{\rm p}}
\def\iC{{\scriptscriptstyle \C}}
\def\iR{{\scriptscriptstyle \R}}

\def\cA{{\cal A}}
\def\cB{{\cal B}}
\def\cC{{\cal C}}
\def\cD{{\cal D}}
\def\cE{{\cal E}}
\def\cF{{\cal F}}
\def\cH{{\cal H}}

\def\cL{{\cal L}}
\def\cM{{\cal M}}

\def\cO{{\cal O}}
\def\cP{{\cal P}}
\def\cR{{\cal R}}
\def\cS{{\cal S}}
\def\cT{{\cal T}}

\def\cW{{\cal W}}

\mathcode`:="603A
\def\col{\,:\,}

\def\um{{\rm um}}
\def\dual{{\rm d}}
\def\temp{{\rm temp}}
\def\rmp{{\rm p}}
\def\ep{{\rm ep}}
\def\ev{{\rm ev}}
\def\hyp{{\rm hyp}}

\def\mer{{\rm mer}}
\def\rmT{{\rm T}}

\def\glob{{\rm glob}}
\def\hglob{{\rm hglob}}

\def\bp{{}^\backprime}
\def\bs{\backslash}

\def\asmid{\mid}

\def\bfP{{\bf P}}

\def\bfe{{\bf e}}

\def\cC{{\cal C}}
\def\Hyp{{\cal H}}

\def\specinp{b}
\def\stoneR{{}^{*1\!}R}

\def\start{{{}^*t}}

\def\DP{\gD(P)}
\def\DrP{\gD_r(P)}
\def\DrQ{\Delta_r(Q)}
\def\gDr{\gD_r}
\def\gSp{\Sigma^+}
\def\gSrF{\gS_r(F)}

\def\gSr{\Sigma_r}

\def\umuP{\umu_P}
\def\umu{\underline \mu}

\def\spX{{\rm X}}
\def\spXp{\spX_+}
\def\spXoneP{\spX_{1P}}
\def\spXFv{\spX_{F,v}}
\def\spXFvp{\spX_{F,v,+}}
\def\spXQv{\spX_{Q,v}}
\def\spXzerov{\spX_{0,v}}

\def\spXP{\spX_P}
\def\spXoneQv{\spX_{1Q, v}}
\def\spXoneQvp{\spX_{1Q,v,+}}

\def\spXQu{\spX_{Q,u}}
\def\spXQup{\spX_{Q,u,+}}
\def\spXPvp{\spX_{P,v,+}}
\def\spXPv{\spX_{P,v}}
\def\spXFup{\spX_{F,u,+}}

\def\spXQe{\spX_{Q,e}}
\def\spXQvp{\spX_{Q,v,+}}

\def\spXonePv{\spX_{1P, v}}
\def\spXoneQe{\spX_{1Q,e}}

\def\AsubF{\cA_F^{\start}}
\def\cAFv{\cA^\start(\spXFv\col \tau_F)}

\def\cAtwoQ{\cA_{2,Q}}

\def\cAFtwo{\cA_{2,F}}

\def\cAt{\cA^t}
\def\cAstart{\cA^{{}^*t}}
\def\cAtwoP{\cA_{2,P}}
\def\cAFt{\cA^{{}^* t}_F}
\def\cAtemp{\cA_{\rm temp}}

\def\CephypQY{C^{\ep, \hyp}_{Q,Y}(\spXp\col \tau)}
\def\Cephyp{C^{\ep, \hyp}_0(\spXp\col \tau)}

\def\cEhypQYgdglob{\cE^\hyp_{Q,Y}(\spXp\col \tau\col \gd)_\glob}

\def\vanfamF{\cE^{\rm hyp}_F(\spX\col \tau)}
\def\cEhypPallX{\cE^\hyp_P(\spX\col \tau)}

\def\cEhypQYgd{\cE^\hyp_{Q,Y}(\spXp\col \tau\col \gd)}

\def\cEhypQallX{\cE^\hyp_Q(\spX\col \tau)}

\def\Mer{\cM}

\def\Ci{C^\infty}

\def\Cci{C_c^\infty}

\def\oCFv{{}^\circ \cC_{F,v}}
\def\oC{{}^\circ\cC}

\def\DX{{\msy D}(\spX)}

\def\laur{{\rm laur}}
\def\Ad{{\rm Ad}}
\def\Wave{{\cal J}}
\def\tWave{\widetilde{\Wave}}
\def\Fou{\cF}
\def\tFou{\widetilde{\Fou}}

\def\iFv{{\rm i}_{F,v}}
\def\Exp{{\rm Exp}\,}
\def\Res{{\rm Res}}
\def\uRes{\underline{\rm Res}}

\def\Ltwod{L^2_{\rm d}}

\def\faq{\fa_\iq}

\def\faqd{\fa_{\iq}^*}
\def\faqdc{\fa_{\iq\iC}^*}
\def\faFqdc{\fa_{F\iq\iC}^*}
\def\faFq{\fa_{F\iq}}
\def\faQq{\fa_{Q\iq}}
\def\faQqdc{\fa_{Q\iq\iC}^*}
\def\faQqd{\fa_{Q\iq}^*}
\def\faQqp{\fa_{Q\iq}^+}
\def\faFprq{\fa_{F'\iq}}
\def\faFqp{\fa_{F\iq}^+}
\def\faFprqp{\fa_{F'\iq}^+}

\def\faPqdc{\fa_{P\iq\iC}^*}
\def\faPqp{\fa_{P\iq}^+}
\def\faqp{\fa_\iq^+}

\def\faP{\fa_P}
\def\faqreg{\fa_{\iq}^{\rm reg}}
\def\staPqd{{}^*\fa_{P\iq}^*}
\def\staFqdc{{}^*\fa_{F\iq\iC}^*}
\def\staFqd{{}^*\fa_{F\iq}^*}
\def\staPq{{}^*\fa_{P\iq}}
\def\staQq{{}^*\fa_{Q\iq}}
\def\faFqd{\fa_{F\iq}^*}
\def\staFq{{}^*\fa_{F\iq}}
\def\staQqdc{{}^*\fa_{Q\iq\iC}^*}
\def\staQqd{{}^*\fa_{Q\iq}^*}
\def\faPq{\fa_{P\iq}}
\def\faPqd{\fa_{P\iq}^*}
\def\staPqdc{{}^*\fa_{P\iq\iC}^*}

\def\bfaqd{\bar\faqd}
\def\bfaQqd{\bar\faQqd}
\def\bfaPqd{\bar\faPqd}
\def\bfaFqd{\bar\fa_{F\iq}^*}
\def\fbkd{\fb_{\ik}^*}
\def\fbk{\fb_\ik}
\def\fbp{\fb_\ip}
\def\stfbP{{}^*\fb_{P}}

\def\starfb{{}^*\fb}

\def\fpdual{{\fp^\dual}}
\def\fkdual{{\fk^\dual}}

\def\fgdual{\fg^\dual}

\def\fbdc{\fb_{\iC}^*}

\def\fv{{\mathfrak v}}
\def\fd{{\mathfrak d}}
\def\dfv{{}^\dual\fv}

\def\dfb{{}^\dual \fb}

\def\stfbPd{{}^*\fb_P^*}

\def\stbQdc{{}^*\fb_{Q\iC}^*}

\def\stfbPdc{{}^*\fb_{P\iC}^*}
\def\stfbPv{{}^*\fb_{P,v}}
\def\stfbQdc{{}^*\fb_{Q\iC}^*}

\def\stfb{{}^*\fb}

\def\fgc{\fg_\iC}

\def\AQq{A_{Q\iq}}
\def\Aqp{A_\iq^+}
\def\Aq{A_\iq}
\def\AQqp{A_{Q\iq}^+}
\def\AFq{A_{F\iq}}
\def\APq{A_{P\iq}}
\def\APqp{A_{P\iq}^+}
\def\Aqreg{A_\iq^\reg}
\def\stAQq{{}^*A_{Q\iq}}

\def\KM{K_\iM}

\def\MQgs{M_{Q\gs}}

\def\MoneQ{M_{1Q}}
\def\GL{{\rm GL}}

\def\minparabs{\cP_\gs^{\rm min}}
\def\allparabs{\cP_\gs}
\def\rmP{{\rm P}}
\def\parone{\cP^1_\gs}
\def\repparabs{{{\bf P}_\gs}}

\def\repPmin{{\rm P}^{\min}}

\def\staroneP{{}^{*\!1\!}P}

\def\Igdgl{I_{\gd, \gl}}

\def\WT{{\rm WT}}

\def\ik{{\rm k}}
\def\iT{{\rm T}}
\def\mc{0}

\def\DQmaps{{\rm D}_Q}

\def\DPmaps{{\rm D}_P}

\def\rmY{{\rm Y}}

\def\rmL{{\rm L}}

\def\dens#1{{}^\circ\!\!\wedge\!(#1)}

\def\genEFv{E^\circ_{F,v}}
\def\genEF{E^\circ_F}
\def\dE{E^*}
\def\nE{{E^\circ}}

\def\Eps{E_{+,s}}
\def\nEQu{E^\circ_{Q,u}}
\def\dEF{E_F^*}

\def\nEQv{E^\circ_{Q,v}}

\def\nCQP{C^\circ_{Q|P}}
\def\nC{C^\circ}

\def\WKH{W_{K \cap H}}
\def\NKaq{N_K(\faq)}
\def\NKQaq{N_{K_Q}(\faq)}
\def\FcW{{}^F\cW}
\def\QcW{{}^Q\cW}
\def\NKHaq{N_{K\cap H}(\faq)}
\def\PcW{{}^P\cW}

\def\NKPaq{N_{K_P}(\faq)}

\def\gL{\Lambda}
\def\dotvar{\,\cdot\,}

\def\pr{{\rm pr}}

\def\Cartan{\theta}

\def\prFv{\pr_{F,v}}

\def\prQv{\pr_{Q,v}}
\def\rmiQv{\rmi_{Q,v}}
\def\rmiFv{\rmi_{F,v}}
\def\rmi{{\rm i}}

\def\gsdual{\gs^\dual}

\def\lengthX{l_\spX}
\def\lspX{l_\spX}

\def\sing{{\rm sing}\,}
\def\cl{{\rm cl}\,}
\def\span{{\rm span}\,}
\def\supp{{\rm supp}\,}
\def\image{{\rm im}\,}
\def\dega{{\rm deg}_a}
\def\order{{\rm order}}
\def\rega{{\rm reg}_a}
\def\reg{{\rm reg}\,}

\def\ressup{t}

\def\Int{{\rm Int}}

\def\simPQ{\sim_{P|Q}}
\def\simQP{\sim_{Q|P}}

\def\Laur{{\rm Laur}\,}
\def\Laustar{\cL_*}
\def\Lau{{\cal L}}

\def\tauM{\tau_\iM}

\def\dK{{\widehat K}}
\def\Vtau{V_\tau}
\def\types{\vartheta}

\def\Vtaud{V_\tau^*}

\def\Ind{{\rm Ind}}

\def\oCtau{\oC(\tau)}

\begin{document}
\title{The Plancherel decomposition\\
for a reductive symmetric space \\{I. Spherical functions}}
\author{E.~P.~van den Ban and H.~Schlichtkrull}
\date{}
\maketitle
\begin{abstract}
We prove the Plancherel formula for spherical Schwartz functions
on a reductive symmetric space. Our starting point is an inversion formula
for spherical smooth compactly supported functions. The latter
formula was earlier obtained from the most continuous part of the Plancherel formula
by means of a residue calculus.
In the course of the present paper we also obtain new proofs of the uniform
tempered estimates for normalized Eisenstein integrals
and of the Maass--Selberg relations satisfied by the associated C-functions.
\end{abstract}
\tableofcontents
\section{Introduction}\eqnreset
In this paper and its sequel \bib{BSpl2}
we determine the Plancherel decomposition for a reductive symmetric
space $\spX = G/H.$ Here $G$ is a real reductive Lie group of Harish-Chandra's class
and  $H$ is an open subgroup of the group $G^\gs$ of fixed points for an involution
$\gs$ of $G.$ In the present paper we establish the Plancherel formula for $K$-finite (spherical)
Schwartz functions on $X,$ with $K$ a  $\gs$-invariant maximal compact subgroup
of $G.$ In \bib{BSpl2} we shall derive the Plancherel decomposition, in the sense
of representation theory, from it.

The space $\spX$ carries a $G$-invariant measure $dx;$ accordingly
the regular representation $L$ of $G$ in $L^2(\spX)$ is unitary. The Plancherel
decomposition amounts to an explicit decomposition of $L$ as a direct integral
of irreducible unitary representations. The reductive group $G$ is a symmetric space
of its own right, for the left times right action of $G \times G.$ In this `case
of the group',
the explicit Plancherel decomposition was obtained in the sixties
and early seventies in the  work of Harish-Chandra, see \bib{HC1}, \bib{HC2},
\bib{HC3}.
His ideas, in particular those on the role of Eisenstein integrals and
the Maass--Selberg relations satisfied by them,
have been a major influence in our work. On the other
hand, our approach to the Plancherel formula is via a residue calculus, and thus in a
sense closer in spirit to the work of
R.P.\ Langlands on the spectral decomposition in the theory of automorphic forms, see
\bib{LEis}.

The results of this paper and \bib{BSpl2}
were found and announced in the fall of 1995
when both authors were visitors of the
Mittag--Leffler Institute in Djursholm, Sweden.
At the same time P.\ Delorme announced his proof of the Plancherel theorem.
His results
have appeared in a series of papers, partly in collaboration
with J.\ Carmona,
\bib{CDn}, \bib{Dtr}, \bib{Dpl}.
At the time of the announcement
we relied on the theorem of Carmona and Delorme on the
Maass--Selberg relations for Eisenstein integrals, \bib{CDn}, Thm.~2,
which in turn relied on Delorme's paper \bib{Dtr}.
On the other hand, we also announced the proof of
a Paley--Wiener theorem for reductive symmetric spaces, generalizing
Arthur's theorem \bib{Arthur} for the case of the group.
The proof of the Paley--Wiener theorem has now appeared in \bib{BSpw},
which is independent of the present paper and \bib{BSpl2}.
The present paper as well as \bib{BSpl2} and \bib{BSpw} rely
on \bib{BSfi} and \bib{BSanfam}.

Since the time of announcement
we have been able to derive the
Maass--Selberg relations from those associated with
the most continuous part of the decomposition; these had been obtained earlier in
\bib{Bps2},\bib{Bsd}. The resulting proof of the Plancherel theorem in the present
paper and \bib{BSpl2} is independent of the one in \bib{Dpl}; moreover, it follows
a completely different approach. Finally, we mention that T.~Oshima has announced
a Plancherel formula in \bib{Opl}, p.\ 604, but the details have not appeared.
For a more extensive survey of recent developments, we refer the reader to
\bib{BSecm}.

For Riemannian symmetric spaces, the Plancherel formula is due
to Harish-Chandra \bib{HCsf} and \bib{HCds2}, p.~48.
Later, it became incorporated in the above mentioned
formula for the group. For further results in harmonic analysis on Riemannian
symmetric spaces, we refer the reader to the references given in \bib{Hel}.

For reductive symmetric spaces of type $G_\iC/G_\iR,$ the Plancherel formula
is due to P.~Harinck, \bib{Hpl}. Earlier, the Plancherel decomposition had been determined for
specific non-Riemannian spaces of rank one, by V.F.~Molchanov, J.~Faraut, G.~van Dijk and others;
see \bib{HS}, p.~185, for references.

We first give a rough
outline of the contents of this paper and its sequel \bib{BSpl2}.
The following global picture should be kept in mind.
We first concentrate on the Plancherel formula
for $K$-finite functions, with $K\subset G$ a
maximal compact subgroup that is chosen to be $\gs$-stable. The latter condition
is equivalent to the condition that the Cartan involution $\Cartan$ determined by
$K$ commutes with $\gs.$ The  building blocks of the formula will
be discrete series representations of $\spX$ and generalized
principal series of the form $\Ind_P^G(\xi \otimes \nu \otimes 1),$
with $P = M_P A_P N_P$ a $\gs\Cartan$-stable parabolic subgroup of $G$
with the indicated Langlands decomposition,
$\xi$ a discrete series representation of $\spX_P:=M_P/M_P\cap H $
and $\nu$ contained in the space $i\faPqd$ of unitary characters of $A_P/A_P \cap H.$
For purposes of exposition this introduction is written under
the simplifying assumption that the number
of open $H$-orbits on $P\backslash G$ is one. In general the
open orbits are parametrized by a finite set $\PcW$ of representatives,
and then one should take for $\xi$
the discrete series representations of all the spaces
$\spX_{P,v}: = M_P/M_P \cap v H v^{-1},$
$v \in \PcW.$

In \bib{BSmc} we obtained the most continuous part of the Plancherel
decomposition; this is the part built up from representations obtained
by induction from a minimal $\gs\Cartan$-stable
parabolic subgroup $P_\mc = M_\mc A_\mc N_\mc.$
Here $M_\mc/M_\mc\cap H$ is compact, so the theory
of the discrete series did not critically enter at this stage. On the level of $K$-finite
functions the most continuous part of the formula is described
via a Fourier transform $\Fou_\mc,$
which in turn is defined in terms of Eisenstein integrals
$\nE(P_\mc\col \gl).$ The latter are essentially matrix coefficients of the principal
series induced from $P_\mc$ and behave finitely under the action of
the algebra $\DX$ of invariant differential operators
on $\spX.$

{}From the most continuous part of the
Plancherel decomposition we derived
in \bib{BSfi} a Fourier inversion formula for functions in $\Cci(\spX\col \tau),$
the space of smooth compactly supported $\tau$-spherical functions on $\spX,$
with $\tau$ a finite dimensional unitary representation of $K.$ This formula
expresses a function $f$ in terms of the meromorphic extension of $\Fou_\mc f$
to the complexification of $i\fa_{P_0 \iq}^*.$

The strategy of the present paper is to put the inversion formula in a form
that makes it valid for functions in the Schwartz space $\cC(\spX\col\tau).$
This requires the introduction of Eisenstein integrals $\nE(P\col \nu),$
for imaginary $\nu \in i\faPqd,$
via residues of the Eisenstein integrals $\nE(P_\mc\col \dotvar).$
To show that these residual Eisenstein
integrals define Fourier transforms on the Schwartz space we need
the Maass--Selberg relations.
It is here that the theory of the discrete series, initiated  by M.\ Flensted-Jensen
in \bib{FJds} and further developed by
T.\ Oshima and T.\ Matsuki in \bib{OMds}, enters. In our proofs
we do not need the full
classification of the discrete series. However, for
the theory of the constant term developed in \bib{Cn} to apply,
both the necessity of the rank condition
and the fact that the infinitesimal $\DX$-characters of discrete series representations
are real and regular
(see Theorem \refer{t: infinitesimal characters L two d real and regular new},
due to \bib{OMds}),
play a crucial role.

The resulting inversion formula for Schwartz functions is called
the spherical Plancherel formula, see Theorems
\refer{t: inversion for Schwartz} and \refer{t: spherical Planch on Schwartz}.
It naturally leads to the spherical Plancherel formula for $L^2$-functions, Theorem
\refer{t: spherical Planch on Ltwo}.
The present paper finishes at this point, where it is not yet
clear that the residual Eisenstein integrals are related to induced
representations. This fact will be established in the second paper \bib{BSpl2}
by using the vanishing theorem of \bib{BSanfam}.
The contributions of all $K$-types can then be collected and lead to the
representation theoretic Plancherel formula.
At the end of the second paper it will also be shown that
the residual Eisenstein integrals $\nE(P\col \nu)$
equal the normalized Eisenstein integrals introduced
in Delorme's paper
\bib{D1n}. The idea is to use the automatic continuity
theorem of W.\ Casselman and N.R.\ Wallach (\bib{Cas}, \bib{Wal2})
to show that the residual Eisenstein integrals are
matrix coefficients. An asymptotic analysis
then completes the identification.

We shall now give a more detailed outline of the present paper.
The first few sections concern preliminaries. In particular, in
Section \refer{s: normalization} we specify the normalizations
of the residual operators and the measures used in the rest of the paper.
In Section \refer{s: vanishing theorem} we
give a formulation of the vanishing theorem of \bib{BSanfam},
in a form suitable for this paper. Let $\faq$ be a maximal abelian subspace
of $\fp \cap \fq,$ the intersection of the $-1$ eigenspaces in $\fg$  for $\Cartan$ and $\gs,$
respectively. Let $\allparabs$ denote the set of $\gs\Cartan$-stable
parabolic subgroups of $G$ containing $\Aq:= \exp \faq.$
For each $Q \in \allparabs$ we
introduce a space $\cE^\hyp_Q(\spX\col \tau)$
of families $\{f_\nu\}$ of spherical generalized eigenfunctions
on $\spX,$ depending meromorphically on the parameter $\nu \in \faQqdc.$
Here $\faQq:= \fa_Q \cap \fq.$
The vanishing theorem asserts that $f_\nu = 0$ for all $\nu,$
as soon as the coefficient of $e^{\nu - \rho_Q}$ in the asymptotic expansion along $Q$
vanishes, for all $\nu$ in a non-empty open subset of $\faQqdc.$

In Section \refer{s: Fourier inversion} we recall the inversion formula of
\bib{BSfi}. Let $\gS$ be the root system of $\faq$ in $\fg$ and let $W$ be the associated
Weyl group. Let $\gD$ be the system of simple roots associated with the minimal
element $P_\mc$ from $\allparabs.$ For each $F \subset \gD,$ let $P_F = M_FA_FN_F$
denote the associated standard parabolic subgroup in $\allparabs.$
Then the inversion formula is of the form
$$
f = \sum_{F \subset \gD} t(\faFqp)\, T^t_F f,
$$
where
\begin{equation}
\naam{e: formula for T t F with ge F}
T^t_Ff(x) = t(\faFqp)\, |W|\,
\int_{i \faFq + \geps_F} \int_\spX K^t_F(\nu\col x\col y)\, f(y) \; dy\; d\mu_F(\nu).
\end{equation}
Here $t$ is a choice of $W$-invariant even residue weight on $\gS$ and $\geps_F$ is an element
of $\faFqp,$ sufficiently close to zero (if $F =\gD,$ we may take
$\geps_F = 0$).
Moreover, $d\mu_F$ is a suitable choice of Lebesgue measure on $i\faFqd + \geps_F.$
The kernel functions $K^t_F(\nu\col x \col y) \in \End(\Vtau)$
are obtained from residual operators acting on a combination
of normalized and partial Eisenstein integrals for
$P_\mc,$ see \bib{BSfi} for
details.
They are meromorphic in the variable $\nu \in \faQqdc$ and smooth spherical
and $\DX$-finite in both of the variables $x, y \in \spX.$
Essentially, the idea is that
the kernel $K^t_F$ determines the projection onto the part of $L^2(\spX\col \tau)$
determined by the induction from the standard parabolic subgroup $P_F.$

To make the above formula valid for Schwartz functions
it is necessary to establish it with  $\geps_F = 0,$ for
every $F.$ This can be achieved
by using Cauchy's formula, once we have established the
regularity at $i\faFqd$ of the kernel functions $K_F^t(\nu, x, y)$
in the variable $\nu.$ In addition to this we need estimates
that are tempered in the variables $x,y$ with uniformity in $\nu \in i\faQqd.$
All this is taken care of by a long inductive argument, that ranges over the
Sections
\refer{s: spaces of residue type} --
\refer{s: Fourier inversion for Schwartz functions}.
We shall describe the structure of the argument,
which goes by induction on the $\gs$-split rank of
 $G,$
at a later stage in this introduction.

In Section \refer{s: the generalized Eisenstein integral}
we recall the definition of the generalized Eisenstein integral $\genEF(\nu),$
for $\nu \in \faFqdc.$
In \bib{BSfi}, Def.~10.7, see also \bib{BSanfam}, Remark 16.12,
this Eisenstein integral
was defined by means of a linear combination of
residual operators (a so called Laurent functional)
applied to the Eisenstein integral $\nE(P_\mc\col \gl)$ with respect to the variable $\gl.$
As a family in the parameter $\nu,$ the generalized Eisenstein integral
belongs to $\cE^\hyp_{P_F}(\spX\col \tau).$
Hence, in view of the vanishing theorem, it can be characterized uniquely in terms of its
asymptotic behavior along $P_F;$ see Theorem \refer{t: char of genEis by asymp new}.

In Corollary \refer{c: nE tempered}
we show that the generalized
Eisenstein integral is tempered for regular imaginary values of $\nu.$ This fact
can be derived from a limitation on the asymptotic exponents,
see Theorem \refer{t: exponents of genEF},
caused by the support of the residual operators. Here Thm.~3.15 of \bib{BSres} is crucial.
In the next section, in Proposition
\refer{p: init estimate genEis}, we establish
uniformly moderate estimates for the generalized Eisenstein integrals.
These estimates come from similar estimates for $\nE(P_\mc\col \dotvar),$
which survive the application of the residual operators.

In Section \refer{s: spaces of residue type} we start with the preparation of the long
inductive argument mentioned above.
The reductive symmetric
pair $(G,H)$
is said to be of residue type if the following two conditions
 are fulfilled,
(a)  $G$ has compact center
modulo $H$ and (b) the operator $T^t_\gD$ is the projection
onto the discrete part $L^2_d(\spX\col \tau)$ of $L^2(\spX\col \tau)$
(which may be trivial).
Condition
 (b) implies that the mentioned operator,
and hence the associated `discrete' kernel
$K^t_\gD,$ is independent of the particular choice of residue weight $t.$
Moreover, from condition (b) it follows  that $L^2_d(\spX\col \tau)$
is finite dimensional,
a result known for all pairs $(G,H)$ as a consequence of \bib{OMds},
see also Remark \refer{r: finite dimensional Ltwod}.

We proceed with the induction in Section \refer{s: the normalized
Eisenstein integral}. A parabolic subgroup $P \in \allparabs$ is
said to be of residue type if the pair $( M_P , M_P \cap H)$ is.
A subset $F \subset \gD$ is said to be of residue type, if
the associated standard parabolic subgroup $P_F$ is.
In the course of the inductive argument, many results in Sections
\refer{s: spaces of residue type} --
\refer{s: Fourier inversion for Schwartz functions}
will initially
be proved under the assumption that $(G,H)$ or a parabolic subgroup from
$\allparabs$ is of residue type.
Such results will always be indicated with the abbreviation {(RT)}
following their declaration. The chain of results marked (RT)
is needed in the induction step of Theorem \refer{t: inversion for Schwartz},
where by induction on the $\gs$-split rank of $G$
it is shown that  all groups from $\allparabs$
are of residue type. In particular, also all
pairs $(G,H),$ with $G$ having compact center modulo $H,$ are of residue type.
It thus follows that the results marked (RT)
are valid in the generality stated (see also Remark \remRT).

The kernel $K^t_F(\nu \col \dotvar\col \dotvar)$ is determined by its asymptotic
expansion along $P_F\times P_F,$
in view of the vanishing theorem. The coefficient of $e^{\nu -\rho_F}\otimes e^{-\nu - \rho_F}$
in this expansion
is essentially the discrete kernel of $M_F/M_F\cap H.$ If $F$ is of residue type,
then the discrete kernel,
and hence $K_F^t,$  is independent of the particular choice of $t.$ Therefore,
so is the generalized Eisenstein integral.
{}From then on we call this Eisenstein integral the normalized one
and denote it by $\nE(P_F\col \nu).$ It is a meromorphic function
of $\nu \in \faFqdc,$ with values in $\Ci(\spX) \otimes \Hom(\cA_{2,F}, \Vtau),$
where $\cA_{2,F} = L^2_d(\spX_F \col \tau_F).$ (Without the simplifying assumption
mentioned above, the latter space is replaced by a suitable direct sum over $\FcW.$)
The unique characterization
of the normalized Eisenstein integral by means of the vanishing theorem then
allows us to define it for $P_F$ replaced by any parabolic subgroup $P \in \allparabs$
of residue type.
In the case of the group,
the characterization allows us to identify the normalized
Eisenstein integral with
Harish-Chandra's, as defined in \bib{HCeis}, \S 6, Thm.\ 6; see
Remark \refer{r: comparison with nE of HC}.

The definition of the normalized Eisenstein integral in turn
allows us to define a kernel function
$ K_P$ for arbitrary $P \in \allparabs$ of residue type, generalizing the kernels
for standard parabolic subgroups of residue type. In terms of the normalized Eisenstein
integrals, the kernel is given by the formula
\begin{equation}
\naam{e: formula kernel}
K_P(\nu \col x \col y) = |W_P|^{-1} \nE(P\col \nu \col x)\dE(P\col \nu \col y).
\end{equation}
Here $W_P$ is the subgroup of $W$ that corresponds to the Weyl group of $\spX_P,$
and $\dE(P\col \nu \col \dotvar)$ is the dualized Eisenstein integral. The latter is
the  function in
$\Ci(\spX)\otimes \Hom(\Vtau, \cA_{2,P})$
defined by
$$
\dE(P\col \nu \col y): = \nE(P \col -\bar \nu \col y)^*.
$$
Two parabolic subgroups $P,Q \in \allparabs$ are said to be associated, notation $P\sim Q,$
if their $\gs$-split components $\faPq$ and $\faQq$ are conjugate under
$W.$ If $P, Q \in \allparabs$ are associated, and if one of them is of residue type,
then so is the other. Moreover, the set
$$
W(\faQq\mid \faPq) = \{s|_{\faPq} \mid s \in W,\;\;s(\faPq) \subset \faQq\}
$$
is non-empty. The main result of the section is
Theorem \refer{t: invariance kernel}, which asserts
that the kernel function is unchanged if $P, \nu$ are replaced by $Q, s\nu,$
with $Q \sim P$ and $s \in W(\faQq\mid \faPq).$
For $P$ minimal, this result is
a consequence of the Maass--Selberg relations for the Eisenstein integral
$\nE(P\col \dotvar),$
in view of (\refer{e: formula kernel}). For arbitrary $P$ of
residue type the result follows from the minimal
case by $W$-equivariance properties of the
residue calculus.

In Section \refer{s: eigenvalues} we describe
the action
of $\DX$ on the normalized Eisenstein integral $\nE(P\col \nu),$ for $P$ of residue type.
The diagonalization of $\D(\spX_P)$ on
$L^2_d(\spX_P\col \tau_P),$ where $\tau_P:= \tau|_{K \cap P},$
induces a simultaneous diagonalization of the action
of $\DX$ on the Eisenstein integral, in view of the vanishing theorem, see
Corollary \refer{c: deco cAPtwo in eigenspaces}.
In the next section this result is used to show
that the uniform moderate estimates of the Eisenstein integral can be improved
to uniform tempered estimates, exploiting a technique that goes back to \bib{Bps2} and \bib{Wal1}

In Section \refer{s: infinitesimal characters} we recall the mentioned result of \bib{OMds}
on the $\DX$-characters of the discrete series in
Theorem \refer{t: infinitesimal characters L two d real and regular new}.
This result is of crucial importance for Section \refer{s: the constant term},
where we determine the constant term of the normalized Eisenstein integral.
In addition we  use the theory of the constant term as
developed in \bib{Cn}, see  also \bib{HC1}.

The constant term of $\nE(P\col \nu)$ along a $Q \in \allparabs$ with
$Q\sim P$ describes $e^{\rho_Q}$ times the top order asymptotic behavior along $Q;$ it is given by
\begin{equation}
\naam{e: constant term of nE along Q}
E^\circ_Q(P\col \nu \col ma)\psi = \sum_{s \in W(\faQq\mid\faPq)}
a^{s \nu}\,[\nC_{Q\mid P} (s \col \nu) \psi](m),
\end{equation}
for $\nu \in \faPqd$ generic,  $a \in A_{Q\iq}$ and $m\in \spX_Q.$
Here $\nC_{Q|P}(s\col \nu),$ the normalized $C$-function,
is a meromorphic $\Hom(\cAtwoP, \cAtwoQ)$-valued function
of $\nu \in \faPqdc.$ For the description of the constant term without the
simplifying assumption mentioned above, see Corollary
\refer{c: constant term and C}.

In Section \refer{s: the Maass--Selberg relations}
we derive the Maass--Selberg relations from
the invariance property of the kernel (\refer{e: formula kernel}) mentioned
above. They assert the following identity of meromorphic functions in the variable
$\nu \in \faPqd,$
$$
\nC_{Q|P}(s\col - \bar \nu)^*\nC_{Q|P}(s\col \nu)= I,
$$
for $P,Q\in \allparabs$ associated and
of residue type, and for  $s \in W(\faQq\mid \faPq).$
In particular, it follows that the normalized $C$-functions are unitary for imaginary $\nu.$
This in turn shows that the constant term
(\refer{e: constant term of nE along Q}) is regular for imaginary $\nu.$ By
a result from  \bib{BCDn}
this implies that the Eisenstein integral $\nE(P\col \nu)$ is regular for imaginary
$\nu,$ see Theorem \refer{t: regularity nE}.
Because of the uniform tempered estimates formulated in Corollary \refer{c: second estimate dE},
it becomes possible to define a spherical Fourier transform $\Fou_P$
in the next section by the formula
\begin{equation}
\naam{e: spherical Fou in intro}
\Fou_P f (\nu) = \int_\spX \dE(P\col \nu\col x) f(x)\, dx,
\end{equation}
for $f \in \cC(\spX\col \tau),$ the space of $\tau$-spherical Schwartz functions on $\spX.$
Proposition \refer{p: Fou continuous on Schwartz space} asserts that $\Fou_P$ is a
continuous linear
map into the Euclidean Schwartz space $\cS(i\faPqd) \otimes \cAtwoP,$
if $P$ is of residue type.
In Section \refer{s: Wave packet} it is shown, using a result from \bib{BCDn},
that the adjoint wave packet transform, given by the
formula
\begin{equation}
\naam{e: wave packet in intro}
\Wave_P \gf(x) = \int_{i\faPqd} \nE(P\col \nu\col x) \gf(\nu)\, d\mu_P(\nu)
\end{equation}
is a continuous linear map from $\cS(i\faPqd) \otimes \cAtwoP$ into $\cC(\spX\col \tau),$
see Theorem \refer{t: Wave continuous on Schwartz space}. Here $d\mu_P$ is
Lebesgue measure on $i\faPqd,$  normalized as in Section \refer{s: normalization}.

In Section \refer{s: Fourier inversion for Schwartz functions} the long inductive
argument is completed as follows. In the proof
of Theorem \refer{t: inversion for Schwartz},
it is shown that every $P \in \allparabs$ is of residue
type, by induction on the $\gs$-split rank of $G.$ The hypothesis
of the induction step implies that one may assume that $G$ has compact center
modulo $H$  and that every $F \subsetneq \gD$ is of residue type.
In view of the regularity of the normalized Eisenstein integrals, hence of the
kernels $K_F(\nu \col \dotvar\col \dotvar),$ for $F \subsetneq\gD$ and  $\nu \in i\faFqd,$
the formula (\refer{e: formula for T t F with ge F}) becomes valid with $\geps_F = 0$
for every subset $F \subset \gD$ (recall
that $\geps_\gD = 0$). Moreover, by the definition of the transforms $\Fou_F := \Fou_{P_F}$
 and $\Wave_F := \Wave_{P_F},$ it takes
the form
$$
f = T_\gD^t f + \sum_{F \subsetneq \gD} [W:W_F] \,t(P_F)\, \Wave_F\Fou_F f.
$$
{}From this one reads off that $T_\gD^t$ maps $\Cci(\spX\col \tau)$ into $\cC(\spX\col \tau),$
from which it readily follows that $T_\gD^t$
is the restriction to $\Cci(\spX\col \tau)$
of the
orthogonal projection onto $L^2_d(\spX\col \tau).$
This argument completes the induction step; moreover, at the same time it shows that
\begin{equation}
\naam{e: I as sum Wave Fou}
I = \sum_{F \subset \gD} [W:W_F] \,t(P_F)\,\Wave_F \Fou_F
\end{equation}
on $\Cci(\spX\col \tau)$ and hence on $\cC(\spX\col \tau)$ by continuity and density.
It is at this point that the role of the residue weight in the harmonic analysis becomes
clear. Define the equivalence relation $\sim$ on the powerset $2^\gD$ by $F\sim F'\iff
P_F \sim P_{F'}.$ Then by the Maass--Selberg relations the composed transform $\Wave_F \Fou_F$
only depends on the class of $F$ in $2^\gD/\!\!\sim.$ Collecting the terms in
(\refer{e: I as sum Wave Fou})
according to such classes we obtain, by an easy counting argument, the following
Fourier inversion formula which is independent of the choice of residue weight
\begin{equation}
\naam{e: last inversion on Schwartz}
I = \sum_{[F] \in 2^\gD/\sim} [W:W^*_F]\, \Wave_F \Fou_F,
\end{equation}
with $W^*_F$ the normalizer of $\faFq$ in $W.$
In other words, the residue weight describes
the weight by which an element in the class
of $F$ contributes to the above inversion formula.

In Section \refer{s: properties of the Fourier transforms}
we give a precise description of the
kernels and images of the Fourier
transforms and their adjoints. This leads to the spherical Plancherel theorem
for Schwartz functions, Theorem \refer{t: spherical Planch on Schwartz}
and the similar theorem for $L^2$-functions,
Theorem \refer{t: spherical Planch on Ltwo}.
In particular, the summands in (\refer{e: last inversion on Schwartz})
extend to  $L^2(\spX\col \tau)$ and become orthogonal projections onto mutually orthogonal
subspaces.

We are grateful to the organizers of the special year at the
Mittag--Leffler Institute, during which these results were found,
and we thank Mogens Flensted-Jensen for several helpful discussions.
We are also grateful to Jacques Carmona and Patrick Delorme for
informing us about their results on the Maass--Selberg relations,
which played a crucial role for us, as mentioned above.

\section{Notation and preliminaries}\eqnreset
\naam{s: notation}
Throughout this paper, $G$ will be a real reductive group of Harish-Chandra's
class, $\gs$ an involution of $G$ and $H$ an open subgroup
of $G^\gs,$ the set of fixed points for $\gs.$
The associated reductive symmetric space is denoted by
$$
\spX = G/H.
$$
The algebra of $G$-invariant differential operators on $\spX$ is denoted by $\DX.$

We fix a Cartan involution $\Cartan$ of $G$ that commutes with
$\gs;$ thus, the associated maximal compact subgroup $K = G^\Cartan$
is $\gs$-invariant. We follow the convention to
denote Lie groups with roman capitals, and their Lie algebras
with the corresponding gothic lower cases. In particular, $\fg$ denotes
the Lie algebra of $G.$ The infinitesimal involutions of $\fg$
associated with $\Cartan$ and $\gs$ are denoted by the same symbols.

We equip $\fg$ with a $G$-invariant
non-degenerate bilinear form $B$ that restricts to the Killing form on $[\fg, \fg],$
that is positive definite on $\fp,$ negative definite on $\fk,$ and
for which $\gs$ is symmetric. Then $\inp{X}{Y}: = - B(X, \Cartan Y)$
defines a positive definite inner product on $\fg$ for which the involutions
$\Cartan$ and $\gs$ are symmetric.
Accordingly, the decompositions $\fg = \fk \oplus \fp = \fh \oplus \fq$
into the $+1$ and $-1$ eigenspaces of these involutions,
respectively, are orthogonal. If $\fv \subset \fg$ is a linear subspace,
we agree to identify $\fv^*$ with a linear subspace of $\fg$
via the inner product $\inp{\cdot}{\cdot}.$ Finally, we equip the linear dual $\fg^*$
of $\fg$ with the dual inner product, and the complexified dual $\fg^*_\iC$ with
its  complex bilinear extension.

We fix a maximal abelian subspace $\faq$ of $\fp \cap \fq$ and
denote the associated system of restricted roots
by $\gS = \gS(\fg, \faq).$ This is a, possibly non-reduced, root system;
the associated Weyl group is denoted by $W = W(\gS).$ We recall that
$W \simeq N_K(\faq)/Z_K(\faq),$ naturally. Accordingly, the natural
image of $N_{K\cap H}(\faq)$ in $W$ is denoted by $W_{K\cap H}.$

If $P$ is a parabolic subgroup of $G,$ we denote its Langlands decomposition by
$P = M_P A_P N_P$ and put $M_{1P}:= M_P A_P.$
A parabolic subgroup that is invariant under the composed involution $\gs\Cartan$
is called a $\gs$-parabolic subgroup. The set of $\gs$-parabolic subgroups containing
$\Aq:= \exp \faq$ is finite and denoted by $\allparabs.$

We shall briefly recall the structure
of the parabolic subgroups from $\allparabs$, meanwhile fixing notation.
For details we refer to \bib{BSanfam}, \S 2.
If $P \in \allparabs,$ then $M_P$ and $A_P$ are $\gs$-invariant, and
$\gs N_P = \Cartan N_P = \bar N_P.$
The algebra $\faP$ is $\gs$-invariant, hence decomposes as $\fa_P = \fa_{P\ih} \oplus \faPq,$
the vector sum of the intersections of $\fa_P$ with
$\fh$ and $\fq,$ respectively.
We put $\APq:= \exp \faPq$ and $M_{P\gs}: = M_P (A_P \cap H)$ and call
$
P = M_{P\gs} \APq N_P
$
the $\gs$-Langlands decomposition of $P.$

As usual, we denote by $\rho_P$ the linear functional
$\frac12 \tr [\ad(\dotvar)|_{\fn_P}] \in \fa_P^*.$
The following lemma is of importance in
the theory of induced representations.
\begin{lemma}
Let $P \in \allparabs.$ Then $\rho_P$ belongs to $\faPqd.$
\end{lemma}
\proof
The algebras $\fa_P$ and $\fn_P$ are $\gs\Cartan$-invariant, hence $- \gs \rho_P =
\gs \Cartan \rho_P = \rho_P.$ This implies that $\rho_P$ vanishes on $\fa_{P\ih},$ hence
belongs to $\faPqd.$
\qed

The space $\faPq$ is contained in $\faq.$ Let $\gS_P$ denote
the collection of roots from $\gS$ that vanish on $\faPq.$ Then
$$
\faPq = \cap_{\ga \in \gS_P} \ker \ga.
$$
The subgroup of $W$ generated by the reflections in the roots of $\gS_P$
is denoted by $W_P.$ It equals the centralizer of $\faPq$ in $W.$

Let $\gS(P)$ be the collection
of roots from $\gS$ that occur in $\fn_P$ as an $\faq$-weight.
Then $\fn_P$ is the vectorial direct sum of the root spaces $\fg_\ga,$ for $\ga \in \gS(P).$
We put
$$
\gS_r(P) := \gS(\fn_P, \faPq) = \{ \ga|_{\faPq}  \mid \ga \in \gS(P)\}.
$$
The set
$$
\faPqp := \{ X \in \faPq\mid \ga(X) > 0\text{for all} \ga \in \gS(P)\}
$$
is non-empty. Therefore, the elements of $\gS_r(P)$ are non-zero
linear functionals on $\faPq.$ Moreover,
$\faPqp$ is a connected component of the complement $\fa_{P\iq}^\reg$
of the union of their null spaces.
We put $\APqp:= \exp \faPqp.$

The collection of weights in $\gSr(P)$ that cannot be expressed as a sum of two
elements of $\gSr(P)$ is denoted by $\gDr(P).$ We recall from \bib{BSanfam}, beginning
of \S 3, that the set $\gDr(P)$ is linearly independent over $\R$ and spans $\gS_r(P)$
over $\N.$

If $X \in \faq,$ then $X \in \fp,$ hence $\ad X$ diagonalizes with real eigenvalues.
It is well known that the sum  of the eigenspaces for the non-negative eigenvalues
is a parabolic subalgebra of $\fg.$ Its $\Cartan$-stable Levi component $\fm_{1X}$
and its nilpotent radical $\fn_X$ are given by
$$
\fm_{1X} = \ker \ad X,\qquad \fn_X = \oplus_{\ga\in \gS, \,\ga(X) > 0}\; \;\;\fg_\ga.
$$
The associated parabolic subgroup of $G$ is denoted by $P_X.$
If $P \in \allparabs$ and $X \in \faPqp,$
then it follows from \bib{BSanfam}, Eqn.\ (2.4), that $P = P_X.$
{}From $\gs\Cartan X = X$ it follows that $P_X \in \allparabs.$

Let $\sim$ be the relation of parabolic equivalence on $\faq,$ with respect to
the root system $\gS.$ Thus, $X \sim Y$ if and only if for each $\ga \in \gS$
we have $\ga(X)> 0 \iff \ga(Y) > 0.$ It readily follows from the definition
given above that $X \sim Y \iff P_X = P_Y.$

\begin{lemma}
\naam{l: parabolic equivalence classes}
The map $P \mapsto \faPqp$
is a bijection from $\allparabs$ onto the set $\faq/\!\sim$ of parabolic
equivalence classes.
\end{lemma}

\proof
If $P \in \allparabs$ and $X \in \faPqp$ then $P = P_X,$ as said above. Hence, the map
$X \mapsto P_X$ is a surjection from $\faq$ onto $\allparabs.$
By the last assertion before Lemma \refer{l: parabolic equivalence classes}, the map factors
to a bijection from $\faq/\!\sim$ onto $\allparabs.$
If $X \in \faq,$ let $P = P_X.$ Then  $\gS(P) = \{\ga \in \gS \mid \ga(X) > 0\},$
hence $\gS_P = \{\ga \in \gS\mid \ga(X) = 0\},$ and we see that $[X] = \faPqp.$
Thus, $P \mapsto \faPqp$ is the inverse to $[X] \mapsto P_X.$
\qed

It follows from the description in Lemma \refer{l: parabolic equivalence classes}
that the Weyl group $W$ acts on the finite set
$\allparabs.$
We recall from \bib{BSres}, Def.\ 3.2, that a residue weight on $\gS$ is
a map $\faqd/\!\!\sim \:\to [0,1]$ such that for every $Q \in \allparabs,$
\begin{equation}
\naam{e: WT sum one}
\sum_{P \in \allparabs,\; \faPq = \faQq} t(\faPqp) = 1.
\end{equation}
The collection of residue weights on $\gS$ is denoted by $\WT(\gS).$
Via the bijection of Lemma \refer{l: parabolic equivalence classes},
a weight $t \in \WT(\gS)$ will also be viewed as a map $t: \allparabs \to [0,1].$
A residue weight $t \in \WT(\gS)$ is said to be $W$-invariant
if $t(w \faPqp) = t(\faPqp)$ for all $P \in \allparabs$ and $w \in W,$
and even if $t(\faPqp) = t(-\faPqp)$
for all $P \in \allparabs.$

Let $\minparabs$ the collection of minimal elements in $\allparabs.$
Then $P \mapsto \faPqp$ is a bijection from $\minparabs$ onto the collection
of open chambers for $\gS$ in $\faq.$ To emphasize this, we shall also write
$\faqp(P):= \faPqp$ and $\Aqp(P): = \APqp$ for $P \in \minparabs.$
Accordingly, $W$ acts simply
transitively on $\minparabs.$ Note that for $P \in \minparabs,$ $\gS(P) = \gS_r(P)$
is  a positive system for $\gS$
and $\gD(P):= \gD_r(P)$ the associated collection of simple roots.

We fix a system $\gS^+$ of positive roots for $\gS;$ let $\gD$ be the associated
collection of simple roots. Given $F \subset \gD$ we define
$$
\faFq: = \cap_{\ga \in F} \ker \ga
$$
and denote by $\faFqp$ the subset of elements $X \in \faFq$ such that $\gb(X) > 0$ for
$\gb \in \gD\setminus F.$ Then $\faFqp$ is a parabolic equivalence class.
The associated parabolic subgroup $P_F$ is called the standard parabolic
subgroup determined by $F.$ We adopt the convention to replace an index or superscript
$P_F$ by $F.$
In particular, the Langlands decomposition of $P_F$ is denoted by $P_F = M_F A_F N_F$
and the centralizer of $\faFq$ in $W$ by $W_F.$
Let
$$
W^F: = \{s \in W\mid s(F) \subset \gS^+\}.
$$
Then the canonical map $W \to W/W_F$
induces a bijection $W^F \to W/W_F.$

We write $P_0$ for $P_\emptyset,$ $P_0 = MAN_0$ for its Langlands
decomposition and
$M_1 :=MA.$
Then $P = M_1 N_P$ for every $P \in \minparabs.$

If $P \in \allparabs$ and $v \in \NKaq,$ we define
\begin{equation}
\naam{e: defi spXPv}
\spXPv: = M_P / M_P \cap vHv^{-1}.
\end{equation}
Here $M_P$ is a real reductive group of Harish-Chandra's class and $M_P \cap vHv^{-1}$
is an open subgroup of the group of fixed points for the involution $\gs^v: M_P \to M_P$
defined by $\gs^v(m) = v \gs( v^{-1} m v) v^{-1}.$
Thus, the space in (\refer{e: defi spXPv}) is a
reductive symmetric space in the class under consideration.
Moreover, $\Cartan|_{M_P}$ is a Cartan involution of $M_P$ that commutes with
$\gs^v;$ the associated maximal compact subgroup is $K_P := K \cap M_P.$

Note that as (an isomorphism class of) a $M_P$-homogeneous space,
the symmetric space $\spXPv$ depends on $v$ through
its class in the double coset space
$W_P\backslash W / W_{K\cap H}.$ Throughout this paper,
$\PcW$ will denote
a choice of representatives in $N_K(\faq)$ of
$W_P\backslash W/W_{K\cap H}.$ In general, if $f$ is a surjective map
from a set $A$ onto a set $B,$ then by a choice of representatives
for $B$ in $A,$ we mean a subset $\cB \subset A$ such that $f|_{\cB}: \cB \to B$
is a bijection.

Let $\staPq$ denote the orthocomplement of $\faPq$ in $\faq.$ Then
$$
\staPq = \fm_P \cap \faq.
$$
Moreover, for every $v \in N_K(\faq),$ this space is the analogue
of $\faq$ for the triple $(M_P, K_P, M_P \cap v H v^{-1}),$
see \bib{BSanfam}, text following (3.4).

In analogy with (\refer{e: defi spXPv}), we define
$\spXonePv:= M_{1P}/M_{1P} \cap v H v^{-1},$
for $P \in \allparabs$ and $v \in \NKaq.$ The multiplication
map $M_P \times \APq \to M_{1P}$ induces a diffeomorphism
\begin{equation}
\naam{e: deco spXonePv}
\spXonePv \simeq \spXPv \times \APq.
\end{equation}
If $v = e,$ we agree to omit $v$ in the notation of the spaces in this product,
so that $\spX_{1P} = M_{1P}/M_{1P}\cap H \simeq M_P/M_P\cap H \times \APq.$

We end this section with collecting some basic facts about
$\Cartan$-stable Cartan subspaces of $\fq,$ meanwhile
fixing notation. We
define the dual real form $\fg^d$ of $\fg$
as the real form of $\fg_\iC$ given by
$
\fg^d = \ker(\gs\Cartan - I) \oplus i \ker (\gs \Cartan  + I).
$
Let $\Cartan_\iC$ and $\gs_\iC$
be the complex linear extensions of $\Cartan$ and $\gs,$ respectively.
Then $\Cartan^d:= \gs_\iC|_{\fg^d}$ is a Cartan involution of $\fg^d$
and $\gs^d:= \Cartan_\iC|_{\fg^d}$ is an involution of $\fg^d$ commuting
with $\Cartan^d.$

If $\fv$ is any $\gs$- and  $\Cartan$-stable  subspace of $\fg,$
then $\dfv: = \fv_\iC \cap \fgdual$ is a
$\gs^\dual$- and  $\Cartan^\dual$-stable subspace of $\fgdual,$ whose
complexification equals that of $\fv.$

If $\fb$ is a $\Cartan$-stable Cartan subspace of $\fq,$ then
$\fb = \fbk \oplus \fbp,$ where $\fbk:= \fb \cap \fk$ and $\fbp: = \fb \cap \fp.$
Moreover,
$$
\dfb: = i \fbk \oplus \fbp
$$
is a $\gsdual$-stable maximal abelian subspace
of $\fpdual.$
We denote by
$\gS(\fb)$ the root system of $\dfb$ in $\fgdual,$
by $W(\fb)$ the associated Weyl group and by
$I(\fb)$ the space of $W(\fb)$-invariants
in $S(\fb),$ the symmetric algebra of $\fb_\iC.$
Moreover, we denote the
associated Harish-Chandra isomorphism by
\begin{equation}
\naam{e: ggfd}
\gg_{\dfb}: \;\;\;U(\fgdual)^\fkdual/U(\fgdual)^\fkdual \cap U(\fgdual) \fkdual
\to I(\fb).
\end{equation}
As usual, if $\fl$ is a real Lie algebra, we denote by $U(\fl)$ the universal
algebra of its complexification.
Via the natural isomorphism
\begin{equation}
\naam{e: iso invariant diff ops}
\DX \simeq U(\fg)^\fh/U(\fg)^\fh \cap U(\fg)\fh =
U(\fgdual)^\fkdual/U(\fgdual)^\fkdual \cap U(\fgdual)\fkdual,
\end{equation}
see \bib{Bps2}, Lemma 2.1,
we shall identify the algebra $\DX$ with the algebra
on the left-hand side of (\refer{e: ggfd}) and thus view
the Harish-Chandra isomorphism $\gg_{\dfb}$
as an algebra isomorphism from $\DX$ onto $I(\fb);$ as such it is denoted
by $\gg = \gg_\fb.$

If $P \in \allparabs$ and $\fb$ a $\Cartan$-stable Cartan subspace of $\fq$
containing $\faPq,$
we agree to write $\starfb_P: = \fb \cap \fm_P.$ Then $\starfb_P$ is
a $\Cartan$-stable Cartan subspace of $\fm_{P} \cap \fq$ and
\begin{equation}
\naam{e: deco fb with star part}
\fb = \starfb_P \oplus \faPq,
\end{equation}
with orthogonal summands. If $P$ is minimal, then $\fb$ is maximally split,
and we suppress the index $P,$ so
that $\fb = \starfb \oplus \faq.$
We shall write $W(\starfb_P)$ for
the Weyl group of the pair $(\fm_{P\iC}, \starfb_P).$ Via the decomposition
(\refer{e: deco fb with star part}) this Weyl group
is naturally identified with  $W_P(\fb),$  the centralizer of $\faPq$ in $W(\fb).$

\section{Weyl groups}\eqnreset
\naam{s: Weyl groups}
In this section we discuss a straightforward
generalization of well known results on Weyl groups,
see \bib{HC1}, p.\ 111.

If $\fa_1$ and $\fa_2$ are abelian subspaces of $\fp,$
then following \bib{HC1}, p.\ 112, we define the set
$$
W(\fa_2\mid\fa_1) :=
 \{s \in \Hom(\fa_1, \fa_2) \mid \exists \; g \in G:\;\; s = \Ad(g)|_{\fa_1}\}.
$$
{}From the definition it is obvious that the set $W(\fa_2\mid \fa_1)$ consists of
injective linear maps. In particular, if $\dim \fa_1 = \dim \fa_2,$ it consists
of linear isomorphisms. Finally, if $\fa_1 = \fa_2,$ the set is a subgroup of
$\GL(\fa_1) = \GL(\fa_2).$
We note that $W(\fa_1\mid \fa_1)$ naturally acts from the right on $W(\fa_2\mid \fa_1),$
whereas $W(\fa_2\mid \fa_2)$ naturally acts from the left. If $\dim \fa_1 = \dim \fa_2,$ then
both of these actions are transitive and free. If $\fa_\fp$ is a maximal abelian subspace
of $\fp,$
then by $W(\fg, \fa_\fp)$ we denote the Weyl group of  $\gS(\fg, \fa_\fp), $
the root system of $\fa_\fp$ in $\fg.$

\begin{lemma}
\naam{l: first lemma on Weyl groups}
Let $\fa_1$ and $\fa_2$ be abelian subspaces of $\fp.$
\begin{enumerate}
\itema
The set $W(\fa_1\mid \fa_2)$ is finite.
\itemb
If $\gf \in \Int(\fgc)$ maps $\fa_1$ into $\fa_2$ then $\gf|_{\fa_1} \in W(\fa_2\mid \fa_1).$
\itemc
If $s \in W(\fa_2\mid\fa_1),$ then there exists a $k \in K_e$ such that $s = \Ad(k)|_{\fa_1}.$
\itemd
Assume that $\fa_1$ and $\fa_2$ are contained in a maximal abelian subspace
$\fa_\fp$ of $\fp.$
Then
$$
W(\fa_2\mid \fa_1) = \{ t \in \Hom(\fa_1, \fa_2) \mid \exists s \in W(\fg, \fa_\fp):\;\;
 t =s|_{\fa_1}\}.
$$
\end{enumerate}
\end{lemma}

\proof
All assertions are immediate consequences of Corollaries 1, 2 and 3 of
\bib{HC1}, p.\ 112.
\qed

\begin{cor}{\ }
\naam{c: W is W faq faq}
\begin{enumerate}
\itema
$W = W(\faq\mid \faq).$
\itemb
Let $\fa_\fp$ be a maximal abelian subspace of $\fp,$ containing $\faq.$
Then the map $k \mapsto \Ad(k)|_{\faq}$ is a surjection
from $N_{K_e}(\faq)  \cap N_{K_e}(\fa_\fp)$ onto $W.$
\end{enumerate}
\end{cor}

\proof
The map $k \mapsto \Ad(k)|_{\faq}$ induces a natural isomorphism
$\NKaq/Z_K(\faq) \simeq W,$ see, e.g., \bib{Bps1}, Lemma 1.2.
Hence, $W \subset W(\faq\mid \faq).$ For the converse inclusion, select a maximal
abelian subspace $\fa_\fp$ of $\fp,$ containing $\faq.$ Then by
Lemma \refer{l: first lemma on Weyl groups} (d), any element $t \in W(\faq \mid \faq)$
is the restriction of an element $s \in W(\fg, \fa_\fp).$ There
exists a $k \in N_{K_e}(\fa_\fp)$ such that $s = \Ad(k)|_{\fa_\fp}.$ The element
$k$ necessarily normalizes $\faq.$ Thus, we obtain the converse inclusion
and also the validity of assertion (b).
\qed

The following lemma generalizes Lemma 1 of \bib{HC1}, p.\ 111.
Let $\faq$ be a maximal abelian subspace of $\fp \cap \fq.$ Let $\fa_\fp$ be
a maximal abelian subspace of $\fp$ containing $\faq,$ and
$\fj$ a Cartan subalgebra of $\fg$ containing $\fa_\fp.$ We denote
by $W(\fg_\iC, \fj_\iC)$ the Weyl group of the root system of $\fj_\iC$ in $\fg_\iC.$

\begin{lemma}
\naam{l: second lemma on Weyl groups}
Two elements of $\faq$ are conjugate under $\Int(\fg_\iC)$ if and only if they are
conjugate under any one of the following groups
$$
W(\fgc , \fj_\iC), \;\;W(\fg, \fa_\fp),\;\; W = W(\fg, \faq),
\;\;N_{K_e}(\fa_\fp) \cap N_{K_e}(\faq).
$$
Moreover, given $P \in \minparabs$ and $H \in \faq,$ there is a unique
element $H_0 \in \cl\faqp(P)$ which is conjugate to $H$ under $W.$
\end{lemma}

\proof
If $P \in \minparabs,$ then $\faqp(P)$ is the open positive chamber for
the positive system $\Sigma(P)$ of the root system $\gS.$ Also,
$W$ is the Weyl group of $\Sigma.$
Thus, the final assertion follows by a well known property of Weyl groups.
We turn to the assertions about equivalence of conjugation.

For the first two listed groups the equivalence follows from Lemma 1 in \bib{HC1}, p.\ 111.
For the equivalence for the third group,  let
$H_1, H_2 \in \faq$ and assume that $H_2 = \gf(H_1)$ for some $\gf \in \Int(\fgc).$
We may fix $P \in \minparabs$ such that $H_1 \in \cl(\faqp(P)).$ There exists
a $s \in W$ such that $s^{-1}(H_2) \in \cl(\faqp(P)).$
Fix a choice
$\Sigma^+(\fg, \fa_\fp)$ of positive roots for $\Sigma(\fg, \fa_\fp)$ that is compatible
with $\Sigma(P),$ and let $\fa_\fp^+$ be the associated positive chamber.
Then $\cl(\faqp(P)) \subset \cl(\fa_\fp^+).$
Since $W$ is naturally isomorphic to $\NKaq/Z_K(\faq),$
the elements $s^{-1}(H_2)$ and $H_1$ are conjugate under $\Int(\fgc).$ Hence, they
are already conjugate under $W(\fg, \fa_\fp).$ Being both contained in $\cl(\fa_\fp^+),$
the elements must be equal and we conclude that $H_2 = s(H_1).$ The equivalence
for the third group now
follows.

Using  Corollary \refer{c: W is W faq faq} (b),
we immediately obtain the equivalence for the fourth group
from the one for the third.
\qed

\begin{lemma}
Let $\fa$ be a linear subspace of $\faq$ and assume that
$\gf \in \Int(\fgc)$ maps $\fa$ into $\faq.$
Then there exists a $s \in W$ such that $s|_\fa = \gf|_\fa.$
\end{lemma}

\proof
The proof is identical to the proof of Cor.~2 in \bib{HC1}, p.~112, with use of
Lemma \refer{l: second lemma on Weyl groups} instead of \bib{HC1}, Lemma 1.
\qed

\begin{cor}
\naam{c: W a1 a2 in q}
Let $\fa_1, \fa_2$ be linear subspaces of $\faq,$ then
$$
W(\fa_2 \mid \fa_1) =
 \{t \in \Hom(\fa_1, \fa_2) \mid \;\exists\; s \in W:\;\;\; t = s|_{\fa_1}\;\}.
$$
\end{cor}

We briefly interrupt our discussion of Weyl groups to collect
some useful facts about conjugacy classes
of the parabolic subgroups from $\allparabs.$

\begin{lemma}
\naam{l: on parabolics}
Let $\fa_\fp$ be a maximal abelian subspace of $\fp,$ containing $\faq,$ and
let $Q \in \allparabs.$
\begin{enumerate}
\itema
There exists a $k \in N_{K_e}(\faq) \cap N_{K_e}(\fa_\fp)$
such that  $kQk^{-1}$ is standard.
\itemb
If $F, F' \subset \gD$ are such that $P_F$ and $P_{F'}$ are conjugate under
$G,$ then $F = F'.$
\itemc
There exists a unique subset $F \subset \gD$ such that $Q$ is conjugate to $P_F$ under
$G.$
\itemd
If $P \in \allparabs$ is conjugate to $Q$ under $G,$ then it is already conjugate
to $Q$
under
$N_{K_e}(\faq) \cap N_{K_e}(\fa_\fp).$
\end{enumerate}
\end{lemma}

\proof
There exists a $s \in W$ such that the parabolic equivalence class
$s(\faQqp)$ is contained in $\cl \faqp(P_0),$ hence equals $\faFqp,$ for
some $F \subset \gD.$ It follows that $sQs^{-1} = P_F,$ see Section
\refer{s: notation}. Now apply Corollary \refer{c: W is W faq faq} (b)
to obtain (a).

For (b), we note that $P_F$ and $P_{F'}$ both contain the minimal standard $\gs$-parabolic
subgroup $P_0.$ Hence, $P_F = P_{F'},$ by \bib{HC1}, p.~111, Lemma 2.
This implies that $F = F',$ see Section \refer{s: notation}.
Assertions (c) and (d) both follow from combining (a) and (b).
\qed

We end this section with a discussion of automorphisms
connecting $\Cartan$-stable Cartan subspaces of $\fq;$ see Section \refer{s: notation}
for basic notation.

If $\fb_1$ and $\fb_2$ are two $\Cartan$-stable abelian
subspaces of $\fq,$ then we define
\begin{equation}
\naam{e: defi W fb one two}
W(\fb_2\mid \fb_1):=
\{ \gf |_{\fb_{1\iC}} \;\mid
         \; \gf \in \Int(\fgc),\,\,\gf({}^\dual \fb_1) \subset {}^\dual\fb_2 \}.
\end{equation}
Note that ${}^\dual \fb_2$ and ${}^\dual \fb_1$
are abelian subspaces of $\fpdual.$ Using the notation
of the first part of this section, relative to the algebra
$\fgdual = \fkdual \oplus \fpdual$
with the indicated Cartan decomposition,
we see that complex linear extension induces a natural isomorphism
$$
W({}^\dual \fb_2\mid {}^\dual \fb_1) \simeq W(\fb_2\mid \fb_1).
$$
In particular, it follows from this that  the set in
(\refer{e: defi W fb one two}) is finite.
Moreover, if $\fb_1$ and $\fb_2$ are contained in $\fp \cap \fq,$
the notation (\refer{e: defi W fb one two})
is consistent with the notation introduced earlier
in this section.
\begin{lemma}
\naam{l: conjugacy in Cartan}
Let $\fb_1$ and $\fb_2$ two $\Cartan$-stable subspaces
of a fixed $\Cartan$-stable
Cartan subspace $\fd$ of $\fq.$ Then
$$
W(\fb_2\mid \fb_1) =
\{ s |_{\fb_{1\iC}} \;\mid \; s\in W(\fd),\;\;  s(\dfb_1) \subset \dfb_2 \}.
$$
\end{lemma}

\proof
This follows from Corollary \refer{c: W a1 a2 in q}
applied with $\fa_1 = {}^d\fb_1,$ $\fa_2 = {}^d\fb_2$
and $\fa_\fp = {}^d\fd.$
\qed

If $s \in W(\fb_2\mid \fb_1),$ then by $s^*$ we denote the map
$\fb_{2\iC}^* \to \fb_{1\iC}^*$ given by pull-back, i.e.,
$$
s^* \nu := \nu \after s, \qquad (\nu \in \fb_{2\iC}^* ).
$$
If $\fb_1, \fb_2\subset \fq$ are two $\Cartan$-stable Cartan subspaces,
then $\dfb_1$ and $\dfb_2$ are conjugate under an interior
automorphism of $\fg^d$ that commutes with $\Cartan^\dual;$
hence, the set $W(\fb_2\mid \fb_1)$ is non-empty and consists of isomorphisms.
If $s$ is any isomorphism from this set,
we denote its natural extension to the symmetric algebras
by $s$ as well. This extension
maps $I(\fb_1)$ into $I(\fb_2).$

\begin{lemma}
\naam{l: gg for diff css}
Let $\fb_1, \fb_2$ be $\Cartan$-stable Cartan subspaces of $\fq.$
Then $W(\fb_2\mid \fb_1)\neq \emptyset.$ Moreover, if
$s \in W(\fb_2\mid \fb_1),$ then
$$
s\after \gg_{\fb_1} = \gg_{\fb_2}.
$$
\end{lemma}

\proof
The first assertion follows from the discussion preceding
the lemma.

Let $K^d$ be the analytic subgroup of  $\Int(\fg_\iC)$ generated
by $e^{\ad \fk^d}.$ Then by Lemma
\refer{l: first lemma on Weyl groups}, applied to $\dfb_1,\dfb_2 \subset \fp^\dual,$
there exists an element $k \in K^d$
such that $s =k|_{{}^d\fb_1}.$
The action of $k$ induces the identity on $U(\fg^d)^{\fk^d} / U(\fk^d)^{\fk^d} \cap
U(\fg^d)\fk^d.$ Hence, if $D$ belongs to the latter algebra, then
$s \gg_{\fb_1}(D) = k \gg_{\fb_1}(k^{-1} \cdot D) = \gg_{\fb_2}(D).$
\qed

\section{Laurent functionals and operators}\eqnreset
\naam{s: Laurent functionals}
In this section we briefly recall the concept of Laurent functional,
introduced in \bib{BSanfam}; meanwhile, we fix notation that will be used in the rest
of the paper.
For details we refer to Sections 10 and 11 of \bib{BSanfam}.

Let $V$ be a finite dimensional real linear space, equipped with a positive
definite inner product $\inp{\dotvar}{\dotvar}.$ Its complexification
$V_\iC$ is equipped with the complex bilinear extension of the inner product.
We write $P(V)$ for the algebra of polynomial functions $V_\iC \to \C,$
and $S(V)$ for the symmetric algebra of $V_\iC.$ We identify the latter algebra
with the algebra of translation invariant differential operators on $V,$
which in turn is identified with the algebra of translation invariant
holomorphic differential operators on $V_\iC.$
In both settings, $u \in V$ is identified with the differential operator
$f \mapsto df(\dotvar)u.$

Let $X$ be a finite subset of non-zero elements of $V.$
By an $X$-hyperplane
in $V_\iC$ we mean an affine hyperplane of the form $H = a + \xi^\perp_\iC,$
with $a \in V_\iC$ and $\xi \in X.$ The hyperplane $H$ is said to be
real if $a$ may be chosen in $V.$
By
\begin{equation}
\naam{e: defi Pi X V}
\Pi_{X}(V)
\end{equation}
we denote the collection of polynomial functions $p \in P(V_\iC)$
with zero locus $p^{-1}(0)$ equal to a finite union of $X$-hyperplanes.
The subset consisting of $p$ with zero locus a finite union
of real $X$-hyperplanes is denoted by $\Pi_{X, \R}(V).$
Note that $\Pi_X(V)$ consists of all polynomial functions
that may be written as a non-zero multiple of a product of
factors of the form $\inp{\xi}{\dotvar} - c,$ with $\xi \in X$ and $c \in \C.$
The subset $\Pi_{X,\R}(V)$ consists of such products with $c \in \R$ in all
factors.

By an $X$-configuration in $V_\iC$ we
mean a locally finite collection of $X$-hyperplanes in $V_\iC.$ The configuration
is said to be real if all its hyperplanes are real. If $a \in V_\iC,$
then by $\Mer(V_\iC, a, X)$ we denote the space of germs at $a$ of meromorphic
functions with singular locus contained in the union of the hyperplanes
$a + \xi^\perp_\iC,$ for $\xi \in X.$ Let $\N^X$ denote the space of maps
$X \to \N.$ For $d \in \N^X$ we define the polynomial function
$\pi_{a,d} \in \Pi_X(V)$ by
$$
\pi_{a,d}(z) = \prod_{\xi \in X} \inp{\xi}{z-a}^d.
$$
Let $\cO_a(V_\iC)$ denote the space of germs of holomorphic functions at $a.$
Then $\Mer(V_\iC, a, X)$ is the union of the spaces $\pi_{a, d}^{-1} \cO_a(V_\iC)$
for $d \in \N^X.$
The space $\Mer(V_\iC, a, X)^*_\laur$ of $X$-Laurent functionals
at $a$ is defined as the subspace of  $\Mer(V_\iC, a, X)^*$ consisting of $\Lau$
with the property that for every $d \in \N^X$ there exists a $u_{\cL,d} \in S(V)$ such that
$$
\Lau \gf = [u_{\cL,d} \pi_{a, d} \gf](a), \;\; \text{for all}\;
\gf \in \pi_{a,d}^{-1} \cO_a(V_\iC).
$$
The element $u_{\cL}$ belongs to a projective limit space $S_{\leftarrow}(V,X)$
whose definition is suggested by the above,
see \bib{BSanfam}, Sect.~10, for more details.
Moreover, the map $\Lau \mapsto u_\Lau$ defines a linear isomorphism
\begin{equation}
\naam{e: iso Laurent functional to projlim}
\Mer(V_\iC, a, X)^*_\laur \; {\buildrel \simeq \over \longrightarrow} \; S_{\leftarrow}(V,X),
\end{equation}
see \bib{BSanfam}, Lemma 10.4.

The space on the left-hand side of the above isomorphism
only depends on $X$ through its proportionality class.
More precisely, a finite set  $X' \subset V\setminus \{0\}$ is said
to be proportional to $X$ if every element of one of the sets $X, X'$
is proportional to an element of the other set. If $X$ and $X'$ are proportional
sets, then $\Mer(V_\iC, a, X)^*_\laur = \Mer(V_\iC, a, X')^*_\laur,$
see \bib{BSanfam}, Lemma 10.3.

If $\Omega \subset V_\iC$ is open and $E$ a complete locally convex space,
then a (densely defined) $E$-valued function $f$ on $\Omega$ is
said to be meromorphic if for every $z_0 \in \Omega$
there exists an open neighborhood $\Omega_0$ of $z_0$ and
a holomorphic function $g\in \cO(\Omega_0)$ such that $gf|_{\Omega_0}$ is
a holomorphic $E$-valued function on $\Omega.$ The space of $E$-valued meromorphic
functions on $\Omega$ is denoted by $\Mer(\Omega, E).$
A point $z \in \Omega$ is said to be a regular
point of $f \in \Mer(\Omega, E)$ if $f$ is holomorphic in a neighborhood of $z.$
The collection of regular points of $f$ is denoted by $\reg f.$

Let $\Hyp$ be an $X$-configuration in $V_\iC.$
By $\Mer(V_\iC, \Hyp, E)$ we denote the space of meromorphic functions
$V_\iC \to E$ with singular locus contained in $\cup \Hyp.$
We agree to write $\Mer(V_\iC, \Hyp): = \Mer(V_\iC, \Hyp,\C).$

The space $\Mer(V_\iC, \Hyp, E)$
is topologized as follows.
Let $X^0 \subset X$
be minimal  subject to the condition that $X^0$ and $X$ are proportional.
For each $X$-hyperplane $H \subset
V_\iC$ there exists a unique $\ga_H \in X^0$ and a unique first order polynomial
function $l_H$ of the form $z \mapsto \inp{\ga_H}{z} - c,$ with $c \in \C,$
such that $H = l_H^{-1}(0).$

We
denote by $\N^\Hyp$ the space of maps $\Hyp \to \N.$
For $d \in \N^\Hyp$ and $\omega \subset V_\iC$ a bounded
subset, we define the polynomial function
$\pi_{\omega,d} \in \Pi_X(V)$
by
\begin{equation}
\naam{e: defi pi omega d}
\pi_{\omega, d} = \prod_{H \in \Hyp\atop H \cap \omega \neq \emptyset} l_H^{d(H)}.
\end{equation}
A change of choice of $X^0$ only causes a change of this polynomial
by a non-zero factor. If $E$ is a complete locally convex space, we define
$\Mer(V_\iC, \Hyp, d, E)$ to be the space of meromorphic
$E$-valued functions $\gf$ on $V_\iC$ with the property that
$\pi_{\omega, d} \gf$ is holomorphic on $\omega,$ for every bounded
open subset $\omega \subset V_\iC.$ This space is equipped with the
weakest locally convex topology that makes every map
$\gf \mapsto \pi_{\omega,d} \gf|_\omega$ continuous from $\Mer(V_\iC, \Hyp, d, E)$
to $\cO(\omega, E).$ This topology is complete; it is Fr\'echet if $E$
is. We equip $\N^\Hyp$ with the partial ordering $\preceq$ defined
by $d \preceq d' \iff \forall H \in \Hyp:\; d(H) \leq d'(H).$
We now have
$$
\Mer(V_\iC, \Hyp, E) = \cup_{d \in N^\Hyp} \;\Mer(V_\iC, \Hyp, d, E).
$$
Accordingly, we equip the space on the left-hand side with the
direct limit locally convex topology.

Any non-empty intersection of $X$-hyperplanes in $V_\iC$ is called an $X$-subspace of $V_\iC.$
An $X$-subspace $L$ of $V_\iC$ may be written as $L = a + V_{L\iC},$
with $V_L\subset V$ a real linear subspace.
Let $V^\perp_{L\iC}$ denote the complexification
of its orthocomplement. The intersection $V^\perp_{L\iC} \cap L$ consists of a single
point $c(L),$ called the central point of $L.$ Via the translation
$x \mapsto c(L) + x$ from $V_{L\iC}$ onto $L$ we equip $L$ with the structure of a
complex linear space, together with a real form with a positive definite
inner product on it.

If $\Hyp$ is an $X$-configuration and
$L \subset V$ an $X$-affine subspace, we define
$\Hyp_L$ to be the collection of affine hyperplanes in $L$
of the form $H \cap L,$ with $H \in \Hyp$ a hyperplane that properly intersects $L.$
Let $X(L): = X \cap V_L^\perp$ and let $X_r \subset V_L$ be the image  of $X\setminus X(L)$
under the orthogonal projection onto $V_L.$  The image of $X_r$ in $L$ under translation
by $c(L)$ is denoted by $X_L.$ Thus, $(L, X_L)$ is the analogue of $(V_\iC, X).$
The collection $\Hyp_L$ is a $X_L$-configuration in $L.$

If $\cL$ is a Laurent functional in $\Mer(V_{L\iC}^\perp, 0, X(L))_\laur^*,$
then $\cL$ induces a continuous linear map
$$
\cL_* : \Mer(V_\iC, \Hyp) \to \Mer(L, \Hyp_L),
$$
given by the formula
\begin{equation}
\naam{e: defi Laustar}
\cL_*\gf(\nu) := \cL [ \gf(\dotvar + \nu)],
\end{equation}
for $\gf \in \Mer(V_\iC, \Hyp)$
and generic $\nu \in L.$ The map $\cL_*$ belongs to
the space
$
\Laur(V_\iC, L, \Hyp)
$
of Laurent operators
$\Mer(V_\iC, \Hyp) \to \Mer(L, \Hyp_L),$ as defined in \bib{BSres}, Sect.~1.3,
see also \bib{BSanfam}, Sect.~11.
It follows from the definition of Laurent operator
combined with the isomorphism
$\Mer(V^\perp_{L\iC}, 0, X(L))^*_\laur \simeq S_{\leftarrow}(V^\perp_L, X(L))$
given by (\refer{e: iso Laurent functional to projlim})
that the map $\cL \mapsto \cL_*$
defines a linear surjection
\begin{equation}
\naam{e: surjection Laurent functionals to operators}
\Mer(V_L, X(L), 0)^*_\laur \to
\Laur(V_\iC, L, \Hyp) \to 0.
\end{equation}
Accordingly, a Laurent operator may alternatively be defined as any
continuous linear operator $\Mer(V_\iC, \Hyp) \to \Mer(L, \Hyp_L)$ of the form
$\Lau_*$ with $\Lau$ a Laurent functional from the space on the left-hand side of
(\refer{e: surjection Laurent functionals to operators}).

More generally, if $\Lau \in \Mer(V_L, X(L), 0)^*_\laur$ and if $E$ is a complete locally convex space
then the algebraic tensor product $\Lau_* \otimes I_E$ has a unique extension
to a continuous linear map $\Mer(V_\iC, \Hyp, E) \to \Mer(L, \Hyp_L, E)$
that we briefly denote by $\Lau_*$ again.

The concept of Laurent functional may be extended as follows, see
\bib{BSanfam}, Def.~10.8.
Let $\Mer(V_\iC, *,X)^*_\laur$ be the disjoint union of the spaces
$\Mer(V_\iC, a, X)^*_\laur,$ for $a \in V_\iC.$ An $X$-Laurent functional on $V_\iC$
is defined to be a finitely supported section of $\Mer(V_\iC, *,X)^*_\laur,$
i.e., a map $\Lau : V_\iC \to \Mer(V_\iC, *,X)^*_\laur$ with
$\Lau_a \in \Mer(V_\iC, a ,X)^*_\laur$ for every $a \in V_\iC$ and with
$\supp \Lau := \{a \in V_\iC\mid \Lau_a \neq 0\}$  a finite
set. The set of all $X$-Laurent functionals on $V_\iC$ naturally forms a complex linear
space, denoted $\Mer(V_\iC, X)_\laur^*.$

Let $\Hyp, L$ be as before. If $a \in V_\iC$ we denote by $\Hyp_L(a)$ the
$X_L$-configuration consisting of all hyperplanes $H'$ in $L$ for which there exists a $H \in\Hyp$
such that $H' = L \cap [(-a) + H].$
If $S \subset V_{L\iC}^\perp$ is a finite
subset, we put
\begin{equation}
\naam{e: defi Hyp L S}
\Hyp_L(S) = \cup_{a \in S} \Hyp_L(a).
\end{equation}
Let now $\Lau \in \Mer(V^\perp_\iC, X)^*_\laur$
and put $S = \supp \Lau.$ Then from the above discussion it follows in a straightforward
manner that the formula (\refer{e: defi Laustar}) defines a continuous linear map
$
\Lau_*: \Mer(V_\iC, \Hyp) \to \Mer(L, \Hyp_L(S)).
$
As above, if $E$ is a complete locally convex space, then the tensor
product map $\Lau_* \otimes I_E$ has a unique extension
to a continuous linear map
\begin{equation}
\naam{e: Laustar with E}
\Lau_*: \Mer(V_\iC, \Hyp, E) \to \Mer(L, \Hyp_L(S), E);
\end{equation}
see \bib{BSanfam}, Cor.~11.6.

By $\Mer(V_\iC, X, E)$ we denote the space of meromorphic $E$-valued functions on $V_\iC$
with singular locus contained in the union of an $X$-configuration. Every Laurent
functional $\Lau \in \Mer(V_{L\iC}^\perp, X(L))^*_\laur$ determines
a unique continuous linear map
$
\Lau_*: \Mer(V_\iC, X, E) \to \Mer(L, X_L, E)
$
such that $\Lau_*$ restricts to the map (\refer{e: Laustar with E}) for
every $X$-configuration $\Hyp$ in $V_\iC.$ See \bib{BSanfam}, Lemma 11.8, for details.

\section{Normalization of residues and measures}\eqnreset
\naam{s: normalization}
For the explicit determination of the constants in the Plancherel formula,
it is of importance to specify the precise normalizations of residual operators
and measures that will be used in the rest of this paper.

Let $P_0$ be the standard parabolic subgroup in $\minparabs$ and let
$t \in \WT(\Sigma)$ be a
$W$-invariant residue weight, see \S~\refer{s: notation}.
Let $\specinp$ be a $W$-invariant positive definite
inner product on $\faqd.$
Associated with the data $ \Sigma^+,  t, \specinp,$
we defined in \bib{BSres}, beginning of \S~3.4,
for each subset $F\subset \gD$
and every element $\gl \in \staFqd,$
a universal residue operator
\begin{equation}
\naam{e: universal residue operator}
\Res^t_{\gl + \faFqd}:= \Res^{P_0,t}_{\gl + \faFqd},
\end{equation}
which encodes
the procedure of taking a residue along the affine subspace
$\gl + \faFqdc$ of $\faqdc.$
In \bib{BSres}, text below Eqn.~(3.6),
this  residue operator is introduced as an element of a project limit space
$S_{\leftarrow}(\staFqd, \bar \gS_F^+),$ defined in \bib{BSres}, \S~1.3;
here $\bar \gS_F^+$ denotes the collection of indivisible roots
in $\gS^+ \cap \staFqd.$
However, to make the residue operator into an object as canonical as possible,
we shall prefer to view it as a Laurent functional.

Applying the results of Section \refer{s: Laurent functionals}
with $V = \staFqd,$ $X = \gS_F,$ $X^0 = \bar \gS_F^+$
and $a =0,$ we obtain an isomorphism
\begin{equation}
\naam{e: iso Laurent functionals with proj lim}
\Mer(\staFqdc, 0, \Sigma_F)^*_\laur \simeq S_{\leftarrow}(\staFqd, \bar \gS_F^+).
\end{equation}
In the present paper, the universal residue operator
(\refer{e: universal residue operator})
is accordingly viewed as an element of the space $\Mer(\staFqdc, 0, \Sigma_F)^*_\laur$
of $\gS_F$-Laurent functionals at the origin in $\staFqdc.$

By Eqn.~(\refer{e: surjection Laurent functionals to operators}), with $V = \faqd,$ $X = \gS$
and $L= \gl + \faFqd,$ we see that $\Lau \mapsto \Laustar$
induces a surjective linear map from $\Mer(\staFqdc, \gS_F)^*_\laur$
onto $\Laur(\faqdc, \gl + \faFqdc, \Hyp).$
In this context we  omit the star in the notation, and use
the notation (\refer{e: universal residue operator}) also for the Laurent operator
defined by the universal residue operator. Thus, (\refer{e: defi Laustar}) becomes
$$
(\Res_{\gl + \faFqd}^\ressup \gf)(\nu)
= \Res_{\gl + \faFqd}^\ressup [\gf(\dotvar + \nu)],
$$
for every $\gf \in \Mer(\faqdc, \Hyp)$ and $\nu \in \faFqdc$ generic.
In this way the notation becomes compatible with
the notation of \bib{BSres}.

All definitions in \bib{BSres} are given with reference to the fixed
$W$-invariant inner product $b$ on $\faqd,$ denoted $\inp{\dotvar}{\dotvar}$
in \bib{BSres}, so that a priori the
universal residue operator depends on the particular choice
of the inner product. However, as we will show, the dependence
is through certain measures determined by $b.$ To explain this,
it is convenient to first introduce some general terminology.

If $\fv$ is a real finite dimensional vector space, let $\dens{\fv}$ denote the one
dimensional real linear space of densities on $\fv,$ i.e., the space of maps
$\omega: \fv^n \to \R,$ where $n = \dim \fv,$
transforming according to the rule $\omega\after A^n  =
|\det A| \,\omega,$ for all $A \in \End(\fv).$ Evaluation at the origin induces
a natural isomorphism
from the space of translation invariant densities on $\fv,$ where $\fv$ is
viewed as a manifold,
onto $\dens{\fv};$ we shall identify accordingly. Consequently, via
integration the space $\dens{\fv}$ may be identified with the space of Radon measures
$\R d\gl,$
where $d\gl$ is a choice of Lebesgue measure on $\fv.$ If $\fv$ is equipped
with a positive definite inner product, then by the normalized density on $\fv$
we mean the unique element $\omega \in \dens{\fv}$ such that
$\omega(e_1, \ldots, e_n) = 1,$ for every orthonormal basis $e_1, \ldots, e_n$ of $\fv.$

We shall often encounter the situation
that $\fv = i \fb$ with $\fb$ a subspace of a real linear space $V;$ here multiplication
by $i$ is defined in the complexification $V_\iC$ of $V.$ If $V$ comes equipped with
a positive definite inner product $\inp{\dotvar}{\dotvar},$
we extend it to $V_\iC$ by complex bilinearity, and equip $i\fb$ with the positive
definite inner product $-\inp{\dotvar}{\dotvar}.$ Accordingly, in this setting
it makes sense to speak
of the normalized density on $i\fb.$

Let $d\mu \in \dens{i \fb}.$ If $\gl\in V_\iC,$
then we shall
adopt the convention to also denote
by $d\mu$ the density on the real affine subspace $\gl +  i \fb$ of $V_\iC,$
obtained by transportation under
the translation $X \mapsto \gl + X,$ $i \fb \to \gl + i \fb.$
Accordingly, by unoriented integration, the density $d \mu$ determines
a real Radon measure on $\gl + i\fb,$ which we shall denote by the same symbol.

We now return to the dependence of the residue operator
(\refer{e: universal residue operator}) on the choice of $b.$
For every $\ga \in \gS,$ the orthogonal reflection $s_\ga$
is independent of the particular
choice of $b,$ and therefore so are the
root hyperplanes
$\ga^\perp,$ and, more generally, the root spaces $\faFqd$ and their orthocomplements
$\staFqd,$ for $F \subset \gD.$
Combining this observation with the uniqueness statement of \bib{BSres},
Theorem 1.13,
it follows that the residue operator $\Res^\ressup_{\gl + \faFqd}$ can be completely
defined in terms of the data $\gS^+, t, \gl +\faFqd$ and
$b;$ moreover, it depends on the latter datum through
the quotient measure on $i \staFqd$ of the normalized
Lebesgue measures $d\gl$ on ${i\faqd}$ and $d\mu_F$ on ${i\faFqd}.$

To keep track of constants coming from comparing residue operators related
to different choices
of the mentioned inner product we shall introduce  a version of the residue operator
(\refer{e: universal residue operator}) that is independent
of the choice of $b.$

The unnormalized residue operator is defined
as the unique Laurent functional
\begin{equation}
\naam{e: true universal residue operator}
\uRes^\ressup_{\gl + \faFqd} \in
\Mer(\staFqdc, 0, \Sigma_F)^*_\laur \otimes \Hom(\dens{i \faqd}, \dens{i \faFqd})
\end{equation}
satisfying the following requirement. Let (\refer{e: universal residue operator}) be
defined relative to the given choice of $b$
and let  $d\gl \in \dens{i \faqd}$ and
$d\mu_F \in \dens{i \faFqd}$ be the normalized
densities associated with $b.$ Then the requirement is that
\begin{equation}
\naam{e: formula universal residue operator}
\uRes^\ressup_{\gl + \faFqd}(\gf)(d\gl)
= \Res^\ressup_{\gl + \faFqd} (\gf)\, d\mu_F,
\end{equation}
for $\gf \in \Mer(\staFqdc, 0, \Sigma_F).$
{}From the above mentioned dependence of (\refer{e: universal residue operator}) on
$b$ through the quotient density on $i\staFqd$ of $d\gl$ and $d\mu_F,$ it follows that
the residue operator in
(\refer{e: true universal residue operator})
only depends on
the data $ \Sigma^+, t, F, \gl$ and not on the choice of $b.$
In other words, if (\refer{e: universal residue operator}), $d\gl$ and $d\mu_F$ had been defined
relative to an arbitrary $W$-invariant inner product on $\faq,$ then formula
(\refer{e: formula universal residue operator}) would be valid as well.

Suppose now that for each $F \subset \gD$
a non-trivial density $d\mu_F\in \dens{i \faFqd}$ is given. In particular,
$d\gl:= d\mu_\emptyset$ is given. Then we can use the formula
(\refer{e: formula universal residue operator}) to define residual operators
$\Res^\ressup_{\gl + \faFqd}.$
With this definition, the integral formula of \bib{BSres}, Thm.~3.16, is valid.

In the rest of this  paper, we fix a choice $dx$ of invariant measure on $\spX.$
In the rest of this section we will describe how this choice determines all
other choices of normalization of measures, and, by the preceding discussion,
all choices of normalization of residual operators.

As in \bib{BSmc}, \S~3, the measure $dx$ determines a choice $da$ of Haar measure
on $\Aq,$ and a choice $d\gl$ of Lebesgue measure on $i\faqd.$
A change of the measure $dx$ by
a multiplication by a positive factor $c$ causes a change
of $da$ by the same factor. This in turn causes a change of the measure
$d\gl$  by the  factor $c^{-1}.$
It follows that the product measure $dx \, d\gl $ on $\spX \times i\faqd$
is independent of the particular choice of the measure $dx.$

In order to be able to use the formula (5.2) of \bib{BSfi},
we normalize the Lebesgue measures $d\mu_F$ of $i \faFqd,$ for $F \subset \gD,$ as
in the text following the mentioned formula. We describe this normalization in a somewhat more
general setting, in terms of the above terminology.
Let $B$ be the bilinear form of $\fg,$ fixed in the beginning of
Section \refer{s: notation}.
Via restriction and dualization, $B$ induces a positive definite inner product
on $\faqd,$ which we
denote by $B$ as well.
Let $c > 0$ be the positive constant such
that $d\gl$ corresponds to the density on $i\faqd,$ normalized
relative to the inner product
$c B.$

If $P \in \allparabs,$ then $d\mu_P$ denotes
the Lebesgue measure on $i\faPqd$ normalized with respect to $c B.$
If
$\fa_{G\iq}^* = \{0\},$ which occurs if and only if $G$ has a compact center modulo $H,$
we agree that $d \mu_G$ has total volume $1,$ in accordance with
\bib{BSfi}, text following (5.2).

The residual operators $\Res_{\gl + \faFqd}^\ressup$ are now normalized by
(\refer{e: formula universal residue operator}) and with respect to the
choices of normalizations of measures made.
All results of \bib{BSres} and \bib{BSfi} needed in this paper are valid
with the normalization of measures and residual operators just described.

If $P \in \allparabs,$ we denote by $d\gl_P$ the choice of Lebesgue measure
on $i \staPqd$ for which
\begin{equation}
\naam{e: product measure for P}
d \gl = d\gl_P\, d\mu_P.
\end{equation}
If $v \in \PcW,$ then by the above discussion of the normalization of $d\gl,$ applied to
the space $\spXPv,$ a choice $dx_{P,v}$ of invariant measure on
$\spXPv$ corresponds in one-to-one
fashion with a choice
of Lebesgue measure $d\gl_{P,v}$ on $i \staPqd.$ Throughout this paper we agree to select
 $dx_{P,v}$
so that $d\gl_{P,v} = d\gl_P.$

We end this section with the observation that for $P,Q \in \allparabs$ with $\faPq$ and
$\faQq$ conjugate under $W,$ the measures $d\mu_P$ and $d\mu_Q$ are $W$-conjugate. Indeed,
this follows from the $W$-invariance of the
inner product $B.$ {}From (\refer{e: product measure for P})
we see that the measures $d\gl_P$ and $d\gl_Q$ are $W$-conjugate as well.

\section{A vanishing theorem}\eqnreset
\naam{s: vanishing theorem}
Let $Q \in \allparabs.$ Throughout
this paper, we assume $(\tau, \Vtau)$ to be a finite dimensional
unitary representation of $K.$
In this section we introduce a space $\cE^\hyp_Q(\spX\col \tau)$
of meromorphic
families of  $\DX$-finite $\tau$-spherical functions
and show that the vanishing theorem of \bib{BSanfam} applies to it.

Let $\faqreg$ denote the set of regular elements in $\faq$ for
the root system $\gS$ and put $\Aqreg: = \exp \faqreg.$
We define a subset of $\spX$ by
$$
\spXp: = K \Aqreg H.
$$
According to \bib{BSanfam}, Sect.~2, this set is open dense in $\spX.$
Let $\cW \subset \NKaq$ be a choice of representatives
for $W/\WKH.$ Then, for each $P \in \minparabs,$
$$
\spXp = \cup_{v \in \cW}\;\; K A_{P\iq}^+  v H\qquad\text{(disjoint union).}
$$
By $\Ci(\spXp \col \tau)$ we denote the space of smooth
functions $f: \spXp \to \Vtau$ that are $\tau$-spherical,
i.e.,
$$
f(kx) = \tau(k)f(x) \qquad (x \in \spXp,\, k\in K).
$$
By $\cA(\spXp\col \tau)$ we denote the subspace of $f \in \Ci(\spXp\col \tau)$
that behave finitely under the action of the algebra $\DX.$
Moreover, we denote the subspaces of these spaces consisting
of functions that extend smoothly from $\spXp$ to all of $\spX$ by
$\Ci(\spX\col \tau)$ and $\cA(\spX\col \tau),$ respectively.

Let $P \in \allparabs$ and $v \in \NKaq.$ We put $K_P: = K \cap M_P$
and define $\spXPvp$
to be the analogue of the set $\spXp$ for the triple $(M_P, K_P, M_P \cap vH v^{-1}).$
In particular, $\spXPvp$ is an open dense subset of $\spXPv.$

We define the function $R_{P,v}: M_{1P} \to ]\,0,\infty\,[$  as in
\bib{BSanfam}, Section 3. According to \bib{BSanfam}, Lemma 3.2,
this function is left $K_P$- and right $M_{1P}\cap v H v^{-1}$-invariant. Moreover,
if $P \neq G,$ $a \in \Aq$ and $u\in \NKPaq,$ then
$$
R_{P,v}(au) = \max_{\ga \in \gS(P)} a^{-\ga}.
$$
Finally, $R_{P,v} \geq 1$ on $\spXPv.$ The function $R_{P,v}$
is of importance for the description of a domain of convergence
for the series expansion that describes the asymptotic behavior
of a function from $\cA(\spXp\col \tau)$ along $(P,v).$ To be more precise,
we define, for $0 <  r < 1,$
$$
\APqp(r) : = \{a \in \APq \mid \forall \ga \in \DrP: \;\; a^{-\ga} < r \}.
$$
Then the following property, see \bib{BSanfam}, Lemma 3.3,
is relevant for the mentioned
description of the domain of convergence.
For $m \in \spXPv$ and $a \in \APq,$
$$
m \in \spXPvp, \;\;a \in \APqp(R_{P,v}(m)^{-1}) \;\;\implies \;\;mavH \in \spXp.
$$
We can now describe the mentioned series expansion along $(P,v),$
together
with a domain of convergence.
According to \bib{BSanfam}, Lemma 5.3 and Thm.~3.4,
a function $f \in \cA(\spXp \col \tau)$ admits a
converging series expansion of the form
\begin{equation}
\naam{e: series for f in q}
f(mav) = \sum_{\xi \in E} a^{\xi} q_\xi (P,v\asmid f, \log a, m),
\end{equation}
for $m \in \spXPvp$ and $a \in \APqp(R_{P,v}(m)^{-1}).$
The set $E$ in (\refer{e: series for f in q}) is a subset of
$\faPqdc$ contained in a set of the
form $E_0 - \N\DrP:= E_0 + (-\N \DrP),$ with $E_0\subset \faPqdc$ finite.
In addition, there exists a $k \in \N$ such that, for every $\xi \in E,$
the expression $q_\xi(P,v\asmid f)$
belongs to $P_k(\faPq, \Ci(\spXPvp\col \tau_P)),$
the space of polynomial functions $\faPq \to \Ci(\spXPvp\col \tau_P)$
of degree at
most $k.$ Here $\tau_P$ stands for $\tau|_{K_P}.$

The series on the right-hand side of
(\refer{e: series for f in q}) converges neatly in the sense of \bib{BSanfam},
Def.~1.2,
for each $m \in \spXPvp,$ and for $a$ in the indicated range (depending on $m$).
The functions $q_\xi$ are uniquely determined by these properties.

The set of $\xi \in E$ for which $q_\xi(P,v\asmid f) \neq 0$ is called
the set of exponents of $f$ along $(P,v),$ and denoted by $\Exp(P,v\asmid f).$
We agree to write $q_\xi = 0$ for $\xi \in \faPqdc\setminus E.$

Using the above terminology we shall introduce the space
$\cE^\hyp_Q(\spX \col \tau)$ in a number of steps. First,
following \bib{BSanfam}, Def.~12.1, we introduce a suitable space
of meromorphic families of spherical functions. We agree to write
$P_0$ for the standard minimal $\gs$-parabolic subgroup. An index
or superscript $P_0$ will be replaced by $0.$ In particular,
$\spX_{0,v} = \spX_{P_0, v}$ and $\tau_0 = \tau_{P_0}.$
Note that $\spX_{0,v} = \spX_{P, v}$ and $\tau_0 = \tau_P$ for every
$P \in \minparabs.$

\begin{defi}
\naam{d: Cephyp}
Let $Q \in \allparabs$ and let $Y \subset \staQqdc$ be a finite subset.
We define
\begin{equation}
\naam{e: CephypQY}
\CephypQY
\end{equation}
to be the space of functions
$f:\faQqdc \times \spXp \to \Vtau,$ meromorphic
in the first variable, for which there exist a constant $k \in \N,$
a $\gS_r(Q)$-hyperplane configuration $\Hyp$ in $\faQqdc$ and a function
$d: \cH \to \N$ such that
the following conditions are fulfilled.
\begin{enumerate}
\itema
The function $\gl \mapsto f_\gl$ belongs to
$\Mer(\faQqdc, \Hyp, d, \Ci(\spXp\col \tau)).$
\itemb
For every $P \in \minparabs$ and $v \in \NKaq$ there exist
(necessarily unique)
functions
$
q_{s,\xi}(P,v\asmid f)$ in $ P_k(\faq) \otimes
\Mer(\faQqdc, \Hyp, d, \Ci(\spXzerov\col \tau_0)),$
for $s \in W/W_Q$ and $\xi \in -sW_Q Y + \N \Delta(P),$ with the
following property.
For all $\gl \in \faQqdc\setminus \cup \Hyp,$  $m \in \spXzerov$ and $a \in \Aqp(P),$
\begin{equation}
\naam{e: expansion Cephyp}
f_\gl(mav) = \sum_{s \in W/W_Q} a^{s\gl - \rho_P} \sum_{\xi  \in -sW_Q Y + \N \gD(P)}
a^{-\xi}\, q_{s,\xi}(P,v\asmid f, \log a)( \gl, m),
\end{equation}
where the $\DP$-exponential polynomial series with coefficients in $\Vtau$
converges neatly on $\Aqp(P).$
\itemc
For every $P \in \minparabs,$ $v \in \NKaq$ and $s\in W/W_Q,$ the series
$$
\sum_{\xi \in -sW_Q Y + \N\DP} a^{-\xi} q_{s, \xi}(P,v\asmid f,\log a)
$$
converges neatly on $\Aqp(P),$ as an exponential polynomial series
with coefficients in the space
$\Mer(\faQqdc, \Hyp, d, \Ci(\spXzerov \col \tauM)).$
\end{enumerate}
Finally, we define
\begin{equation}
\naam{e: Cephyp}
\Cephyp: = C^{\ep,\hyp}_{P_0, \{0\}}(\spXp\col \tau).
\end{equation}
\end{defi}

\begin{rem}
\naam{r: C ep hyp}
If $Q' \in \allparabs$ and $\fa_{Q'\iq} = \faQq,$ then
$\gSr(Q) \subset \gSr(Q') \cup [-\gSr(Q')].$ Hence, the notions
of $\gSr(Q)$- and $\gSr(Q')$-configuration coincide. It follows that the space
(\refer{e: CephypQY})
 depends on $Q$ through its
$\gs$-split component $\faQq.$

It is sufficient to require conditions (b) and (c) either for all $P \in \minparabs$
and a fixed $v,$ or for a fixed $P \in \minparabs$ and all $v$ in a choice
of representatives for $W/W_{K\cap H}$ in $\NKaq;$ see \bib{BSanfam}, Rem.~7.2, for
details.
\end{rem}

If $f \in \CephypQY,$ then following \bib{BSanfam}, Def.~12.4, we
define the asymptotic degree of $f,$ denoted  $\dega(f),$ to be
the smallest integer $k$ for which there
exist $\Hyp, d$ such that the conditions of Definition \refer{d: Cephyp}
are fulfilled.
Moreover, we denote by $\Hyp_f$ the smallest $\gS_r(Q)$-configuration
in $\faQqdc$
such that the conditions of
Definition \refer{d: Cephyp} are fulfilled with $k = \dega(f)$ and for
some $d: \Hyp_f \to \N.$
We denote by $d_f$ the $\preceq$-minimal
map $\Hyp_f \to \N$ for which the conditions of the definition are fulfilled
with $\Hyp = \Hyp_f$ and $k = \dega f.$
Finally, we put $\rega(f) := \faQqdc\setminus \cup\Hyp_f.$

We extend $\faq$ to a Cartan subspace $\fb$ of $\fq;$ clearly, $\fb$ is $\Cartan$-stable.
If $\mu \in \fbdc,$ then $I_\mu := \ker \gg(\dotvar \col \mu)$ is
an ideal of codimension one in $\DX.$

If $P \in \allparabs,$ we
write $\DPmaps$ for the set of finitely supported maps $\gd: \stfbPdc \to \N,$ see
(\refer{e: deco fb with star part}).
For $\gd \in \DPmaps$
and $\gl \in \faPqdc$ we define the ideal $I_{\gd, \gl}$ in $\DX$ by
$$
I_{\gd,\gl}: = \prod_{\gL \in \supp \gd} I_{\gL + \gl}^{\gd(\gL)}.
$$
This ideal is cofinite, since it is
a product of finitely  generated cofinite ideals.

Following \bib{BSanfam}, Def.~12.8, we introduce the following space
of meromorphic families of
$\DX$-finite functions.

\begin{defi}
\naam{d: defi cEhyp Q Y gd}
Let $Q \in \allparabs$ and $\gd \in \DQmaps.$ Then for $Y \subset \staQqdc$ a finite
subset we define
$$
\cEhypQYgd
$$
to be the space of functions $f \in \CephypQY$ such that,
for all $\gl \in \faQqdc\setminus \Hyp_f,$
the function $f_\gl: x \mapsto f(\gl, x)$ is annihilated by the cofinite ideal
$\Igdgl.$

Finally, we define
$$\cE^\hyp_0(\spXp\col \tau \col \gd):=
\cE^\hyp_{P_0,\{0\}}(\spXp\col \tau\col\gd).$$
\end{defi}

Following \bib{BSanfam}, Def.~12.8,  we introduce
the following subspace of meromorphic families of $\DX$-finite functions in
$\Ci(\spXp\col \tau)$
satisfying a certain additional assumption. Let $\parone$ be the
collection of parabolic subgroups $P \in \allparabs$ whose $\gs$-split component
$\faPq$ has codimension one in $\faq.$

\begin{defi}
\naam{d: cEhypQgdglob}
Let $Q\in \allparabs.$  Then for $\gd \in \DQmaps$ and
$Y\subset \staQqdc$ a finite subset
we define $\cEhypQYgdglob$ to be the space of
families $f \in \cEhypQYgd$ satisfying the following
condition.
\begin{enumerate}
\item[{\ }{\ }]
For every $s \in W,$ every $P \in \parone$ with $s\faQq\not\subset \faPq$
and all $v \in \NKaq,$ there exists an open dense subset $\Omega$ of $\rega f$
with the following property. For every $\gl \in \Omega,$
every $\xi \in s\gl|_{\faPq} + W_P s W_Q Y|_{\faPq} - \rho_{P}-  \N\gDr(P)$
and all $X \in \faPq,$ the function $m \mapsto q_\xi(P, v\mid f_{\gl}, X, m),$
originally defined on $\spXPvp,$ extends smoothly to all of $\spXPv.$
\end{enumerate}
\end{defi}

\begin{rem}
In
\bib{BSanfam}, Def.~9.5 and Def.~8.4,
a family $f$ satisfying the property stated in the above display was said
to be $[s]$-global along $(P, v).$
\end{rem}

\begin{defi}
\naam{d: cEhypQallX}
Let $Q \in \allparabs.$ For $\gd \in \DQmaps$ and
$Y \subset \staQqd$ a finite subset, we
define $\cE^\hyp_{Q,Y}(\spX\col\tau\col\gd)$
to be the space of families $f \in \cEhypQYgdglob$ satisfying the following condition.
\begin{enumerate}
\item[]
For all $\nu$ in the complement of a locally finite union
of analytic null-sets, the function  $f_\nu,$ initially defined
on $\spXp,$ extends to a smooth function on all of $\spX.$
\end{enumerate}
Finally, we define
$$
\cE^\hyp_{Q,Y}(\spX\col \tau):= \cup_{\gd}\;
\cE^\hyp_{Q,Y}(\spX\col\tau\col\gd)
\text{and}
\cE^\hyp_Q(\spX\col \tau):= \cup_{Y}\;
\cE^\hyp_{Q,Y}(\spX\col\tau),
$$
where $\gd$ runs over $\DQmaps$ and $Y$ over the collection of finite subsets
of $\staQqd.$
\end{defi}

\begin{rem}
\naam{r: dependence cEhypQY on faQq}
One readily verifies that the space $\cE^\hyp_{Q,Y}(\spX\col \tau\col \gd)$
depends on $Q$ through its split component $\faQq.$ See also Remark
\refer{r: C ep hyp}.
\end{rem}

\begin{rem}
\naam{r: Cephyp for Q is G}
If $G$ has compact center modulo $H,$ then $\fa_{G\iq} = \{0\}.$
Therefore, the map $f \mapsto f_0$ defines an embedding
of $\cE^\hyp_G(\spX \col \gd)$ into $\cA(\spX\col \tau),$ the space
of $\DX$-finite functions in $\Ci(\spX\col \tau).$ We claim that
this map is in fact a bijection. To see this, let $g \in \cA(\spX\col \tau).$
Then by \bib{Bas}, Thm.\ 5.3, the family
$f: \fa_{G\iq\iC}^* \times \spX \to \Vtau$
defined by $f_0 = g$ belongs to $C^{\ep, \hyp}_{G,Y}(\spXp \col \tau)$
for some finite subset $Y \subset {}^*\fa_{G\iq\iC}^* = \faqdc.$
Moreover, $\DX g$ is a finite dimensional space. In the notation introduced
before Definition \refer{d: defi cEhyp Q Y gd}, let $S$ be the set
of $\gL \in {}^*\fb_{G\iC}^* = \fb_\iC^*$ such that $\DX g$ has a non-trivial subspace
on which $\DX$ acts by the character $\gg(\dotvar\col \gL).$
Then there exists a map $\gd: {}^*\fb_{P\iC}^* \to \N$ supported by $S,$
such that $g$ is annihilated by $I_{\gd, 0}.$ It is now readily seen that
the family $f$ belongs to $\cE^\hyp_{G,Y}(\spX \col \tau\col \gd).$
\end{rem}

\begin{lemma}
\naam{l: Weyl on cEhypQY}
Let $f \in \cE^\hyp_Q(\spX \col \tau)$ and put $\Hyp = \Hyp_f$ and $d =d_f.$
Then $\nu \mapsto f_\nu$ is a meromorphic $\Ci(\spX\col \tau)$-valued function
in the space $\Mer(\faQqdc, \Hyp, d, \Ci(\spX\col \tau)).$
\end{lemma}

\proof
This follows by using condition (a) of Definition \refer{d: Cephyp}
and applying Cor.~18.2 of \bib{BSanfam}.
\qed

\begin{rem}
Let $\gd \in \DQmaps$ and let $Y \subset \staQqdc$ be a finite subset.
It can be shown that every family $f \in \cEhypQYgd$ that satisfies the
displayed condition of  Definition \refer{d: cEhypQallX} automatically belongs
to $ \cEhypQYgdglob,$ hence to $\cEhypQallX.$
In case $\max \gd \leq 1,$ this follows from
\bib{Bps2}, Thm.~12.8.
For general $\gd \in \DQmaps$ one may proceed along similar lines, see also
\bib{Wal1}, Chapter 4. However, we shall not need such a result in the present paper,
since in all cases where we could apply it, the property of
Definition \refer{d: cEhypQgdglob} has already been
established in \bib{BSanfam} for the functions involved.
The present remark justifies the notation used.
\end{rem}

The following special case of the vanishing
theorem of \bib{BSanfam} will play an important role in the rest of this paper.

\begin{thm}
\naam{t: vanishing thm new}
{\rm (Vanishing theorem)\ }
Let $Q \in \allparabs$ be standard and let
$\QcW\subset \NKaq$ be a choice
of representatives for $W_Q\bs W/ \WKH.$
Let $f \in \cEhypQallX$ and assume that there exists a non-empty open
subset $\Omega \subset \rega f$ such that, for each $v \in \QcW,$
$$
q_{\nu - \rho_P}(P,v\asmid f_\nu) = 0, \qquad (\nu \in \Omega).
$$
Then $f = 0.$
\end{thm}

\proof
This is a special case of
\bib{BSanfam}, Thm.~12.10.
\qed

We shall often use the vanishing theorem in combination with the
following lemma to relate families of eigenfunctions.

\begin{lemma}
\naam{l: Weyl group on cEhyp}
Let $P,Q \in \allparabs,$  let $s \in W$ be such that $s\faPq = \faQq$
and assume that $f \in \cE^\hyp_{P,Y}(\spX \col \tau),$
with $Y \subset \staPqdc$ a finite subset.
Then the family
$f^s : \faQqdc \times \spX \to \Vtau,$
defined by  $f^s(\nu, x) = f(s^{-1}\nu, x),$
belongs to $\cE^\hyp_{Q, sY}(\spX\col \tau).$
\end{lemma}

\proof
In view of Remark \refer{r: dependence cEhypQY on faQq}
we may assume that $Q = sP s^{-1}.$
There exists a $\gd \in \DPmaps$ such that
$f \in \cE^\hyp_{P,Y}(\spX \col \tau \col \gd).$
By Lemma \refer{l: conjugacy in Cartan}
there exists a $\tilde s \in W(\fb)$
such that $s  =\tilde s|_{\faq}.$ The element $\tilde s$
maps $\stfb_{P\iC}$ onto $\stfb_{Q\iC}.$ Let $\tilde s^*$ denote
its transpose $\stfbQdc \to \stfbPdc.$ Then $\gd \after \tilde s^* \in \DQmaps.$
Moreover, one readily checks from the definitions
that $f^s \in \cE^\hyp_{Q,sY}(\spX\col \tau\col \gd \after \tilde s^*).$
\qed

\section{Meromorphy of asymptotic expansions}\eqnreset
Let $f$ belong to the space $\cE^\hyp_Q(\spX\col \tau),$ defined in
Definition \refer{d: cEhypQallX}. The mentioned definition refers back to Definition
\refer{d: Cephyp}, according to which,
for $v \in \NKaq$ and for $P$ a minimal group in $\allparabs,$ the function
$f_\nu$ admits an expansion along  $(P,v)$ that depends
meromorphically on the parameter $\nu \in \faQqdc$ in a well defined
sense. It follows from \bib{BSanfam} that an analogous
result holds for arbitrary $P \in \allparabs.$ For its formulation,
we need a particular type of subset of the symmetric space $\spXPv.$
For $1\leq R \leq \infty,$ we define the set
$$
\spXPvp[R] := \{x \in \spXPvp\mid R_{P,v}(x) < R \};
$$
see \bib{BSanfam}, Eqn.\ (3.7) for details.

We also need an equivalence relation $\simPQ$ on $W$  to describe
asymptotic exponents along $(P,v)$ without redundance.
The relation is defined by
\begin{equation}
\naam{e: defi sim P Q}
s \simPQ t \iff \forall  \nu \in \faQqd:\;\; s \nu|_{\faPq} = t \nu|_{\faPq}.
\end{equation}
If $Y\subset \staQqdc$ is a finite subset and  $\gs \in W/\!\!\simPQ,$ we put
$$
\gs \cdot  Y: = \{ s \eta|_{\faPq} \mid s \in \gs, \eta \in Y\};
$$
see \bib{BSanfam}, \S~6, for details. We recall from \bib{BSanfam}, Lemma 6.5,
that $W/\!\!\sim_{P|Q} \simeq W/W_Q,$ if $P\in \minparabs$ and $Q\in \allparabs.$

\begin{prop}
\naam{p: meromorphy of asymp along walls}
Let $Q \in \allparabs,$ $\gd \in \DQmaps$ and $Y \subset \staQqd$
a finite subset.
Let $f$ be a family in $\CephypQY$ and put $k = \dega f.$

Let $P \in \allparabs$ and $v \in \NKaq.$
Then $\Exp(P,v\mid f_\nu) \subset W(\nu + Y)|_{\faPq} - \rho_P - \N \gDr(P),$
for every $\nu \in \rega f.$ Moreover,
there exist unique
functions
$$
q_{\gs, \xi}(P,v\mid f) \in P_k(\faPq) \otimes
\Mer(\faQqdc, \Hyp_f , d_f, \Ci(\spXPvp\col \tau_P)),
$$
for $\gs \in W/\!\!\simPQ$ and $\xi \in - \gs \cdot Y + \N\DrP,$
such that, for every $\nu \in \rega f,$
$$
f_\nu(mav) = \sum_{\gs \in W/\simPQ} a^{\gs \nu - \rho_P}
\sum_{\xi \in - \gs \cdot Y  + \N\DrP}
a^{-\xi} q_\xi(P,v\mid f, \log a)(\nu, m),
$$
for all $m \in \spXPvp$ and $a \in \APqp(R_{P,v}(m)^{-1}),$
where the exponential polynomial series in the variable $a$
with coefficients in $\Vtau$
is neatly convergent in the indicated range. In particular, for all $\nu$ in
an open dense subset of $\faQqdc$ and all
$\gs \in W/\!\!\simPQ$ and $\xi \in - \gs \cdot Y + \N\DrP,$
\begin{equation}
\naam{e: q gs xi equal to q}
q_{\gs, \xi}(P,v\mid f)(X,\nu) =
q_{\gs\nu|_{\faPq} - \rho_P -\xi}(P,v\mid f_\nu, X), \qquad (X \in \faPq).
\end{equation}
Finally, for each $\gs \in W/\!\!\simPQ$ and every $R > 1,$ the series
$$
\sum_{\xi \in -\gs\cdot Y + \N \gDr(P)} a^{-\xi} q_{\gs, \xi}(P,v\mid f, \log a)
$$
converges neatly on $\APqp(R^{-1})$ as a $\gDr(P)$-exponential
polynomial series with coefficients in $\Mer(\faQqdc, \Ci(\spXPvp[R]\col \tau_P).$
\end{prop}

\proof
This follows from \bib{BSanfam}, Thm.~7.7 and Lemma 12.7.
\qed

The following result
is based on the meromorphic nature of the series in the above proposition.
It may be considered a natural companion to \bib{BSanfam}, Lemma 7.9.

\begin{lemma}
\naam{l: restriction of exponents by holomorphy new}
Let $Q \in \allparabs,$ $Y \subset \staQqdc$ a finite subset
and $f \in \CephypQY.$
Let $P \in \allparabs, v\in \NKaq.$
Assume that for every
$\gs \in W/\!\!\simPQ$ a set $E_\gs \subset - \gs\cdot Y + \N \DrP$ is given
such that, for $\nu $ in a non-empty open subset $\Omega$ of $\rega f,$
\begin{equation}
\naam{e: inclusion exponents with E gs new}
\Exp(P, v\mid f_\nu) \subset \bigcup_{\gs\in W/\simPQ} \gs\nu|_{\faPq} - \rho_P - E_\gs.
\end{equation}
Then (\refer{e: inclusion exponents with E gs new}) holds for every $\nu \in \reg f.$
\end{lemma}

\proof
Let $\gs_0 \in W/\!\!\simPQ$ and $\xi \in -\gs_0 \cdot Y + \N \DrP$ be such that
$q_{\gs_0, \xi}(P, v\mid f) \neq 0.$ By
Proposition \refer{p: meromorphy of asymp along walls} there exists an open dense subset
$\Omega \subset \rega f$ such that, for all $\nu \in \Omega,$
(\refer{e: q gs xi equal to q}) is valid. In particular, it follows that,
for $\nu \in \Omega,$
$\gs_0\nu |_{\faPq} - \rho_P - \xi \in \Exp(P,v\mid f_\nu);$ hence,
$\gs_0\nu |_{\faPq} - \rho_P - \xi$ belongs to the union on the right-hand
side of (\refer{e: inclusion exponents with E gs new}).
By \bib{BSanfam},
Lemma 6.2, the sets
$\gs\nu|_{\faPq} + \gs\cdot Y - \N \DrP,$ for $\gs \in W/\!\!\simPQ,$
are mutually disjoint, for $\nu$ in an open dense subset $\Omega'$ of $\Omega.$
It follows that $\xi \in E_{\gs_0}.$

{}From the above and Proposition \refer{p: meromorphy of asymp along walls}
we conclude that for all $\gs \in W/\!\!\simPQ$ and
$\xi \in (- \gs\cdot Y +\N \DrP)\setminus E_\gs,$
the meromorphic function $\nu \mapsto q_{\gs, \xi}(P,v\mid f, \nu)$ is zero.
Hence, for $\nu \in \rega f,$
\begin{equation}
f_\nu(ma v)= \sum_{\gs \in W/\simPQ} a^{\gs \nu - \rho_P}
\sum_{\xi \in E_\gs} a^{-\xi}\,q_{\gs, \xi}(P,v\mid f)(\log a, \nu , m),
\end{equation}
for $m \in \spX_{P,v,+}$ and $a \in \APqp(R_{P,v}(m)^{-1}).$ Thus, the inclusion
(\refer{e: inclusion exponents with E gs new}) holds for $\nu \in \rega f.$ It remains
to extend the domain of its validity to all of $\reg f.$

Let $R \geq 1.$ Then for each $\gs \in W/\!\!\simPQ,$ the series
\begin{equation}
\naam{e: series F gs}
F_\gs(a) = \sum_{\xi \in E_\gs} a^{-\xi}\,q_{\gs, \xi}(P,v\mid f)(\log a)
\end{equation}
converges as a $\DrP$-exponential polynomial series in $a \in \APqp(R^{-1}), $
with coefficients in the space
$\Mer(\faQqdc, \Hyp, d, \Ci(\spX_{P,v,+}[R]\col \tau_P));$
here $\Hyp = \Hyp_f$ and $d = d_f.$

Let $\nu_0 \in \faQqdc$ be such that the meromorphic $\Ci(\spXp \col \tau)$-valued
function
$\nu \mapsto f_\nu$ is regular at $\nu_0.$
Let $\Hyp(\nu_0)$ be the collection of $H \in \Hyp$ that contain $\nu_0.$
Since $\Hyp$ is a $\gS_r(Q)$-configuration in $\faQqdc,$
the collection $\Hyp(\nu_0)$ is finite and there exists a bounded  open neighborhood
$\omega$ of $\nu_0$ in $\reg(f)$ such
that $\Hyp(\nu_0) = \{H \in \Hyp \mid H \cap \omega \neq \emptyset\}.$
Put $\pi = \pi_{\omega,d},$ see (\refer{e: defi pi omega d}).
Then $\gf \mapsto \pi \gf|_{\omega}$ defines a continuous linear map
from $\Mer(\faQqdc, \Hyp, d, \Ci(\spXPvp[P]\col \tau_P))$
into $\cO(\faQqdc, \Ci(\spXPvp[R]\col \tau_P)).$
In particular, the series (\refer{e: series F gs}) multiplied with $\pi$
converges neatly on $\APqp(R^{-1})$ as an exponential polynomial series
with values in the space $\cO(\omega, \Ci(\spXPvp[R]\col\tau_P)).$

It follows from \bib{BSanfam}, Lemma 10.5,
that there exists a $u \in S(\faQqd),$ such that
$\gf(\nu_0) = u(\pi\gf)(\nu_0),$ for $\gf \in \cO(\omega).$
We apply this to $\nu \mapsto f_\nu.$
Then,
\begin{eqnarray*}
f(\nu_0, ma v) & = & u[\pi(\dotvar) f(\dotvar, ma v)](\nu_0)\\
& = &
\sum_{\gs \in W/\simPQ} a^{\gs \nu_0 - \rho_P}
\sum_{j = 1}^k p_{\gs, j}(\log a) U_{\gs, j} [\pi(\dotvar) F_\gs(a)(\dotvar, m)](\nu_0),
\end{eqnarray*}
with finitely many $p_{\gs, j} \in P(\faPq)$ and $U_{\gs, j} \in S(\faQqd),$
as is readily seen by application of the Leibniz rule. Note that
$\deg p_{\gs, j} + \order(U_{\gs, j}) \leq l:= \order(u).$

{}From \bib{BSanfam}, Lemma 1.10, we obtain that the elements $U_{\gs, j} \in S(\faQqd),$
viewed as constant coefficient differential
operators in the variable $\nu,$ may be applied termwise to the series
for $\pi(\dotvar) F_\gs(a)(\dotvar, m),$
without disturbing the nature of the convergence of the series. This leads to
the existence of polynomial functions
$\bar q_{\gs, \xi} \in P_{k + l}(\faPq) \otimes \Ci(\spXPvp \col \tau_P),$
for $\gs \in W/\!\!\simPQ, \xi \in E_\gs,$ such that
$$
f_{\nu_0}(ma v) =
\sum_{\gs \in W/\simPQ} a^{\gs \nu_0 - \rho_P}
\sum_{\xi \in E_\gs}
a^{-\xi} \bar q_{\gs, \xi}(\log a, m),
$$
for $a \in \APqp(R^{-1})$ and $m \in \spXPvp[R].$
The series on the right-hand side converges neatly as a $\DrP$-exponential polynomial
series on $\APqp(R^{-1}),$ with coefficients in $\Ci(\spXPvp[R]\col \tau_P).$
It follows that the inclusion
(\refer{e: inclusion exponents with E gs new}) is valid for $\nu = \nu_0.$
\qed

\section{Fourier inversion}\eqnreset
\naam{s: Fourier inversion}
In this section we recall the Fourier inversion formula
from \bib{BSfi} that will be the starting
point for the derivation of the spherical Plancherel formula.

Let $\cW\subset \NKaq$ be a choice  of representatives for $W/\WKH.$
This choice determines the space $\oC = \oCtau,$ defined as the formal
Hilbert direct sum of finite dimensional Hilbert spaces
\begin{equation}
\naam{e: defi oC}
\oC := \bigoplus_{v \in \cW}\;\Ci(M/M\cap v H v^{-1} \col \tauM),
\end{equation}
where $\tauM$ denotes the restriction of $\tau$ to $\KM:= K \cap M.$
Given $P \in \minparabs,$ and $\psi \in \oC,$
we define the normalized Eisenstein integral $\nE(P\col \psi)$
as in \bib{BSft}, \S~5.  This Eisenstein integral is a meromorphic
$\Ci(\spX\col \tau)$-valued
function on $\faqdc,$ which depends linearly on $\psi.$
It naturally arises in representation theory,
essentially as a sum of matrix coefficients of the minimal principal series
of $\spX.$ However, it can also be characterized by some of
its analytic properties, as follows.

\begin{prop}
\naam{p: characterization min Eis new}
Let $P \in \minparabs$ and $\psi \in \oC.$
The function $\nu \mapsto \nE(P\col \psi\col \nu)$ is
the unique  family in $\cEhypPallX$ with the following property.
For each  $v \in \cW$ and  for $\nu$ in a dense open subset of $\faqdc,$
$$
q_{\nu - \rho_P}(P,v\asmid \nE(P\col \psi\col \nu),\dotvar, m) = \psi_v(m).
$$
Here $q_{\nu -\rho_P}$ is the coefficient in the expansion
(\refer{e: series for f in q}) for $f = \nE(P\col \psi\col \nu).$
\end{prop}

\begin{rem}
It is implicit in the above formulation that the function
on the left-hand side of the above equation is constant as a polynomial
function of the variable indicated by the dot.
It is known that for $\nu$ in a dense open subset of $\faqdc,$ there are no terms with $\log a$
in the expansion (\refer{e: series for f in q}) with $f = \nE(P\col \psi\col \nu);$
see \bib{BSmc}, Thm.~7.5.
\end{rem}

\proof
Uniqueness follows from the vanishing theorem, Thm.~\refer{t: vanishing thm new}.
Thus, it suffices to show that, for $\psi \in \oC,$
the $\Ci(\spX\col \tau)$-valued meromorphic function
$f: \nu \mapsto \nE(\psi\col \nu)$ on $\faqdc$ satisfies the properties
mentioned above. By linearity with respect to $\psi$ we may assume
that $\psi \in \oC[\gL],$ for some $\gL \in i\fbkd = i\stfbPd,$
in the notation of \bib{BSanfam}, text preceding Lemma 14.3.
It now follows from \bib{BSanfam}, Lemma 14.5, that
 $f \in \cE^\hyp_P(\spXp\col \tau\col \gd_\gL)_\glob,$
where $\gd_\gL$ is the characteristic function of the
subset $\{\gL\}$ of $\stfbPdc.$ Finally, it follows from
\bib{BSft}, beginning of Sect.~5, that $\nu \mapsto f_\nu$ is a
$\Ci(\spX\col \tau)$-valued meromorphic
function on $\faqdc.$ Therefore, $f \in \cEhypPallX.$

Combining \bib{BSft}, Eqns.~(49), (45) and the display following the latter equation,
we finally obtain that
$$
q_{\nu - \rho_P}(P,v\asmid f_\nu\col \dotvar\col m) =
[p_{P|P, 0}(1 \col \nu) C_{P|P}(1 \col \nu)^{-1} \psi]_v(m)
=
\psi_v(m),
$$
for each $v \in \cW$ and all $m \in M.$
\qed

\begin{lemma}
\naam{l: Hyp f real for minEis}
Let $P \in \minparabs$ and $\psi \in \oC.$ Let $f$ be the family in
$\cE^\hyp_P(\spX\col \tau)$ defined by $f_\nu = \nE(P\col \psi\col \nu),$
for $\nu \in \faqdc.$
Then the $\gS(P)$-configuration $\Hyp_f$,
defined as in the text  preceding
Definition \refer{d: defi cEhyp Q Y gd}, is real.
\end{lemma}

\proof
In view of \bib{BSanfam}, Eqn.~(14.12), it suffices to prove a similar
statement for the function $f$ of \bib{BSanfam}, Lemma 14.3.
We observe that in the proof of the mentioned lemma, $\Hyp_f$ is shown
to be contained in a $\gS(P)$-configuration $\Hyp''.$ This configuration
is given as $\Hyp'' = t^{-1} \Hyp \cup \Hyp'$ with $t \in W.$ Now $\Hyp$
is real by \bib{BSanfam}, Cor.~14.2. It therefore remains
to show that $\Hyp'$ is real.
For this it suffices to show that the singular locus
of the normalized $C$-function $\nC_{P|P}(t\col \dotvar)$
is the union of a real $\gS$-configuration. By \bib{BSanfam}, Eqn.~(14.6),
it suffices to show that the singular loci of $C_{P|P}(1\col \dotvar)^{-1}$
and $C_{P|P}(t\col \dotvar)$ are likewise. In view of \bib{Bps2}, Cor.~15.5,
it suffices to show that the functions $C_{Q|P}(1\col \dotvar)^{\pm 1},$ for
$Q \in \minparabs,$
all have a singular locus equal to the union of a real $\gS$-configuration.
The latter assertion follows by the argument
following the proof of Lemma 3.2 in \bib{BSfi}.
\qed

We write $\nE(P\col \dotvar)$ for the meromorphic
$\Ci(\spX, \Hom(\oC, \Vtau))$-valued function on $\faqdc$ given
by
$$
\nE(P\col \nu\col x)\psi =  \nE(P\col \psi \col \nu \col x),
$$
for $x \in \spX,$ $\psi \in \oC$ and generic $\nu \in \faqdc.$
Following
\bib{BSfi}, Eqn.~(2.3), we define the dualized Eisenstein integral as the
$\Ci(\spX, \Hom(\Vtau, \oC))$-valued meromorphic function on $\faqdc$ given
by
$$
\dE(P\col \nu \col x): = \nE(P\col - \bar \nu \col x)^*.
$$
Finally, we introduce the partial Eisenstein integrals $\Eps(P\col \dotvar),$
for $s \in W,$ as in \bib{BSfi}, Eqn.~(2.9), see also \bib{BSanfam}, Eqn.~(14.11).
Let $1 \otimes \tau$ denote the natural representation
of $K$ in $\Hom(\oC, \Vtau) \simeq \oC^* \otimes \Vtau .$ Then the partial
Eisenstein integrals are meromorphic $\Ci(\spXp\col 1 \otimes \tau)$-valued
functions on $\faqdc.$  Moreover, for $\nu\in \faqdc$ a regular point, the
partial Eisenstein integral $\Eps(P\col \nu)$ is $\DX$-finite.
By \bib{BSanfam}, Lemma 5.3,
it therefore has converging expansions like (\refer{e: series for f in q}), but with
$q_\xi(Q,v\asmid \Eps(P\col \nu))$ a $\Ci(\spXQvp\col 1 \otimes \tau)$-valued
polynomial function on $\faQq;$ see \bib{BSanfam}, Def.~2.1,
for details. The exponents
of the partial Eisenstein integrals are restricted by
$$
\Exp(P,v\asmid \Eps(P\col \nu))\subset s\gl - \rho_P -\N\DP,
$$
for all $v \in \cW.$ Finally, according to \bib{BSanfam}, Eqn.~(165),
$$
\nE(P\col \nu) = \sum_{s \in W} \Eps(P\col \nu) \text{on} \spXp.
$$
The mentioned properties determine the partial Eisenstein integrals
completely, see \bib{BSanfam}, Lemma 2.2.

We shall now investigate the dependence of the Eisenstein integrals on the
choice of $\cW.$ To this purpose, let $\bp \cW$ be a second choice
of representatives for $W/\WKH$ in $\NKaq.$
We denote by $\bp\oC$ the associated space define by (\refer{e: defi oC}),
with $\bp \cW$ in place of $\cW.$ The associated Eisenstein and partial Eisenstein
integrals are similarly indicated with a backprime.

\begin{lemma}
\naam{l: dependence Eis on cW}
There exists a unique linear map $R: \oC \to \bp \oC$ such that
\begin{equation}
\naam{e: Eis and bp Eis}
\nE(P\col \nu \col x) = \bp\nE(P\col \nu \col x) \after R,
\end{equation}
for all $x \in \spX$ and generic $\nu \in \faqdc.$
The map $R$ is an isometric isomorphism.
\end{lemma}

\proof
For every $w \in \cW,$ let $\bp w$ denote the unique element of $\bp \cW$ that
represents the same class in $W/\WKH.$ Then for every $w \in \cW$ we may select an element
$l_w \in \KM$ such that $\bp w \in  l_w w \NKHaq.$ The right regular action
$R_{l_w}: \Ci(M) \to \Ci(M)$ induces a linear isomorphism $R_w$ from
$\Ci(M/M \cap  wH w^{-1}\col \tauM)$
onto $\Ci(M/M \cap\bp wH\bp w^{-1}\col \tauM).$ Let $R: \oC \to \bp\oC$
be the direct sum of
the isomorphisms $R_w,$ for $w \in \cW.$ Then obviously $R$ is an isometry.
Let $\psi \in \oC.$ Then by Proposition \refer{p: characterization min Eis new}
the map $g: \nu \mapsto \bp \nE(P\col \nu \col x)\after R \psi$
belongs to $\cEhypPallX.$ Moreover, it follows from the same proposition that
\begin{equation}
\naam{e: expression q of g}
q_{\nu - \rho_P}(P, \bp w\asmid g_\nu, \dotvar, m) = (R\psi)_{\bp w}(m),
\end{equation}
for $\nu \in \rega(g),\, w \in \cW$ and $m \in M.$ Also, by the definition of $R,$
\begin{equation}
\naam{e: expression for R psi}
(R\psi)_{\bp w}(m) =  \psi_w(m l_w),
\end{equation}
for each $w \in \cW$ and all $m \in M.$

On the other hand, still by Proposition \refer{p: characterization min Eis new},
the map $f: \nu \mapsto \nE(P\col \nu\col x)\psi$
belongs to $\cEhypPallX$ as well, and for $\nu \in \rega(f),$ $w \in \cW$ and $m \in M,$
$$
q_{\nu - \rho_P}(P, w \asmid f_\nu, \dotvar, m) = \psi_w(m).
$$
This equation remains valid if we
replace $w$ in the expression on the left-hand side by any element $\bar w$ of $w \NKHaq,$
see \bib{BSanfam}, Lemma 3.7. Now $\bp w = l_w \bar w,$ for some $\bar w \in w\NKHaq;$
hence, applying \bib{BSanfam}, Lemma 3.6, with $v = \bar w$
and $u = l_w,$ we obtain that
\begin{equation}
\naam{e: expression  q of f}
q_{\nu - \rho_P}(P, \bp w \asmid f_\nu, \dotvar, m) =
\psi_w(m l_w).
\end{equation}
Comparing (\refer{e: expression q of f}) with (\refer{e: expression q of g}) and
(\refer{e: expression for R psi}) and applying the uniqueness
statement of Proposition \refer{p: characterization min Eis new}, we infer that
$f =g.$ On the other hand, if $R: \oC \to \bp \oC$ is a linear map such that
(\refer{e: Eis and bp Eis})
is valid, then $f =g,$ hence (\refer{e: expression q of g}) and (\refer{e: expression q of f})
are equal. This implies (\refer{e: expression for R psi}) and shows that
$R$ is uniquely determined by the requirement (\refer{e: Eis and bp Eis}).
\qed

\begin{lemma}
\naam{l: partial minimal kernel}
Let $P \in \minparabs$ and $s \in W.$ Then, for all $x \in \spXp, y \in \spX$ and
generic $\nu \in \faqdc,$ the element
\begin{equation}
\naam{e: partial minimal kernel}
\Eps(P\col \nu\col x) \dE(P\col \nu \col y) \in \End(\Vtau)
\end{equation}
does not depend on the choice of $\cW$ made in the text preceding
(\refer{e: defi oC}).
\end{lemma}

\proof
Let $R$ be the isometry of
Lemma \refer{l: dependence Eis on cW}.
{}From the fact that
the partial Eisenstein integrals are uniquely determined by the properties
mentioned in the text above Lemma \refer{l: dependence Eis on cW},
it follows that they satisfy the
transformation property (\refer{e: Eis and bp Eis})
with on both sides $\nE$ replaced by $\Eps,$ for
$s \in W.$ See also \bib{BSanfam}, Lemma 2.2.

On the other hand, taking adjoints of the homomorphisms on both sides of
(\refer{e: Eis and bp Eis}), and substituting $-\bar \nu$ for $\nu,$ we obtain that
$$
\dE(P\col \nu \col x) = R^* \after \bp\dE(P\col \nu \col x),
$$
for all $x \in \spX$ and generic $\nu \in \faqdc.$
{}From the unitarity of $R$ it now follows that the endomorphism
(\refer{e: partial minimal kernel}) does not change if we replace
$\Eps$ and $\dE$ by $\bp \Eps$ and $\bp\dE,$ respectively.
\qed

In the following we consider Eisenstein integrals associated
with the standard parabolic subgroup $P_0 \in \minparabs,$ but suppress the
symbol $P_0$ in the notation.
Moreover, we agree to write
$E_+$ for $E_{+,1}.$ We recall from \bib{BSfi}, p.\ 41, that
the meromorphic functions $\nu \mapsto \nE(\nu)$ and $\nu \mapsto E_+(\nu)$
have singular sets that are locally finite unions of real $\gS$-hyperplanes in
$\faqdc.$
Let $\Hyp$ be the collection of the singular hyperplanes for $\nE(\dotvar)$
and $E_{+}(\dotvar).$ Let $t \in \WT(\Sigma)$ be a $W$-invariant even residue weight,
see the text following (\refer{e: WT sum one}).
Associated with the data $\Sigma^+, t, $
we define, for each subset $F\subset \gD$ and every element $\gl \in \staFqdc,$
the residue operator $\Res^\ressup_{\gl + \faFqd}$ as in
\S~\refer{s: normalization}, with respect to the normalization
of the Lebesgue measures $d\gl$ of $\faqd$ and $d\mu_F$ of $\faFqd$
given at the end of that section.

The data $\Sigma^+, F$ and $\Hyp$ determine a finite subset
$\gL(F)$ of $ -\R_+ F,$ the negative of the closed cone spanned by $F,$
see \bib{BSfi}, Eqn.~(5.1), where between the set brackets the requirement `for some $t$'
should be added. We now recall, from \bib{BSfi}, Eqn.~(5.7),
the definition of the kernel $K_F^t(\nu\col x \col y) \in \End(\Vtau),$
for $(x,y) \in \spXp\times \spX$ and generic $\nu \in \faFqdc,$
by
\begin{equation}
\naam{e: kernel K t F}
K_F^t(\nu\col x\col y) =
\sum_{\gl \in \gL(F)}
\Res^\ressup_{\gl + \faFqd}
\left(
\sum_{s \in W^F}
E_+(s\dotvar\col x)\after \dE(s\dotvar\col y)\right) (\gl + \nu).
\end{equation}
{}From the definition it follows that $\nu \mapsto K_F^t(\nu \col \dotvar \col \dotvar)$
is a meromorphic
function on $\faFqdc$ with values in the space
$\Ci(\spXp \times \spX\col \tau\otimes \tau^*)$ and with
singularities along the hyperplanes of a real
$\gSr(F)$-configuration. Here $\tau\otimes \tau^*$ denotes the tensor product
representation of $K \times K$ in $\End(V_\tau) \simeq V_\tau \otimes V^*_\tau.$

The  residue operators in (\refer{e: kernel K t F})
depend on the choices of $d\gl$ and $d\mu_F,$
see the discussion in \S~\refer{s: normalization};
therefore, so does the kernel $K_F^t.$

\begin{lemma}
\naam{l: data determining K F t}
Let the data $(G,H,K, \tau, \faq, \Sigma^+)$ be fixed as above.
Let $F\subset \gD$ and let $t \in \WT(\Sigma)$ be a $W$-invariant even residue weight.
These data completely determine  $K_F^t \,d \mu_F \,dy,$ the product   of the kernel $K_F^t,$
defined by (\refer{e: kernel K t F}),
with the product measure $d\mu_F \, dy$ on $i\faFqd \times \spX.$

In particular, $K_F^t \,d\mu_F \,dy$ is independent of the particular choice
of $\cW,$ made in the text preceding (\refer{e: defi oC}), and of the choices of
$B$ and $dy,$ made in Sections \refer{s: notation} and \refer{s: normalization}.
\end{lemma}

\proof
Put
$$
k_F(\dotvar) = \sum_{s \in W^F}
E_+(s\dotvar \col x)\after \dE(s\dotvar\col y).
$$
It follows from
Lemma \refer{l: partial minimal kernel} that
$k_F$ depends only on the data mentioned, and not on $\cW, B, dy.$
Moreover, from  (\refer{e: kernel K t F}) and
(\refer{e: formula universal residue operator}),
it follows
that
$$
K_F^t(\nu \col x\col y) d\mu_F d y =
\sum_{\gl \in \gL(F)}
\uRes^\ressup_{\gl + \faFqd} (k_F(\cdot + \nu))(d\gl) dy.
$$
The occurring residue operators only depend on the data mentioned, and the product measure
$d\gl\, dy$ only depends on the choice of $(G,H, K, \faq),$
by the discussion in \S~\refer{s: normalization}.
All assertions now follow.
\qed

\begin{rem}
Since $d\mu_F$ and $dy$ do not depend on the choice of $\cW,$
it follows that the same holds for the kernel $K_F^t.$
This fact
has already silently been exploited
in \bib{BSfi}, text below Lemma 8.1, where the choice
of $\cW$ is adapted to the set $F.$
\end{rem}

According to \bib{BSfi}, Cor.\ 10.10,
the kernel $K_F^t(\nu\col x\col y)$ extends smoothly to all of $\spX$ in the variable
$x;$ more precisely, $K^t_F$ is a meromorphic $\Ci(\spX \times \spX, \End(\Vtau))$-valued
function on $\faFqdc$ with singularities along real $\gSr(F)$-hyperplanes.

{}From \bib{BSfi}, Eqn.\ (5.5) and (5.8),
we recall the definition of the linear operator
$T_F^t$ from $\Cci(\spX\col \tau)$ to $\Ci(\spXp\col \tau)$ by
\begin{equation}
\naam{e: defi T F t}
T_F^t f(x) = |W| \, t(\faFqp)\; \int_{\geps_F+ i\faFqd}\int_X K_F^t(\nu\col x\col y)
f(y)\; dy\; d\mu_F(\nu),
\end{equation}
for $f\in \Cci(\spX\col\tau)$ and $x \in \spXp.$
Here $d\mu_F$ is the translate by $\geps_F$ of the Lebesgue measure on
$i\faFqd$ normalized as in Section \refer{s: normalization}.
Moreover, $\geps_F$ is a point in the chamber $\faFqp,$
arbitrary but
sufficiently close to the origin.

\begin{rem}
\naam{r: on T gD}
If $G$ has compact center modulo $H,$ then $\fa_{\gD\iq} = \{0\}$ and
$t(\fa_{\gD\iq}^+) =1.$ Moreover, the above
is to be understood so that integration relative to $d\mu_\gD$ means
evaluation in $0;$ also, $\geps_\gD = 0.$
In this case we agree to write $K_\gD^t(x\col y) = K_\gD^t(0\col x\col y),$
so that the formula for $T_\gD^t$ becomes
\begin{equation}
\naam{e: formula for T gD}
T_\gD^t f(x) = |W| \; \int_X K_\gD^t(x\col y)\,f(y)\; dy
\end{equation}
for $f \in \Cci(\spX\col\tau)$ and $x \in \spXp.$
\end{rem}

In \bib{BSfi}, Cor.\ 10.11, it is shown that in fact $T_F^t$ maps into
$\Ci(\spX\col \tau)$ and defines a continuous
linear operator $\Cci(\spX\col \tau)\to \Ci(\spX\col \tau).$ Moreover,
by \bib{BSfi}, Thm.~1.2, it follows that
\begin{equation}
\naam{e: inversion formula with T t F}
I = \sum_{F \subset \gD} T^t_F \text{on} \Cci(\spX\col \tau).
\end{equation}

\begin{lemma}
\naam{l: T F t independent of dx}
Let $(G,H,K,\tau, \faq, \gS^+)$ and $(F,t)$ be data as
in Lemma \refer{l: data determining K F t}.
These data determine the operator $T_F^t,$ defined by (\refer{e: defi T F t}), completely.
\end{lemma}

\proof
This follows from Lemma \refer{l: data determining K F t}.
\qed

We finish this section with a discussion of how the kernels $K_F^t$ and
the operators $T_F^t$ behave under isomorphisms of reductive symmetric spaces.

Let $\gf: G \to \bp G$ be an isomorphism of reductive groups
of Harish-Chandra's class, and put $\bp H = \gf(H),$ $\bp K =\gf(K),$
$\bp \tau:= \tau \after \gf^{-1},$ $\bp\faq = \gf(\faq).$
Let $\bp \gS$
be the root system of $\bp\faq$ in $\bp \fg$ and let $\bp W$ denote the associated
Weyl group.
The isomorphism $\gf$ naturally induces the linear isomorphism $\faqdc \to \bp \faqdc$
given by
$$
\nu \mapsto \bp \nu := \nu \after \gf^{-1}|_{\bp \faq}.
$$
This isomorphism restricts to an isomorphism of root systems
$\gS \to  \bp\gS.$ Let $\gS^+$ be a positive system
for $\gS$ and $\bp\gS^+$
the corresponding positive system for $\bp\gS.$ Let $\gD$ and $\bp \gD$
be the collections of simple roots for $\gS^+$ and $\bp\gS^+,$
respectively. We denote by $\bp F$
the image in $\bp \gD$ of a subset $F \subset \gD.$

The map $\gf$ also naturally induces a bijection  from  $\WT(\gS)$ onto $\WT(\bp \gS),$
which we denote by $t \mapsto \bp t.$ If $t \in \WT(\gS)$ is $W$-invariant and even,
then $\bp t \in \WT(\bp \gS)$ is $\bp W$-invariant and even.

We put $\spX = G/H$ and $\bp \spX = \bp G /\bp H.$ Then the map
$\gf$ factors to a diffeomorphism
$\bar\gf:  \spX \to \bp\spX.$ This diffeomorphism induces the isomorphism
$\gf_* : \Ci(\spX\col \tau) \to \Ci(\bp\spX\col \bp\tau),$
given by $f \mapsto f \after \bar \gf^{-1}.$ It maps
$\Cci(\spX\col \tau)$ onto $\Cci(\bp\spX\col \bp\tau).$
We select invariant measures $dx$ and $\bp d x$ on $\spX$ and $\bp\spX,$
respectively. As in Section \refer{s: normalization},
this choice determines Lebesgue measures $d\gl$ and $\bp d \gl$
on $i \faqd$ and $i\bp\faqd,$ respectively. As in Section \refer{s: notation}
we fix bilinear forms
$B$ and $\bp B$ on $\fg$ and $\bp \fg,$ respectively. These choices determine
densities $d \mu_F$ and $d\mu_{\bp F}$ on $i \faFqd$ and $i \bp\fa_{\bp F\iq}^*,$
respectively.

Let $K^{\bp t}_{\bp F}$ be the analogue of the kernel
$K_F^t$ for the data
$(\bp G, \bp H, \bp K, \bp \tau, \bp \faq, \bp \gS^+, \\
\bp F,  \bp t, \bp d x, \bp B)$ in place of
$(G,H,K,\tau, \faq, \gS^+, F, t, dx, B).$ Moreover, let $T^{\bp t}_{\bp F}$ be the
associated analogue of the operator $T^t_F.$

\begin{lemma}
\naam{l: conjugation T gD t}
Let notation be as above. Then
\begin{equation}
\naam{e: K t F and iso}
K_{\bp F}^{\bp t} (\bp \nu \col \bar\gf(x) \col \bar\gf(y))
               \gf^*(d\mu_{\bp F})\,\bar\gf^* (\bp dy)
 = K_F^t(\nu \col x \col y)\, d\mu_F dy,
\end{equation}
for $x,y \in \spX$ and generic $\nu \in \faqdc.$
Moreover, the corresponding operators are related by
\begin{equation}
\naam{e: T t F and iso}
T^{\bp t}_{\bp F} \after  \gf_* =  \gf_* \after T^t_F \;\;\text{on}\;\;
\Cci(\spX\col \tau).
\end{equation}

\end{lemma}
\proof
In view of Lemma \refer{l: data determining K F t} it suffices to
prove the identity (\refer{e: K t F and iso}) in case $\bp \cW, \bp dx$ and $\bp B$ are
compatible with $\cW, dx$ and $B,$ via $\gf.$
It then follows from the definition
of the kernels that
$$
K^{\bp t}_{\bp F}(\bp \nu\col \bar\gf(x)\col \bar \gf (y)) = K^t_F(\nu \col x \col y),
$$
whence (\refer{e: K t F and iso}).

Equation (\refer{e: T t F and iso}) now
follows by combining (\refer{e: K t F and iso}) with (\refer{e: defi T F t})
and using the relations induced by $\gf$ between the data associated
with $G$ and $\bp G.$
\qed

\section{The generalized Eisenstein integral}\eqnreset
\naam{s: the generalized Eisenstein integral}
In this section we shall use the vanishing theorem
to give an alternative
characterization of the generalized Eisenstein
integral  defined in \bib{BSfi}, Def.\ 10.7.
This characterization, which is in the spirit of Proposition
\refer{p: characterization min Eis new},
will be used throughout the paper.

For the moment we assume the $G$ has compact center modulo $H.$
Then, with notation as in Remark \refer{r: on T gD},
we define  the space
\begin{equation}
\naam{e: defi cA t as span}
\cAt(\spX\col \tau) := \span \{ K^t_\gD( \dotvar \col y)v\mid
y \in \spXp,\, v \in \Vtau\}.
\end{equation}
This space equals the space $\cC_\nu$ of \bib{BSfi}, Eqn.\ (10.2),
with $F = \gD$ and $\nu = 0.$
It is  finite
dimensional and consists of $\DX$-finite functions in $C^\omega(\spX\col \tau),$
see \bib{BSfi}, Lemma 10.1. We adopt the new notation (\refer{e: defi cA t as span})
in stead of $\cC_\nu$
to avoid confusion with
the space
defined by (\refer{e: defi oC}).

\begin{lemma}
\naam{l: cA t is image T t gD}
Let $G$ have compact center modulo $H$,
and let $t \in \WT(\Sigma)$ be a $W$-invariant even residue
weight.
Then the space $\cAt(\spX \col \tau)$ equals the image $T^t_\gD( \Cci(\spX\col \tau)).$
\end{lemma}

\proof
{}From \bib{BSfi}, Lemma 10.2,
it follows that $K_\gD^t$ extends to a real analytic  function
$\spX \times \spX \to \End(\Vtau) \simeq \Vtau \otimes \Vtaud$ which
is $\tau \otimes \tau^*$-spherical.

By density of $\spXp$ in $\spX,$ continuity of $K^t_\gD$ and
finite dimensionality of the space $\cAt(\spX\col \tau),$ it follows that the latter
contains the function $K_\gD^t(\dotvar\col y)v$
for every $y \in \spX$ and $v \in \Vtau.$
This implies that
$T^t_\gD$ maps $\Cci(\spX\col \tau)$ into $\cAt(\spX\col \tau).$

To see that the converse inclusion holds, suppose that $\xi$ is a linear
functional of $\cAt(\spX\col \tau),$ vanishing on $\image(T^t_\gD).$
Then it suffices to show that $\xi = 0.$

For every $x\in \spX$ and $\eta \in \Vtaud,$
let $\xi_{x,\eta}$ denote the linear functional $\gf \mapsto \eta(\gf(x))$
on $\cAt(\spX\col \tau).$ The intersection of the kernels
of these linear functionals, as $x \in \spX$ and $\eta \in \Vtaud,$ is zero.
Therefore, these linear functionals span the dual of $\cAt(\spX\col \tau),$
and we see that there exist $n\geq 1,$
$x_1, \ldots, x_n\in \spX$ and $\eta_1, \ldots, \eta_n \in \Vtaud,$
such that $\xi = \sum_{j=1}^n \xi_{x_j, \eta_j}.$
In view of (\refer{e: formula for T gD}), the fact that $\xi$
vanishes  on $\image(T^t_\gD)$
implies that
$$
\sum_{j=1}^n \int_\spX \eta_j K_\gD^t(x_j\col y) f(y) \; dy = 0
$$
for all $f \in \Cci(\spX\col \tau).$ By sphericality of $K_\gD^t$ in the second
variable, the above integral also
vanishes for all functions $f\in \Cci(\spX, \Vtau).$ Hence,
$\sum_j \eta_j K_\gD^t(x_j\col \dotvar ) = 0$ as a function in $\Ci(\spX,\Vtaud).$
It follows that $\xi(K_\gD(\dotvar\col y)v) = 0,$ for all $y \in \spXp$ and $v \in \Vtau.$
In view of (\refer{e: defi cA t as span}), this implies that $\xi = 0.$
\qed

We now assume that $G$ is arbitrary again.
Let $F \subset \gD$ and let $\FcW\subset \NKaq$ be a choice of representatives
for $W_F\bs W /\WKH.$
If $t \in \WT(\gS)$ we denote by $\start$ the induced residue weight of $\gS_F,$
see \bib{BSres}, Eqn.~(3.16).
 Let $t$ be $W$-invariant and even; then $\start$ is $W_F$-invariant and even.

If $v\in \FcW,$ let $K_F^{\start}(\spXFv\col m \col m'),$
for $m,m' \in \spXFv,$ denote the analogue of $K_\gD^t$ for the symmetric
space $\spXFv.$ Note that $M_F$ has a compact center, so that
the discussion of the beginning of this section applies to $M_F$ in stead of $G.$
In particular, the data $(M_F, H_F, K_F, \tau_F, \staFq, \gS_F^+, \start)$
determine the finite dimensional space
$$
\cA^{\start}(\spXFv\col \tau_F) = \span\{K_F^{\start}(\spXFv\col \dotvar \col m')
\mid m' \in \spXFvp\}.
$$
Note that this space was denoted $\cC_{F,v}$ in \bib{BSfi}, Eqn.\ (10.7).
We define
\begin{equation}
\naam{e: defi AsubF}
\AsubF =\oplus_{v \in \FcW}\;\;\;\cA^{\start}(\spXFv\col \tau_F);
\end{equation}
this formal direct sum
was denoted $\cC_F$ in \bib{BSfi}, Def.\ 10.7.
The natural projections and embeddings associated with
the above direct sum are denoted by
$$
\prFv:\;\; \AsubF \to \cA^{\start}(\spXFv\col \tau_F) \text{and}
\rmiFv: \;\; \cA^{\start}(\spXFv\col \tau_F)\to \AsubF,
$$
for $v \in \FcW.$ Given $\psi \in \AsubF$ we shall also write $\psi_v:= \prFv \psi.$

The generalized Eisenstein integral $\genEF(\psi \col \nu),$
defined in \bib{BSfi}, Def.\ 10.7, is a function
in $\Ci(\spX\col\tau)$ that depends linearly on $\psi\in \AsubF$ and meromorphically
on $\nu \in \faFqdc.$ We shall not repeat the definition here,
but instead give a characterization based on the vanishing theorem,
Theorem \refer{t: vanishing thm new}.
The following result will allow us to
show that $\genEF(\psi \col \dotvar)$ belongs to
the space of families $\cE_F^\hyp(\spX\col \tau)$ introduced
in Definition \refer{d: cEhypQallX} with $Q = P_F.$ For its formulation,
we recall some notation from \bib{BSfi}, \S~8.

In the rest of this section we write $\nE(\gl\col x): = \nE(P_0\col \gl \col x).$
Similarly, if $v \in \FcW,$
we write $\nE(\spX_{1F,v}\col \nu \col m)$ for the
normalized Eisenstein integral of $\spX_{1F,v},$ associated with
the minimal parabolic subgroup $M_{1F} \cap P_0.$ The analogue of the space
$\oC$ for the latter Eisenstein integral is
\begin{equation}
\naam{e: defi oC F v}
\oC_{F,v}: = \oplus_{w \in \cW_{F,v}}\;\;\; \Ci(M/M\cap wv H v^{-1} w^{-1}\col \tauM).
\end{equation}
Here $\cW_{F,v} \subset N_{M_F \cap K}(\faq)$ is a choice
of representatives for $W_F/W_{K_F \cap vHv^{-1}};$ see \bib{BSfi}, Eqn.~(8.2).
Adapting the set $\cW$ if necessary, we may assume that $\cW_{F,v} \subset \cW.$
Then $\cW$ is the disjoint union
of the sets $\cW_{F,v},$ for $v \in \FcW,$ see \bib{BSfi}, Lemma 8.1.
Accordingly,  $\rmi_{F,v}$ denotes the natural inclusion $\oC_{F,v} \to \oC$, defined
as the identity on each component of (\refer{e: defi oC F v}). Moreover,
\begin{equation}
\naam{e: deco oC over FcW}
\oC = \oplus_{v \in \FcW}\;\;\; \rmiFv(\oC_{F,v}).
\end{equation}
We denote the associated projection operator $\oC \to \oC_{F,v}$ by $\pr_{F,v},$
for $v \in \FcW.$

\begin{lemma}
\naam{l: Lau on nE new}
Let $\cL$ be a Laurent functional in
$\Mer(\staFqdc, \Sigma_F)^*_\laur\otimes \oC.$
Then the family $g:  \faFqdc \to \Ci(\spX\col \tau),$  defined by
$$
g(\nu, x): = \Lau[\nE(\nu + \dotvar \col x)]
$$
belongs to $\cE_{F,Y}^\hyp(\spX\col \tau),$ with $Y = \supp \Lau.$
Moreover, if $v \in \FcW,$ then
for $\nu$ in a dense  open subset of $\faFqdc,$
\begin{equation}
\naam{e: q of family g is Lau Eis}
q_{\nu - \rho_F}(P_F, v \asmid g_\nu , X, m) = \Lau[ \nE(\spX_{1F,v}\col \dotvar + \nu  \col m)
\after \pr_{F,v}],
\end{equation}
for all $X \in \faFq$ and $m \in \spXFvp.$ Here $\pr_{F,v}$
 denotes
the projection associated with (\refer{e: deco oC over FcW}).
\end{lemma}

\proof
It suffices to prove the result for a Laurent
functional of the form $\Lau = \Lau' \otimes \psi,$
with $\Lau' \in \Mer(\staFqdc, \gS_F)^*_\laur$
and $\psi \in \oC.$ Define the family $f$ by
$f(\nu \col x): = \nE(\nu \col x)\psi.$

It follows from \bib{BSft}, p.~52, Lemma 14, that there exists a locally finite
collection $\Hyp$ of $\gS$-hyperplanes in $\faqdc$ and a map $d: \Hyp \to \N$
such that $f$ belongs
to $\Mer(\faqdc, \Hyp, d, \Ci(\spX\col \tau)).$
{}From \bib{BSanfam}, Lemma 13.1, applied with
$Q = P_F,$ it follows
that $g$ is a meromorphic
function on $\faFqdc$ with values in $\Ci(\spX\col \tau).$

It follows from \bib{BSanfam}, Lemma 14.5, that there exists a $\gd \in D_P$
such that $f$ belongs to $\cE^\hyp_0(\spXp \col \tau \col \gd)_{\hglob},$
see \bib{BSanfam}, Def.~13.10, for the definition of the latter space.
According to \bib{BSanfam}, Thm.~13.12, this implies that
$g \in \cE^\hyp_{F,Y}(\spXp\col \tau)_\glob.$
We conclude that $g \in \cE^\hyp_{F,Y}(\spX\col \tau).$

The family $f$ equals the family $f_W$ defined in \bib{BSanfam}, Prop.~15.4.
It follows from that proposition, applied with $Q = P_F$ and with $\Lau'$ in place
of $\Lau,$ that (\refer{e: q of family g is Lau Eis}) holds for each $v \in \FcW,$
generic $\nu \in \faFqdc$ and all $X \in \faFqd$ and $m \in \spXFvp.$
Combining this with \bib{BSanfam}, Theorem 7.7, Eqn.~(7.14), we see that
(\refer{e: q of family g is Lau Eis}) holds for all $\nu$ in a dense  open
subset of $\faFqdc,$ every $v \in \FcW,$ $X \in \faFqd$ and all
$m \in \spXFvp.$
\qed

\begin{thm}
\naam{t: char of genEis by asymp new}
Let $\psi \in \AsubF.$ Then $g: \nu \mapsto \genEF(\psi \col \nu)$
is the unique family in $\cE_F^\hyp(\spX\col \tau)$ with the following property.
For all $\nu$ in some non-empty open subset of $\faFqdc$ and each $u \in \FcW,$
\begin{equation}
\naam{e: q of genEis is psi}
q_{\nu - \rho_F}(P_F, u \asmid g_\nu)(X, m) = \psi_u(m),
\qquad (X \in \faFq,\, m \in \spX_{F,v,+}).
\end{equation}
\end{thm}

\begin{rem}
\naam{r: Eisenstein for P is G}
If $F = \gD$ and $G$ has compact center modulo $H,$ then
$\faFq = \{0\}$ and $f \mapsto f_0$ defines a bijection from
$\cE^\hyp_F(\spX\col \tau)$ onto $\cA(\spX\col \tau),$ the space of $\DX$-finite
functions in $\Ci(\spX\col \tau),$ see
Remark \refer{r: Cephyp for Q is G}.
Moreover, ${}^F\cW$ consists of one element
which one may take to be $1,$ $M_F/M_F\cap H \simeq \spX,$
and $\cA_F^{\start} \simeq \cA^t(\spX\col \tau).$ Finally, with notation as in the
above theorem, $g_0 = \psi,$ so that $\psi \mapsto E_F^\circ(\psi \col 0)$
is the inclusion map $\cA^t(\spX\col \tau) \to \cA(\spX\col \tau).$
\end{rem}

\begin{rem}
\naam{r: intro gL XFv F}
In the proof of Theorem \refer{t: char of genEis by asymp new}, we will encounter the set
\begin{equation}
\naam{e: intro gL XFv F}
\gL(\spXFv, F) \subset - \R_+ F,
\end{equation}
which is defined to be the analogue of the set $\gL(\gD)$ of (\refer{e: kernel K t F}),
for the data $(\spXFv, \staFq, \gS^+_F)$ in place of $(\spX, \faq, \gS^+).$
\end{rem}

\proof
Uniqueness follows from the vanishing theorem, Theorem \refer{t: vanishing thm new};
hence, it suffices to prove existence.

In view of (\refer{e: defi AsubF}),
we may assume that
$\psi \in \cAFv$
for some $v\in \FcW.$
According to \bib{BSfi}, Eqn.\ (10.9),
we may then express $\psi$ as
follows
\begin{equation}
\naam{e: defi psi with res}
\psi(m) = \sum_{\gl \in \gL(\spXFv, F)}
\Res_\gl^{\start}[\nE(\spXFv\col - \dotvar\col m)\Phi(\dotvar)],
\qquad(m \in \spXFv),
\end{equation}
where
$$
\Phi(\gl) =
\sum_{j =1}^k E^*_+(\spXFv\col -\gl\col m_j) v_j \in \oCFv, \qquad
(\gl \in \staFqdc),
$$
with $\{m_1, \ldots, m_k\}$ a finite subset of $\spXFvp,$ and
$\{v_1, \ldots, v_k\}$ a finite subset of $\Vtau.$
We now note that
\begin{equation}
\naam{e: defi cR F}
\cR_F := \sum_{\gl \in \gL(\spXFv, F)}
\Res_\gl^{\start}
\end{equation}
is a Laurent functional in $\Mer(\staFqdc, \Sigma_F)^*_\laur;$
moreover, according to
\bib{BSfi}, Def.\ 10.7,
the generalized Eisenstein integral is given by
$$
g(\nu,x) =
\cR_F [ \nE(\nu - \dotvar\col x)\after \iFv \Phi(\dotvar)].
$$
Define the Laurent functional
$\Lau_0 \in \Mer(\staFqdc, \Sigma_F)_\laur^* \otimes \oC$ by
\begin{equation}
\naam{e: defi Lau zero}
\Lau_0 \gf : = \cR_F[ \gf(- \dotvar) \after \iFv \Phi(\dotvar)],
\end{equation}
for $\gf \in \Mer(\staFqdc, \Sigma_F) \otimes \oC^*.$
Then the generalized Eisenstein integral is given by
$$
g(\nu, x) = \Lau_0[\nE(\nu + \dotvar \col x)].
$$
It now follows from Lemma \refer{l: Lau on nE new} that
$g \in \cE_F^\hyp(\spX\col \tau)$
 and that,
for $\nu$ in an open dense subset of $\faFqdc$ and all $X \in \faFq$ and $m \in \spXFvp,$
\begin{eqnarray}
q_{\nu - \rho_F}(P_F,u\asmid g_\nu, X, m) &=&
\Lau_0[ \nE(\spXFv  \col \nu + \dotvar \col m) \after \pr_{F,u}]\nonumber \\
&=&
\cR_F[ \nE(\spXFv \col \nu -\dotvar \col m)
\after \pr_{F,u} \after \iFv \Phi(\dotvar)].
\naam{e: cR F applied to Eis}
\end{eqnarray}
If $u \neq v,$ then the latter expression equals $0.$ Since also $\psi_u = 0,$
the identity (\refer{e: q of genEis is psi}) then follows.
On the other hand, if $u = v,$ then
$\pr_{F,u} \after \iFv \Phi(\dotvar) = \Phi(\dotvar);$
hence, (\refer{e: cR F applied to Eis}) equals the expression on the right-hand side of
(\refer{e: defi psi with res}), and since $\psi_v = \psi,$ the identity
(\refer{e: q of genEis is psi}) follows.
\qed

\begin{cor}
\naam{c: cor to thm 4.3}
There exists a locally finite collection $\Hyp_0$ of hyperplanes in $\faFqdc$ such
that the following holds. Let $\psi \in \cA^{\start}_F$ and let $g$ be defined
as in Theorem \refer{t: char of genEis by asymp new}.
Then the meromorphic
function $\nu \mapsto g_\nu$ is regular on the
complement of $\cup\Hyp_0.$ Moreover, for every
$u \in \FcW,$ $X \in \faFq$ and $m \in \spXFvp,$ formula (\refer{e: q of genEis is psi})
holds for
all $\nu \in \faFqdc\setminus \cup \Hyp_0.$
\end{cor}
\proof
Let $\bar 1$ denote the image of $1$ in $W/\!\!\sim_{P_F|P_F},$ see
(\refer{e: defi sim P Q}).
Let
\begin{equation}
\naam{e: intro q bar one zero}
\nu \mapsto q_{\bar 1, 0}(P_F, u \mid g)(\nu, X, m)
\end{equation}
be the function in $\Mer(\faFqdc, \gSr(F), \Vtau),$ defined
as in Proposition \refer{p: meromorphy of asymp along walls}.
By \bib{BSanfam}, Thm.~7.7, there
exists a locally finite collection $\Hyp_0$ of hyperplanes in $\faFqdc$ such that
$\nu \mapsto g_\nu$ is regular on $\faFqdc\setminus \cup \Hyp_0$ and for
every $u \in \FcW,$  all $X \in \faFq$ and $m \in \spXFup,$
\begin{equation}
\naam{e: equality q bar one zero and q nu}
q_{\bar 1, 0}(P_F, u \mid g)(\nu, X, m) = q_{\nu - \rho_F}(P_F,u\mid g_\nu , X, m),
\end{equation}
for all $\nu \in \faFqdc\setminus \cup \Hyp_0.$
By linearity in $\psi$ and finite dimensionality of the space $\cA^{\start}_F,$
the collection $\Hyp_0$ can be chosen independent of $\psi.$
Combination of  (\refer{e: equality q bar one zero and q nu}) and
(\refer{e: q of genEis is psi}) gives that the meromorphic
function
(\refer{e: intro q bar one zero}) is constant and equal to $\psi_u(m).$
In view of (\refer{e: equality q bar one zero and q nu})
it now follows that (\refer{e: q of genEis is psi})
holds for all $X \in \faFq,$ $m \in \spXFup$
and $\nu \in \faFqdc\setminus \cup\Hyp_0 .$
\qed

{}From the uniqueness statement in Theorem \refer{t: char of genEis by asymp new}
it follows that the generalized Eisenstein integral
$\genEF(\psi)\in \cE_F^\hyp(\spX\col \tau)$
depends linearly on $\psi \in \cAFt.$ We agree to write
$\genEF(\nu\col x)\psi: = \genEF(\psi\col \nu\col x),$
for $x \in \spX$ and generic $\nu \in \faFqdc.$ Accordingly, we view the generalized
Eisenstein integral as a meromorphic function on $\faFqdc$
with values in $\Ci(\spX \col \tau \otimes 1);$
here $\tau \otimes 1$ denotes the tensor product representation
in $\Hom(\cAFt, \Vtau) \simeq \Vtau \otimes (\cAFt)^*.$
In accordance with
\bib{BSanfam}, Def.~10.7, we put
\begin{equation}
\naam{e: genEFv}
\genEFv(\nu \col x): = \genEF(\nu\col x) \after \iFv \in \Hom(\cAFv, \Vtau),
\end{equation}
for $v \in \FcW,$ $x \in \spX$ and generic $\nu \in \faFqdc.$

\begin{lemma}
\naam{l: psi as Lau of nE}
Let $v\in \FcW$ and let  $\psi \in \cAFv.$
\begin{enumerate}
\itema
There exists a Laurent functional
$\Lau \in \Mer(\staFqdc, \Sigma_F)^*_\laur \otimes \oC_{F,v}$
such that
\begin{equation}
\naam{e: psi as Lau of nE}
\psi(m) = \Lau [ \nE(\spXFv \col \dotvar \col m)], \qquad (m \in \spXFv).
\end{equation}
\itemb
There exists a functional as in (a) with support contained in $\gL(\spXFv, F),$
the set introduced in Remark \refer{r: intro gL XFv F}. In particular, the support
of this functional is real.
\itemc
If $\Lau$ is any Laurent functional as in (a), then, for all $x \in \spX,$
$$
\genEFv(\nu: x) \psi = \Lau[ \nE(\nu + \dotvar : x)\after \iFv]
$$
as an identity of meromorphic functions in $\nu \in \faFqdc.$
\end{enumerate}
\end{lemma}

\proof
As in the proof of Theorem \refer{t: char of genEis by asymp new} we
may express $\psi$ by (\refer{e: defi psi with res}).
Let $\Lau_0$ be defined as in (\refer{e: defi Lau zero}) and let the Laurent functional
$\Lau \in \Mer(\staFqdc, \Sigma_F)^*_\laur \otimes \oC_{F,v}$ be defined
by $\Lau \gf = \Lau_0 (\gf(\dotvar) \after \pr_{F,v} ),$ for
$\gf \in \Mer(\staFqdc, \Sigma_F) \otimes \oC_{F,v}^*.$
Since $\pr_{F,v}\after \iFv = I$ on $\oCFv,$ it follows
from (\refer{e: defi Lau zero}), (\refer{e: defi cR F}) and
(\refer{e: defi psi with res}) that $\Lau$ satisfies (\refer{e: psi as Lau of nE}).
We observe that $\supp \Lau \subset \supp \cR_F \subset \gL(X_{F,v}, F);$
in particular, $\Lau$ has support contained in $\staFqd.$ This establishes
(a) and (b).

Now assume that $\Lau$ is a Laurent functional as in (a).
Let $\Lau' $ be the Laurent functional in $\Mer(\staFqdc, \gS_F)^*_\laur \otimes \oC$
defined by $\Lau' \gf = \Lau (\gf(\dotvar) \after \iFv).$
Then it follows from Lemma \refer{l: Lau on nE new}
that the family $g: \faFqdc \times \spX \to \Vtau$
defined by
$$
g(\nu, x) = \Lau'[\nE(\nu + \dotvar \col x)]
$$
belongs to $\cE_F^\hyp(\spX\col \tau)$
and satisfies, for $u \in \FcW,$  $\nu $ in an open
dense subset of $\faFqdc$
and all $X \in \faFq$ and $m \in \spX_{F,u},$
\begin{eqnarray*}
q_{\nu - \rho_F}(P_F, u \asmid g_\nu) &=&  \Lau'[\nE(\spXFv\col \dotvar) \after \pr_{F,u}]\\
&=&\Lau[\nE(\spXFv\col \dotvar )\after\pr_{F,u} \after \iFv ]\\
&=&\pr_{F,u}\after \iFv \psi.
\end{eqnarray*}
Here we note that the last equality is obvious for $u \not= v,$ since then $\pr_{F,u} \after \iFv = 0.$
On the other hand, if $u = v,$ then $\pr_{F,u}\after \iFv =I$ on $\cAFv$
and the equality follows from the assumption on $\psi.$
It now follows from Theorem \refer{t: char of genEis by asymp new} that
$
g(\nu, x) = \genEF(\nu \col x)\after \iFv \psi.
$
\qed
Combining the above result with  Lemma \refer{l: Hyp f real for minEis},
we obtain
the following information on the asymptotic coefficients of the generalized
Eisenstein integral. We put
\begin{equation}
\naam{e: defi rmY F}
\rmY(F):= \cup_{v \in \FcW}\;\; \gL(\spXFv, F).
\end{equation}
This is a finite subset of $-\R_+F,$ which in turn is contained in $\staFqd.$
\begin{lemma}
\naam{l: globality q of genEF newer}
Let $F \subset \gD$ and
$\psi \in \cA^{\start}_F.$
The family $f: (\nu, x) \mapsto \genEF(\nu\col x)\psi$ belongs
to $\cE^\hyp_{F,\rmY(F)}(\spX\col \tau).$ Moreover, the $\gSr(F)$-configuration
$\Hyp_f,$ defined in the text preceding Definition \refer{d: defi cEhyp Q Y gd}, is real.

Put  $k = \dega f$ and let
$Q \in \allparabs,$  $u \in \NKaq.$ Then, for every
$\gs \in W/\!\!\sim_{Q|P_F}$ and all $\xi \in - \gs\cdot\rmY(F) + \DrQ,$
\begin{equation}
\naam{e: q gs xi in Mer}
q_{\gs,\xi}(Q,u\mid f) \in P_k(\faQq) \otimes
\Mer(\faFqdc, \Hyp_f, d_f, \Ci(\spX_{Q,u} \col \tau_Q)).
\end{equation}
\end{lemma}

\proof
{}From Theorem \refer{t: char of genEis by asymp new}
it follows that $f \in \cE^\hyp_{F,Y}(\spX\col \tau),$
with $Y$ a finite subset of $\staFqdc.$

For the first two assertions we may assume that $f_\nu = \genEFv(\nu)\psi,$
with $\psi \in \cA^{\start}(\spXFv\col \tau_F).$
Select $\Lau \in \Mer(\staFqdc, \gS_F)^*_\laur \otimes \oCFv$
as in Lemma \refer{l: psi as Lau of nE} (b). According to
Lemma \refer{l: Hyp f real for minEis} there exists
a real $\gS$-configuration $\Hyp$ in $\faqdc$ such that for every
$\psi' \in \oCFv,$ the family $g: \gl \mapsto \nE(\gl)\rmi_{F,v}\psi',$ which belongs
to $\cE^\hyp_0(\spX\col \tau),$ satisfies $\Hyp_g \subset \Hyp.$
It now follows from Lemma \refer{l: psi as Lau of nE} (c), combined with
\bib{BSanfam}, Prop.~13.2, that $Y \subset Y(F)$ and that
$\Hyp_f \subset \Hyp_F(Y(F)),$ with
the latter set defined as in \bib{BSanfam}, Eqn.~(11.6).
It follows from the mentioned definition
and the fact that $\Hyp$ and $Y(F)$ are real, that
$\Hyp_F(Y(F))$ and hence $\Hyp_f$ are real. It remains
to establish (\refer{e: q gs xi in Mer}).

Let $Q,u, \gs$ be as asserted.
With a reasoning as above, it follows from
Lemma \refer{l: psi as Lau of nE}, combined with \bib{BSanfam}, Lemma 14.5
and Proposition 13.9, that $f$ is holomorphically $\gs$-global
along $(Q,u)$ (see \bib{BSanfam}, Definition 13.6).
Let $\xi \in -\gs \cdot Y + \N \DrQ.$ Then (\refer{e: q gs xi in Mer}) follows
by application of
\bib{BSanfam}, Proposition 13.8.
\qed

\section{Temperedness of the Eisenstein integral}\eqnreset
\naam{s: temperedness}
In this section we show that the generalized Eisenstein integral
$\genEF(\nu)\psi,$ defined in the previous section, is tempered
for regular values of $\nu$ in $i\faFqd.$

Let us first recall the notion of temperedness.
Following \bib{Bps2}, p.\ 415, we define the function $\Theta: \spX \to \R$ by
$$
\Theta(x) = \sqrt{\Xi(x\gs(x)^{-1})},
$$
where $\Xi$ is the elementary spherical function $\gf_0$ associated
with the Riemannian symmetric space $G/K.$

Moreover, we define the function $\lspX: \spX \to [0,\infty\,[,$
denoted $\tau$ in \bib{Bps2}, by
\begin{equation}
\naam{e: defi lspX}
\lspX(kah) = |\log a|,\qquad (k \in K, \, a \in \Aq,\, h \in H).
\end{equation}

\begin{defi}
\naam{d: temperedness}
A $\DX$-finite function $f$ in $\Ci(\spX\col \tau)$
is said to be tempered if there exists
a $d \in \N$ such that
\begin{equation}
\naam{e: tempered estimate}
\sup_{\spX} (1 + \lspX)^{-d} \Theta^{-1} \| f\| < \infty.
\end{equation}
The space of these functions is denoted by
$
\cAtemp(\spX\col \tau).
$
\end{defi}

The following lemma gives a criterion for temperedness in terms of
exponents. We assume that
$\repPmin\subset \minparabs$ is a choice of representatives for
$\minparabs/\WKH$ and that $\cW \subset \NKaq$ is a choice of representatives
for $W/\WKH.$ We also assume that $P_1$ is a fixed element of $\minparabs.$

\begin{lemma}
\naam{l: criterion for temperedness}
Let $f \in \Ci(\spX\col \tau)$ be a $\DX$-finite function.
Then the following conditions are equivalent.
\begin{enumerate}
\itema
$f \in \cAtemp(\spX\col \tau).$
\itemb
For each $ P \in \repPmin$  and
every $\xi \in \Exp(P, e \asmid f)$ the estimate
$\Re \xi + \rho_P \leq 0$ holds
on
$\faqp(P).$
\itemc
For each $v \in \cW$ and every $\xi \in \Exp(P_1,v\mid f)$
the estimate
$\Re \xi + \rho_{P_1} \leq 0$ holds on $\faqp(P_1).$
\end{enumerate}
\end{lemma}
\proof
By sphericality and the decomposition $G ={\rm cl}\, \cup_{P \in \repPmin} K\Aqp(P)H,$
see \bib{Bas}, Cor.~1.4 and top of p.~232,
the estimate (\refer{e: tempered estimate}) is equivalent to the requirement
that, for each $P \in \repPmin,$
$$
\sup_{a \in \Aqp(P)} (1 + |\log a|)^{-d} \Theta(a)^{-1}\|f(a)\| < \infty.
$$
By \bib{Bps2}, Prop.~17.2,
there exist constants $C > 0$ and $N \in \N$ such that, for each $P \in \repPmin,$
$$
a^{-\rho_P} \leq\Theta (a) \leq
C\,(1 + |\log a|)^N a^{-\rho_P}\qquad (a \in \Aqp(P)).
$$
Therefore, (\refer{e: tempered estimate}) is equivalent
to the existence of a constant $d' \in \N$ such that, for
each $P \in \repPmin,$
$$
\sup_{a \in \Aqp(P)} (1 + |\log a|)^{-d'} a^{\rho_P}\|f(a)\| < \infty.
$$
According to \bib{Bas}, Thm.~6.1, this condition is in turn equivalent
to (b). This establishes the equivalence of (a) and (b), for any
choice of $\repPmin.$ The equivalence of (b) and (c) follows
from the observation that $\{v^{-1} P_1 v  \mid v \in \cW\} $ is a choice
of representatives for $\minparabs/\WKH$ combined with the fact
that $\Exp(v^{-1}P_1 v, e \mid f) = v^{-1}\Exp(P_1, v\mid f),$ for $v \in \cW,$
by \bib{BSanfam}, Lemma 3.6.
\qed

If  $G$ has compact center modulo $H$ and  $t \in \WT(\gS)$ is
a $W$-invariant even residue weight, then according to
\bib{BSfi}, Lemma 10.3,
there exists, for every choice of Hilbert structure on the space
$\cAt(\spX\col \tau),$
a unique endomorphism $\ga = \ga^t$ of this space, such that
\begin{equation}
\naam{e: K t gD and ga}
K_\gD^t( x\col y) = \bfe(x)\after \ga \after \bfe(y)^*,
\end{equation}
for $x,y \in \spXp.$ Here the map
$\bfe(x): \cAt(\spX\col \tau) \to \Vtau$ is defined by
$\gf \mapsto \gf(x).$ The corresponding function $\bfe,$ with values in
$\Hom(\cAt(\spX\col \tau) ,\Vtau) \simeq  \Vtau \otimes \cAt(\spX\col \tau)^*,$
is a $\tau \otimes 1$-spherical
real analytic function on $\spX.$ We recall from \bib{BSfi}, Lemma 10.3,
that $\ga$ is self-adjoint and bijective.

In the following we assume that
$t \in \WT(\gS)$ is a $W$-invariant even
residue weight and that $F \subset \gD.$
We equip each finite dimensional space $\cA^{\start}(\spXFv\col \tau_F),$
for $v\in \FcW,$ with a positive definite inner product.
Moreover, we equip the direct sum space $\cA^{\start}_F,$ defined by
(\refer{e: defi AsubF}), with the direct sum inner product, denoted $\hinp{\dotvar}{\dotvar}.$
Here and in the following, we use a bar in the notation of an inner product
to indicate its sesquilinearity. Moreover, all such inner products
will be antilinear in the second variable.

Let $\ga_{F,v}^\start = \ga_{F,v}$ be the analogue of the
endomorphism $\ga$ for $(\spXFv, \tau_F),$
and let $\ga_F^\start= \ga_F \in \End(\cA^{\start})$
be the direct sum of the $\ga_{F,v},$ for $v \in \FcW.$
Then $\ga_F$ is self-adjoint and bijective.
Moreover, according to \bib{BSfi}, Prop.~10.9, see also Lemma~10.2,
we have,
for $x,y \in \spX,$
\begin{equation}
\naam{e: K t F as prod Eis}
K_F^t(\nu \col x \col y) = \genEF(\nu\col x)\after \ga_F \after \dEF(\nu \col y)
\end{equation}
as an identity of $\End(\Vtau)$-meromorphic functions
in the variable $\nu \in \faFqdc.$  Here $\dEF$ denotes the dual
generalized Eisenstein integral, defined by
\begin{equation}
\naam{e: defi dEF}
\dEF(\nu \col y):= \genEF(-\bar \nu \col y)^* \in \Hom(\Vtau, \cA_F^{\start}),
\end{equation}
for $y \in \spX$ and generic $\nu \in \faFqdc.$

\begin{lemma}
\naam{l: dEis spans cA t F}
There exists a locally finite collection $\Hyp_1$ of affine hyperplanes
in $\faFqdc,$ such that $\nu \mapsto E^*_F(\nu)$ is regular
on $\faFqdc\setminus\cup\Hyp_1,$ and such that the following holds.
For every $\nu \in \faFqdc\setminus\cup\Hyp_1,$
\begin{equation}
\naam{e: cA F start as span dEF}
{\rm span}\{E_F^*(\nu \col y)v\mid y \in \spXp, v \in \Vtau\} = \cA_F^{\start}.
\end{equation}
\end{lemma}

\proof
Let $\Hyp_0$ be the collection of hyperplanes of Corollary \refer{c: cor to thm 4.3}
and let $\Hyp_1$ be the image of $\Hyp_0$ under the map $\nu \mapsto -\bar\nu.$
In view of (\refer{e: defi dEF}) the function $\nu \mapsto E_F^*(\nu)$ is
regular on the complement of $\cup \Hyp_1.$

Let $\nu \in \faFqdc\setminus\cup \Hyp_1$ and let $\psi \in \cA^{\start}_F.$
Assume that
$$
\hinp{\psi}{E_F^*(\nu \col y)v} = 0 \text{for all} y \in \spXp, v \in \Vtau.
$$
Using (\refer{e: defi dEF}) we see that $\genEF(-\bar \nu) \psi = 0.$
It now follows from Corollary \refer{c: cor to thm 4.3} that
$\psi_u = 0$ for each $u \in \FcW.$ Hence, $\psi = 0$ and
(\refer{e: cA F start as span dEF}) follows.
\qed

In the following we write $\rho = \rho_{P_0},$ where
$P_0$ denotes the standard parabolic subgroup
in $\minparabs.$

\begin{lemma}
\naam{l: exponents of K F t}
Let $v \in \Vtau$ and $y \in \spXp.$ Then the family $f: (\nu,x) \mapsto
K_F^{t}(\nu \col x \col y) v$ belongs to $\vanfamF.$ Moreover,
for every $\nu \in  \reg f$ and each $u \in \cW,$
\begin{equation}
\naam{e: exponents K F t in temp set}
\Exp(P_0,u\mid f_\nu) \subset W^F(\nu + \gL(F)) - \rho -\N \gD,
\end{equation}
where $\gL(F)$ denotes the finite subset of $-\R_+F$
introduced in (\refer{e: kernel K t F}).
\end{lemma}

\proof
The first assertion follows from (\refer{e: K t F as prod Eis})
and Theorem \refer{t: char of genEis by asymp new}.

According to \bib{BSfi}, Prop.~3.1,
the function $\gl \mapsto \dE(P_0\col \gl\col y)v$ belongs to the
space $\Mer(\faqdc, \gS) \otimes \oCtau.$
Combining this with \bib{BSanfam}, Lemma 14.3, we deduce that the family
$h: \faqdc \times \spXp \to \Vtau,$ defined by
$$
h(\gl, x) = \sum_{s \in W^F} E_{+,s}(P_0\col \gl \col x) \dE(P_0\col \gl \col y)v
$$
belongs to $\cE_0^\hyp(\spXp \col \tau),$ hence to
$C^{\ep, \hyp}_0(\spXp\col \tau).$ Using \bib{BSanfam}, Eqn.~(14.13),
we see that if  $s \in W,$ $\mu \in \N\gD$ and $u\in \NKaq,$ then
\begin{equation}
\naam{e: implication when q of h non trivial}
q_{s, \mu}(P_0,u\mid h) \neq 0 \;\;\implies\;\; s \in W^F.
\end{equation}
In the notation of \S~\refer{s: normalization},
define the Laurent functional $\Lau \in \Mer(\staFqdc, \gS_F)^*_\laur$ by
$$
\Lau = \sum_{\gl \in \gL(F)} \Res_{\gl + \faFqd}^{t};
$$
then $\supp \Lau \subset \gL(F).$
It follows from (\refer{e: kernel K t F}) that $f = \Laustar h.$
{}From \bib{BSanfam}, Prop.~13.2 (b), it now follows that there exists
an open dense subset $\Omega \subset \faFqdc$ such that, for $\nu \in \Omega,$
$$
\Exp(P_0,u\mid f_\nu) \subset \{ s(\nu +\gl) - \rho - \mu\mid
s \in W,\gl \in \gL(F), \mu \in \N \gD,\, q_{s,\mu}(P_0, v\mid h) \neq 0\}.
$$
In view of (\refer{e: implication when q of h non trivial})
this implies that the inclusion
(\refer{e: exponents K F t in temp set}) holds for $\nu \in \Omega.$
{}From  $f \in \cE^\hyp_F(\spX\col \tau)$ it follows in particular that there exists a
finite subset $Y \subset \staFqdc$
such that $f \in C^{\ep, \hyp}_{F,Y}(\spXp\col \tau).$
The canonical map $W \to W/W_F$ restricts to a bijection $s \mapsto \bar s$ from
$W^F$ onto $W/W_F.$ For $s \in W^F$ we put $E_{\bar s} = -s\gL(F) +\N \gD.$
We now apply Lemma \refer{l: restriction of exponents by holomorphy new},
with $Q = P_F,$ $P = P_0,$ so that $W/\!\!\simPQ \simeq W/W_F,$
and with  $E_\gs$ as just defined, for $\gs \in W/W_F.$  Then it
follows that the inclusion (\refer{e: exponents K F t in temp set}) holds for
$\nu \in \reg f.$
\qed

\begin{thm}
\naam{t: exponents of genEF}
Let $\psi \in \cA_F^{\start}$ and $p \in \Pi_{\gSrF}(\faFqd).$
Then $g: (\nu, x) \mapsto p(\nu) \genEF(\nu\col x)\psi$
defines a family in
$\vanfamF.$ Moreover, for each $v \in \cW$
and every $\nu \in \reg g,$
\begin{equation}
\naam{e: inclusion exponents genEF}
\Exp(P_0,v\mid g_\nu) \subset W^F(\nu + \gL(F)) - \rho -\N \gD,
\end{equation}
where $\gL(F)$ denotes the finite subset of $-\R_+F$ introduced
above (\refer{e: kernel K t F}).
\end{thm}

\proof
The first assertion follows from
Theorem \refer{t: char of genEis by asymp new}.
By Lemma \refer{l: dEis spans cA t F} there exists a $\nu_0 \in \faFqdc$
and elements $y_j \in \spXp$ and $v_j \in \Vtau,$ for $1 \leq j \leq r,$
such that $\ga_F\after\dEF(\nu_0\col y_j)v_j, \, 1 \leq j \leq r$
is a basis for
$\cA^{\start}_F.$
Define meromorphic $\cA^{\start}_F$-valued functions
on $\faFqdc$ by
$\psi_j: = \ga_F\after \dEF(\dotvar \col y_j)v_j,$ for $1 \leq j \leq r.$
By standard arguments involving analyticity and linear algebra it follows
that $(\psi_j(\nu)\mid 1 \leq j \leq r)$ is a basis for $\cA^{\start}_F,$
for $\nu$ in an open dense subset of $\faFqdc.$
Moreover, $\psi \in \cA^{\start}_F$ may be expressed as a linear combination
$\psi = \sum_{1 \leq j \leq r} c_j(\nu) \psi_j(\nu),$
with  meromorphic functions  $c_j: \faFqdc \to \C.$ Using
(\refer{e: K t F as prod Eis})
we
now deduce that
$$
g(\nu\col x) =  \sum_{j=1}^r c_j(\nu) p(\nu) K^t_F(\nu \col x \col y_j) v_j,
$$
as an identity of meromorphic functions in the variable $\nu \in \faFqdc.$
{}From Lemma \refer{l: exponents of K F t} it follows
that there exists a  dense  open subset
$\Omega \subset \faFqdc$ such that $\nu \mapsto g_\nu$ is regular on $\Omega$
and such that, for $\nu \in \Omega,$ the inclusion
(\refer{e: inclusion exponents genEF}) is valid.
{}From  $g \in \cE^\hyp_F(\spX\col \tau)$ it follows that there exists a
finite subset $Y \subset \staFqdc$
such that $g \in C^{\ep, \hyp}_{F,Y}(\spXp\col \tau).$
By the same argument as at the end of the proof of Lemma
\refer{l: exponents of K F t} we now conclude that the inclusion
(\refer{e: inclusion exponents genEF})
is valid for every $\nu \in \reg g.$
\qed

\begin{cor}
\naam{c: nE tempered}
Let notation be as in Theorem \refer{t: exponents of genEF}.
Then, for each $\nu \in i\faFqd \cap \reg g,$
$$
g_\nu \in  \cA_\temp(\spX\col \tau).
$$
\end{cor}

\proof
Let $\nu \in i\faFqd \cap \reg g.$ Then from (\refer{e: inclusion exponents genEF})
it follows that every $(P_0,v)$-exponent of $g_\nu$ is
of the form $\xi = s(\nu + \eta) - \rho - \mu,$ with $s \in W^F, \,\eta \in \gL(F)$
and $\mu \in \N\gD.$ Now $\gL(F) \subset - \R_+ F,$ hence $s \eta \in - \R_+ \gS^+.$
It follows that $\Re\xi + \rho  = s \eta - \mu \in - \R_+ \gD,$
hence $\Re \xi + \rho \leq 0$ on $\Aqp(P_0).$ In view of
Lemma \refer{l: criterion for temperedness},
this implies that $g_\nu \in \cA_\temp(\spX\col \tau).$
\qed

\section{Initial uniform estimates}\eqnreset
In this section we shall derive estimates for the generalized
Eisenstein integrals $\genEF(\nu),$ with uniformity in the
parameter $\nu \in \faFqdc,$ from similar estimates for the
normalized Eisenstein integral $\nE(\gl)= \nE(P_0\col \gl).$ The
idea is that estimates of the latter survive the application
of certain Laurent functionals.

We start with an investigation of the type of estimates involved.
For $Q \in \allparabs$ and $R \in \R,$ we define
$$
\faQqd(Q,R): =
\{\nu \in \faQqdc \mid \Re \inp{\nu}{\ga} < R, \;\;\; \forall \ga \in \gS_r(Q) \}.
$$
The closure of this set is denoted by $\bfaQqd(Q,R).$ It is readily seen
to consist of all elements $\nu \in \faQqdc$ with
$\Re \inp{\nu}{\ga} \leq R$ for all $\ga \in \gSr(Q).$

In the following lemma we assume that $S$ is a finite subset
of $\faQqd\setminus \{0\}$ and we use the notation of
Section \refer{s: Laurent functionals}.
\begin{lemma}
\naam{l: estimate preservation}
Let $R \in \R,$ $p \in \Pi_S(\faQqd),$ $u \in S(\faQqd)$ and $n \in \N.$
Then for every real number $R_- < R$ and every $\gd >0$
there exists a constant $C>0$ with the following
property.

Assume that $V$ is any  complete locally convex space,
$s$ a continuous seminorm on $V$ and $b > 0$ a constant. Moreover,
let $f: \faQqd(Q,R) \to V$ be
a holomorphic function satisfying the estimate
$$
s( p(\nu)  f(\nu) ) \leq (1 + |\nu|)^n e^{b |\Re \nu|},
$$
for all $\nu \in \faQqd(Q, R).$ Then
$$
s(u f(\nu))\leq C (1 + |\nu|)^n e^{b \gd}e^{b |\Re \nu|},
$$
for all $\nu \in \faQqd(Q, R_-).$
\end{lemma}

\proof
It suffices to prove this on the one hand for $u=1$ and $p$ arbitrary
and on the other hand for $p=1$ and $u$ arbitrary.
In the first case the proof is essentially the same as
that of Lemma 6.1 in \bib{Bps2}, which is based on an application
of Cauchy's integral formula.

In the second case the proof relies on a
straightforward application of Cauchy's integral
formula.
\qed

Let $\Hyp$ be a $\gS$-configuration in $\faqdc.$ For $Y \subset \staQqdc$
a finite subset, we define the $\gS_r(Q)$-configuration
$\Hyp_Q(Y)= \Hyp_{\faQqdc}(Y)$ in $\faQqdc$ as in (\refer{e: defi Hyp L S}),
with $L = \faQqdc$
and $S = Y,$  see also
\bib{BSanfam}, text preceding Cor.~11.6.
Thus, for $\nu \in \faQqdc,$ we have
$$
\nu \in \faQqdc\setminus \cup\Hyp_Q(Y)
\iff
\{\forall \gl \in Y\, \forall H\in \Hyp :\quad \gl + \nu \in H \implies \gl + \faQqdc \subset H\}.
$$
Let now $\Lau \in \Mer(\staQqdc, \Sigma_Q)^*_\laur$ have support contained in
the finite subset $Y$ of $\staQqdc.$ For any locally convex space $V$
we have an associated continuous linear operator $\Laustar$ as in
(\refer{e: Laustar with E}). The following result
expresses the continuity with uniformity in the space $V.$

\begin{lemma}
\naam{l: continuity Laustar}
Let $\Hyp, Y, \Lau$ be as above, and let $d: \Hyp \to \N$ be a map.
Then there exists a map $d': \Hyp_Q(Y) \to \N$
with the following property. For any locally convex space $V,$ the
prescription
$$
\Lau_* f(\nu) = \Lau [f(\dotvar + \nu)]
$$
defines a continuous linear operator
$$
\Lau_*: \Mer(\faqdc, \Hyp, d, V) \to \Mer(\faQqdc, \Hyp_Q(Y), d', V).
$$

\end{lemma}
\proof
This is Cor.~11.6 of \bib{BSanfam}.
\qed

A real $\gSr(Q)$-configuration $\Hyp'$ in $\faQqdc$
consists of hyperplanes of the form
$$
H_{\ga,s} := \{\nu \in \faQqdc\mid \inp{\ga}{\nu} =s \},
$$
with $\ga \in \gSr(Q)$ and $s \in \R.$ The configuration $\Hyp'$
is called $Q$-bounded if there exists a constant $s_0 \in \R$ such that
$
H_{\ga, s} \in \Hyp' \implies s \geq s_0,
$
for all $\ga \in \gSr(Q), s\in \R.$
See \bib{BSres}, text before Lemma 3.1, for the similar notion for $Q$ minimal.

\begin{lemma}
\naam{l: Hyp Q Y is Q bounded}
Let $Q\in \allparabs,$ $P \in \minparabs$ and $P \subset Q.$
Let $Y \subset \staQqd$ be a finite subset.
\begin{enumerate}
\itema
If $\Hyp$ is a $P$-bounded real $\gS$-configuration in $\faqdc,$
then $\Hyp_Q(Y)$ is a $Q$-bounded real $\gS_r(Q)$-configuration in $\faQqdc.$
\itemb
If $\Hyp'$ is a $Q$-bounded real $\gS_r(Q)$-configuration in $\faQqdc,$
then for every $R \in \R$ the collection
$
\{H \in \Hyp' \mid H \cap \bfaQqd(Q,R) \neq \emptyset \}
$
is finite.
\itemc
If $\Hyp'$ is as in (b), then for every $R \in \R$ there exists
a constant $R_+ >  R$ such that
$H \cap \bfaQqd(Q,R) \neq \emptyset \iff H \cap \faQqd(Q,R_+) \neq \emptyset,$
for every $H \in \Hyp'.$
\end{enumerate}
\end{lemma}

\proof
There exists $t_0 \in \R$ such that
the hyperplanes in $\Hyp$ are all of the form $H_{\gb, t},$ with
$\gb \in \gS(P)$ and $t\in [t_0, \infty\, [.$  Let $\eta \in Y$
and assume that $-\eta + H_{\gb, t}$ intersects $\faQqdc$ in a proper hyperplane $H'.$
Then it follows that the restriction $\ga = \gb|_{\faQq}$ is non-zero, hence belongs
to $\gSr(Q).$ Moreover, $H' = H_{\ga, s}$  with $s = t - \inp{\gb}{\eta}.$ Let $m$ be the
maximum of the numbers $\inp{\gb}{\eta},$ for $\gb\in \gS(P)\setminus \gS_Q$ and $\eta \in Y$
and put $s_0 = t_0 - m.$ Then it follows that every hyperplane
from $\Hyp_Q(Y)$ is of the form $H_{\ga, s},$ with $\ga \in \gSr(Q)$ and $s \geq s_0.$
This establishes (a).

To prove (b), fix $\ga \in \gSr(Q)$ and put $I_{\ga, R} =
\{s \in \R \mid H_{\ga,s} \in \Hyp,\;
H_{\ga, s} \cap \bfaQqd(Q,R) \neq \emptyset\},$ for every $R > 0.$
Then it suffices to show that $I_{\ga,R}$ is finite.
Since $\Hyp$ is locally finite, the set $I_{\ga,R}$ is discrete, and since $\Hyp$ is
$Q$-bounded, the set $I_{\ga,R}$ is bounded from below.
If $h \in H_{\ga, s} \cap \bfaQqd(Q,R),$
then $s = \inp{\ga}{h} \leq R.$ It follows that the set
$I_{\ga, R}$ is bounded from above by $R.$
Hence, $I_{\ga, R}$ is finite.

For (c) we observe that $R \leq R' \implies I_{\ga,R} \subset I_{\ga,R'}.$
Fix $R' > R.$ Using that
$I_{\ga,R'}$ is discrete, we see that we may choose $R_+ \in \;]\, R, R'\, [\,$
sufficiently close to $R$ so that $I_{\ga, R_+} = I_{\ga, R}$ for all $\ga \in \gSr(Q).$
The constant $R_+$ has the required property.
\qed

If $\Hyp'$ is a $Q$-bounded real $\gSr(Q)$-configuration in $\faQqdc,$
and $d': \Hyp' \to \N$ a map, then,
for $R \in \R,$ we define the polynomial function $\pi_{Q,R,d'}$
on $\faQqdc$ in analogy with (\refer{e: defi pi omega d}) by
$$
\pi_{Q,R,d'}: = \prod_{H}
l_H^{d'(H)},
$$
where the product is taken over the collection of $H \in \Hyp'$
whose intersection with $\bfaQqd(Q,R)$ is non-empty; this collection
is finite by Lemma \refer{l: Hyp Q Y is Q bounded} (b).
It follows from Lemma \refer{l: Hyp Q Y is Q bounded} (c)
that $\pi_{Q,R,d'} = \pi_{Q,R_+, d'},$
for $R_+ > R$ sufficiently close to $R.$

\begin{prop}
\naam{p: induced estimate Laustar}
Let $Q\in \allparabs,$ $P \in \minparabs$ and $P \subset Q.$
Let $Y \subset \staQqd$ be a finite subset and
let $\Lau \in \Mer(\staQqdc, \gS_Q)^*_\laur$ be
a Laurent functional with $\supp \Lau \subset Y.$
Let $\Hyp$ be a $P$-bounded $\gS$-configuration in $\faqdc,$
$d: \Hyp \to \N$ a map, and let $d': \Hyp_Q(Y) \to \N$ be associated with
the above data as in Lemma \refer{l: continuity Laustar}.
Let $M > \max_{\eta \in Y} |\Re \eta|$ and assume that
$R,R' \in \R$ are constants with $Y + \faQqd(Q, R') \subset \faqd(P,R_-)$
for some $R_- < R.$

There exists a constant $k \in \N$
and for every $n \in \N$ a constant $C > 0$ with the following
property.

If $V$ is a complete locally convex space, $s$ a continuous seminorm on $V,$
$b > 0$ a positive constant and
$\gf$ a function in $ \Mer(\faqdc, \Hyp, d, V)$ satisfying the estimate
$$
s( \pi_{P,R,d}(\gl)   \gf (\gl) ) \leq (1 + |\gl|)^n e^{b |\Re \gl|},
$$
for all $\gl \in \faqd(P, R),$
then the function $\Lau_* \gf \in \Mer(\faQqdc, \Hyp_Q(Y), d', V)$
satisfies the estimate
$$
s( \pi_{Q,R',d'}(\nu) \Lau_* \gf (\nu )) \leq C (1 + |\nu|)^{n+k}
e^{b M} e^{b |\Re \nu|}
$$
for all $\nu \in \faQqd(Q,R').$
\end{prop}

\proof
It suffices to prove this for the case that $\supp \Lau$ consists of
a single point $\gl_0 \in Y.$
Let $\Hyp_0$ be the collection of $H \in \Hyp$ containing
$\gl_0 + \faQqdc.$ Then for every $H \in \Hyp_0$ there is a unique
indivisible root $\ga_H \in \Sigma_Q \cap \Sigma(P)$ such that
$H=\gl_0 + (\ga_H^\perp)_{\C}.$ We define
the affine function $l_H: \faqdc \to \C$ by $l_H(\gl) = \inp{\gl - \gl_0}{\ga_H}.$
Then $H = l_H^{-1}(0).$ We define the polynomial function
$q_0: \faqdc \to \C$  by
$$
q_0 = \prod_{H \in \Hyp_0} l_H^{d(H)}.
$$
{}From the definition of the space of Laurent functionals
in $\Mer(\staQqdc, \gS_Q)^*_\laur$
supported at $\gl_0,$ see \S \refer{s: Laurent functionals},
it follows that there
exists a $u \in S(\staQqd)$ such that on
a function $f \in (q_0|_{\staQqdc})^{-1}\cO_{\gl_0}$ the action of the Laurent
functional is given by $\Lau f = u(q_0|_{\staQqdc}f)(\gl_0).$
It follows that for $\gf \in \Mer(\faqdc, \Hyp, d,V)$ and
$\nu \in \faQqdc \setminus \cup \Hyp_Q(Y),$
$$
\Lau_* \gf(\nu) = \Lau(\gf(\dotvar + \nu)) = u(q_0(\dotvar)\gf(\dotvar + \nu))(\gl_0)
=
u(q_0\gf)(\gl_0 + \nu).
$$
As in the proof of \bib{BSres}, Lemma 1.2, we
infer that there exist a polynomial function $\pi \in \Pi_{\gS_r(Q)}(\faQqd)$
and finitely many $q_j \in P(\faQqd)$ and $u_j \in S(\faqd),$
all independent of $V,s,b$ and $\gf,$ such that
\begin{equation}
\naam{e: Laurent estimate 1}
\pi(\nu) \, u (q_0 \gf) (\gl_0 + \nu)
=
\sum_j q_j(\nu) u_j (\pi_{P, R, d} \gf)(\gl_0 + \nu),
\end{equation}
for $\nu \in \faQqd(Q, R').$
Multiplying both sides of (\refer{e: Laurent estimate 1})
with a suitable
polynomial function we see that we may as well assume that $\pi = \pi_0
\pi_{Q,R', d'},$ for some $\pi_0\in \Pi_{\gSr(Q)}(\faQqd).$
We obtain
\begin{equation}
\naam{e: formula with pi zero psi}
\pi_0(\nu) \psi(\nu) = \sum_j q_j(\nu) u_j(\pi_{P,R,d}\gf)(\gl_0 + \nu),
\end{equation}
where we have written $\psi = \pi_{Q,R', d'} \Lau_* \gf.$

Let $k$ be the maximum of the degrees of
the polynomials $q_j.$ Then there exists a constant $D > 0,$
independent of $V,s,b$ and $\gf,$
 such that
for every $j,$
\begin{equation}
\naam{e: estimate q j}
|q_j(\nu)| \leq D (1 + |\nu|)^k, \qquad (\nu \in \faQqdc).
\end{equation}
Put $m = \max_{\eta \in Y} |\Re \eta|$ and fix $\gd >0$
such that $m + 2 \gd < M.$
We may select  constants $R_+' > R'$ and $R_- < R$ such that
$Y + \bfaQqd(Q,R_+') \subset \faqd(P, R_-).$ Adapting $R_+'$ if necessary,
we may in addition assume that
\begin{equation}
\naam{e: equality pi s}
\pi_{Q, R_+', d'} = \pi_{Q, R', d'},
\end{equation}
see the text preceding the proposition.

Let now $n \in \N$ and $b >0,$
and assume that $\gf$ satisfies the hypotheses of the proposition.
It follows from Lemma \refer{l: estimate preservation}, applied with $P$ in place of $Q,$
that there exist constants
$C_j > 0,$ independent of $V,s,b$ and $\gf,$
such that
\begin{equation}
\naam{e: estimate action Lau 2}
s(u_j(\pi_{P,R,d}\,\gf)(\gl) ) \leq C_j (1 + |\gl|)^n e^{b \gd} e^{b |\Re \gl|},
\end{equation}
for $\gl \in \faqd(P, R_-).$

Using the estimate
$1 + |\gl_0 + \nu| \leq (1 + |\gl_0|)(1 + |\nu|)$
and combining (\refer{e: formula with pi zero psi}),
(\refer{e: estimate q j}) and (\refer{e: estimate action Lau 2}),
we obtain
$$
s(\pi_0(\nu) \psi(\nu)) \leq
C' (1 + |\nu|)^{n +k} e^{b \gd} e^{b (|\Re\nu| + m)},\qquad (\nu \in \faQqd(Q, R_+')),
$$
with
$$
C' = D (\sum_{j} C_j) (1 + |\gl_0|)^n.
$$
{}From (\refer{e: equality pi s}) we see
that the function $\psi$ is holomorphic on  $\faQqd(Q, R_+').$
We may therefore apply Lemma \refer{l: estimate preservation}
with $\psi,$ $\pi_0,$ $C'$ and $[C'e^{b\gd + b m}]^{-1} s$ in place
of $f, p, u$ and $s,$ and with $R_+', R'$ in place of $R, R_-$, respectively.
Using that $m + 2\gd < M,$
we obtain the desired estimate, with $C>0$ a constant that is independent of $V,s,b$ and $\gf.$
\qed

In the rest of the section we shall apply the above results to Eisenstein integrals.
We start with a suitable estimate for Eisenstein integrals associated with minimal
$\gs$-parabolic subgroups.

\begin{lemma}
\naam{l: estimate minimal Eis}
Let $P \in \minparabs.$ Then there exists a $\bar P$-bounded real
$\gS$-configuration $\Hyp$ in $\faqdc$ and a map $d: \Hyp \to \N$
such that the function $\gl \mapsto \nE(P\col \gl)$
belongs to the space $\Mer(\faqdc, \Hyp, d, \Ci(\spX) \otimes \Hom(\oC, \Vtau)).$

Let $R \in \R$ and let $p \in \Pi_{\gS}(\faqd)$ be a polynomial
such that the function $\gl \mapsto p(\gl) \nE(P\col \gl)$
is holomorphic on a neighborhood of $\bfaqd(\bar P, R).$ Then there exists a constant
$r > 0$ and for every $u \in U(\fg)$ constants $n \in \N$ and $C>0$
such that
\begin{equation}
\naam{e: estimate minimal Eis}
\| p(\gl) \nE(P\col \gl \col u ; x)\| \leq C (1 + |\gl|)^n e^{(r + |\Re \gl|)\lengthX(x)},
\end{equation}
for all $\gl \in \bfaqd(\bar P, R)$ and $x \in \spX.$
(See (\refer{e: defi lspX}) for the definition of the function $\lspX.$)
\end{lemma}

\proof
First assume that $\tau = \tau_\types,$ defined as in \bib{BSft}, text after Eqn.~(28),
with $\types \subset \dK$ a finite subset.
Then for  $x \in \Aq$ the estimate (\refer{e: estimate minimal Eis})
follows from \bib{Bps2}, Corollary 16.2 and Proposition 10.3, combined with
the fact that $\nE(P\col \gl) = E^1(\bar P\col \gl),$ see
\bib{BSft}, Eqn.\ (52). In view of the decomposition $\spX = K\Aq (eH),$ the estimate
now follows for general $x \in \spX$ by sphericality of the Eisenstein
integral. Finally, for general $\tau$ the estimate follows
by application of the `functorial' dependence of the Eisenstein integral on
$\tau,$ see \bib{BSft}, Eqn.~(32).
\qed

We can now prove the following analogous result
for the generalized Eisenstein integral.

\begin{prop}
\naam{p: init estimate genEis}
Let $F \subset \gD,$ $v \in \FcW$
and let $t \in \WT(\gS)$ be a $W$-invariant even residue weight.

There exists a $\bar P_F$-bounded, real $\gSr(F)$-hyperplane
configuration $\Hyp_F$ in $\faFqdc$ and a map $d_F: \Hyp_F \to \N$ such
that $\nu \mapsto \genEF(\nu)$ belongs to the space
$$
\Mer(\faQqdc, \Hyp_F, d_F, \Ci(\spX) \otimes \Hom(\cAFt, \Vtau)).
$$
Moreover, if $R'\in \R$ and if
$p$ is any polynomial in $\Pi_{\gSr(F)}(\faFqd)$
such that $\nu \mapsto p(\nu)\genEFv(\nu)$ is holomorphic
on a neighborhood of $\bfaFqd(\bar P_F,R'),$ then
there exist a constant $r> 0$ and for every $u \in U(\fg)$ constants
$n \in \N$ and $C>0,$ such that
\begin{equation}
\naam{e: estimate for nEFv}
\|p(\nu) \genEFv(\nu\col u ; x)\| \leq C (1 + |\nu |)^n e^{(r  + |\Re \nu |)\lengthX(x)},
\end{equation}
for all $\nu \in \bfaFqd(\bar P_F,R')$ and $x \in \spX.$
\end{prop}

\begin{rem}
This result is a sharpening of the estimate
given in \bib{BSfi}, Lemma 10.8.
\end{rem}

\proof
According
to Lemma \refer{l: psi as Lau of nE} (b,c), the
generalized Eisenstein integral may be expressed as
\begin{equation}
\naam{e: genEFv as Lau on nE}
\genEFv(\nu\col x) \psi = \Lau_*[\nE(P_0\col \dotvar \col x) \after \iFv](\nu),
\end{equation}
with $\Lau \in \Mer(\faFqdc, \gS_F)^*_\laur \otimes \oCFv$
a Laurent functional whose support is contained in a finite
subset $Y \subset \staQqd.$
Let $\Hyp,d$ be associated with $P = P_0$ as in Lemma \refer{l: estimate minimal Eis}.
Let $\Hyp_F:= \Hyp_F(Y)$ and $d_F: =d'$ be associated with the
data $P=P_0,$ $Q = P_F,$ $\Hyp,$ $d,$ $\Lau$ as in Lemma \refer{l: continuity Laustar},
then $\Hyp_F$ is a real $\bar P_F$-bounded $\gSr(F)$-configuration,
by Lemma \refer{l: Hyp Q Y is Q bounded}.
The first assertion of the proposition
follows  by application of Lemma \refer{l: continuity Laustar}.

Fix $R \in \R$ such that $\faFqd(\bar P_F, R') + Y \subset \faqd(\bar P_0, R_-) $
for some $R_- < R.$  Let $r$ be the constant
of Lemma \refer{l: estimate minimal Eis}
applied with $P_0, R$ and $\pi_{\bar P_0, R, d}$ in place of $P, R$ and
$p,$ respectively.
Fix  $u \in U(\fg).$
Then according to Lemma \refer{l: estimate minimal Eis} there exist constants
$n_0\in \N$ and $C_0>0,$
such that for all $x \in \spX,$
$$
\| \pi_{\bar P_0, R, d}(\gl) \nE(P_0\col \gl \col u;x)\after \iFv \|
\leq C_0 (1 + |\gl|)^{n_0} e^{(r + |\Re \gl|)\lengthX(x)},
$$
for all $\gl \in \bfaqd(\bar P_0, R).$
Let $x \in \spX.$
We apply Proposition \refer{p: induced estimate Laustar} to the function $\gf = \gf_x$
from  $\Mer(\faqdc, \Hyp, d, \Hom(\oCFv, \Vtau)),$ given by
$$
\gf_x = \nE(P_0\col \dotvar \col u ; x)\after \iFv,
$$
with the constant $b = l_\spX(x)$ and the seminorm $s = C_0^{-1} e^{-rb} \|\dotvar\|.$
Let $M > \max_{\eta \in Y}|\Re \eta|.$ Then we obtain the estimate
$$
\| \pi_{\bar P_F, R', d_F}(\nu) \Lau_*(\gf_x)(\nu) \| \leq C_1 (1 + |\nu|)^{n_0 +k}
e^{(r + |\Re \nu| + M)\lengthX(x)},
$$
for $\nu \in \faFqd(\bar P_F, R'),$
with constants $k \in \N$ and $C_1 >0$ that are independent of $x, \nu.$
We now note that
\begin{eqnarray*}
\Lau_*(\gf_x) & = & \Lau_*(\nE(P_0\col \dotvar\col u;x)\after\iFv)(\nu) \\
&=&
\Lau_*(\nE(P_0\col \dotvar)\after\iFv)(\nu)(u;x)\\
&=&
\genEFv(\nu\col u;x),
\end{eqnarray*}
as a meromorphic identity in $\nu \in \faQqdc.$
The second of the above identities involves the interchange of $u$ and $\Laustar,$ which
is allowed by the continuity of $\Lau_*,$ see
Lemma \refer{l: continuity Laustar}. The third identity
is obtained by application of (\refer{e: genEFv as Lau on nE}).
Thus, we obtain, for all $x \in \spX$ and all $\nu \in \faFqd(\bar P_F, R'),$
the estimate
$$
\| \pi_{\bar P_F, R', d_F}(\nu)
\genEFv(\nu\col u;x) \| \leq C_1 (1 + |\nu|)^{n_0 +k}
e^{(r + |\Re \nu| + M)\lengthX(x)}.
$$
This proves the result for the particular polynomial $p = \pi_{\bar P_F,R', d_F}.$
For $p$ equal to a multiple of $\pi_{\bar P_F,R', d_F}$ the result
now also follows, since any polynomial from $P_d(\faFqd),$ $d \in \N,$
can be estimated from above by a function of the form $C (1+ |\nu|)^d.$

Let $p$ now be an arbitrary element of $\Pi_{\gSr(F)}(\faFqd)$
satisfying the hypothesis. Fix $R_+' > R'.$
Then an estimate of type (\refer{e: estimate for nEFv}) holds on $\faFqd(\bar P_F, R_+')$
with $p \pi_{\bar P_F, R'_+ , d_F}$ in place of $p.$
By application of Lemma \refer{l: estimate preservation} this
implies an estimate of the form (\refer{e: estimate for nEFv}), with the required
dependences of the constants.
\qed

\section{Symmetric pairs of residue type}\eqnreset
\naam{s: spaces of residue type}
By
$\Ltwod(X)$ we denote the discrete part of $L^2(\spX),$
i.e., the closed span in $L^2(\spX)$ of all the irreducible
closed subspaces of $L^2(\spX).$
Accordingly we define
$$
\Ltwod(\spX\col \tau): = ( \Ltwod(\spX)\otimes V_\tau)  \cap L^2(\spX\col \tau).
$$
For the following definition we recall from Section \refer{s: Fourier inversion}
that the data $(G,H,K,\tau, \faq, \gS^+)$ together with a
$W$-invariant even residue
weight $t\in \WT(\gS)$  determine the continuous linear operator
$T^t_\gD: \Cci(\spX\col \tau) \to \Ci(\spX\col \tau).$
If $G$ has a compact center modulo $H,$
then this operator is given by the formula (\refer{e: formula for T gD}).

\begin{defi}
\naam{d: good residue type new}
The
reductive symmetric pair $(G,H)$ is said to be
of residue type if the following conditions are fulfilled.
\begin{enumerate}
\itema
The group $G$ has a compact center modulo $H.$
\itemb
For any choice of the data $(K,\faq),$  the following
requirement is fulfilled.
For every finite dimensional unitary representation $\tau$ of $K,$
every choice $\gS^+$ of positive roots  for $\gS$
and every $W$-invariant even residue weight
$t \in \WT(\Sigma),$
the operator $T_\gD^{t}: \Cci(\spX\col \tau) \to \Ci(\spX \col \tau)$
is the restriction of the orthogonal projection
$L^2(\spX \col \tau) \to \Ltwod(\spX \col \tau).$
\end{enumerate}
\end{defi}

\begin{rem}
\naam{r: technical reasons residue type}
The above definition is given for technical reasons. Together with
Definition \refer{d: parabolic of residue type}, where the notion
of residue type for a parabolic subgroup from $\allparabs$ will be introduced,
it plays a role in a long chain of reasoning that will be used
in an induction step in the proof of Theorem
\refer{t: inversion for Schwartz}. (The induction goes by induction on
the $\gs$-split rank of $G.$)
As part of the mentioned theorem
it is asserted that in fact every pair $(G,H)$ with $G$ having compact center modulo
$H$ and every group from $\allparabs$ is of residue type.

In the course of the chain of reasoning,
many results will first be proved under the assumption that an
involved parabolic subgroup from $\allparabs$ or an involved
reductive  symmetric pair $(G,H)$ is of residue type.
Such results will be marked with (RT) after their
declaration. The additional hypotheses will be clearly stated at
the beginning of their proofs. Within the chain of reasoning,
until Theorem \refer{t: inversion for Schwartz},
the results marked (RT) will only be used if these additional
hypotheses are assumed to be fulfilled. The mentioned
theorem implies that the additional hypotheses are in fact always fulfilled
so that in the end the results marked (RT) are valid as stated.

Within the chain of reasoning, it also happens that definitions need
extra hypotheses concerning residue type in order to be valid. These definitions
will be marked (RT) as well. The extra hypotheses are stated in a subsequent remark.
Within the chain of reasoning such definitions will only
be used if the extra hypotheses are assumed to be fulfilled.
In the end Theorem \refer{t: inversion for Schwartz} implies that
the extra hypotheses are always fulfilled,
so that the definitions marked (RT) are valid as stated.
\end{rem}

\begin{rem}
\naam{r: compact is residue type}
If $\faq = \{0\},$ then $\spX$ is compact and
the operator $T_\gD$ is understood to be the identity
operator of $\Cci(\spX\col \tau).$ Thus, conditions (a) and (b) of the definition
are fulfilled and in this case $(G,H)$ is of residue type.
\end{rem}

\begin{rem}
\naam{r: stability residue type under iso}
It follows from Lemma \refer{l: conjugation T gD t}
that the notion of residue type is
stable under isomorphisms of reductive
symmetric pairs.
\end{rem}

\begin{rem}
\naam{r: one K suffices for residue type}
Condition (b) of the definition is valid as soon as
a particular choice of the data $(K,\faq)$ satisfies the mentioned requirement.
Indeed, assume that  $(K,\faq)$ satisfies the requirement and let
$\bp K \subset G$ be a second $\gs$-invariant maximal compact subgroup,
$\fg = \bp \fk \oplus \bp\fp$ the associated Cartan decomposition,
and $\bp \faq \subset \fq\cap \bp \fp$ an associated maximal abelian subspace.
Then there exists a (unique) $g \in \exp(\fh \cap \fp)$ such that
$gKg^{-1} = \bp K,$
see \bib{Loos1}, p.~153.
Now $\Ad(g) \faq$ is maximal abelian in $\bp \fp \cap \fq,$ hence there
exists an element $k \in \bp K_e\cap H_e$ such that $\Ad(kg) \faq = \bp \faq.$
Let $\gf: G \to G$ be conjugation by $kg,$ then $\gf$ maps the data
$(G,H, K, \faq)$ onto $(G,H, \bp K, \bp \faq).$ In view of
Lemma \refer{l: conjugation T gD t}
it follows that the requirement in (b)
for the pair $(K,\faq)$ is equivalent
to the similar requirement for the pair $(\bp K, \bp \faq).$
\end{rem}

We recall from \bib{Bps2}, \S~17, that the Schwartz space
$\cC(\spX\col \tau)$ is defined to be the space of functions
$f \in \Ci(\spX\col \tau)$ such that, for every $u \in U(\fg)$
and $n \in \N,$
\begin{equation}
\naam{e: Schwartz estimate}
s_{u,n}(f):= \sup_{x \in \spX} (1 + \lspX)^{n} \Theta(x)^{-1} \|uf(x)\| < \infty;
\end{equation}
see also the beginning of \S~\refer{s: temperedness}.
The Schwartz space is equipped with the Fr\'echet topology
determined by the seminorms $s_{u,n}.$ By \bib{Bas}, Lemma 7.2, the operators
from $\DX$ act on $\cC(\spX\col \tau)$ by continuous linear endomorphisms.
We define  $\cA_2(\spX\col \tau)$ to be the space of $\DX$-finite
functions in $\cC(\spX\col \tau).$

\begin{lemma}\restypeone
\naam{l: cAt equals cAtwo}
Let $G$ have compact center modulo $H.$
Then, for every $W$-invariant
even residue weight $t \in \WT(\gS),$
\begin{equation}
\naam{e: equalities cAt Ltwod and cAtwo}
\cA^t(\spX\col \tau) = \Ltwod(\spX\col \tau) = \cA_2(\spX\col \tau).
\end{equation}
In particular, the space $\cA_2(\spX\col \tau)$ is finite dimensional.
\end{lemma}

\begin{rem}
\naam{r: finite dimensional Ltwod}
The fact that $L^2_d(\spX\col \tau)$ is finite dimensional is
also an immediate consequence of the classification of the discrete series
in \bib{OMds}. In the present paper it would not have been advantageous
to use this known fact. Only the spectral
properties of the discrete series as formulated in Theorem
\refer{t: infinitesimal characters L two d real and regular new}
are needed. The mentioned finite dimensionality
naturally follows from the finite dimensionality of
$\cA^t(\spX\col \tau)$, by the nature of our argument.
\end{rem}

\proof We give the proof
under the assumption that $(G,H)$ is of residue
type, see Remark \remRT.
It follows from
 Definition \refer{d: good residue type new} (b)
that $T_\gD^t(\Cci(\spX\col \tau))$ is dense in $\Ltwod(\spX\col \tau).$
By Lemma \refer{l: cA t is image T t gD} it follows that
$\cA^t(\spX\col \tau)$ is dense in  $\Ltwod(\spX\col \tau).$
By finite dimensionality of the first of these spaces, the first equality
in (\refer{e: equalities cAt Ltwod and cAtwo}) follows.
In particular, it  follows that the space $\Ltwod(\spX\col \tau)$
consists of smooth $\DX$-finite functions; by \bib{Bas}, Thm.~7.3 with $p=2$
it is therefore contained
in $\cA_2(\spX\col \tau).$ Conversely, if $f \in \cA_2(\spX\col \tau),$
then $f$ is $K$-finite and $\DX$-finite. Hence, by a well known result of
Harish-Chandra its $(\fg, K)$-span in $\cC(\spX, \Vtau)$ is a $(\fg, K)$-module
of finite length; see \bib{Var}, p.~312, Thm.~12 and \bib{Wal1}, p.~112, Thm.~4.2.1.
The  closure of this span in $L^2(\spX) \otimes \Vtau$
is therefore a finite direct
sum of irreducible representations.
The mentioned closure contains $f;$
hence, $f \in L^2_d(\spX \col \tau).$
\qed

Assume that $G$ has compact center modulo $H$ and that $t \in \WT(\gS)$ is
a $W$-invariant even residue weight. We recall that
a choice of Hilbert structure
on the space $\cAt(\spX\col \tau)$ uniquely determines an endomorphism
$\ga = \ga^t \in \End(\cAt(\spX\col \tau))$
such that (\refer{e: K t gD and ga}) holds.

\begin{lemma}\restypeone
\naam{l: endomorphism ga}
Assume that $G$ has compact center modulo $H$
and let $t \in \WT(\gS)$ be
a $W$-invariant even residue weight.
\begin{enumerate}
\itema
Assume that the space
$\cAt(\spX\col \tau)$ is equipped with the restriction
of the inner product from $L^2(\spX) \otimes \Vtau,$ see
(\refer{e: equalities cAt Ltwod and cAtwo}). Then the
endomorphism $\ga,$ determined by (\refer{e: K t gD and ga}),
equals $|W|^{-1}$ times
the identity operator of $\cAt(\spX \col \tau).$
\itemb
The kernel $K^t_\gD$ is independent of the residue weight $t.$
\end{enumerate}
\end{lemma}

\proof We give the proof under the assumption that $(G,H)$ is of residue
type, see Remark \remRT.
It follows from (\refer{e: equalities cAt Ltwod and cAtwo}) that the real analytic
$\tau \otimes 1$-spherical
function  $\bfe$ on $\spX$ attains its values in
$\Hom(\Ltwod(\spX \col \tau), \Vtau) \simeq \Vtau \otimes \Ltwod(\spX\col \tau)^*.$
Hence, $\bfe^*: y\mapsto \bfe(y)^*$ is a real analytic $1 \otimes \tau^*$-spherical
function on $\spX,$ with values in $\Hom(\Vtau, \Ltwod(\spX\col \tau)).$

We define the continuous linear operator $P: \Cci(\spX\col \tau) \to \Ltwod(\spX\col \tau)$
by
$$
P f = \int_X \bfe(y)^* f(y) \; dy.
$$
Then one readily verifies that $\hinp{Pf}{\psi} = \hinp{f}{\psi}$ for every
$\psi \in \Ltwod(\spX\col \tau).$ It follows that $P$ equals the restriction
to $\Cci(\spX\col \tau)$ of the orthogonal projection
$L^2(\spX\col \tau) \to \Ltwod(\spX\col \tau).$
Hence, $P = T_\gD^t.$ Combining this with (\refer{e: formula for T gD})
we obtain that, for all $x \in \spX$ and
all $f \in \Cci(\spX\col \tau),$
\begin{eqnarray*}
|W| \int_X K_\gD^t (x\col y) f(y)\; dy &=& T_\gD^t f (x)\;\;
= \;\;\bfe(x) (Pf)\\
&=&
\int_X \bfe(x)\after \bfe(y)^* f(y) \; dy.
\end{eqnarray*}
Since $\bfe(x)\after \bfe(\dotvar)^*$ and $K_\gD^t(x\col \dotvar)$ are smooth and
$\tau^* \otimes 1$-spherical
functions on $\spX,$ with values in $\End(\Vtau),$
it follows from the above identities that
$|W| K_\gD^t(x\col y) = \bfe(x) \after \bfe(y)^*$ for all $x,y \in \spX.$
This implies (\refer{e: K t gD and ga})
with $\ga$ equal to $|W|^{-1} I_{\Ltwod(\spX\col \tau)}.$ Hence, (a) holds.
Assertion (b) is now immediate.
\qed

\section{The normalized Eisenstein integral}\eqnreset
\naam{s: the normalized Eisenstein integral}
In this section we shall define the normalized Eisenstein integral,
initially for the class
of parabolic subgroups introduced in the following definition, see
Remark \remRT.
\begin{defi}
\naam{d: parabolic of residue type}
A parabolic subgroup $P \in \allparabs$ is said to be of
residue type (relative to $H$)
if for every
$v \in \NKaq$ the pair $(M_P, M_P \cap vHv^{-1})$ is of
residue type.
A subset $F \subset \gD$ is said to be of residue type if the associated
standard $\gs$-parabolic subgroup
$P_F$ is of residue type.
\end{defi}

\begin{rem}
In view of Remark \refer{r: stability residue type under iso}
it suffices to require the above condition
for $v$ in a choice of representatives $\PcW\subset \NKaq$ of $W_P\bs W/ \WKH.$
\end{rem}

\begin{rem}
\naam{r: G of residue type}
We write ${}^\circ G$ for $M_G,$ the Langlands $M$-component of $G$ viewed as a parabolic
subgroup. Thus, ${}^\circ G$ equals the intersection of the kernels $\ker |\chi|,$
for $\chi: G \to \C^*$ a character.
By the previous remark, $G$ (viewed as an element from $\allparabs$)
is of residue type relative
to $H$ if and only if the pair $({}^\circ G, {}^\circ G \cap H)$
is of residue type. If $G$ has compact center modulo $H,$
then the pair $({}^\circ G, {}^\circ G \cap H)$ is of residue type
if and only if $(G,H)$ is.
\end{rem}

\begin{defi}
\naam{d: associated parabolics}
Two parabolic subgroups $P,Q \in \allparabs$ are said to be associated
if their  $\gs$-split components
$\faPq$ and $\faQq$ are conjugate under $W.$ The equivalence relation
of associatedness is denoted by $\sim.$
\end{defi}

\begin{lemma}
\naam{l: residue type and sim}
Assume that $P \in \allparabs$ is of residue type.
Then every $Q\in \allparabs$ with $Q\sim P$ has the same property.
\end{lemma}
\proof
If $Q \sim P,$ there exists a $k \in \NKaq$ such that $\faQq = \Ad(k) \faPq.$ {}From
this it follows that $M_Q = k M_P k^{-1}.$ If $v \in \NKaq,$ then $M_Q \cap vHv^{-1}$
equals $k(M_P \cap wH w^{-1}) k^{-1},$ with $w = k^{-1}v \in \NKaq.$
The pair $(M_P, M_P \cap w H w^{-1})$ is of residue type, and by Remark
\refer{r: stability residue type under iso} we conclude
that $(M_Q, M_Q \cap vHv^{-1})$ is of residue type as well.
\qed
Let $P \in \allparabs.$
We equip the space $\spXPv,$ for $v \in \PcW,$ with the invariant measure
$dx_{P,v}$ specified at the end of Section \refer{s: normalization}.
The space $\cA_2(\spXPv\col \tau_P)$ is equipped with the
inner product
from $L^2(\spX_{P,v}, \Vtau),$ for $v \in \PcW.$
Moreover, the
space $\cA_{2,P} = \cA_{2,P, \PcW}$
is defined to be the formal direct  sum
\begin{equation}
\naam{e: defi cAtwoP}
\cA_{2,P}: = \bigoplus_{v \in \PcW} \cA_2(\spXPv\col \tau_P),
\end{equation}
equipped
with the direct sum inner product. The space $\cA_{2,P}$
is finite dimensional by Lemma \refer{l: cAt equals cAtwo}. Application
of this lemma requires $P$ to be of residue type, see Remark \remRT.

We agree to denote by $\pr_{P,v}: \cA_{2,P}\to  \cA_2(\spXPv\col \tau_P)$
the natural projection operator, for $v \in \PcW,$
and by $\rmi_{P,v}$ the associated
natural embedding operator.

In the following we shall use the characterization of the generalized
Eisenstein integral by its asymptotic behavior, see Theorem
\refer{t: char of genEis by asymp new},
to define an Eisenstein integral for
arbitrary parabolic subgroups.

\begin{prop}\restypetwo{\ }
\naam{p: intro nE P psi}%
\begin{enumerate}
\itema
Let $P\in \allparabs.$
For every $\psi \in \cA_{2,P}$ there exists
a unique family $\nE(P\col \psi) \in \cE_{P}^\hyp(\spX\col \tau)$
with the following property.
For all $\nu$ in a non-empty open subset of
$\faPqdc,$ each $v \in \PcW,$ every $X \in \faPq$ and every $m \in \spXPvp,$
$$
q_{\nu -\rho_P} (P, v\asmid \nE(P\col \psi\col \nu), X , m) = \psi_v(m).
$$
\itemb
Let $F \subset \gD.$
Then, for every $W$-invariant even residue weight $t \in \WT(\gS),$ the space
$\cA^\start_F$ equals $\cA_{2, F}: = \cA_{2, P_F};$  here we assume that in the definition
of both spaces the same set $\FcW$ has been used.
Moreover, for every $\psi \in \cA_{2,F},$
$$
\nE(P_F \col \psi \col \nu \col x) = \genEF(\nu\col x)\psi,
$$
for all $x \in \spX$ and generic $\nu \in \faFqdc.$
\end{enumerate}
\end{prop}

\proof
Here we prove (a) under the assumption that $P$ is of residue type
and (b) under the assumption that $F$ is of residue type, see Remark \remRT.

Uniqueness follows from Theorem \refer{t: vanishing thm new}.
Thus, it suffices to establish existence.
We will first do this for $P$ of residue type and equal to a standard parabolic subgroup
$P_F,$ with $F \subset \gD.$ Let $t \in \WT(\gS)$ be any
$W$-invariant even residue
weight. Let ${}^*t$ be the induced residue weight of $\gS_F.$
Then it follows from Lemma \refer{l: cAt equals cAtwo},
applied to the pair $(M_F, M_F \cap vHv^{-1}),$
that $\cAstart(\spXFv\col \tau_F) = \cA_2(\spXFv\col \tau_F).$
Moreover, if $\psi \in \cA_2(\spXFv\col \tau_F),$ then
$\nE(P_F\col \psi) := \genEF(\dotvar)\psi$ satisfies the desired property, by
Theorem \refer{t: char of genEis by asymp new}.

Now assume that $P$ is general and of residue type, let $\PcW\subset \NKaq$
be a choice of representatives for $W_P\bs W/\WKH$ and let $\psi \in \cAtwoP.$
There exists a $u \in \NKaq$ such that
$u^{-1}P u = P_F,$ with $F \subset \gD.$ Moreover, $P_F$ is of
residue type, by Lemma \refer{l: residue type and sim}. The set
$\FcW = u^{-1} (\PcW)$ is a choice
of representatives for $W_F\bs W/\WKH$ in $\NKaq.$
For $v \in \PcW,$ let $\rho_{\tau,u}$ be the linear isomorphism
from $\Ci(\spX_{F,u^{-1}v,+}\col \tau_F)$ onto $\Ci(\spX_{P, v,+}\col \tau_P),$
defined as in \bib{BSanfam}, Eqn.\ (3.24).
We define the function
$\psi' \in \cAFtwo$ by $\psi_{u^{-1}v}' = \rho_{\tau,u}^{-1} \psi_{v},$
for $v \in \PcW.$ Define the meromorphic family $f: \faPqdc  \to \Ci(\spX\col \tau)$
by
$$
f_\nu( x) = \nE(P_F\col \psi' \col u^{-1} \nu \col x),
$$
for $x \in \spX$ and generic $\nu \in \faPqdc.$ By
Theorem \refer{t: char of genEis by asymp new}
and Lemma \refer{l: Weyl group on cEhyp}, the family $f$
belongs to $\cE_P^\hyp(\spX\col \tau).$ Moreover, it follows from
\bib{BSanfam}, Lemma 3.6, that, for $v \in \PcW$ and  $\nu$ in a dense  open subset
of $\faPqdc,$
\begin{eqnarray*}
q_{\nu - \rho_P}(P, v\asmid f_\nu) &=&
[\Ad(u^{-1})^* \otimes \rho_{\tau,u}]\,
q_{u^{-1}\nu - \rho_F}(P_F, u^{-1}v \asmid f_\nu)\\
&=& \rho_{\tau,u} \psi_{u^{-1}v}' = \psi_v.
\end{eqnarray*}
This establishes the result with $\nE(P\col \psi) =f.$
\qed
{}From the uniqueness assertion in Proposition \refer{p: intro nE P psi}
it follows that the meromorphic function $\nE(P\col \psi):
\faPqdc \to \Ci(\spX\col \tau)$
depends linearly on $\psi.$

\begin{defi}\restypetwo
\naam{d: dual Eisenstein new}
Let $P \in \allparabs.$
For $\psi \in \cAtwoP,$ let
$\nE(P\col \psi)$
denote the unique family in $\cE^\hyp_P(\spX\col \tau)$
of Proposition \refer{p: intro nE P psi} (a).

The meromorphic $\Ci(\spX, \Hom(\cA_{2,P}, \Vtau))$-valued function
$\nE(P\col \dotvar) = \nE(\spX\col P\col \dotvar)$ on $\faPqdc,$ defined by
$$
\nE(P\col \nu\col x)\psi = \nE(P\col \psi\col  \nu \col x),
$$
for  $\psi \in \cAtwoP,$ $x \in \spX$ and generic $\nu \in \faQqdc,$
is called the normalized Eisenstein integral associated with
the parabolic subgroup $P$ (and the choice $\PcW$).

The meromorphic $\Ci(\spX, \Hom( \Vtau, \cA_{2,P}))$-valued function
$E^*(P\col \dotvar)$ on $\faPqdc$ defined by
$$
E^*(P\col \nu \col x) = \nE(P\col -\bar \nu \col x)^*,
$$
for $\psi \in \cAtwoP,$
$x \in \spX,$ and generic $\nu \in \faPqdc,$
is called the dual Eisenstein integral associated with $P.$
\end{defi}

\begin{rem}
The above definition requires the validity of Proposition \refer{p: intro nE P psi},
which inside the chain of reasoning leading up to
Theorem \thmendRT,
requires $P$ to be of residue type, see Remark \remRT.
\end{rem}

\begin{rem}
\naam{r: comparison with nE of HC}
In the case of the group, the normalized Eisenstein integral
defined above is essentially equal to
the one defined by Harish-Chandra \bib{HCeis}, \S 6, Thm.\ 6. This is seen
as follows. Let $\bp G$ be a real reductive group of Harish-Chandra's class,
let $G = \bp G \times \bp G,$ let $\gs: G \to G$ be the involution
given by $(\bp x,\bp y) \mapsto (\bp y, \bp x)$
and let $H = G^\gs$ be the diagonal subgroup.
Then $(\bp x, \bp y) \mapsto \bp x(\bp y)^{-1}$
induces a $G$-diffeomorphism $\Phi: G/H \to \bp G.$
Let $\bp\fg = \bp \fk \oplus \bp \fp$ be a Cartan decomposition,
$\bp \Cartan$ the associated involution and put $\Cartan = \bp \Cartan \times \bp \Cartan.$
Then $\Cartan$ is a Cartan involution commuting with $\gs.$
Let $\bp K$ and $K = \bp K \times \bp K$ be the associated maximal
compact subgroups of $\bp G$ and $G,$ respectively,
and let $(\tau, \Vtau)$ be a finite
dimensional unitary representation of $K.$
Let $\tau_1, \tau_2$
be the unitary representations of $\bp K$ in $\Vtau$ defined by
$\tau_1(\bp k) = \tau(\bp k, 1)$ and $\tau_2(\bp k) = (1, \bp k).$
Let $\bp \tau$ denote the pair of commuting representations $(\tau_1, \tau_2).$
Then pull-back by $\Phi$ induces a linear isomorphism $\Phi^*$
from the space
$\Ci(\bp G \col \bp \tau)$ of smooth $\bp \tau$-spherical functions on $\bp G,$
onto $\Ci(G/H\col\tau).$

Let $\fa_{\bp \fp}$ be maximal abelian in $\bp \fp,$ then
$\faq: = \{(X,-X)\mid X \in \fa_{\bp \fp}\}$ is maximal abelian in $\fp\cap \fq.$
Let $\bp P \subset \bp G$ be a
parabolic subgroup containing $A_{\bp \fp}.$ Then $P: = \bp P \times \bp \bar P$ belongs
to $\allparabs.$ Moreover, the map $(X,Y) \mapsto X - Y$ is surjective from
$\fa_P = \fa_{\bp P} \times \fa_{\bp P}$ onto $\fa_{\bp P}$ and induces a linear
isomorphism from $\fa_{P\iq} = \fa_P \cap \faq$ onto $\fa_{\bp P},$ mapping
$\fa_{P\iq}^+$ onto $\fa_{\bp P}^+.$ The complexified
adjoint map $\gf^*$ is the linear isomorphism
from $\fa_{\bp P\iC}^*$ onto $\fa_{P\iq\iC}^*$ given by
$\bp \nu \mapsto (\bp \nu, - \bp \nu).$

We observe that $M_P =  M_{\bp P} \times M_{\bp P},$ so that
$\spX_P = M_{\bp P} \times M_{\bp P} /H\cap M_P,$ which is $M_P$-diffeomorphic
to $M_{\bp P}$ under the map $\Phi_P$ induced by restricting $\Phi.$
It is readily seen that $\Phi_P^*$ restricts
to a linear bijection from the finite dimensional space
$L^2_d(M_{\bp P} \col \bp \tau_{M_{\bp P}})$ onto $ L^2_d(\spXP \col \tau_P)= \cA_{2,P}.$
Let $\psi \in L^2_d(\spXP \col \tau_P)$ and consider the family $f: \faqdc \times G \to \Vtau$
defined by $f_\nu  = \Phi^* \nE(\bp P \col \Phi_P^{*-1} \psi \col \gf^{*-1}
\nu/i),$
where the normalized Eisenstein integral is Harish-Chandra's.
By holomorphy
of Harish-Chandra's unnormalized Eisenstein integral
combined with meromorphy and the product structure of
Harish-Chandra's $C$-function $C_{\bp P| \bp P}(1 \col \bp \nu),$
it follows that the family $f_\nu \in \Ci(\spX\col \tau)$
satisfies condition (a)
of Definition \refer{d: Cephyp}.
Via the subrepresentation
theorem of \bib{CM}, combined with induction by stages, Harish-Chandra's
Eisenstein integral can be realized by matrix coefficients
of the minimal principal series of $\bp G.$ Using this information
it can be deduced that $f$ satisfies the
remaining conditions (b), (c) of Definition \refer{d: Cephyp}.
In fact, it is now readily checked that $f$ belongs to the space
$C^{\hyp}_P(\spX \col \tau),$
defined in Definition \refer{d: cEhypQallX}.
Moreover, from the information on the constant term
of Harish-Chandra's Eisenstein integral, see \bib{HCeis}, \S 5, Thm.~5,
it follows that
$q_{\nu - \rho_P}(P,1\mid f_\nu, X, m) = \psi(m),$ in the notation of
Proposition \refer{p: intro nE P psi}, for generic $\nu \in i\faPqd.$
By meromorphy it now follows
that the family $f$ satisfies the condition of Proposition
\refer{p: intro nE P psi} (a)
(note that we may take ${}^P\cW = \{1\}$ here).
Hence,
$$
\Phi^* \nE(\bp P \col \Phi_P^{*-1} \psi \col \gf^{*-1} \nu/i) = \nE(P \col \nu)\psi.
$$
\end{rem}

\begin{rem}
\naam{r: Eis and Delorme}
At the end of the sequel to this paper, \bib{BSpl2},
we will show that the normalized Eisenstein integral introduced above
coincides (up to a change from $\nu$ to $-\nu$)
with the one introduced by J.\ Carmona and P.\ Delorme in \bib{CDn}.
\end{rem}

\begin{rem}
\naam{r: normalized Eisenstein for P is G}
If $G$ has compact center
modulo $H,$ then $A_{G\iq} = \{0\}$ and $\cA_{2,G} = \cA_2(G/H
\col \tau).$ In view of Lemma \refer{l: cAt equals cAtwo}, it follows
from Remark \refer{r: Eisenstein for P is G} that
 $\nE(G\col 0 \col x)$ equals the evaluation map
$\ev_x: \cA_2(\spX \col \tau) \to V_\tau,$ $\psi \mapsto \psi(x).$
Accordingly,
$\dE(G\col 0\col x) = \ev_x^* \in \Hom(\Vtau, \cA_2(\spX\col \tau)).$
\end{rem}

The following result describes the dependence of the normalized Eisenstein
integral on a member $P$ of a class in $\allparabs/W,$ as well
as on the choice of $\PcW.$

\begin{lemma}\restypetwo
\naam{l: Eis and R P}
Let $P \in \allparabs,$
let $s \in W,$ and let
$Q:= sPs^{-1}.$  Let $\PcW$ and $\QcW$ be choices of representatives
in $\NKaq,$ for $W_P\bs W/\WKH$ and $W_Q\bs W/\WKH,$ respectively.
Then there exists a unique linear map $R_P(s): \cAtwoP \to \cAtwoQ$
such that
\begin{equation}
\naam{e: Eis and R P}
\nE(Q\col s \nu \col x)\after R_P(s) = \nE(P\col \nu \col x),
\end{equation}
for $x \in \spX$ and generic $\nu \in \faPqdc.$  The map $R_P(s)$ is bijective
and unitary.
\end{lemma}

\proof
We  give the proof under the assumption that
$P$ is of residue type, see Remark \remRT.
Left multiplication by $s$ induces a bijective map $W_P\bs W/\WKH \to W_Q\bs W/\WKH.$
Via the natural bijections $\PcW \to W_P\bs W/\WKH$ and $\QcW \to W_Q\bs W/\WKH$
we transfer the induced map to a bijection $\bar s: \PcW \to \QcW.$

Let $\psi \in \cAtwoP$ and define the meromorphic family
$f: \faQqdc \to \Ci(\spX\col \tau)$
by
\begin{equation}
\naam{e: f as Eis P}
f_\gl(x) = \nE(P\col s^{-1}\gl \col x)\psi,
\end{equation}
for $x \in \spX$ and generic
$\gl \in \faQqdc.$ Then $f \in \cE^\hyp_Q(\spX\col \tau),$ by Lemma
\refer{l: Weyl group on cEhyp}.
Select $v \in \PcW.$ We may select a representative
$u_s$ in $\NKaq$ of a Weyl group element from $s W_P$ such that
$u_s v =  \bar s(v) w_s$ for some $w_s \in N_{K\cap H}(\faq).$
Note that $\spX_{Q, u_s v} = \spX_{Q, \bar s(v)}.$ Hence, we may define
the bijective linear map
$$
\rho_{\tau, u_s}:
\cA_2(\spX_{P, v}\col \tau_P) \to \cA_2(\spX_{Q, \bar s(v)}\col \tau_Q)
$$
as in \bib{BSanfam}, Eqn.\ (3.24).  This map is
unitary in view of the choice of invariant measures on
$\spX_{P, v}$ and $\spX_{Q, \bar s(v)},$ specified at the end of Section
\refer{s: normalization}.

It follows from \bib{BSanfam}, Lemmas 3.7 and 3.6,
that
\begin{eqnarray}
q_{\gl - \rho_Q}(Q, \bar s(v) \asmid f_\gl)
&=&
q_{\gl - \rho_Q}(Q, u_s v \asmid f_\gl)\nonumber\\
&=&  [\Ad(u_s^{-1})^*\otimes \rho_{\tau, u_s}]
q_{s^{-1}\gl - \rho_P}(P, v \asmid f_\gl)\nonumber\\
&=&
\naam{e: q s nu of f nu}
\rho_{\tau, u_s} \psi_v,
\end{eqnarray}
for generic $\gl \in \faQqdc.$ Hence, by Definition \refer{d: dual Eisenstein new},
\begin{equation}
\naam{e: f as Eis Q}
f_\gl(x) = \nE(Q\col \gl \col x)\psi',
\end{equation}
with $\psi_{\bar s (v)}' = \rho_{\tau, u_s} \psi_v$ for $v \in \PcW.$
We define the bijective linear map $R_P(s): \cAtwoP \to \cAtwoQ$ by
\begin{equation}
\naam{e: formula for R P}
(R_P(s)\psi)_{\bar s(v)} = \rho_{\tau, u_s} \psi_v.
\end{equation}
Then $\psi' = R_P(s) \psi$ and (\refer{e: Eis and R P}) follows from
(\refer{e: f as Eis P}) and (\refer{e: f as Eis Q})
by substituting  $s \nu$ for $\gl.$  {}From the definition
it follows that $R_P(s)$ is unitary.

To establish uniqueness, let $R_P(s): \cAtwoP \to \cAtwoQ$
be a linear map.
Let $\psi \in \cAtwoP,$ define $f$ as above, and
define the meromorphic family $g: \faQqdc \to \Ci(\spX\col \tau)$
by $g_\gl = \nE(Q\col \gl \col \dotvar)R_P(s)\psi.$ Then
\begin{equation}
\naam{e: q equals R P s psi}
q_{\gl - \rho_Q}(Q, \bar s(v)\mid g_\gl) = (R_P(s)\psi)_{\bar s(v)}.
\end{equation}
Now assume that (\refer{e: Eis and R P}) holds. Then $g = f;$ combining
(\refer{e: q s nu of f nu}) and (\refer{e: q equals R P s psi})
we  obtain
(\refer{e: formula for R P}).
\qed

Let $P\in \allparabs.$
Then for all  $x,y \in \spX,$ the meromorphic $\End(\Vtau)$-valued function
on $\faPqdc$ given by
\begin{equation}
\naam{e: nu mapsto nE after dE}
\nu \mapsto \nE(P\col \nu \col x) \dE(P\col \nu \col y)
\end{equation}
depends a priori on the choice of the set
$\PcW.$

\begin{cor}\restypetwo
\naam{c: independence kernel of PcW}
Let $P \in \allparabs.$
Then for every $x,y \in \spX$ the function
(\refer{e: nu mapsto nE after dE}) is independent of the particular choice of $\PcW.$
\end{cor}

\proof
Here  we assume that $P$ is of residue type, see Remark \remRT.
The result then
follows from application of
Lemma \refer{l: Eis and R P}, with $s =1.$
\qed

\begin{prop}\restypetwo
\naam{p: init estimate normalized Eis}
Let $P \in \allparabs.$
There exists a $\bar P$-bounded, real $\gS_r(P)$-hyperplane
configuration $\Hyp = \Hyp_P$ in $\faPqdc$ and a map $d = d_P: \Hyp \to \N$ such
that $\nu \mapsto \nE(P\col \nu)$ belongs to the space
$$
\Mer(\faPqdc, \Hyp, d, \Ci(\spX) \otimes \Hom(\cA_{2,P}, \Vtau)).
$$
Moreover, if $R \in \R$ and if
$p$ is any polynomial in $\Pi_{\gS_r(P)}(\faPqd)$
such that $\nu \mapsto p(\nu)\nE(P\col \nu)$ is holomorphic
on a neighborhood $\bfaPqd(\bar P,R),$
 then
there exist a constant $r> 0$ and for every $u \in U(\fg)$ constants
$n \in \N$ and $C>0,$ such that
\begin{equation}
\naam{e: estimate for nE P}
\|p(\nu) \nE(P\col \nu\col u ; x)\| \leq C (1 + |\nu |)^n
 e^{(r  + |\Re \nu |)\lengthX(x)},
\end{equation}
for all $\nu \in \bfaPqd(\bar P,R)$ and $x \in \spX.$
\end{prop}

\proof
Here we prove the result for $P$ of residue type, see Remark \remRT.

We first assume that $P = P_F$ with $F \subset \gD$ of residue type.
In this case,
$\nE(P\col \nu) = \genEF(\nu),$ by Proposition \refer{p: intro nE P psi}.
Hence, the result follows
from Proposition \refer{p: init estimate genEis} by summation
over $\FcW,$ see  (\refer{e: genEFv}) and (\refer{e: defi AsubF}).

Next, let $P \in \allparabs$ be a general parabolic subgroup
of residue type.
There exists a $s \in W$ such that $sPs^{-1} =  P_F ,$
by Lemma \refer{l: on parabolics} (c). Since $P$ is of residue
type, $P_F$ is of residue type as well, see Lemma \refer{l: residue type
and sim}. By Lemma \refer{l: Eis and R P} and
Proposition \refer{p: intro nE P psi}
there exists a unitary map
$R_P(s): \cA_{2,P} \to \cA_{F,2}$ such that
$$
\nE(P_F \col s\nu \col x)\after R_P(s) = \nE(P \col \nu \col x)
$$
for all $x$ and generic $\nu\in \faPqdc.$ The result
now follows by application of the first part of the proof.
\qed

The following result limits the exponents of the normalized
Eisenstein integral along a minimal
parabolic subgroup. To formulate it we need the following notation.
Let $Q \in \allparabs$ and  let $P \in \minparabs$ be such that
$P \subset Q.$ We put
\begin{eqnarray*}
\gD_Q(P) & := & \{\ga \in \gD(P) \mid \ga|_{\faQq} = 0\};\\
W^{P|Q}&:= &\{t \in W \mid t(\gD_Q(P))\subset \gS(P)\}.
\end{eqnarray*}
Let $s \in W$ be the unique element such that $P = sP_0s^{-1}.$
Then $s^{-1}Qs$ contains $P_0$ hence equals $P_F$ for some
subset $F \subset \gD;$ note that $F$ is uniquely determined
by $Q$ in view of Lemma \refer{l: on parabolics}.
We define
$$
\gL(P|Q):= s \gL(F),
$$
where $\gL(F)$ is the finite subset of $-\R_+F$ introduced
in (\refer{e: kernel K t F}). We note that $\gL(P|Q)$ is a finite
subset of $-\R_+ \gD_Q(P).$

\begin{prop}\restypetwo
\naam{p: exponents normalized Eis}
Let $Q \in \allparabs$
and let $P \in \minparabs$ be contained in $Q.$
Let $\psi \in \cA_{2,Q}$ and $q \in \Pi_{\gS_r(Q)}(\faQqd).$
Then $f: (\nu, x) \mapsto q(\nu) \nE(Q\col \nu\col x)\psi$
defines a family in $\cE_Q^\hyp(\spX\col \tau).$ Moreover,
for each $v \in \NKaq$ and every $\nu \in \reg f,$
\begin{equation}
\naam{e: inclusion exponents and gL P Q}
\Exp(P,v\mid f_\nu) \subset W^{P|Q}(\nu + \gL(P|Q)) - \rho_P - \N \gD(P).
\end{equation}
In particular, $f_\nu \in \cA_\temp(\spX\col \tau),$ for every
$\nu \in i\faQqd \cap \reg f.$
\end{prop}

\begin{rem}
For $P$ minimal, the assertion about temperedness is due to \bib{Bps2}, Thm.~19.2, in view of
\bib{BSft}, Eqn.~(52). For general $P$ the assertion about temperedness
is due to
\bib{D1n}, Thm.~1, in view of
Remark
\refer{r: Eis and Delorme}.
\end{rem}

\proof
We give the proof under the assumption that $Q$ is of
residue type, see Remark \remRT.
Let $s \in W$ and $F \subset \gD$ be as in the text preceding the corollary.
Let the polynomial function $p: \faFqdc \to \C$ be defined by
$p(\nu)= q(s\nu).$ Then $p \in \Pi_{\gSr(F)}(\faFqd).$ It follows
from Lemma \refer{l: Eis and R P} with $P_F$ in place of $P$
that $f(\nu, x) = g(s^{-1}\nu, x),$ for $x \in \spX$ and generic $\nu \in \faQqdc,$
where $g: (\gl,x)  \mapsto p(\gl) \genEF(\gl \col x) R_{P_F}(s)^{-1} \psi.$
By Theorem \refer{t: exponents of genEF} the family $g$ belongs
to $\cE_F^\hyp(\spX\col \tau).$ By Lemma \refer{l: Weyl group on cEhyp}
it follows that  $f \in \cE^\hyp_Q(\spX\col \tau).$
Moreover, let $\nu \in \reg f;$ then $s^{-1}\nu \in \reg g$ and
by the last mentioned  theorem it follows that,
for every $u \in \NKaq,$
\begin{eqnarray*}
\Exp(P_0, u \mid f_\nu) & = &
\Exp(P_0, u \mid g_{s^{-1}\nu}) \\
&  \subset &  W^F(s^{-1} \nu + \gL(F)) - \rho - \N \gD.
\end{eqnarray*}
On the other hand, by \bib{BSanfam}, Lemma 3.6, it follows that,
for $v \in \NKaq,$
$$
\Exp(P,v \mid f_\nu) = s\Exp(P_0, \bar s^{-1} v \mid f_\nu),
$$
where $\bar s$ is any representative of $s$ in $\NKaq.$
We conclude that
$$
\Exp(P, v\mid f_\nu) \subset sW^F s^{-1}(\nu + s\gL(F)) - \rho_P - \N\gD(P).
$$
Now $s\gL(F) = \gL(P|Q)$ by definition. Moreover, one
readily verifies that $s W^F s^{-1} = W^{P|Q}.$ Hence,
(\refer{e: inclusion exponents and gL P Q}) follows.
The final assertion follows from the similar assertion for $g,$ which in turn follows by application
of Corollary \refer{c: nE tempered}.
\qed

In the theory of the constant term, we shall need the following
result on the coefficients of the asymptotic expansions of the Eisenstein
integral.

\begin{lemma}\restypetwo
\naam{l: q belongs to Mer newer}
Let $P \in \allparabs$
and let $\psi \in \cAtwoP.$
The family $f: (\nu, x) \mapsto \nE(P\col \nu\col x)\psi$ belongs to
$\cE^\hyp_{P,Y}(\spX\col \tau),$ for a suitable finite subset $Y \subset \staPqd.$
Moreover, the $\gSr(P)$-configuration
$\Hyp_f,$ defined as in the text before Definition
\refer{d: defi cEhyp Q Y gd}, is real.

Let $k = \dega f.$
Then for every $Q \in \allparabs,$ $v \in \NKaq,$ each $\gs \in W/\!\!\sim_{Q|P}$
and all $\xi \in - \gs\cdot Y + \N \DrQ,$
\begin{equation}
\naam{e: q gs xi in Mer for nEP}
q_{\gs,\xi}(Q,v\mid f) \in P_k(\faQq) \otimes
\Mer(\faPqdc, \Hyp_f, d_f, \Ci(\spXQv \col \tau_Q)).
\end{equation}
\end{lemma}

\proof
We
give the proof under the assumption that
$P$ is of residue type, see Remark \remRT.
In view of Lemma \refer{l: on parabolics} (a), there
exist $s \in W$ and $F \subset \gD$ such that $P = sP_Fs^{-1}.$
In view of Lemma \refer{l: Eis and R P} with $P$ and $P_F$ in place
of $Q$ and $P,$ respectively,
we may as well assume that $P = P_F$ for some $F \subset \gD.$
In this case we have $\nE(P\col \nu) = \genEF(\nu)$
by Proposition \refer{p: intro nE P psi}.
Hence, the result follows from Lemma \refer{l: globality q of genEF newer}.
\qed

\begin{lemma}\restypetwo
\naam{l: formula kernel F}
Let $F \subset \gD$
and
let  $t \in \WT(\gS)$ be a $W$-invariant
even residue weight.
Then, for all $x,y \in \spX,$
\begin{equation}
\naam{e: K F t in terms of Eis}
K^t_F(\nu \col x\col y) = |W_F|^{-1}
\nE(P_F \col \nu \col x)\after \dE(P_F \col \nu \col y),
\end{equation}
as an identity of $\End(\Vtau)$-valued meromorphic
functions in the variable $\nu \in \faFqdc.$

In particular, the function $K_F^t$  does not depend on the residue
weight $t,$ nor on the choice of $\FcW.$
\end{lemma}

\proof
We give the proof under the assumption that
$F$ is of residue type, see Remark \remRT.
{}From Proposition \refer{p: intro nE P psi} (b) we recall
that $\cAFt = \cAFtwo.$
Accordingly, we equip the space $\cAFt$ with the inner product
described in the text preceding (\refer{e: defi cAtwoP}). As in the text
preceding (\refer{e: K t F as prod Eis}),
this choice of inner product determines an endomorphism
$\ga_{F,v} \in\End(\cA^{\start}(\spXFv\col \tau_F)),$ for each $v \in \FcW.$
The endomorphism  $\ga_{F,v}$ is the analogue for the space $\spXFv$ of the endomorphism
$\ga,$ described in (\refer{e: K t gD and ga}). Thus, $\ga_{F,v} = |W_F|^{-1} I,$
by Lemma \refer{l: endomorphism ga}.
Let $\ga_F \in \End(\cAFt)$ be the direct sum of the $\ga_{F,v},$ for $v \in \FcW.$
Then from (\refer{e: K t F as prod Eis}) we obtain that
$$
K_F^t(\nu\col x\col y) = |W_F|^{-1}\, \genEF(\nu \col x) \genEF(-\bar\nu\col y)^*,
$$
for all $x,y \in \spX$ and generic $\nu \in \faFqdc.$ Now use Definition
\refer{d: dual Eisenstein new} and Proposition \refer{p: intro nE P psi}
to conclude the validity
of (\refer{e: K F t in terms of Eis}). It is now obvious
that $K_F^t$ does not depend on $t;$ it follows by application
of Corollary
\refer{c: independence kernel of PcW} that it does not depend on $\FcW$ either.
\qed

\begin{rem}
\naam{r: drop t in K F t}
In view of Lemma \refer{l: formula kernel F} we agree to omit $t$ in the notation $K_F^t.$
\end{rem}

\begin{defi}
\restypetwo
\naam{d: K P}
Let $P \in \allparabs.$
We define the meromorphic
function $K_P: \faPqdc \to \Ci(\spX\times \spX , \End(\Vtau))$
by
$$
K_P(\nu \col  x \col y) = |W_P|^{-1}\nE(P\col \nu \col x)\dE(P\col \nu \col y),
$$
for $x,y \in \spX$ and generic $\nu \in \faPqdc.$
\end{defi}

In the chain of reasoning leading up to Theorem \thmendRT{}
this definition requires $P$ to be of residue type, since only then
the Eisenstein integral is well-defined, see Remark \remRT.

\begin{rem}
\naam{r: K P F is K F}
If $P = P_F,$
for  $F \subset \gD,$ then
$K_{P} = K_F,$ in view of Lemma \refer{l: formula kernel F}.
\end{rem}

\begin{lemma}\restypetwo
\naam{l: action W on K}
Let $P \in \allparabs.$
Then, for every $s \in W$
and all $x,y \in\spX,$
$$
K_P(\nu\col  x\col y) = K_{sPs^{-1}} (s\nu \col  x \col  y)
$$
as a meromorphic identity in $\nu \in \faPqdc.$
\end{lemma}

\proof
We give the proof under the assumption
that $P$ is of residue type,
see Remark \remRT.
Put  $Q =s P s^{-1}.$
Since the inner product $cB,$ specified in Section \refer{s: normalization},
is $W$-invariant,
the normalized measures $d\mu_P$ and $d\mu_Q$ are $s$-conjugate.
Moreover, since $W_Q = sW_P s^{-1},$
we have $|W_Q| = |W_P|.$ The result now follows from combining
Definition \refer{d: K P} and Lemma \refer{l: Eis and R P}.
\qed

\begin{thm}\restypetwo
\naam{t: invariance kernel}
Let $P,Q\in \allparabs$ be associated parabolic subgroups.
Then for every
$s \in W(\faQq\mid \faPq),$ and all $x,y \in \spX,$
\begin{equation}
\naam{e: invariance kernel}
K_Q(s\nu \col x\col y) = K_P(\nu \col x\col y),
\end{equation}
as an identity of $\End(\Vtau)$-valued
meromorphic functions of the variable $\nu \in \faPqdc.$
\end{thm}

\proof
We give the proof under the assumption that $P$ and $Q$ are
of residue type, see Remark \remRT.
Since $P$ and $Q$ are associated, $\dim \faQq = \dim\faPq;$
hence, $s$ is a linear bijection from $\faPq$ onto $\faQq.$
{}From Corollary \refer{c: W a1 a2 in q} it
follows that there exists $w \in W$ such that $w(\faPq) = \faQq$
and $w|_{\faPq} = s.$

Assume first that $P = P_F$ and $Q = P_{F'},$ with $F,F' \subset \gD.$
Then the result follows from
\bib{BSfi}, Lemma 6.2, since $K_F^t = K_{P_F}$ and $K_{F'}^t = K_{P_{F'}},$
for any $W$-invariant even residue weight $t,$ see
Remarks \refer{r: drop t in K F t}
and \refer{r: K P F is K F}.

Next assume that
$P$ and $Q$ are arbitrary.
Then the result follows by using
Lemma \refer{l: on parabolics} (a) and
Lemma \refer{l: action W on K}.
\qed

\section{Eigenvalues for the Eisenstein integral}\eqnreset
\naam{s: eigenvalues}
In this section we  investigate the action of $\DX$
on the normalized Eisenstein integral.

Let $P \in \allparabs.$ We define the algebra homomorphism
$\mu_P: \DX \to \D(\spXoneP)$ as in \bib{Bps2}, text following (20).
Here $\spXoneP := \spX_{1P,e} = M_{1P}/M_{1P} \cap H.$
Let $\fb\subset \fq$ be a $\Cartan$-stable Cartan subspace containing
$\faPq$ and let $\gg_\fb$ be the associated Harish-Chandra isomorphism
from $\DX$ onto $I(\fb).$ Let $W_P(\fb)$ denote the centralizer
of $\faPq$ in $W(\fb),$ and $I_P(\fb)$ the ring of $W_P(\fb)$-invariants in
$S(\fb).$ Moreover, let $\gg^{\spXoneP}_\fb$ denote the associated Harish-Chandra
isomorphism $\D(\spXoneP) \to I_P(\fb).$
Then we recall from \bib{Bps2}, Eqn.\ (21), that
\begin{equation}
\naam{e: composition gg and mu}
\gg^{\spXoneP}_\fb \after \mu_P = \gg_\fb.
\end{equation}
If $v \in \NKaq,$ then following \bib{BSanfam}, text above Lemma 4.12,
we define the algebra homomorphism
$\mu_P^v: \D(G/vHv^{-1}) \to \D(\spX_{1P, v})$ as $\mu_P$ for the
triple $(G, vHv^{-1}, P)$ instead of $(G, H, P).$ Moreover, we define
the algebra homomorphism $\mu_{P,v}: \DX \to \D(\spX_{1P, v})$ by
\begin{equation}
\naam{e: defi mu P v}
\mu_{P,v}= \mu_P^v \after \Ad(v),
\end{equation}
where $\Ad(v)$ denotes the isomorphism $\DX \to \D(G/vHv^{-1})$
induced by the adjoint action by $v.$
Since $\APq$ is central in $M_{1P},$ it follows
from (\refer{e: deco spXonePv})  that
\begin{equation}
\naam{e: deco DspXonePv}
\D(\spX_{1P, v}) \simeq \D(\spX_{P,v}) \otimes S(\faPq).
\end{equation}
Accordingly, if $D \in \DX,$  we shall write $\mu_{P,v}(D\col \dotvar)$
for $\mu_{P,v}(D),$ viewed as a
$\D(\spX_{P,v})$-valued polynomial function on $\faPqdc.$
If $D \in \DX,$ $\nu \in \faPqdc,$ and $v \in \PcW,$ then
$\mu_{P,v}(D\col \nu) \in \D(\spXPv)$ acts on the space $\cA_2(\spX_{P,v}\col \tau_P)$
(see the text preceding Lemma \refer{l: cAt equals cAtwo})
by an endomorphism that we denote by $\umu_{P,v}(D\col \nu).$
The direct sum of these endomorphisms, for $v \in \PcW,$
is an endomorphism of the space $\cA_{2,P},$ denoted $\umu_P(D \col \nu).$

\begin{lemma}\restypetwo
\naam{l: action of DX on nE}
Let $P\in \allparabs.$
Then
$$
D \nE(P\col \nu)  = \nE(P\col \nu)\, \umuP (D\col \nu), \qquad (D \in \DX).
$$
\end{lemma}

\proof
We give the proof under the assumption
that $P$ is of residue type, see Remark \remRT.
Let $\psi \in \cAtwoP.$  Then the
family $f:\faPqdc \times \spX \to \Vtau,$ defined by
$$
f(\nu, x) = \nE(P\col \nu\col x)\psi,
$$
belongs to $\cE_{P}^\hyp(\spX\col \tau),$
by Proposition \refer{p: intro nE P psi}.
Let $D \in \DX.$ The
family $Df: (\nu , x) \mapsto Df_\nu(x)$ belongs
to $\cE_{P}^\hyp(\spX\col \tau)$ as well, by Definition \refer{d: cEhypQallX}
and \bib{BSanfam}, Lemma 9.8.
Moreover, by \bib{BSanfam}, Lemma 6.2, there exists a dense  open subset $\Omega$ of
$\faPqdc$  such that,
for $\nu \in  \Omega,$ the element
$\nu - \rho_P$ is a leading exponent of $f_\nu$ along $(P, v).$ Hence, by
\bib{BSanfam}, Lemma 4.12,
it follows that, for $\nu \in \Omega,$ $X \in \faPq$ and $m \in \spXPvp,$
$$
q_{\nu - \rho_P}(P, v\mid D f)(X, \nu, m) =
\mu_{P,v}(D) \gf_\nu(\exp X m),
$$
where the function $\gf_\nu: \spX_{1P,v, +} \to \Vtau$ is defined
by
\begin{equation}
\naam{e: gf nu defined as q}
\gf_\nu(ma) =
a^{\nu} q_{\nu - \rho_P}(P, v\mid f_\nu, \log a, \nu, m),
\end{equation}
for $a \in \APq$ and $m \in \spXPvp.$
By Proposition \refer{p: intro nE P psi}, the expression on the right-hand side of
(\refer{e: gf nu defined as q}) equals $a^\nu \psi_v(m),$
and we see that
\begin{equation}
\naam{e: q of Df as psi}
q_{\nu - \rho_P}(P, v\mid D f)(X, \nu, m) = \mu_{P,v}(D\col \nu) \psi_v(m),
\end{equation}
for $\nu$ in a dense  open subset of $\faPqdc,$ $m \in \spX_{P,v,+}$ and $X \in \faPq.$

On the other hand, $\nu \mapsto \umuP(D\col \nu)\psi$
is a polynomial $\cAtwoP$-valued function on $\faPqdc.$
It readily follows that the family
$$
g: (\nu, x) \mapsto \nE(P\col \nu \col x)\, \umuP(D\col \nu) \psi
$$
belongs to $\cE_P^\hyp(\spX\col \tau).$
Moreover, by Proposition \refer{p: intro nE P psi}, it follows that
\begin{equation}
\naam{e: q of g as psi}
q_{\nu - \rho_P}(P, v\mid g)(X, \nu, m) = \mu_{P,v}(D\col \nu) \psi_v(m),
\end{equation}
for each $v \in \PcW,$ $\nu$ in a dense  open subset of $\faPqdc,$ all $X \in \faPq$ and all
$m \in \spX_{P,v,+}.$
It follows from (\refer{e: q of Df as psi}) and (\refer{e: q of g as psi}) that the family
$Df -g \in \cE_P^\hyp(\spX\col \tau)$
satisfies all hypotheses of  Theorem \refer{t: vanishing thm new}.
Therefore, $Df = g.$
\qed

In the rest of this section we shall study the eigenvalues of the endomorphism
$\umu_P(D\col \nu)$ of $\cAtwoP.$
For a start, we collect some facts about the action of $\DX$ on $\cA_2(\spX\col \tau).$

Let $L^2_d(\spX)$ be the discrete part of $L^2(\spX),$ defined as in the beginning
of Section \refer{s: spaces of residue type}.
It follows from \bib{Bfm}, Thm.~1.5, that the space $L^2_d(\spX)$ admits a
decomposition as an orthogonal direct sum
of closed $G$-invariant subspaces on each of which $\DX$ acts by scalars
(in the distribution sense).
Let $\fb$ be a $\Cartan$-stable Cartan subspace of $\fq.$
We denote by $\rmL(\spX, \fb)$ the collection of infinitesimal characters
$\gL \in \fbdc$
for which the associated
character $\gg(\dotvar \col \gL) = \gg_\fb(\dotvar\col \gL)$ of $\DX$
occurs as
a simultaneous eigenvalue in the decomposition mentioned.

The
elements of the $\DX$-module $\cA_2(\spX\col \tau)$
are $\DX$-finite and belong to $L^2_d(\spX) \otimes \Vtau.$
It follows that $\cA_2(\spX\col \tau)$ splits as an algebraic
direct sum of $\DX$-submodules on which the action of $\DX$
is by infinitesimal characters from
$\rmL(\spX, \fb).$
More precisely, for $\gL \in \fbdc$ we put
$$
\cA_2(\spX\col \tau\col \gL) := \{ f\in \cA_2(\spX\col \tau)\mid
 Df = \gg(D\col \gL) f,\;\;\forall D \in \DX
\}.
$$
This space is finite dimensional by \bib{Bfm}, Lemma 3.9.
It
depends on $\gL$ through its class $[\gL]$ in $\fbdc /W(\fb);$
we therefore also denote it with $[\gL]$ in place of $\gL.$
Let $\rmL(\fb, \tau) = \rmL(\spX, \fb, \tau)$ denote the collection
of $\gL \in \fbdc$ for which $\cA_2(\spX\col \tau\col \gL) \neq 0.$
Then $\rmL(\fb,\tau)$ is a
$W(\fb)$-invariant subset of $\rmL(\spX,\fb)$ and we have
the following algebraic direct sum decomposition into joint
eigenspaces for $\DX,$
\begin{equation}
\naam{e: deco cA two in eigenspaces}
\cA_2(\spX\col \tau)= \bigoplus_{\gL \in \rmL(\fb, \tau)/W(\fb)}
\cA_2(\spX\col \tau\col \gL).
\end{equation}
The summands in this decomposition are finite dimensional and
mutually orthogonal
with respect to the inner product from $L^2(\spX\col \tau).$
Moreover,
the decomposition
is finite by Lemma \refer{l: cAt equals cAtwo}. In the chain of reasoning
leading up to Theorem \thmendRT{},
finiteness of the decomposition requires $(G,H)$ to be of residue type,
see Remark \remRT.

\begin{lemma}
\naam{l: conjugacy of rmL fb tau}
Let $\fb_1, \fb_2 \subset \fq$ be two $\Cartan$-stable Cartan subspaces.
Each element $s$ from $W(\fb_2\mid \fb_1),$ which set is non-empty by
Lemma \refer{l: gg for diff css},
maps $\rmL( \fb_1, \tau)$ bijectively onto $\rmL( \fb_2,  \tau).$
\end{lemma}

\proof
This follows by application of Lemma \refer{l: gg for diff css}.
\qed

\begin{lemma}
\naam{l: faP cap v fq}
Let $P \in \allparabs$ and $v \in \NKaq.$ Then $\faP \cap \Ad(v) \fq = \faPq.$
\end{lemma}

\proof
$\faP \cap \Ad(v) \fq = \Ad(v) ( \fa_{v^{-1} P v} \cap \fq) =
\Ad(v) (\fa_{v^{-1}P v} \cap \faq) = \faP \cap \faq = \faPq.$
\qed

Let now $\fb$ be a Cartan subspace of $\fq$ containing $\faq$
and let $v \in \NKaq.$
Then $\fb^v: = \Ad(v) \fb$ is a Cartan subspace of $\Ad(v)\fq,$ which
contains $\faq.$ In particular, $\fb^v$ contains
$\faPq,$ hence is contained in the latter's centralizer
$\fm_{1P}.$
We write $\stfbPv : = \fb^v \cap \fm_P.$  Then
$$
\fb^v = \stfbPv \oplus \faPq.
$$
In view of Lemma \refer{l: faP cap v fq} this is the analogue
of the decomposition (\refer{e: deco fb with star part}) for the
Cartan subspace $\fb^v$ related to symmetric pair
$(\fm_{1P}, \fm_{1P} \cap \Ad(v) \fh).$
The restriction of $\Ad(v)$ to $\fb$ determines an element
of $\Hom(\fb, \fb^v)$ that we denote by $\bar v.$ The restriction $\bar v|_{\faq}$ is
an element of $W.$ The latter set equals $W(\faq\mid \faq),$
by Corollary \refer{c: W a1 a2 in q};
hence, by Lemma \refer{l: conjugacy in Cartan}, applied with $\fb, \faq, \faq$ in place
of $\fd, \fb_1, \fb_2,$ there exists
an element $s \in W(\fb)$ such that $s = \bar v$ on $\faq.$
It readily follows that $\bar v \after s^{-1} \in \Hom(\fb,\fb^v)$
equals the identity on $\faPq,$ hence
maps $\stfbP$ isomorphically onto $\stfbPv.$ Note that this map
maps $W(\stfbP)$-orbits onto $W(\stfbPv)$-orbits. The
induced map from
$
\stfb_{P\iC}^*/W(\stfbP)
$
to
$(\starfb_{P, v})^*_\iC/W(\stfbPv) $
is bijective and depends on $v,$ but is independent of the particular choice of $s.$
Given $\gL \in \stfbPdc,$ we define
\begin{equation}
\naam{e: defi cA two XPv gL}
\cA_2(\spXPv\col \tau_P\col \gL) := \cA_2(\spXPv\col \tau_P\col \bar v \after s^{-1}\gL).
\end{equation}
Moreover, we define $\rmL_{P,v}(\fb, \tau)$ to be the set of $\gL \in \stfbPdc$
for which the above space is non-trivial. Then
\begin{equation}
\naam{e: image rmLPv under bar v s minus one}
\bar v\after s^{-1} \; \rmL_{P,v}(\fb, \tau) = \rmL(\spXPv,  \stfbPv, \tau_P).
\end{equation}
Thus, $\rmL_{P,v}(\fb, \tau)$ is a
$W(\stfbP)$-invariant
subset of $\stfbPdc.$

\begin{cor}\restypetwo
\naam{c: deco cAPtwo in eigenspaces}
Let  $P\in \allparabs$
and let $\fb \subset \fq$
be a $\Cartan$-stable Cartan subspace containing $\faq.$
Then
\begin{equation}
\naam{e: deco cA two P in eigenspaces}
\cA_{2,P}
= \bigoplus_{v \in \PcW}\; \bigoplus_{\gL \in \rmL_{P,v}(\fb, \tau)/W(\stfbP)}
\rmi_{P,v}\cA_2(\spXPv \col \tau_P \col \gL),
\end{equation}
with a finite orthogonal direct sum of finite dimensional spaces.
If $D \in \DX$ and $\nu \in \faPqdc,$ then for every
$v \in \PcW$ and $\gL \in \rmL_{P,v}(\fb, \tau),$
$$
\umu_{P,v}(D\col \nu) =
\gg_\fb(D\col \gL + \nu) I \text{on} \rmi_{P,v} \cA_2(\spXPv\col \tau_P\col \gL).
$$
\end{cor}

\proof
We give the proof under the assumption that $P$ is of residue type,
see Remark \remRT.
By (\refer{e: defi cAtwoP}) the space $\cA_{2,P}$ is the orthogonal direct sum
of the spaces $\cA_2(\spXPv\col \tau_P),$ as $v \in \PcW.$ Fix $v \in \PcW.$
By the assumption on $P,$ the pair
$(M_P, M_P \cap vPv^{-1})$
is of residue type,
hence $\cA_2(\spXPv\col \tau_P)$ is finite dimensional and by
(\refer{e: deco cA two in eigenspaces})
it is the orthogonal direct sum of the spaces
$\cA_2(\spXPv\col \tau_P\col\gL'),$
with $\gL'\in \rmL(\spXPv, \starfb_{P,v}, \tau_P)/W(\starfb_{P,v}).$
It now follows from (\refer{e: defi cA two XPv gL}) and
(\refer{e: image rmLPv under bar v s minus one}) that $\cA_2(\spXPv\col \tau_P)$
is the orthogonal
direct sum of the spaces $\cA_2(\spXPv\col \tau_P\col\gL),$ for
$\gL \in \rmL_{P,v}(\fb,\tau)/W(\stfb_P);$ moreover, the sum is finite and
the summands are finite dimensional. This establishes
(\refer{e: deco cA two P in eigenspaces}), with the asserted properties.

Let $\gL \in \rmL_{P,v}(\fb, \tau).$
Then by (\refer{e: image rmLPv under bar v s minus one}),
$\gL' := \bar v s^{-1}  \gL$ belongs to
$\rmL(\spX_{P,v}, \stfb_{P,v}, \tau_P).$
Let now $\psi \in \cA_2(\spXPv\col \tau_P\col \gL).$
Then, writing $D^v = \Ad(v)D$ for $D \in \DX,$
\begin{eqnarray*}
\umu_{P,v}(D\col \nu) \psi
&=& \mu_P^v(D^v\col \nu) \psi \\
&=& \gg_{\starfb_{P,v}}^{\spXPv}(\mu_P^v(D^v \col \nu)\col \gL') \psi \\
&=&  \gg_{\fb^v}^{\spX_{1P,v}}(\mu_P^v(D^v)\col \gL' + \nu ) \psi.
\end{eqnarray*}
In the last equation we have used that
$\gg_{\fb^v}^{\spX_{1P, v}} = \gg_{\starfb_{P,v}}^{\spXPv} \otimes I,$
in accordance with (\refer{e: deco DspXonePv}).
Combining (\refer{e: composition gg and mu}) for the triple $(G/vHv^{-1}, \fb^v, P)$
in place of
$(G/H, \fb, P)$ with (\refer{e: defi mu P v}), we obtain that
\begin{eqnarray*}
\umu_{P,v}(D\col \nu) \psi &=& \gg_{\fb^v}^{G/vHv^{-1}}(D^v \col \gL' + \nu ) \psi
\\
&=&
\gg_\fb(D\col \Ad(v)^{-1}(\gL' + \nu ))\psi \\
&=&
\gg_\fb(D\col s\after \Ad(v)^{-1}(\gL' + \nu ))\psi \\
&=&
\gg_\fb(D\col \gL + \nu )\psi.
\end{eqnarray*}
\qed

We define
$\rmL_P( \fb, \tau) \subset \stfbPdc$ to be the union of
the sets $ \rmL_{P,v}(\fb, \tau),$ for $v \in \PcW.$
Moreover, for $\gL$ in this union, we put
$$
\cA_{2,P}(\gL): = \bigoplus_{v \in \PcW}\;\;
\rmi_{P,v}\cA_2(\spXPv \col \tau_P \col \gL).
$$

\begin{cor}\restypetwo
\naam{c: normEis as eigenfunction}
Let
$P\in \allparabs.$
Then
$$
\cA_{2,P} = \bigoplus_{\gL \in \rmL_P(\fb, \tau)/W(\stfbP)} \cA_{2,P}(\gL).
$$
Moreover, if  $\gL \in \rmL_P(\fb, \tau)$ and
$\psi \in \cA_{2,P}(\gL),$ then, for every $D \in \DX,$
$$
D \nE (P\col \nu) \psi =
\gg_\fb(D\col \gL + \nu) \nE(P\col \nu) \psi,
$$
as a meromorphic $\Ci(\spX\col \tau)$-valued identity in $\nu \in \faPqdc.$
\end{cor}

\proof
We give the proof under the assumption that
$P$ is of residue type, see Remark \remRT.
The result follows from Corollary \refer{c: deco cAPtwo in eigenspaces}
combined with Lemma
\refer{l: action of DX on nE}.
\qed

We end this section with a description of the action of $\DX$ on the dualized Eisenstein
integral.
For $D \in \DX$ we define the formal adjoint $D^* \in \DX$
by
\begin{equation}
\naam{e: formal adjoint of D}
\hinp{Df}{g} = \hinp{f}{D^*g},
\end{equation}
for all $f,g \in \Cci(\spX);$ here
$\hinp{\dotvar}{\dotvar}$ denotes the inner product from $L^2(\spX).$
The canonical anti-automorphism $X \mapsto X^\vee$ of $U(\fg)$
induces an anti-automorphism of $U(\fg)^H/U(\fg)^H \cap U(\fg)\fh \simeq \DX,$
which we also denote by $D \mapsto D^\vee.$ If $D \in \DX,$
let $\bar D$ be its complex conjugate, i.e. the differential operator with complex
conjugate coefficients. Then $D^* = \bar D^\vee,$ for every $D \in \DX.$

We recall from \bib{Bas}, Lemma 7.2,
that $D$ restricts to a continuous linear endomorphism
of $\cC(X);$ by density of $\Cci(\spX)$ in $\cC(\spX)$ it follows
that (\refer{e: formal adjoint of D}) also holds for all
$f,g \in \cC(\spX).$

\begin{lemma}
\naam{l: star of mu P v}
Let $P \in \allparabs, $ $v \in \NKaq$ and $D \in \DX.$ Then
\begin{equation}
\naam{e: conjugate of mu P D}
\mu_{P,v}(D^*)  = \mu_{P,v}(D)^*.
\end{equation}
\end{lemma}

\proof
We note that $\mu_{P,e} = \mu_P;$ hence, for $v =e,$ the result
follows by
the same argument as in
\bib{Bps2}, proof of Lemma 19.3. For general $v$ the result follows
by application of (\refer{e: defi mu P v}).
\qed

\begin{lemma}
\restypetwo
\naam{l: conjugate of umuPD}
Let $P \in \allparabs.$
Then, for every $D \in \DX$ and all $\nu \in \faPqdc,$
\begin{equation}
\naam{e: conjugate of umu}
\umu_P(D\col \nu)^* = \umu_P(D^*\col -\bar \nu).
\end{equation}
\end{lemma}

\proof
We give  the proof under the assumption
that $P$ is of residue type, see Remark \remRT.
Let $v\in \PcW.$
The decomposition $\spXonePv \simeq \spXPv \times \APq$
induces an isomorphism
$\D(\spXonePv) \simeq
\D(\spXPv) \otimes S(\faPq).$ Accordingly $(u \otimes p)^* = u^* \otimes p^*,$
for all $u \in \D(\spXPv)$ and $p \in S(\faPq).$ Moreover, $p^*(\nu) =
\overline { p(-\bar \nu)},$ for $\nu \in \faPqdc.$
Hence, $(u \otimes p)^*(\nu) = \overline p(-\bar \nu) u^* = [u\otimes p(-\bar \nu)]^*$
and we see that $u^*(\nu) = u(-\bar \nu)^*$ for $u \in \D(\spXonePv)$
and $\nu \in \faPqdc.$ Applying this to (\refer{e: conjugate of mu P D})
it follows that for $D \in \DX$ and $\nu \in \faPqdc$  we have
$
\mu_{P,v}(D \col \nu)^* = \mu_{P,v}(D^* \col -\bar\nu).
$
By the argument in the text preceding Lemma \refer{l: star of mu P v}, applied to $\spXPv$ in place
of $\spX,$ we infer that
$$
\hinp{\mu_{P,v}(D \col \nu) f}{g} =  \hinp{f}{\mu_{P,v}(D^* \col -\bar\nu) g},
$$
for all $f,g \in \cC(\spXPv\col \tau_P).$ Here $\hinp{\dotvar}{\dotvar}$
denotes the $L^2$-inner product.
In particular, the equation holds for $f,g$
in the subspace $\cA_2(\spX\col \tau_P),$
which is finite dimensional, since $P$ is of residue type.
Hence,
$\umu_{P,v}(D \col \nu)^* = \umu_{P,v}(D^* \col -\bar\nu).$
By orthogonality of the direct sum decomposition
in (\refer{e: defi cAtwoP}),
the result follows by summation over $v \in \PcW.$
\qed

\begin{lemma}\restypetwo
\naam{l: D on dE}
Let $P \in \allparabs.$
Then for every $D \in \DX,$
$$
D\dE(P\col \nu) = \umu_P(D^\vee \col \nu) \dE(P\col \nu),
$$
as a meromorphic identity in $\nu \in\faPqdc.$
\end{lemma}

\proof
We
give the proof under the assumption that $P$ is of residue type,
see Remark \remRT.
By linearity, we may assume that $D$ is real.
It then follows from the definition of the dual Eisenstein integral, see Definition
\refer{d: dual Eisenstein new},
combined with Lemma \refer{l: action of DX on nE}, that
$$
D \dE(P\col \nu) = \umu_P( D \col -\bar \nu)^* \dE(P\col \nu).
$$
The lemma now follows by application of Lemma \refer{l: conjugate of umuPD},
in view of the fact that
$D^* = D^\vee.$
\hbox{\ }\qed

\section{Uniform tempered estimates}\eqnreset
In this section we present straightforward generalizations
of results of \bib{Bps2}, Sect.\ 18, to a setting
involving families $\{f_\nu\}$ of eigenfunctions on $\spX,$
with holomorphic dependence on a parameter $\nu \in \faQqdc,$ where $Q \in \allparabs.$
A similar generalization has been given in \bib{D1n}, Sect.~9. The generalized results
allow us to sharpen uniformly moderate estimates of type (\refer{e: estimate for nE P}) to
uniform tempered estimates. In particular, we obtain
estimates of the latter type for the normalized Eisenstein integral.

We fix $Q \in \allparabs,$ a $\Cartan$-stable Cartan subspace $\fb$ of $\fq$ containing
$\faq$ and an element
$\gL \in \stbQdc,$ cf.~(\refer{e: deco fb with star part}).
For $\geps > 0,$ we put
$$
\faQqd(\geps): = \{ X \in \faQqdc \mid \;\;|\Re X| < \geps \}.
$$
The closure of this set is denoted by $\bar\faQqd(\geps).$
 \begin{defi}
Let $\geps > 0.$
We define $\cE(Q\col \gL\col \geps) = \cE(\spX\col Q\col \gL\col \geps)$ to be the space
of smooth functions $f: \faQqd(\geps) \times \spX \to \C$
satisfying the following conditions.
\begin{enumerate}
\itema
The function $f$ is holomorphic in its first variable.
\itemb
For every $\nu \in \faQqd(\geps),$ the function $f_\nu: x \mapsto f(\nu,x)$
satisfies the following system of differential equations
$$
D f_\nu = \gg(D\col \gL + \nu) f_\nu, \qquad (D \in \DX).
$$
\end{enumerate}
\end{defi}
Note that in this definition it is not required that $f$ is $K$-finite, or spherical,
from the left.
We
also have the following analogue of the space $\Mer(\gL, \geps)$ defined
in \bib{Bps2}, Sect.~18.
\begin{defi}
Let $\geps > 0.$
A function $f \in \cE(Q\col \gL\col  \geps)$ is called uniformly moderate of
exponential rate $r\geq 0$ if for every $u \in U(\fg)$
there exist constants
$n \in \N$ and $C>0$ such that
$$
\|L_u f_\nu(x)\| \leq C (1 + |\nu|)^n e^{r \lengthX(x)},
$$
for all $\nu \in \faQqd(\geps)$ and $x \in \spX.$
The space of all such functions is denoted by $\cE^\um(Q\col \gL\col \geps\col r).$
\end{defi}

\begin{lemma}\restypetwo
\naam{l: normalized Eis in cEum}
Let $Q \in \allparabs$ and let $\geps >0.$ There exists
a polynomial function $p \in \Pi_{\gS_r(Q),\R}(\faQqd)$ such that
the $\Ci(\spX) \otimes \Hom(\cA_{2,Q}, \Vtau)$-valued meromorphic function
$\nu \mapsto p(\nu) \nE(Q\col \nu)$ is regular on
$\bar\faQqd(\geps)$ and such that the following holds.
There exists a constant $r>0$ such that for every
$\gL \in \rmL_Q(\fb,\tau),$ $\psi \in \cA_{2,Q}(\gL)$ and $\eta \in \Vtaud,$
the family $f: (\nu, x) \mapsto \eta(p(\nu) \nE(Q\col \nu\col x))$
belongs to $\cE^\um(Q\col \gL \col \geps \col r).$
\end{lemma}

\proof
We give
the proof under the assumption
that $Q$ is of residue
type,
see Remark \remRT.
Let $R>0$ be such that $\faQqd(\geps) \subset \faQqd(\bar Q,R).$
Then by Proposition
\refer{p: init estimate normalized Eis} there exists a
polynomial function $p \in \Pi_{\gSr(Q)}(\faQqd)$ such that
 the $\Ci(\spX) \otimes \Hom(\cA_{2,Q},\Vtau)$-valued
meromorphic function $F: \nu \mapsto p(\nu) \nE(Q\col \nu) $ is holomorphic
on a neighborhood of $\bar\faQqd(\bar Q, R).$ Moreover,
there exists $r' >0$ and for every $u \in U(\fg)$
constants $n \in \N$
and $C >0$ such that
$$
\| L_uF_\nu(x)\| \leq C(1 + |\nu|)^n e^{(r' + |\Re\nu|)\lspX(x)},
$$
for $x \in \spX,$ $\nu \in \bar\faQqd(\bar Q,R).$ Put $r =r' + \geps.$
Then it follows that $F$ is holomorphic on a neighborhood of
 $\bar\faQqd(\geps)$ and satisfies the
estimates
\begin{equation}
\naam{e: uniform estimate F}
\| L_uF_\nu(x)\| \leq C(1 + |\nu|)^n e^{r\lspX(x)}
\end{equation}
for $x \in \spX$ and $\nu \in \faQqd(\geps).$ Let $f$ be defined as in the lemma.
Then $L_u f_\nu(x) = \eta(L_u F_\nu(x)\psi).$ Hence, it follows from the above
and from Corollary  \refer{c: normEis as eigenfunction} that
$f \in \cE(Q\col \gL\col \geps).$ Finally, it follows from the estimates
(\refer{e: uniform estimate F}) that $f\in \cE^\um(Q\col \gL \col \geps \col r).$
\qed

We also have the following obvious generalization of the notion of
uniformly tempered families; see \bib{Bps2}, Sect.~18.
For $\nu \in \faQqdc$
and $x \in \spX$ we put
$$
|(\nu, x)|: = (1 + |\nu|)(1 + \lengthX(x)).
$$
\medbreak

\begin{defi}
\naam{d: space uniform tempered functions}
Let $\geps > 0.$
A function $f \in \cE(Q\col \gL\col \geps)$ is called uniformly tempered
of scale $s$ if for every $u \in U(\fg)$ there exist constants
$n \in \N$ and $C>0$ such that
$$
|L_u f_\nu(x)| \leq C |(\nu, x)|^n \Theta(x) e^{s|\Re \nu|\lengthX(x)},
$$
for all $\nu \in \faQqd(\geps)$ and all $x \in \spX.$ The space of all
such functions is denoted by
$\cT(Q\col \gL\col \geps\col s).$
\end{defi}

If $f \in \cE^\um(Q\col\gL\col \geps\col r),$ then for every $\nu \in \faQqdc$
the function $f_\nu$ belongs to the space $\cE^\infty_{\gL + \nu, *}(\spX),$
defined in \bib{Bps2}, p.\ 392, see also p.\ 387. If $g$ is any function in
the latter space, then, viewed as a function on $G,$
it has an asymptotic expansion along
every parabolic subgroup $P \in \allparabs$ of the form
$$
g(x\exp tX) \sim \sum_{\xi \in Z - \N\DrP} p_{\xi}(P\asmid g, x, X) e^{t\xi(X)},
$$
as $t \to \infty,$ for $x \in G$ and $X \in \faPqp.$
Here $Z$ is a finite subset of $\faPqdc$ and there exists a $d \in \N$
such that the $p_\xi(P\asmid g)$ are smooth functions $G \to P_d(\faPq),$
for all $\xi.$
Moreover, the functions $p_\xi(P\asmid g)$ are uniquely determined,
see \bib{Bps2}, Theorem 12.8. Accordingly, we may define the set
of exponents of $g$ along the parabolic subgroup $P$ by
\begin{equation}
\naam{e: set of exponents along P}
\Exp(P\asmid g): = \{\xi \in Z -\N\DrP \mid\;\;\; p_\xi(P\asmid g) \neq 0\}.
\end{equation}
We define the partial ordering $\preceq_P$ on $\faPqdc$ by
$$
\gl \preceq_P \mu \iff
\mu - \gl \in \DrP, \qquad (\gl, \mu \in \faPqdc).
$$ The $\preceq_P$-maximal elements in the set
(\refer{e: set of exponents along P}) are called the leading exponents
of $g$ along $P.$ The set of these leading exponents is denoted
by $\Exp_\rmL(P\asmid g).$

\begin{rem}
The above notions of asymptotic coefficients and exponents are related
to the similar notions introduced in \S~\refer{s: vanishing theorem}, as follows.

Let $f\in \cA(\spX\col \tau)$ and assume that every vector component
$\eta \after f,$ for $\eta \in V_\tau^*,$ belongs to $\cE^\infty_{\gL + \nu, *}(\spX).$
For $P \in \allparabs,$ let $\Exp(f\mid P)$ denote the union of the sets $\Exp(\eta\after f\mid P),$
for $\eta \in V_\tau^*;$ by sphericality of $f$ this union equals the union with index $\eta$
ranging over any generating subset of the $K$-module $V_\tau^*.$
If $u \in \NKaq,$ then it readily follows from the definitions
that $\Exp(P,u\mid f) \subset \Exp(P \mid f).$ Moreover, by uniqueness
of asymptotics we have, for
$\xi \in \Exp(P, u\mid f),$ that
$$
\eta(q_\xi(P,u\mid f, X, m)) = p_\xi(P\mid \eta\after f, m u, X),
\qquad (m\in M_P,\; X \in \faPq),
$$
for all $\eta \in V_\tau^*.$
\end{rem}

\begin{lemma}
\naam{l: leading exponent and Wfb}
Let
$\nu \in \faQqdc$ and assume that  $g \in \cE_{\gL + \nu, *}^\infty(\spX).$
Let $P\in \minparabs.$ Then for every $\xi \in \Exp_L(P\asmid g)$
there exists a $s \in W(\fb)$ such that $\xi + \rho_P = s(\nu + \gL)|_{\faq}.$
\end{lemma}

\proof
We recall that $\faq \subset \fb.$
Let $\gSp(\fb)$ be a choice of positive roots for $\gS(\fb)$ that
is compatible with $\gS(P).$ Let $\fg_\iC^+$ be the associated
sum of the positive root spaces and let $\fm_\iC^+$ be its
intersection with $\fm_\iC.$ Let $\gd := \frac12 \tr [\ad(\dotvar)|_{\fg_\iC^+}] \in \fbdc$
and let $\gd_\iiM:= \frac12 \tr [\ad(\dotvar)|_{\fm_\iC^+}] \in i\fbkd.$
Then $\gd = \gd_\iiM + \rho_P.$

Let $\xi$ be a leading exponent along $P.$ Then by
\bib{Bps2}, Cor.\ 13.3 and Lemma 13.1, the function $\gf\in \Ci(M_1)$
defined by
$$
\gf(ma) = a^\xi \, p_\xi(P\asmid g, m , \log a),\qquad (m \in M_\gs, \, a \in \Aq),
$$
is right $M_\gs \cap H$-invariant and satisfies the following system of
differential equations
$$
\mu_P'(D)\gf = \gg(D\col \gL + \nu)\gf, \qquad (D \in \DX).
$$
Here $\mu_P'$ is defined as in \bib{Bps2}, Sect.~2.
Now $M_\gs/M_\gs \cap H \simeq M_{10}/M_{10} \cap H,$ naturally, so that
$\gf$ may be viewed as a function in $\Ci(\spX_{10}).$
By (\refer{e: deco DspXonePv}) with $P = P_0$ and $v =e$ we have
$\D(\spX_{10}) \simeq \D(\spX_0) \otimes S(\faq).$
Since $p_\xi$ is polynomial in $\log a,$ the
second component of the tensor product acts on $\gf$ with a single
generalized eigenvalue $u \mapsto u(\xi).$
On the other hand, we recall from \bib{Bps2}, Lemma 4.8,
that the action of $\D(\spX_0)$ on $C^\infty(\spX_0)_{\KM}$ allows
a simultaneous diagonalization with eigenvalues of the form
$D \mapsto {\gg^{\spX_0}}(D\col \gL_0 + \gd_\iiM),$ with $\gL_0 \in i\fbkd.$
It follows that there exists a $\gL_0 \in i\fbkd$ such that
$$
{\gg^{\spX_0}}(\mu_P'(D\col \xi)\col \gL_0  + \gd_\iiM) =
\gg(D\col \gL + \nu), \qquad (D \in \DX).
$$
The expression on the left-hand side of this expression can
be rewritten as $\gg_P(D\col \gL_0 + \xi +\gd_0 + \rho_P) = \gg(D \col \gL_0 +\xi + \gd),$
from which we conclude that $\gL_0 + \xi + \gd \in W(\fb)(\gL + \nu).$
Since $(\gL_0 +\gd)|_{\faq} = \rho_P,$ it follows that $\xi +\rho_P  = s(\nu + \gL)|_{\faq},$
for some $s \in W(\fb).$
\qed

We can now generalize \bib{Bps2}, Theorem 18.3.
For an appropriate formulation we need the following definition.

\begin{defi}
\naam{d: tempered exponents}
We  say that the exponents of a family
$f \in \cE^\um(Q\col \gL\col \geps\col r)$
are tempered along a minimal $\gs$-parabolic subgroup $P \in \minparabs$
if for every $\nu \in \faQqd(\geps)$ the set of exponents
$\Exp(P\asmid f_\nu)$ satisfies the following
condition. For every $\xi \in \Exp(P\asmid f_\nu),$ there exists
a $s \in W(\fb)$ such that
\begin{enumerate}
\itema
$\Re(s\gL)\leq 0$ on $\faqp(P),$
\itemb
$
\xi \in s(\nu +\gL)|_{\faq} - \rho_P - \N \DP.
$
\end{enumerate}
We denote by $\cE^\um_\iT(Q\col \gL\col \geps\col r)$
the space of functions $f \in \cE^\um(Q\col \gL\col \geps\col r),$ such
that for every $P\in \minparabs$ the exponents of $f$ along $P$
are tempered.
\end{defi}

\begin{rem}
\naam{r: tempered exponents for Q minimal}
If $Q$ is a minimal $\gs$-parabolic subgroup, then it follows
by application of \bib{Bps2}, Thm.~13.7,
that $\cE^\um(Q\col\gL\col \geps\col r) = \cE^\um_\iT (Q\col\gL\col \geps\col r).$
\end{rem}

\begin{thm}
\naam{t: um implies utemp}
Let $Q \in \allparabs$ and let $r>0.$
Then there exists a $s >0$ such that for sufficiently
small $\geps>0,$
$$
\cE^\um_\iT(Q\col\gL\col \geps\col r) \subset \cT(Q\col \gL\col \geps\col s).
$$
\end{thm}

\proof
The proof is a straightforward, but somewhat tedious, adaptation of the proof of
\bib{Bps2}, Theorem 18.3,
with trivial alterations because of the change of the parameter set.
Conditions (a) and (b)
of Definition \refer{d: tempered exponents} are to be used in place of \bib{Bps2}, Theorem 13.7,
see the proof of \bib{Bps2}, Proposition 18.14, to keep track of the exponents occurring
in the asymptotic expansions considered. If $Q$ is minimal, then the mentioned
Theorem 13.7 implies conditions (a) and (b) for any family
$f \in \cE^\um(Q\col\gL\col \geps\col r);$
see Remark \refer{r: tempered exponents for Q minimal}.
\qed

\begin{rem}
Another version of Theorem \refer{t: um implies utemp} is given by
\bib{D1n},  Thm.~3. However, in that paper the
requirement on the exponents in Definition \refer{d: tempered exponents}
is replaced
by the requirement that the function $f_\nu$ is tempered for
every $\nu \in i\faqd.$
By an additional argument it is then shown that this requirement is equivalent
to the one of Definition \refer{d: tempered exponents}, see
\bib{D1n},  Lemma 23.
We shall not need this result,
since by Proposition \refer{p: exponents normalized Eis}
the needed
information on the exponents is known for the normalized Eisenstein integrals
to which Theorem \refer{t: um implies utemp} will be applied.
\end{rem}

\begin{defi}
Let $Q \in \allparabs,$ $ \gL\in \stbQdc,$ $\geps > 0$ and $s >0.$
Then by $\cT(Q,\tau, \gL, \geps, s)$ we denote the space
of smooth functions $f: \faQqd(\geps) \times \spX \to \Vtau$
such that
\begin{enumerate}
\itema
for every $\eta \in \Vtaud$ the family $\eta\after f: (\nu, x) \mapsto \eta(f(\nu,x))$
belongs to $\cT(Q,\gL, \geps , s);$%
\itemb
$f_\nu$ is $\tau$-spherical for every $\nu \in \faQqd(\geps).$
\end{enumerate}
\end{defi}

\begin{thm}\restypetwo
\naam{t: norm Eis uniform tempered}
Let $Q \in \allparabs.$
There exists a polynomial function
$p \in \Pi_{\gSr(Q), \R}(\faQqd)$ and  constants $s > 0$ and $\geps > 0$
such that the meromorphic $\Ci(\spX) \otimes \Hom(\cA_{2,Q}, \Vtau)$-valued
function $\nu \mapsto p(\nu)\nE(Q\col \nu)$ is holomorphic on $\faQqd(\geps),$
and such that the following holds.
For each $\gL \in \rmL_Q(\fb, \tau)$ and every $\psi \in \cAtwoQ(\gL)$
the family $f: (\nu, x) \mapsto p(\nu)\nE(Q\col \nu \col x)\psi$
belongs to $\cT(Q,\tau, \gL, \geps, s).$
\end{thm}

\begin{rem}
For
$Q$ minimal, this result is due
\bib{Bps2}, Thm.~19.2, in view of
\bib{BSft}, Eqn.~(52). For general $Q,$ a similar result for an unnormalized version
of the Eisenstein integral is due to
\bib{D1n}, Thm.~4.
\end{rem}

\proof
We  give the proof under the assumption
that $Q$
is of residue type, see
Remark \remRT.
Fix $\geps >0.$ There exist $p \in \Pi_{\gSr(Q),\R}(\faQqd)$ and $r>0$
as in Lemma \refer{l: normalized Eis in cEum}.
Fix $\gL \in \rmL_Q(\fb, \tau)$ and
$\psi \in \cA_{2,Q}(\gL).$
Define $f: (\nu,x) \mapsto p(\nu)\nE(Q\col \nu \col x)\psi.$
Let $\eta \in \Vtaud$ and define $F: (\nu, x) \mapsto \eta(f(\nu,x)).$
Then by finite dimensionality of $\cA_{2,Q}$ and $\Vtau$ it suffices
to show that there exist $\geps'> 0$ and $s >0$ such that
$F \in \cT(Q\col \gL\col  \geps'\col  s).$

In view of Theorem \refer{t: um implies utemp}
it suffices to show that $F \in \cE_{\rmT}^\um(Q\col \gL\col \geps \col r).$
In view of Lemma \refer{l: normalized Eis in cEum} the function
$F$ belongs to $\cE^\um(Q\col \gL\col \geps \col r).$ Let $P \in \minparabs.$
Then it remains
to be verified that the exponents of $F$ along $P$ are tempered in the
sense of Definition \refer{d: tempered exponents}.

There exists a $v \in \NKaq$ such that $P_1:= v^{-1}Pv \subset Q.$
The meromorphic $\Ci(\spX\col \tau)$-valued function
$\nu \mapsto f_\nu$ is regular on $\faQqd(\geps).$ Moreover,
from Proposition \refer{p: exponents normalized Eis}
and \bib{BSanfam}, Lemma 3.6,
it follows that, for $\nu \in \reg f,$
$$
\Exp(P,e\mid f_\nu) = v\Exp(P_1, v^{-1} \mid f_\nu)
\subset v W^{P_1|Q}(\nu + \gL(P_1|Q)) - \rho_P - \N \gD(P).
$$
Thus, let $\nu_0 \in \faQqd(\geps)$ be fixed, and let
$\xi \in \Exp(P\mid F_{\nu_0}).$
Then we may select $s \in v W^{P_1|Q}$ and $\xi_0 \in - s\gL(P_1|Q) + \N \gD(P)$
such that $\xi = s\nu_0- \rho_P - \xi_0.$ Since $f \in \cE_Q(\spX\col \tau),$
see Definition \refer{d: dual Eisenstein new}, it follows from
Definitions \refer{d: cEhypQallX}, \refer{d: cEhypQgdglob} and \refer{d: defi cEhyp Q Y gd}
that $f \in C^{\ep, \hyp}_{Q,Y}(\spXp\col \tau),$
for a suitable finite subset $Y \subset \staQqdc.$ By Definition \refer{d: Cephyp}
and \bib{BSanfam}, Lemma 6.2, it follows that
$s \nu - \rho_P - \xi_0 \in \Exp(P,e\mid f_\nu),$ for $\nu$ in an open
dense subset of $\faQqdc.$

Let $\xi_1$  be a $\preceq_P$-minimal
element in $- s\gL(P_1|Q) + \N \gD(P)$
with the property that $\xi_1 \preceq_P \xi_0$ and that
$s\nu - \xi_1 - \rho_P \in \Exp(P,e\mid f_\nu)$ for $\nu$ in an open
dense subset of $\faQqdc.$
Then for  $\nu$ in an open dense subset $\Omega$ of $\faQqdc,$
the element $s\nu - \xi_1 - \rho_P$
is a leading exponent of $f_\nu$ along $(P,e).$
By Lemma \refer{l: leading exponent and Wfb}  it follows that
$$
s \nu - \xi_1 - \rho_P \in W(\fb)(\nu +\gL)|_{\faq} - \rho_P,
$$
for $\nu \in \Omega.$
This implies in turn that there exists $t \in W(\fb)$ such
that
$s\nu - \xi_1  = t(\nu + \gL)|_{\faq},$ for all $\nu \in \faQqdc.$ Hence,
$s\nu  = t\nu|_{\faq}$ for all $\nu \in \faQqdc$ and
$-\xi_1 = t\gL|_{\faq}.$ Now $- \xi_1 \in s \gL(P_1 | Q) - \N\DP \subset
- s (\R_+ \gD_Q(P_1)) - \N \DP \subset  - \R_+ \gD(P),$
hence $\Re(t\gL)|_{\faq} = - \xi_1 \leq 0$ on $\faqp(P).$ We complete the proof by
observing that
$$
\xi = s\nu - \rho_P - \xi_0 = s \nu - \xi_1 -\rho_P - (\xi_0 - \xi_1)
\in t(\nu + \gL)|_{\faq} - \rho_P - \N\gD(P).
$$
\qed

\section{Infinitesimal characters of the discrete series}\eqnreset
\naam{s: infinitesimal characters}
In this section we describe a restriction on the set
$\rmL(\spX, \fb)$
of $\DX$-characters
of the discrete series, see the text before
(\refer{e: deco cA two in eigenspaces}). The main result is due to T.\ Oshima
and T.\ Matsuki, \bib{OMds}.

Let $\fb\subset \fq$ be a $\Cartan$-stable Cartan subalgebra.
If $\gL \in \fbdc$ then by $I_\gL$
we denote the kernel of $\gg_\fb(\dotvar \col \gL)$ in $\DX.$
We denote by $\cC(\spX\col \gL)$ the space of $L^2$-Schwartz functions
on $\spX$ annihilated by $I_\gL.$ If $\cC(\spX\col \gL)$ is non-trivial,
then it contains a non-trivial $K$-finite function $f.$
By a well known result of Harish-Chandra, the closed $G$-span of $f$ in
$L^2(\spX)$  is a subrepresentation of finite length;
see \bib{Var}, p.~312, Thm.~12 and \bib{Wal1}, p.~112, Thm.~4.2.1.
Therefore, the mentioned closed $G$-span is contained
in $L^2_d(\spX)$ and we deduce that $\gL \in \rmL( \spX,\fb).$
Conversely, if $\gL \in \rmL( \spX,\fb),$ then there exists
a non-trivial $K$-finite function $f \in L^2_d(\spX)$ that
is annihilated by $I_\gL.$ {}From \bib{Bas}, Thm.~7.3, it follows
that $f$ belongs to $\cC(\spX\col \gL)$ and we see that
the latter space is non-trivial.
We conclude that
\begin{equation}
\naam{e: characterization rmL}
\rmL(\spX,\fb)= \{ \gL \in \fbdc \mid  \cC(\spX\col \gL) \neq 0\}.
\end{equation}

\begin{thm}
\naam{t: infinitesimal characters L two d real and regular new}
Assume that the space $L^2_d(\spX)$ is non-trivial.
Then there exists a compact Cartan subspace $\ft \subset \fq.$
Moreover, each $\gL\in \rmL( \spX,\ft)$ belongs
to $i\ft^*$ and is regular with respect to $\Sigma(\ft).$
\end{thm}

\begin{rem}
\naam{r: rank condition due to OM}
This result, which plays a crucial role in the description
of the constant term of the normalized Eisenstein integral
in Section \refer{s: the constant term},
is essentially due to T.\ Oshima and M.\ Matsuki, \bib{OMds}.
However, we have to be a bit careful here, since in our situation
$G$ is assumed to be of Harish-Chandra's class, whereas
in \bib{OMds} it is assumed that $G$ is semisimple.
\end{rem}

\proof
Fix a Cartan subspace $\fb \subset \fq$ that is fundamental, i.e.,
its compact part $\fb_\ik = \fb \cap \fk$ is of maximal dimension.
Then the assumption that $L^2_d(\spX)$ is non-trivial is equivalent
to the assumption that $\rmL( \spX,\fb )$ is
non-empty. We must show that under this assumption $\fb$ is compact,
and all elements of $\rmL( \spX,\fb )$ belong to $i\fb^*$ and
are regular.

Let $\spX^\circ = G_e / G_e\cap H$ be  the connected component
of the origin in $\spX.$
If $\gL \in \fbdc$ then restriction defines a linear map  $r: \cC(\spX\col \gL)
\to \cC(\spX^\circ\col \gL).$ Conversely, extension by zero defines
a linear embedding $j:  \cC(\spX^\circ\col \gL) \to \cC(\spX\col \gL).$
Now $r\after j = I,$ hence $r$ is surjective.
If the space $\cC(\spX\col \gL)$ is non-trivial, then by $G$-invariance
it follows that $r$ is non-zero, hence its image is non-trivial.
On the other hand, if $\cC(\spX^\circ \col \gL)$ is non-trivial, then
$\cC(\spX \col \gL)$ is non-trivial, by injectivity of $j.$
Thus, from  (\refer{e: characterization rmL}) we see   that
$\rmL(\spX,\fb) = \rmL(\spX^\circ,\fb).$
Therefore, we may as well assume that $G$
is connected.

Let $\fa_{\gD\iq}$ be the intersection in $\faq$ of the root spaces
$\ker \ga,$ $\ga \in \Sigma.$ This space is central in $\fg.$
Hence,  $\fa_{\gD\iq} \subset \fb$ and the group $A_{\gD\iq}: =
\exp \fa_{\gD\iq}$ is central in $G.$

The algebra $U(\fa_{\gD\iq})$
naturally embeds into $\DX$ and into $I(\fb);$ accordingly,
$\gg$ restricts to the identity on $U(\fa_{\gD\iq}).$
Let $\gL \in \rmL(\spX,\fb)$ and let $f$ be a non-trivial
function in $\cC(\spX\col \gL).$ Then it follows that
$R_X f = \gL(X) f$ for all $X \in \fa_{\gD\iq}.$ Let $\gL_0: = \gL|\fa_{\gD\iq}.$
Then it follows that $f(ax) = a^{\gL_0}f(x)$ for all $x \in \spX$
and $a \in A_{\gD\iq}.$ Since $f$ is a non-trivial Schwartz function,
this implies that $\fa_{\gD\iq} = 0.$

Let $\fc$ be the center of $\fg.$ Then it follows that $\fc_\iq:= \fc \cap \fq$
is contained in $\fb \cap \fk.$
Let $\fg_1 := [\fg, \fg].$ Then
$\fb = \fc_\iq \oplus \fb_1,$ with $\fb_1 = \fb \cap \fg_1.$
Accordingly, $I(\fb) = U(\fc_\iq) \otimes I(\fb_1).$

Let $G_1$ be the analytic subgroup of $G$ with Lie algebra $\fg_1$
and let $H_1 = G_1 \cap H.$
The embeddings $\fc_\iq \subset \fg$ and $\fg_1 \subset \fg$ induce
embeddings $U(\fc_\iq) \subset \DX$ and $\D(G_1 / H_1) \subset \DX,$
via which we identify.
Accordingly, $\DX = U(\fc_\iq) \otimes \D(G_1 /  H_1);$
moreover, the map $\gg: \DX \to I(\fb) $ corresponds with the tensor product of
$I_{ U(\fc_\iq)}$ and $\gg_1,$ the Harish-Chandra isomorphism for $(G_1, H_1, \fb_1).$

If $\gL \in \rmL(\spX,\fb), $ let $\gL_\fc := \gL|_{\fc_\iq}$ and $\gL_1 := \gL|_{\fb_1}.$
Then $U(\fc_\iq) \simeq S(\fc_\iq)$ acts by the character $\gL_\fc$ on the non-trivial space
$\cC(\spX\col \gL).$ This character must therefore be an infinitesimal character of
the compact group $\exp(\fc_\iq),$ hence belongs to $i\fc_\iq^*.$
On the other hand, $\D(G_1/H_1)$ acts by the character
$\gg_1(\dotvar \col \gL_1)$ on
$\cC(\spX\col \gL).$ Restriction to $G_1/H_1$ therefore
induces a map $\cC(\spX\col \gL) \to \cC(G_1/H_1, \gL_1),$ which
is non-zero by $G$-invariance and non-triviality of the space $\cC(\spX\col \gL).$
Hence, $\gL_1 \in \rmL( G_1/H_1,\fb_1).$ If $\fb_1$ is contained
in $\fk,$ then so is $\fb$ and if $\gL_1\in i\fb_1^*$ then
$\gL = \gL_\fc + \gL_1 \in i\fb^*;$ finally, if $\gL_1$ is regular, then so is $\gL.$
Therefore, we may as well assume that
$\fg$ is semisimple from the start.

Let $Z(G)$ denote the center of $G$ and put $Z_H:= Z(G) \cap H.$
Since $Z_H$ is discrete and central,
$\bp G := G/Z_H$ is a Lie group with algebra naturally isomorphic
with $\fg.$ The involution $\gs$ factors to an involution
$\bp\gs$ of $\bp G.$ Moreover, $\bp H := H/Z_H,$ viewed as a subgroup of $\bp G,$
is an open subgroup of $\bp G^{\bp \gs}.$  The associated
symmetric space
$\bp \spX: = \bp G / \bp H$  is naturally diffeomorphic
with $\spX$ and it is readily seen that $\rmL(\spX,\fb) = \rmL( \bp \spX, \fb).$
Therefore, it suffices to prove the assertions for $\bp \spX$
and we see that we may as well assume from the start that $Z_H = \{e\}.$

{}From now on we assume that $G$ is connected and semisimple,
and that $Z_H = \{e\}.$ The natural map $\pi: G/H_e \to G/H$ is a finite covering,
hence induces a linear embedding $\pi^*: \cC(G/H) \to \cC(G/H_e)$ by pull-back.
Via the isomorphism (\refer{e: iso invariant diff ops})
we may identify the algebras $\D(G/H)$ and $\D(G/H_e),$
so that $\pi^*(Df) = D\pi^*f,$ for $f \in \cC(G/H).$ Thus, if $\gL \in \rmL(\spX, \fb),$
then the image of $\cC(G/H\col \gL)$ in $\cC(G/H_e)$ is a non-trivial subspace
annihilated by the ideal $I_\gL,$ from which we see that $\gL\in \rmL(G/H_e, \fb).$ It follows
that we may as well assume that $H$ is connected. We will do so from now on.

Let $\fg^d$ be the dual real form of $\fg_\iC$ defined as in
Section \refer{s: notation}.
Via $\ad$ we identify $\fg_\iC$ with the Lie algebra of the complex adjoint group
$G_\iC$  of $G;$ accordingly, we denote by $G^d, K^d$ and $H^d$
the analytic subgroups of $G_\iC$ with Lie algebras $\fg^d, \fk^d$ and $\fh^d,$
respectively.
Via $\Ad$ we may identify $K \cap H$ with a connected subgroup of $G_\iC.$
Accordingly, the map $(k, X)\mapsto k \exp X$ is a diffeomorphism
from $(K \cap H) \times i[\fk \cap \fq]$ onto $H^d.$ Hence,
for every finite dimensional
representation $(\pi, V)$ of $K$ there exists a unique finite dimensional
representation $(\pi^d, V)$ of $H^d$ such that the infinitesimal
representations associated with $\pi$ and $\pi^d$ have the same
complex linear extension to $\fk_\iC.$ It follows that Flensted-Jensen's
dualization procedure, see \bib{FJds}, Thm.~2.3,
defines an injective linear map $f \mapsto f^d$ (denoted $f \mapsto f^\eta$ in \bib{FJds})
from the space $\Ci(G/H)_K$ of $K$-finite smooth functions on $G/H$ into
the space $\Ci(G^d/K^d)_{H^d}$ of $H^d$-finite smooth functions
on $G^d/K^d.$ The map is determined by the property
that, for every $f \in \Ci(G/H)_K$ and all $u \in U(\fk_\iC),$
$$
L_u f|_{\Aq} = L_u f^d|_{\Aq}.
$$
We note that the left $H^d$-types of $f^d$ are all of the form $\pi^d,$
with $\pi$ a finite dimensional irreducible representation of $K.$ We also note that for
$f \in C^\infty(G/H)_K,$
the condition $f \in L^2(G/H)$ can be entirely rephrased in terms of the
function $f^d;$ in fact it is equivalent to the condition that
$L_u f^d|_{\Aq} \in L^2(\Aq, J\,da),$ for all $u \in  U(\fk),$ with $J$
the Jacobian associated with
the decomposition $G = K\Aq H,$ see \bib{BSmc}, (3.1).

Let $D \mapsto D^d$ denote the natural algebra isomorphism
from $\DX$ onto $\D(\spX^d),$ corresponding to (\refer{e: iso invariant diff ops}).
Then  $(Df)^d = D^d f^d,$ for every $f \in \Ci(G/H)_K.$
Moreover, we recall from the text after (\refer{e: iso invariant diff ops})
that $D^d = \gg_{\fa_\fp^d}{}^{-1}\after \gg_\fb(D),$ where we have written
$\fa_\fp^d$ for the maximal abelian subspace $\dfb = \fb_\iC \cap \fg^d$ of $\fp^d.$
Now assume that $\gL \in \rmL(\spX,\fb).$
Then there exists a non-trivial $K$-finite function $f \in \cC(\spX\col \gL).$
It follows that
$f^d \in \Ci(G^d/K^d)_{H^d}$ satisfies the system of differential equations
$Df^d = \gg_{\fa_\fp^d}(D \col \gL)f^d,$ for $D \in \D(G^d/K^d).$

It follows from the above discussion, that the theorem of \bib{OMds},
p.~359, as well as its proof,
can be entirely formulated in terms of the function $f^d,$
and therefore applies without change, see \bib{OMds}, p.~388, note (i) added in proof.
In particular, we may draw the following conclusions.
In the notation of the cited theorem, we may take
$\fa_\fp^d$ as above, and we may select a positive system
$\gS(\fa_\fp^d)^+$ for $\gS(\fg^d, \fa_\fp^d)$ such that $\Re \gL$
is dominant. The hypothesis of part (i) of the cited theorem is fulfilled,
since the non-trivial function $f$ belongs to the space
$\cA_K(G/H, \Mer_\gl) \cap L^2(G/H),$ with $\gl = \gL.$
It follows that $\fb$ is compact,
i.e., is contained in $\fk \cap \fq.$
In the cited theorem we may now take $\ft  =\fb$ and $\fa_\fp' = i\fb.$
Thus, $\fa_\fp^d = \fa_\fp',$ and it follows from part (i) of the cited theorem
that $\gL$ is regular.

We note that $W(\fa_\fp'\mid \fa_\fp^d) = W(\fa_\fp^d),$
so that the elements $\bar x_j = \Ad(x_j)|_{\fa_\fp^d}$ of the cited theorem
belong to $W(\fa_\fp^d).$
It follows from part (iii) of the cited theorem that, for some $j,$ the element
$\bar x_j \gL = \gl^j$ belongs to $\fa_\fp^{d*}.$ This implies
that $\gL \in \fa_\fp^{d*} = i\fb^*.$
 \qed

\begin{cor}
\naam{c: X of grt and inf char fb new}
Let $\fb \subset \fq$ be a $\Cartan$-stable Cartan subspace.
If $\rmL(\spX, \fb) \neq \emptyset$ then there exists
a Cartan subspace $\ft \subset \fq$ with $\ft \subset \fk.$
Moreover, let $t$ be an element of the set $W(\ft\mid \fb),$
which is non-empty by Lemma \refer{l: gg for diff css}.
Then, for every $\gL \in
\rmL(\spX,\fb),$ the element $t\gL$ belongs to $i \ft^*$ and
is regular relative to the root system $\gS(\ft).$
\end{cor}

\proof
Assume that $\rmL(\spX,\fb) \neq \emptyset.$ Then, by definition,
$\Ltwod(\spX)\neq 0.$ By Theorem
\refer{t: infinitesimal characters L two d real and regular new}
there exists a compact Cartan subspace
$\ft \subset \fq.$
Let $t \in W(\fb, \ft).$ Then by Lemma \refer{l: conjugacy of rmL fb tau}
the element $t$ maps
$\rmL(\spX,\ft)$ bijectively onto $\rmL(\spX,\fb).$
The assertion now follows from Theorem
\refer{t: infinitesimal characters L two d real and regular new}.
\qed

In the rest of this section we fix a Cartan subspace $\fb \subset \fq$
containing $\faq.$ If $P \in \allparabs,$ then the $\Cartan$-stable
Cartan subspace $\stfb_P$ of
$\fm_P \cap \fq$ is defined as in the text before (\refer{e: deco fb with star part}).

\begin{lemma}
\naam{l: data associated to gL}
Let $P \in \allparabs,$
$v \in \NKaq$ and assume that
$\rmL_{P,v}(\fb, \tau) \neq \emptyset.$ Then there exist a
Cartan subspace $\hat\fb \subset  \fm_{1v^{-1}Pv}\cap \fq$
and an element $t \in W(\hat \fb\mid \fb)$ with the following properties.
\begin{enumerate}
\itema
${}^*\hat \fb: =\hat \fb \cap  \fm_{v^{-1}Pv} $ is compact, i.e., contained in $\fk;$
\itemb
$t = \Ad(v)^{-1}$ on $\faPq;$
\itemc
the elements of $t\rmL_{P,v}(\fb, \tau)$ belong to   $i\,{}^*\hat\fb^*$
and are regular relative to $\gS(\fm_{v^{-1}P v\iC}, {}^* \hat \fb).$
\end{enumerate}
\end{lemma}

\proof
{}From (\refer{e: image rmLPv under bar v s minus one}) it follows that
$\rmL(\spX_{P,v}, \stfb_{P,v}, \tau_P) \neq \emptyset.$ Hence,
by Corollary \refer{c: X of grt and inf char fb new}
there exists a Cartan subspace $\ft$ of $\fm_P \cap \Ad(v)(\fq)$
that is contained in $\fk.$
Now  $\hat \fb = \Ad(v)^{-1} (\ft \oplus \faPq)$ is a Cartan subspace
of $\fm_{1 v^{-1}P v}\cap \fq$ that satisfies condition (a),
with ${}^*\hat\fb = \Ad(v)^{-1} \ft.$

Fix $t_1 \in W(\hat \fb\mid \fb).$ Then $t_1^{-1} \Ad(v)^{-1}\faPq \subset \dfb.$
Since $\dfb$ is maximal abelian in $\fp^d,$ it follows from
Lemma \refer{l: first lemma on Weyl groups} (d)
that there exists a $t_2 \in W(\fb)$ such that
$t_2 = t_1^{-1} \Ad(v)^{-1}$ on $\faPq.$ It follows that $t = t_1t_2 \in W(\hat\fb\mid \fb)$
satisfies requirement (b).

Finally, let $\gL \in \rmL_{P,v}(\fb, \tau).$ Then, in the notation of
(\refer{e: image rmLPv under bar v s minus one}),
the element $\gL' := \Ad(v)s^{-1} \gL$ belongs to
$\rmL(\spX_{P,v}, \stfb_{P,v}, \tau_P).$
The element $t' = \Ad(v) t s \Ad(v)^{-1}$
belongs to $W(\hat \fb^v \mid \fb^v)$ and equals the identity
on $\faPq,$ hence restricts to an element of  $W(\ft\mid \stfb_{P,v}).$
By Corollary \refer{c: X of grt and inf char fb new} it follows that
$t'\gL'$ belongs to $i \ft^*$ and is regular relative to
$\gS(\fm_{P\iC}, \ft).$ We now observe that  $t'\gL' = \Ad(v) t \gL.$
Hence, $t\gL$ belongs to $i {}^*\hat\fb$ and is regular with
respect to $\gS(\fm_{v^{-1}Pv\iC}, {}^*\hat\fb).$
\qed

\begin{rem}
\naam{r: on description eigenvalue}
Let $P, v, \hat \fb, t$ be as
in Lemma \refer{l: data associated to gL}. Then it
follows from Lemma \refer{l: gg for diff css} that,
for all $\gL \in \rmL_{P,v}(\fb,\tau),$
$\nu \in \faPqdc$
and $D \in \DX,$
$$
\gg_\fb(D \col \gL + \nu) = \gg_{\hat\fb}(D \col t \gL+ t \nu).
$$
\end{rem}

\begin{cor}
\naam{c: gL is real regular}
Let $P \in \allparabs$
and let $\gL \in \rmL_P(\fb, \tau).$
Then $\inp{\gL}{\ga} \in \R \setminus \{0\},$
for each $\ga \in \gS(\fm_{P\iC}, \stfb_P).$
\end{cor}

\proof
Select $v \in \PcW$ such that $\gL \in \rmL_{P,v}(\fb,\tau).$ Let
$\hat \fb, t$ be associated as in Lemma \refer{l: data associated to gL}.
Then by (a) and (b) of the mentioned lemma, $t$ maps $\fb = \stfb_P\oplus \faPq$
onto $\hat\fb = {}^*\hat\fb \oplus \fa_{v^{-1}Pv\iq},$ preserving the decompositions.
The assertion now follows from Lemma \refer{l: data associated to gL} (a) and (c).
\qed

\section{The constant term of the Eisenstein integral}\eqnreset
\naam{s: the constant term}
In this section we describe the constant term of the normalized
Eisenstein integral, introduced in Definition \refer{d: dual Eisenstein new}.
We start by recalling the notion of the constant
term introduced in \bib{Cn}.

If $f \in \cAtemp(\spX\col \tau),$ see Def.~\refer{d: temperedness},
then in particular
$f \in \cA(\spXp\col \tau)$ and $f$ has an expansion of the form
(\refer{e: series for f in q}).
It follows from Lemma \refer{l: criterion for temperedness} (c) combined with
\bib{BSanfam}, Thm.~3.5,
that, for each $Q\in \allparabs$ and every $v\in \NKaq,$
$$
\xi \in \Exp(Q,v \asmid f) \;\;\implies \;\;\Re \xi + \rho_Q \leq 0 \text{on} \faQqp.
$$
We define the function
$f_{Q, v}: \spXoneQvp \to \Vtau$ by
$$
f_{Q, v}(ma) = d_Q(ma)\,\sum_{\xi\in \Exp(Q,v\mid f)\atop
\Re \xi + \rho_Q = 0} a^\xi q_\xi(Q,v\mid f, \log a, m),
$$
for $m\in \spXQvp,$ $a\in \AQq.$
Here $d_Q: \MoneQ \to \R$ is defined by $d_Q(m) = \sqrt{|\det \Ad(m)|_{\fn_Q}}|.$
Note that $d_Q = 1$ on $M_Q$ and on $A_Q\cap H.$
Hence, $d_Q$ factors to a function on
$\spXoneQv;$ in fact, $d_Q(ma) = a^{\rho_Q},$ for $m\in \MQgs$ and $a\in \AQq.$

We note that, for $v\in \NKaq,$ the function $R_v f: x \mapsto f(xv)$
belongs to the space $\cAtemp(G/vHv^{-1} \col \tau).$

\begin{prop}{\ }
\naam{p: characterization constant term}
\begin{enumerate}
\itema
If $u,v\in \NKaq,$ then $(R_v f)_{Q,u} = f_{Q, uv}.$
\itemb
The function $f_{Q,v}$ extends uniquely to smooth
function on $\spXoneQv.$ This extension is the unique function
in $\cAtemp(\spXoneQv\col \tau_Q)$ such that
$$
\lim_{t \to \infty} ( d_Q(m \exp tX) f(m\exp tX v ) - f_{Q, v}(m \exp tX) ) = 0,
$$
for every $m \in M_{1Q}$ and $X \in \faQqp.$
\end{enumerate}
\end{prop}
\proof
The first assertion follows from \bib{BSanfam}, Lemma 3.7.
In view of (a) it suffices to prove the second assertion for $v = e.$
In this case the assertion follows from  \bib{Cn},
proof of Thm.~1.
\qed

Thus, for $v = e,$ the function $f_{Q,v}$ coincides with the constant
term of $f$ along $Q,$ introduced by \bib{Cn},
which in turn generalizes Harish-Chandra's
notion of the constant term for the case of the group, see \bib{HC1},
Sect.~21, Thm.~1.
We shall therefore call $f_{Q,v}$ the constant term of $f$ along $(Q,v).$
The following result, which generalizes a result of Harish-Chandra,
see \bib{HC1}, Sect.~21, Lemma 1, is essentially given in \bib{Cn}, Thm.~1 (b).

\begin{lemma}
\naam{l: transitivity of the constant term}
{\rm (Transitivity of the constant term)\ }
Let $P,Q \in \allparabs$ be such that $P \subset Q.$ Put $\staroneP: = M_{1Q} \cap P.$
Let $v\in \NKaq$ and $u\in \NKQaq.$
Then
$$
(f_{Q,v})_{\staroneP, u} = f_{P, uv}.
$$
\end{lemma}

\proof
For $v = u=e$ the result is equivalent to \bib{Cn}, Thm.~1(b).
Let now $v \in \NKaq$ and $u \in N_{K_Q}(\faq)$ be general. Then right translation by
$u$ defines a linear isomorphism
$
R_u: \cA(\spX_{1Q,v}\col \tau_Q) \to \cA(\spX_{1Q, uv}\col \tau_Q).
$
Hence, applying Proposition \refer{p: characterization constant term} (b) we find that
\begin{equation}
\naam{e: Ru of fQv}
R_u(f_{Q,v}) = f_{Q, uv}.
\end{equation}
Applying Proposition \refer{p: characterization constant term} (a) we see that
\begin{equation}
\naam{e: expression for f Quv}
f_{Q,uv} = (R_{uv} f)_{Q,e}\quad\text{and}\quad (f_{Q,v})_{\staroneP, u} =
(R_u f_{Q,v})_{\staroneP, e}.
\end{equation}
Combining (\refer{e: Ru of fQv}) with (\refer{e: expression for f Quv}), and using the first line
of the proof and Proposition \refer{p: characterization constant term} (a), we finally obtain that
$$
(f_{Q,v})_{\staroneP, u} = ((R_{uv} f)_{Q,e})_{\staroneP, e} = (R_{uv} f)_{P,e} = f_{P,uv}.
$$

\qed
\medbreak\medbreak
The following transformation rule for the constant term will also be useful to
us. If $u,v \in \NKaq,$ we define the map $\rho_{\tau,u}: \Ci(\spX_{1Q,v} \col \tau_Q) \to
\Ci(\spX_{1uQu^{-1}, uv} \col \tau_{uQu^{-1}})$ in accordance with  \bib{BSanfam}, Eqn.~(3.24),
by
$$
\rho_{\tau,u}\gf(m) = \tau(u) \gf(u^{-1} m u),\qquad (m \in M_{1uQu^{-1}}).
$$
One readily checks that $\rho_{\tau,u}$ maps
$\cAtemp(\spX_{1Q,v} \col \tau_Q)$ into $\cAtemp(\spX_{1uQu^{-1}, uv} \col \tau_{uQu^{-1}}).$

\begin{lemma}
\naam{l: second transformation rule constant term}
Let $f \in \cAtemp(\spX\col \tau)$ and let $Q \in \allparabs$ and $u,v \in \NKaq.$
Then
\begin{equation}
\naam{e: second transformation rule constant term}
f_{uQu^{-1}, uv} = \rho_{\tau,u} f_{Q,v}.
\end{equation}
\end{lemma}

\proof
{}From the definition of $d_Q$ one readily verifies that $d_{uQu^{-1}}(u m u^{-1}) =d_Q(m),$
for  $m \in \MoneQ.$ The result now follows by a straightforward application
of Proposition \refer{p: characterization constant term} (b). See
\bib{BSanfam}, Lemma 3.6, for a similar proof.
\qed

Assume that $\Omega \subset i\faPqd$ is open, and $f: \Omega \times \spX \to \Vtau$
a smooth map such that $f_\nu: x \mapsto f(\nu,x)$ belongs
to $\cAtemp(\spX\col \tau)$ for every $\nu \in \Omega.$
If $Q \in \allparabs$ and $v \in \NKaq,$
we shall write $f_{Q,v}$ for the map $\Omega \times \spXoneQv \to \Vtau$
defined by
$$
f_{Q,v}(\nu, m) = (f_\nu)_{Q,v} (m),\qquad (\nu \in \Omega,\; m \in \spXoneQv).
$$
We now  turn our attention to the normalized  Eisenstein integral
$\nE(P\col \dotvar)$ where $P \in \allparabs$
is assumed to be
of residue type.
In the end it will follow that any $P \in \allparabs$ is of this type, see
Remark \remRT, so that this is really no restriction on $P.$
Let $\Omega_P$ be the set of points in
$i\faPqd$ where the function $\nu \mapsto \nE(P\col\nu)$ is regular.
Then for $\nu \in \Omega_P$ and $\psi \in \cAtwoP,$ the function
$\nE(P\col \nu\col \dotvar)\psi$ belongs to $\cAtemp(\spX \col\tau),$ see
Proposition \refer{p: exponents normalized Eis}.
In accordance
with the above, we denote its constant term along $(Q,v),$ for
$Q \in \allparabs$ and $v\ \in \NKaq,$ by $\nEQv(P\col \nu\col \dotvar)\psi.$

\begin{prop}\restypetwo
\naam{p: constant term of the Eisenstein integral}
Let $P,Q \in \allparabs$
and $u \in \NKaq.$
\begin{enumerate}
\itema
The function $\nEQu(P\col \dotvar)$ extends to a meromorphic
$\Ci(\spX_{1Q,u}, \Hom(\cAtwoP, \Vtau))$-valued function on
$\faPqdc,$ with singular set equal to a
locally finite union of real $\gS_r(P)$-hyperplanes.
\itemb
There exists a $\geps>0$ such that, for every $\psi \in \cAtwoP$
and $p \in \Pi_{\gS_r(P)}(\faPqd)$ with the property
that $\nu \mapsto p(\nu) \nE(P\col \nu\col \dotvar)\psi$ is regular on $\faPqd(\geps),$
the function $\nu \mapsto p(\nu) \nEQu(P\col \nu\col \dotvar)\psi$
is regular on $\faPqd(\geps)$ as well.
\itemc
If $\nEQu(P\col \dotvar) \neq 0,$  then $W(\faPq\mid \faQq)$ is non-empty.
\itemd
Let $W(\faPq\mid \faQq)$ be non-empty. Then
there exist unique
meromorphic functions
$E^\circ_{Q,u,s}(P\col \dotvar): \faPqdc \to \Hom(\cAtwoP, \Ci(\spXQu \col \tau_Q)),$
for $s \in W(\faPq\mid \faQq),$
such that, for all $m \in \spXQu$ and $a \in \AQq,$
\begin{equation}
\naam{e: constant term Eis as sum}
\nEQu(P\col \nu\col  ma) = \sum_{s \in W(\faPq\mid \faQq)} a^{s^*\nu }
E^\circ_{Q,u,s}(P\col \nu\col m),
\end{equation}
as an identity of meromorphic functions in the variable
$\nu \in \faQqdc.$ Here $s^*\nu = \nu \after s,$
see \S~\refer{s: Weyl groups}.
 The singular locus of any of the meromorphic functions
$E^\circ_{Q,u,s}(P\col \dotvar),$ for $s \in W(\faPq\mid \faQq),$
is the union of a  locally finite collection
of real $\gS_r(P)$-hyperplanes.
\end{enumerate}
\end{prop}

\proof
We give the proof under the assumption that $P$ is of residue type,
see Remark \remRT.

{\bf (a):\ }
Let $\psi \in \cAtwoP,$ and define $f: (\nu, x) \mapsto \nE(P\col \nu\col x)\psi.$
Then $f \in \cE_{P,Y}^\hyp(\spX\col \tau),$ with $Y \subset \staPqd$ a finite
subset, see Lemma \refer{l: q belongs to Mer newer}.
In particular, it follows that $f \in  C^{\ep, \hyp}_{P,Y}(\spXp\col \tau),$
see Section \refer{s: vanishing theorem}.
The set $\Hyp: = \Hyp_f$ is a real $\gSr(P)$-configuration
in $\faPqdc,$ again by Lemma \refer{l: q belongs to Mer newer}.

Let $\Hyp_0$ be the collection of $H \in \Hyp$ with $H \cap i\faPqd \neq \emptyset.$
Then $\Hyp_0$ is finite, since $\Hyp$ is real.
For every $H \in \Hyp_0$ we select a first degree
polynomial function
$l_H \in P_1(\faPqd)$ with $H = l_H^{-1}(0),$ and put
$$
\pi_0 = \prod_{H \in \Hyp_0} l_H^{d(H)},
$$
with $d = d_f.$
Select $\geps_0>0$ such that
$H \in \Hyp, H \cap \faPqd(\geps_0) \neq \emptyset \implies H \in \Hyp_0.$
Then the family $f^0: (\nu,x) \mapsto \pi_0(\nu)f(\nu, x)$ belongs
to $\cO(\faPqd(\geps_0), \Ci(\spX\col \tau)).$
Moreover, in view of Lemma \refer{l: q belongs to Mer newer},
for every $\gs \in W/\!\!\sim_{Q|P}$ and $\xi \in - \gs\cdot Y +\N\DrQ,$
the function
$$
q^0_{\gs, \xi}(Q,u \mid f) := \pi_0 \,q_{\gs, \xi}(Q,u \mid f)
$$
belongs to $P_k(\faQq) \otimes \cO(\faPqd(\geps_0), \Ci(\spXQu \col \tau_Q));$
here $k = \dega f.$
It follows from \bib{BSanfam}, Lemma 12.7, that
$$
f_\nu^0(mau) = \sum_{\gs \in W/\simQP}a^{\gs \nu - \rho_Q}
\sum_{\xi \in -\gs \cdot Y + \N\DrQ}
a^{-\xi} q^0_{\gs,\xi}(Q,u\mid f, \log a)(\nu, m),
$$
for every $m \in \spXQup,$ and $a \in \AQqp(R_{Q,u}(m)^{-1}),$
where the second series converges neatly in $a.$
For every $\nu \in \Omega_P,$ the function $f^0_\nu$ belongs
to $\cA_\temp(\spX\col \tau),$ see Proposition
\refer{p: exponents normalized Eis}.
Since $Y$ is real, it follows
by uniqueness of asymptotics, for all $m \in \spXQup$
and $a \in \AQq,$ that
\begin{equation}
\naam{e: constant term f nought as sum q}
(f^0_\nu)_{Q,u}(ma) = \sum_{\gs \in W/\simQP \atop 0 \in -\gs \cdot Y + \N\DrQ}
a^{\gs \nu} q_{\gs, 0}^0(Q,u\mid f,\log a)(\nu, m).
\end{equation}
By density and continuity, this expression holds for all $m \in \spXQu$ and $a \in \AQq.$
On the other hand, by the characterization of the constant term in
Proposition \refer{p: characterization constant term} (b),
it follows that, for $\nu \in \Omega_P,$
$$
(f^0_\nu)_{Q,u} = \pi_0(\nu) (f_\nu)_{Q,u}.
$$
Using Lemma \refer{l: q belongs to Mer newer} once again,
we infer from (\refer{e: constant term f nought as sum q}) that
$\nu \mapsto E^\circ_{Q,u}(P\col \nu)\psi = (f_\nu)_{Q,u}$ extends
to a meromorphic $\Ci(\spX_{1Q,u}\col \tau_Q)$-valued
function on $\faPqdc$ with singular set contained in $\cup \Hyp_f.$
This establishes (a).

We will first establish the remaining assertions under the assumption
that $u =e.$

{\bf (b):\ }
Let $\gL \in \rmL_P( \fb,\tau),$ $\psi \in \cAtwoP(\gL)$
and define $f$ as above. For $p \in \Pi_{\gSr(P)}(\faPqd)$
we put $f_p(\nu,x) = p(\nu)f(\nu, x).$

According to Theorem
\refer{t: norm Eis uniform tempered},
there exist $q \in \Pi_{\gSr(\faPqd), \R}(\faPqd)$
and constants $\geps_0>0$ and $s_0>0,$ all independent of $\gL$ and $\psi,$
such that $f_q$ is holomorphic on $\faPqd(\geps_0)$ and
belongs to $\cT(P ,\tau, \gL, \geps_0, s_0).$
Let $\geps_1$ be any constant with
$0 < \geps_1 < \geps_0.$ If $ p \in \Pi_{\gSr(P)}(\faPqd)$ is such
that $f_p$ is holomorphic on $\faPqd(\geps_1),$
then clearly $f_{pq} \in \cT(P ,\tau, \gL, \geps_0, s_0).$
By a repeated application of Cauchy's integral
formula to $f_{pq}(\nu, x),$ with polydiscs  of size
$O((1 +  \lspX(x))^{-1}),$ it now follows that
$ f_p \in \cT(P ,\tau , \gL, \geps_1', s_0),$ for every $\geps_1'$ with
$0< \geps_1' < \geps_1.$ See \bib{Bps2}, Lemma 6.1, for a more detailed indication
of how to use Cauchy's formula.

Let $v \in \PcW$ be such that $\gL \in \rmL_{P,v}(\fb, \tau),$
let $\hat \fb, t$ be as in Lemma \refer{l: data associated to gL}
and put $\hat P = v^{-1}Pv.$ Then, in view of
Remark \refer{r: on description eigenvalue},
the family
$\hat f_p: \fa_{\hat P\iq\iC}^* \times \spX \to \Vtau$
defined by $\hat f_p(\mu, x) = f_p(t^{-1} \mu, x)$ belongs
to $\cT(\hat P ,\tau, t\gL, \geps_0, s_0).$

Since $t\gL \in i{}^*\hat \fb^*,$
by Lemma
\refer{l: data associated to gL} (c), we
may apply  \bib{Cn}, Thm.~3, which in turn is based
on \bib{Bps2}, Thm.~12.9. Let $\geps_1' < \geps_1 < \geps_0$ be as above.
According to the mentioned
theorem there exists a constant $\bar \geps_1' >0$
such that for every $F \in \cT(\hat P ,\tau, t \gL, \geps_1', s_0),$ and all
$m \in \spXoneQe,$ the function $\nu \mapsto (F_\nu)_{Q,e}$
is holomorphic on $\fa_{\hat P \iq}^*(\bar \geps_1').$
{}From the proof of \bib{Cn}, Thm.~3,
it follows that this holds with $\bar \geps_1' = \min(\geps_1', \bar \geps),$
where $\bar \geps >0$ is the constant of \bib{Cn}, Lemma 5.
The latter constant only
depends on $\gL;$ the set $\rmL_P(\fb, \tau)$ is finite, hence we may chose
$\bar \geps$ simultaneously for all $\gL$ under consideration.
We now fix $\geps >0$ such that $\geps < \min(\geps_0, \bar \geps).$
Assume that the hypothesis of part (b) of the theorem is fulfilled.
If we apply the above discussion to the functions $f_p$ and $F = \hat f_p,$ with $\geps_1 = \geps,$
then $\bar \geps_1' = \geps_1',$ and it follows that
the function $\nu \mapsto ((f_p)_\nu)_{Q,e}(m) =
p(\nu)E^\circ_{Q,e}(P\col \nu \col m)\psi$
is holomorphic on $\faPqd(\geps_1'),$ for all $m \in \spX_{Q,e}.$ In view of
part (a) of the theorem, it follows that
$\nu \mapsto p(\nu)E^\circ_{Q,e}(P\col \nu)\psi$ is holomorphic
on $\faPqd(\geps_1')$ as a function with values in $\Ci(\spXoneQe\col \tau_Q).$ This holds
for every $\geps_1' < \geps,$ whence the desired assertion.

{\bf (c):\ }
{}From the hypothesis with $v=e$ it follows that there exists a $\gL \in \rmL_P(\fb, \tau)$
and a $\psi \in \cAtwoP(\gL)$ such that $\nEQu(P\col \cdot)\psi \neq 0.$
Let $\Omega_P'$ be the set of $\nu \in \Omega_P$ such that $\gL + \nu$ is
a $\gS(\fb)$-regular element of $\fbdc.$
It follows from Lemma \refer{l: data associated to gL}
that $\gL$ is regular with respect to $\gS(\fm_{P\iC},\stfb_P).$ Therefore,
$\Omega_P'$ is open dense in $\Omega_P,$ hence
in $i\faPqd.$ We infer that we may select $\nu \in \Omega_P'$
such that $(f_\nu)_{Q,e} \neq 0,$ with notation as introduced in part (a) of this proof.

Fix $v \in \PcW$ such that $\gL \in \rmL_{P,v}(\fb, \tau).$
Let $(\hat\fb, t)$ be  as in Lemma
\refer{l: data associated to gL} and put $\hat P := v^{-1}Pv.$
Then $t \gL$ belongs to $i {}^* \hat\fb^*$ and is
regular relative to $\gS(\fm_{\hat P \iC}, {}^* \hat\fb).$
Now $f_\nu \in \cAtemp(\spX\col \tau)$ and by Remark \refer{r: on description eigenvalue},
\begin{equation}
\naam{e: alternative description eigenequation}
Df_\nu = \gg_{\hat \fb}(D\col t\gL + t\nu) f_\nu,
\end{equation}
for all $D \in \DX.$  Since $t \gL \in i{}^*\hat\fb^*,$ $t\nu \in i\fa_{\hat P \iq}$
and  $t \gL + t \nu$ is regular with respect to $\gS(\hat\fb),$
it follows from  \bib{Cn}, Thm.~2, that
the set $W(\fa_{\hat P \iq}\mid \faQq)$ is non-empty. The map $s \mapsto t^{-1}\after s$
is a bijection from the latter set onto $W(\faPq\mid \faQq),$ which
set is therefore non-empty as well.

{\bf (d):\ }
Uniqueness of the functions $E^\circ_{Q,e,s}(P\col \dotvar)$ is obvious,
by linear independence of the functions $a \mapsto a^{s^* \nu}$
 for generic $\nu.$
We fix $\gL$ and $\psi$ as in part (b) of the proof and
define $f$ as in part (a).
We define the set $\Omega_P'$ as in (c).
Let $\Omega_P''$ be the open dense subset consisting of
$\nu \in \Omega_P'$ with $s^*\nu$ mutually different, for $s \in W(\faPq\mid \faQq).$
Let $v, \hat \fb,t, \hat P$ be as in part (c) of the proof, and fix $\nu \in \Omega_P''.$
In view of (\refer{e: alternative description eigenequation}),
it follows from \bib{Cn}, Thm.~2, that there exists a collection
of functions $f_{t\nu, \hat s, Q} \in \Ci(\spX_{1Q,e}\col \tau_Q),$ for
$\hat s \in W(\fa_{\hat P \iq} \mid \faQq),$
such that
$$
(f_{\nu})_{Q,e}(m)  = \sum_{\hat s \in W(\fa_{\hat P\iq} \mid \faQq)}  f_{t \nu, \hat s, Q}(m),
\qquad (m \in \spX_{1Q,e});
$$
and
$$
f_{t\nu, \hat s,Q}(ma) = a^{\hat s^*t\nu} f_{t\nu, \hat s,Q}(m),
 \qquad (m \in  \spX_{1Q,e},\; a \in \AQq).
$$
Combining these equations,  substituting $ts$ for $\hat s$ and
writing $f_{\nu, s, Q} = f_{t\nu, t s, Q},$ we see that, for all
$m \in \spXQe$ and $a \in \AQq,$
\begin{equation}
\naam{e: deco constant term f}
(f_\nu)_{Q,e}(ma) = \sum_{s \in W(\faPq\mid \faQq)} a^{s^*\nu} f_{\nu, s,Q}(m).
\end{equation}
For every $s \in W(\faPq\mid \faQq)$ there exists an element $\tilde s\in W$ \
such that $s = \tilde s|_{\faQq},$ see Corollary \refer{c: W a1 a2 in q}.
It follows that $s^*\nu = \tilde s^{-1} \nu|_{\faQq},$
for all $\nu \in \faPqdc.$ Using the definition of $\sim_{Q|P}$ we see that
the class of $\tilde s^{-1}$ in $W/\!\!\sim_{Q|P}$ is  uniquely determined
by $s.$ We denote this class by $\gs_s.$ Comparing
(\refer{e: constant term f nought as sum q}) and (\refer{e: deco constant term f}) we see
by uniqueness of asymptotics that $X \mapsto q_{\gs_s, 0}(Q,e\mid f, X)(\nu)$
is constant as a $\Ci(\spX_{Q,e}\col \tau_Q)$-valued function on $\faQq$ and that
$$
f_{\nu, s, Q}(m) = q_{\gs_s, 0}(Q,e\mid f, 0)(\nu, m),
$$
for all $m \in \spXQe.$
We define $E^\circ_{Q,e,s}(P\col \nu\col m)\psi := q_{\gs_s, 0}(Q,e\mid f, 0)(\nu , m).$
Then (\refer{e: constant term Eis as sum}) applied to $\psi$ follows
for $\nu \in \Omega_P''.$
Finally, the assertions on meromorphy  follow from the fact
that $q_{\gs_s, 0}(Q,e\mid f, 0) \in \Mer(\faPqdc, \Hyp_f, d_f, \Ci(\spXQe\col \tau_Q)),$
by Lemma \refer{l: q belongs to Mer newer}.

It remains to establish (b)-(d) under the assumption that $u \in \NKaq$ is arbitrary.
Assertion (b) follows from the already established assertion with $u = e$
by application of Lemma \refer{l: second transformation rule constant term}
with $u,e$ in place of $u, v,$ respectively.

To prove (c), we assume that $E^\circ_{Q,u}(P\col \dotvar) \neq 0$ and
put  $Q' = u^{-1} Q u .$ Using Lemma
\refer{l: second transformation rule constant term}
we infer that $E^\circ_{Q', e}(P\col \dotvar) \neq 0.$
Hence,  from the already established assertion
(c) with $Q',e$ in place of $Q,u$ it follows  that $W(\faPq\mid \fa_{Q'\iq}) \neq
\emptyset.$ Since $s \mapsto \Ad(u)^{-1}\after s\after  \Ad(u)$ induces
a bijection from  $W(\faPq\mid \faQq)$ onto
$W(\faPq\mid \fa_{Q'\iq}),$ it follows that $W(\faPq\mid \faQq) \neq \emptyset.$

Finally, assertion (d) follows from the already established
assertion with $ u^{-1}Qu  ,e$ in place of $Q,u$
by applying Lemma \refer{l: second transformation rule constant term} once more
in a similar fashion as above.
\qed

If $P,Q$ are associated parabolic subgroups in $\allparabs$, see
Def.~\refer{d: associated parabolics},
then $W(\faPq\mid \faQq)$ is a non-empty finite set of isomorphisms
from $\faQq$ onto $\faPq;$ moreover, the natural left action
of $W(\faPq)$ as well as the natural right action of
$W(\faQq)$ on this set is free and transitive.

\begin{prop}\restypetwo
\naam{p: constant term in cAtwo}
Let $P,Q\in \allparabs$ be associated
and let
$v \in \NKaq.$

Then, for each $s \in W(\faPq\mid \faQq)$ and every $\psi \in \cAtwoP,$
the meromorphic $\Ci(\spXQv\col \tau)$-valued
function $\nu \mapsto E^\circ_{Q,v,s}(P\col \nu\col \dotvar)\psi$
on $\faPqdc,$ defined as in (\refer{e: constant term Eis as sum}),
attains its values in the finite dimensional space
$\cA_2(\spXQv\col \tau_Q).$
\end{prop}

\proof
We give the proof under the assumption that
$P$ and $Q$ are of residue type,
see Remark \remRT.
Fix $\psi \in \cAtwoP.$ Let
$\nu \in \Omega_P$ and define the function
$f \in \cAtemp(\spXoneQv\col \tau_Q)$ by
$$
f(m) = \nEQv(P\col \nu \col m) \psi\qquad (m \in \spXQv).
$$
We recall from Section \refer{s: notation} that $\staQq$ is maximal
abelian in $\Ad(v)\fq.$
Let $R_{Q}$ be a proper parabolic subgroup of $M_{Q}$ that
contains ${}^*A_{Q\iq}$ and is stable under the involution
$\gs^v\Cartan.$
In the notation of Section \refer{s: notation},
$R_Q$ is of the form $P_X,$ for some $X \in \staQq,$ relative to $(M_Q, \Cartan)$
in place of $(G, \Cartan).$ Since $X$ is fixed under $\gs\Cartan,$ it
follows that $R_Q$ is $\gs\Cartan$-stable as well. The $\gs^v$-split component
of the Lie algebra of $R_Q$ equals $\fa_{R_Q} \cap \Ad(v)\fq = \fa_{R_Q}\cap \staQq
= \fa_{R_Q} \cap \fq,$ hence equals the $\gs$-split component.
We therefore denote it by $\fa_{R_Q\iq};$ the associated positive
chamber is denoted by $\fa_{R_Q\iq}^+.$

Since $A_Q$ is central in $\MoneQ,$ and stable under both $\gs^v$ and $\gs,$
the group $R_{1Q} = R_Q A_Q$ is a parabolic subgroup
of $\MoneQ$ that contains $\Aq$ and is stable under both  involutions $\gs^v \Cartan $
and $\gs\Cartan.$ The associated
$\gs^v$-split component equals $\Aq = \stAQq \AQq,$ which is also
equal to the $\gs$-split component of $R_{1Q}.$ Accordingly, the
positive chamber is given by $A_{R_{1Q}\iq}^+ = A_{R_Q\iq}^+ \AQq.$

We now claim that every $\xi \in \Exp(R_{1Q}, e\mid f)$ satisfies
\begin{equation}
\naam{e: estimate exponent on pos chamber}
\Re \xi + \rho_{R_{1Q}} \leq 0 \text{\rm on} \fa_{R_{1Q}\iq}^+.
\end{equation}
Indeed, for $R_{1Q}$ minimal this follows from Lemma
\refer{l: criterion for temperedness}. For general $R_{1Q},$ it follows by
application of \bib{BSanfam}, Thm.~3.5.

On the other hand, it is readily seen that
$R = R_{1Q} N_Q$ is an element of $\allparabs$ and that
$R_{1Q} = \stoneR: = R \cap \MoneQ.$
By application of Lemma \refer{l: transitivity of the constant term},
$$
f_{R_{1Q}, e} = E_{R,v}^\circ(P\col \nu)\psi.
$$
{}From $R_Q \subsetneq M_Q$ we infer that $R \subsetneq Q,$
hence $\faQq \subsetneq \fa_{R\iq},$
from which we see that $\dim \fa_{R\iq} > \dim \faPq;$
hence, $W(\faPq \mid \fa_{R\iq}) = \emptyset.$
{}From Proposition \refer{p: constant term of the Eisenstein integral}
(b) it now follows that the function
on the right-hand side of the above equality is zero.
We conclude that $f_{R_{1Q}, e} = 0$ for every $R_Q$ as above.
By definition of the constant term it follows
that every $\xi \in \Exp(R_{1Q}, e \mid f)$ satisfies
$\Re \xi + \rho_{R_{1Q}} \neq 0$ in addition to
(\refer{e: estimate exponent on pos chamber}).
Put $f_s = E^\circ_{Q,v,s}(P\col \nu)\psi.$ Then
$$
f (a m ) = \sum_{s \in W(\faPq\mid \faQq)} a^{s^*\nu} f_s(m),
$$
for $m \in \spXQv$ and $a \in \AQq.$ It follows that every $f_s$ belongs to
$\cAtemp(\spXQv \col \tau_Q).$ Moreover, every
$\xi \in \Exp(R_Q, e\mid f_s)$ satisfies
$\Re \xi + \rho_{R_{Q}} \leq 0$ on $\fa_{R_Q\iq}^+$ and $\Re \xi + \rho_{R_Q} \neq 0.$
In particular, if $R_Q$ is a maximal $\Cartan \gs^v$-stable parabolic subgroup in  $M_Q,$
it follows that every exponent $\xi \in \Exp(R_Q, e\mid f_s)$ satisfies
$\Re\xi + \rho_{R_Q} < 0$ on $\fa_{R_Q\iq}^+.$ This implies that $(f_s)_{R_Q,e} = 0;$
hence, $f_s \in \cA_2(\spXQv\col \tau_Q),$ by \bib{BSft}, Prop.~12.

Let $s \in W(\faPq\mid \faQq).$ We have shown that the function
$\gf: \nu \mapsto E^\circ_{Q,v,s}(P\col \nu ) \psi$ attains
its values in $\cA_2(\spXQv\col \tau_Q)$
for $\nu \in \Omega_P.$ Since $Q$ is of residue type,
$\cA_2(\spXQv\col \tau_Q)$ is a finite
dimensional subspace of $\Ci(\spXQv\col \tau_Q)$
by Lemma \refer{l: cAt equals cAtwo}. By meromorphy it now follows
that $\gf$ is $\cA_2(\spXQv\col \tau_Q)$-valued.
\qed

If $P,Q \in \allparabs$ are associated, then $s \mapsto s^{-1}$ defines a bijection
from $W(\faQq\mid \faPq)$ onto $W(\faPq\mid \faQq).$ In this case we write,
for $s \in W(\faQq\mid \faPq),$
$$
s\nu := (s^{-1})^*\nu = \nu \after s^{-1}, \qquad (\nu \in \faPqdc).
$$

\begin{defi}\restypetwo
\naam{d: c functions as constant term}
Let $P,Q\in \allparabs$ be
be associated.
For each
$s \in W(\faQq\mid \faPq)$ we define the meromorphic
$\Hom(\cAtwoP, \cAtwoQ)$-valued function $\nCQP(s \col \dotvar)$ on $\faPqdc$
by
$$
[\nCQP(s \col \nu) \psi]_v = E^\circ_{Q,v,s^{-1}}(P\col \nu)\psi, \qquad(v \in \QcW).
$$
\end{defi}
In the chain of reasoning leading up to Theorem \thmendRT{},
this definition requires $P$ to be of residue type, since
it depends on the validity of Definition \refer{d: dual Eisenstein new};
see Remark \remRT.

\begin{cor}\restypetwo
\naam{c: constant term and C}
Let $P,Q\in \allparabs.$
For each $s \in W(\faQq\mid \faPq),$
the  $\Hom(\cAtwoP, \cAtwoQ)$-valued meromorphic function
$\nC_{Q|P}(s \col \dotvar)$
on $\faPqdc$ has a singular locus equal to a locally finite
union of
real $\gSr(P)$-hyperplanes.

Let $\nu \in i\faPqd$ be a regular point for $\nE(P\col \dotvar)$
and the $C$-functions $\nC_{Q|P}(s\col \dotvar),$ as $s \in W(\faQq\mid \faPq).$
Let $\psi \in \cAtwoP.$
Then the function $\nE(P\col \nu)\psi,$ which belongs to
$\cAtemp(\spX\col \tau)$ by Proposition \refer{p: exponents normalized Eis},
has the following
constant term along
$(Q,v),$ for $v \in \QcW,$
\begin{equation}
\naam{e: constant term nE with C}
E^\circ_{Q,v}(P\col \nu \col ma )\psi =
\sum_{s \in W(\faQq\mid \faPq)} a^{s \nu} \;
[\nC_{Q|P}(s \col \nu)\psi]_v(m),
\end{equation}
for all $m\in \spXQv$ and $a \in \AQq.$
\end{cor}

\proof
We give the proof under the assumption  that $P$ and
$Q$ are of residue type,
see  Remark \remRT.
The result is an immediate consequence of Propositions
\refer{p: constant term of the Eisenstein integral} and
\refer{p: constant term in cAtwo} combined with Definition
\refer{d: c functions as constant term}.
\qed

\begin{rem}
\naam{r: constant term and c functions}
Formula (\refer{e: constant term nE with C})
above generalizes Harish-Chandra's formula for the constant
term of the normalized Eisenstein integral in \bib{HCeis},
Thm.\ 5-6,
see also \bib{HC2}, Thm.\ 14.1.
Accordingly, the functions $\nCQP(s\col \dotvar),$ for $s \in W(\faQq\mid \faPq),$
will be called normalized $C$-functions.

In the context of reductive symmetric spaces, for minimal $P$ the above result is
due to \bib{Bps2}, Eqn.~(133), in view of \bib{BSft}, Eqn.~(52).
For general $P$ the result is due to \bib{CDn}, Eqn. (5.3). See also Remark \refer{r: Eis and Delorme}.
\end{rem}

\begin{rem}
\naam{r: on C P P}
Note that it follows from the characterization of the normalized Eisenstein integral
in Proposition \refer{p: intro nE P psi} that
$$
\nC_{P|P}(1\col \nu) = I,\qquad (\nu \in \faPqdc).
$$
\end{rem}

\section{The Maass--Selberg relations}\eqnreset
\naam{s: the Maass--Selberg relations}
In this section we derive the Maass--Selberg relations for the normalized
$C$-functions. As a first step we
use the vanishing theorem to prove the following
functional equation for the Eisenstein integral.

\begin{prop}\restypetwo
\naam{p: funct eqn Eis}
Let $P,Q \in \allparabs$ be associated parabolic subgroups.
Then, for each $s \in W(\faQq \mid \faPq)$
and all $x \in \spX,$
$$
\nE(P \col \nu\col x) = \nE(Q \col s \nu\col x)\, \nC_{Q|P}(s\col \nu),
$$
as a meromorphic identity in $\nu \in \faPqdc.$
\end{prop}

\proof
We give the proof under the assumption that
$P$ is of residue type, see Remark \remRT.
Fix $\psi \in \cAtwoP.$
By Corollary \refer{c: constant term and C},
the function $\nu \mapsto C^\circ_{Q|P}(s: \nu)\psi$ belongs
to $\Mer(\faPqdc, \Hyp, \cAtwoQ),$ for some locally finite collection $\Hyp$
of $\gS_r(P)$-hyperplanes in $\faPqdc.$ It follows that the function
$\gl \mapsto C^\circ_{Q|P}(s: s^{-1}\gl)\psi$ belongs
to $\Mer(\faQqdc, s\Hyp, \cAtwoQ).$
Here $s\Hyp$ is a locally finite collection of $\gSr(Q)$-hyperplanes in $\faQqdc.$
By Proposition \refer{p: intro nE P psi},
the family $(\gl, x) \mapsto \nE(Q\col \gl \col x)\psi_Q$
belongs to $\cE^\hyp_Q(\spX\col \tau),$ for every $\psi_Q \in \cAtwoQ.$
We conclude that the family
\begin{equation}
\naam{e: family f defined as nE after C}
f: (\gl, x) \mapsto \nE(Q\col \gl\col x) C^\circ_{Q|P}(s: s^{-1}\gl)\psi
\end{equation}
belongs to $\cE^\hyp_Q(\spX\col \tau)$ as well.
For $\gl$ in the complement $\Omega$  of
a locally finite union of hyperplanes in $i \faQqd,$
the function $f_\gl $ belongs to $\cAtemp(\spX\col \tau),$
and its constant term along $(Q,v),$ for $v \in \QcW,$ is given by
$$
(f_\gl)_{Q,v}(ma) =
\sum_{t \in W(\faQq\mid \faQq)}
      a^{t\gl}\,[\pr_v \nC_{Q|Q}(t\col \gl)\nC_{Q|P}(s \col s^{-1}\gl)\psi](m);
$$
see Corollary \refer{c: constant term and C}.
Taking Remark \refer{r: on C P P} into account, we see that
$q_{\gl - \rho_Q}(Q,v\mid f_\gl)(X,m) =
[\pr_v \nC_{Q|P}(s \col s^{-1}\gl)\psi](m),$ for all $\gl\in \Omega,$
$X \in \faQq$ and
$m \in \spXQvp.$ By application of \bib{BSanfam}, Thm.~7.7, Eqn.~(7.14),
it follows  that
\begin{equation}
\naam{e: q of f and C}
q_{\bar 1, 0}(Q,v\mid f, X)(\gl, m) = [\pr_v \nC_{Q|P}(s \col s^{-1}\gl)\psi](m),
\end{equation}
for generic   $\gl \in i\faQqd,$ $X \in \faQq$ and $m \in \spXQvp.$
By meromorphy this actually holds as an identity of meromorphic
functions in $\gl.$

On the other hand, it follows from Definition \refer{d: dual Eisenstein new}
combined with Lemma \refer{l: Weyl group on cEhyp}
that the family
\begin{equation}
\naam{e: family g defined as nE of sinvgl}
g: (\gl, x) \mapsto \nE(P\col s^{-1}\gl\col x)\psi
\end{equation}
belongs to $\cE^\hyp_Q(\spX\col \tau).$
Moreover, for $\gl$ in the complement of a locally finite union of hyperplanes in
$i \faQqd,$ the function $g_\gl$ belongs
to $\cAtemp(\spX\col \tau),$ and its constant term along $(Q,v)$
is given by
$$
(g_\gl)_{Q,v}(ma) = \sum_{t \in W(\faQq\mid \faPq)} a^{ts^{-1}\gl} \,
[\pr_v\nC_{Q|P}(t\col s^{-1}\gl)\psi](m).
$$
This implies that, for every
$X \in \faQq,$
\begin{equation}
\naam{e: q of g and C}
q_{\bar 1, 0}(Q,v \mid g, X)(\gl, \dotvar) =
\pr_v\nC_{Q|P}(s\col s^{-1}\gl)\psi,
\end{equation}
as a meromorphic identity in $\gl \in \faQqdc.$
{}From (\refer{e: q of f and C}) and (\refer{e: q of g and C})
it follows that the family $f - g \in \cE^\hyp_Q(\spX\col \tau)$
satisfies the hypothesis of the vanishing theorem, Theorem \refer{t: vanishing thm new}.
Hence, $f = g.$ It follows that the meromorphic $\Ci(\spX\col \tau)$-valued
function
$\nu \mapsto f_{s \nu} - g_{s \nu}$ on $\faPqdc$ is zero. This implies the
result, in view of (\refer{e: family f defined as nE after C}) and
(\refer{e: family g defined as nE of sinvgl}).
\qed

\begin{cor}\restypetwo
\naam{c: funct eqn dEis}
Let $P,Q\in \allparabs$ be associated parabolic subgroups.
Then, for each $s \in W(\faQq|\faPq)$ and all $x\in X,$
$$
\dE(P \col \nu\col x) = \nC_{Q|P}(s\col -\bar \nu)^*\,\dE(Q \col s \nu\col x),
$$
as a meromorphic identity in $\nu \in \faQqdc.$
\end{cor}

\proof
We
give the proof under the assumption that  $P$ and $Q$ are of residue type,
see Remark \remRT.
The result follows from Proposition \refer{p: funct eqn Eis} combined
with Definition \refer{d: dual Eisenstein new}.
\qed

We shall now derive the Maass--Selberg
relations for the normalized $C$-functions from the invariance properties
of the kernel $K_P,$ formulated in Theorem \refer{t: invariance kernel}.

\begin{thm}\restypetwo
\naam{t: Maass Selberg}
Let $P,Q \in \allparabs$ be associated.
Then for each $s \in W(\faQq\mid \faPq),$
\begin{equation}
\naam{e: Maass Selberg}
\nCQP(s \col \nu)\nCQP (s \col - \bar \nu)^* = I,
\end{equation}
as an identity of meromorphic $\End(\cAtwoQ)$-valued functions in the variable
$\nu \in \faPqdc;$
\end{thm}

\begin{rem}
For the case of the group the above result was announced by
Harish-Chandra in \bib{HCeis}, Thm.\ 6, with a proof appearing in
\bib{HC3}, Lemma 17.1 (see also Remark \refer{r: comparison with nE of HC}).
For the Riemannian case $H = K,$ which is a special case of that of the
group, the relations were proved in \bib{HelMS}, Thm.\ 6.6.

For general reductive symmetric spaces
and minimal $P$ the result is due to \bib{Bps2}, Thm.~16.3, combined with
\bib{Bsd},
in view of \bib{BSft}, text after Eqn.~(55).
For general $P$ the result is due to J.\ Carmona and P.\ Delorme, \bib{CDn}, Thm.~2 and Prop.~5 (vi).
See also Remark
\refer{r: Eis and Delorme}.
\end{rem}

\proof
We give the proof under the assumption
that $P$ and $Q$ are of
residue type, see Remark \remRT.
It follows from Definition \refer{d: K P} that
\begin{equation}
\naam{e: K Q as Eis and dEis}
|W_Q| \,K_Q(s\nu\col x\col y) = \nE(Q\col s\nu\col x) \dE(Q \col s\nu \col y),
\end{equation}
for all $x,y \in \spX,$ as an identity of meromorphic functions in $\nu \in \faPqdc.$
On the other hand, from the mentioned definition combined with
Proposition \refer{p: funct eqn Eis} and
Corollary \refer{c: funct eqn dEis} it follows that
\begin{eqnarray}
\lefteqn{|W_P|\,K_P(\nu \col x\col y) =
\nE(P\col \nu \col x) \dE(P\col \nu\col y)}
\nonumber\\
&=& \nE(Q\col s\nu \col x)\nCQP(s \col \nu)
\nCQP (s \col - \bar \nu)^* \dE(Q\col s\nu\col y),
\naam{e: KQ to Eis two}
\end{eqnarray}
for $x,y \in \spX,$ and generic $\nu \in \faPqdc.$
Now $|W_P| = |W_Q|$ since $P$ and $Q$ are associated.
Hence, using Theorem \refer{t: invariance kernel}
we  infer that
\begin{eqnarray}
\lefteqn{\nE(Q\col s\nu \col x)\dE(Q\col s\nu\col y)}\nonumber\\
\naam{e: eqality Eis for MS}
&=&\nE(Q\col s\nu \col x)\nCQP(s \col \nu)
\nCQP (s \col - \bar \nu)^* \dE(Q \col  s\nu\col y)
\end{eqnarray}
for all $x,y \in X,$ as identities of meromorphic functions
in the variable $\nu\in \faPqdc.$ Let $v \in \QcW,$
and take $x = ma v$ with $m \in \spXQvp$ arbitrary and
$a$ tending to infinity in $\AQqp.$ Comparing
the coefficients of $a^{s\nu - \rho_Q}$ in the asymptotic expansions
along $(Q,v)$ of the resulting expressions on both sides of
(\refer{e: eqality Eis for MS}), using Definition
\refer{d: dual Eisenstein new}, we obtain
that
$$
\prQv \after  \dE(Q\col s\nu\col y)
=
\prQv \after \nCQP(s \col \nu)
\nCQP (s \col - \bar \nu)^* \dE(Q\col  s \nu\col y),
$$
for all $y \in \spX,$ as an identity of meromorphic functions in the variable
$\nu \in \faPqdc.$ Taking adjoints and substituting $-\bar \nu$ for $\nu$
we now obtain
\begin{equation}
\naam{e: last equality Eis in char MS}
\nE(Q\col s\nu\col y)\after \rmiQv
=
\nE(Q \col s\nu\col y) \nCQP (s \col \nu) \nCQP(s \col -\bar \nu)^*\after  \rmiQv.
\end{equation}
Fix $w \in \QcW,$ and put $y = m a w,$ with $m \in \spX_{Q,w,+}$ arbitrary
and $a \in \AQqp$ tending to infinity. Comparing the coefficients
of $a^{s\nu - \rho_Q}$ in the expansions along $(Q,w)$ of the functions
on both sides of (\refer{e: last equality Eis in char MS}), we obtain
\begin{equation}
\naam{e: last step MS char}
\pr_{Q,w}\after \rmiQv =
\pr_{Q,w}\after \nCQP (s \col \nu) \nCQP(s \col -\bar \nu)^* \after \rmiQv,
\end{equation}
as a meromorphic identity in the variable $\nu \in \faPqdc.$
This holds for arbitrary $v,w \in \QcW;$
in view of the direct sum decomposition (\refer{e: defi cAtwoP})
with $Q$ in place of $P,$
the equality
(\refer{e: last step MS char})
therefore remains valid if the maps $\pr_{Q,w}$ and $\rmiQv$ are replaced by the
identity map of $\cAtwoQ.$
\qed

\begin{rem}
Conversely, if the Maass--Selberg relations
(\refer{e: Maass Selberg}) hold, then the expression
on the right-hand side of (\refer{e: KQ to Eis two})
equals the one on the right-hand side of (\refer{e: K Q as Eis and dEis});
hence (\refer{e: invariance kernel}), the invariance property  of the kernel $K_P,$
follows. Thus, the Maass--Selberg relations are equivalent to the invariance
properties of the kernel.
\end{rem}

\begin{cor}\restypetwo
\naam{c: regularity C}
Let $P,Q \in \allparabs$ be associated parabolic subgroups.
and let $s \in W(\faQq\mid \faPq).$
Then there exists a constant $\geps >0$ such that the meromorphic
$\Hom(\cAtwoP, \cAtwoQ)$-valued function $\nC_{Q|P}(s\col \dotvar)$ is regular
on $\faPqd(\geps).$
\end{cor}

\proof
We give the proof under the assumption that
$P$ and $Q$ are of residue type,
see Remark \remRT.
The corollary is a straightforward consequence of
Corollary \refer{c: constant term and C}
and Theorem \refer{t: Maass Selberg}, combined with the lemma
below.
\qed

\begin{lemma}
Let $P \in \allparabs.$
Let $V$ be a complete locally convex space, and let
$\gf$ be a $V$-valued meromorphic function on $\faPqdc,$
with singular locus $\sing \gf$ contained in a locally finite collection
of real $\gSr(P)$-hyperplanes. Assume that for every $\gl_0 \in i\faPqd$
there exists an open neighborhood $\omega$ of $\gl_0$ in $i\faPqd$
such that $\gf$ is bounded on $\omega \setminus \sing \gf.$
Then there exists a $\geps > 0$ such that $\gf$ is holomorphic on $\faPqd(\geps).$
\end{lemma}

\proof
Since the collection $\Hyp$ of singular hyperplanes of $\gf$ is a real
$\gSr(P)$-configuration, the number of $H \in \Hyp$ with
$H \cap \faPqd(\geps') \cap H \neq \emptyset$
is finite, for every $\geps' > 0.$ Hence, there
exists a $p \in \Pi_{\gSr(P),\R}(\faPqd)$ of minimal degree
such that $\gf_p: \gl \mapsto p(\gl)\gf(\gl)$ is holomorphic
on a neighborhood of $i\faPqd.$ Clearly $\gf_p$
is holomorphic on $\faPqd(\geps),$ for a suitable $\geps >0.$
Assume that $\deg p \geq 1.$
Then there exists a $\ga \in \gSr(P)$ and a constant $c \in \R$
such that $l: \gl \mapsto \inp{\ga}{\dotvar} - c$ is a divisor of $p.$
By minimality of $p$
it follows that $h:= i\faPqd \cap l^{-1}(0)$ is non-empty.
{}From the hypothesis we infer that $\gf_p = 0$ on $h.$ By analytic continuation
it follows that $\gf_p = 0$ on
$h_\iC \cap \faPqd(\geps) = l^{-1}(0) \cap \faPqd(\geps).$
By a straightforward argument involving power series expansion
in the coordinate function $l,$ it now follows
that $l^{-1}\gf_p$ is holomorphic on $\faPqd(\geps).$
This contradicts the minimality of $p.$ Hence, $\deg p = 0$ and the result follows.

\qed

\begin{thm}\restypetwo
\naam{t: regularity nE}
Let $P\in \allparabs.$
Then there exists a constant $\geps >0$ such that $\nu \mapsto \nE(P\col \nu)$
is a holomorphic  $\Ci(\spX, \Hom(\cAtwoP, \Vtau))$-valued function
on $\faPqd(\geps).$
\end{thm}

\begin{rem}
For the group case the above result is due to Harish-Chandra \bib{HC3}.
For general reductive symmetric spaces and for $P$ minimal, the  result is due to
\bib{BSft}, Thm.~2. For non-minimal
$P$ it is due to \bib{CDn}, Thm.~3(i).
\end{rem}

\proof
We give the proof under the assumption that
$P$ is of residue
type, see Remark \remRT.
Let $\gL \in \rmL_P(\fb, \tau)$ and fix $\psi \in \cAtwoP(\gL)$ and $\eta \in \Vtaud.$
Let $\geps_0>0;$ then the family $F: (\gl, x) \mapsto \eta \nE(P\col \gl\col x)\psi$
belongs to the space ${\rm II}_{\mer}(\gL, \geps_0)$
defined in \bib{BCDn}, Def.~3.
{}From Corollary \refer{c: constant term and C} and Corollary \refer{c: regularity C}
it follows that $F$ satisfies the hypotheses of \bib{BCDn}, Thm.~2.
Hence, there exists a $\geps_1 >0$ such that $F$ belongs to the space
${\rm II}'(\gL, \geps_1).$ In particular, this implies that
$\gl \mapsto F_\gl$ is holomorphic on $\faPqd(\geps)$ for some $\geps >0.$
The theorem now follows by linearity and finite dimensionality of $\Vtau$ and $\cAtwoP.$
\qed

\begin{prop}\restypetwo
\naam{p: uniform tempered estimate nE}
Let $P\in \allparabs.$
Then there
exist constants $\geps > 0$ and $s >0$ and for every $u \in U(\fg)$ constants $n \in \N$
and $C > 0,$ such that the function
$\nu \mapsto \nE(P\col \nu)$
is a holomorphic  $\Ci(\spX, \Hom(\cAtwoP, \Vtau))$-valued function
on $\faPqd(\geps)$ satisfying the estimate
$$
\|\nE(P\col \nu \col u;x)\| \leq C |(\nu, x)|^n  \Theta(x) e^{s |\Re \nu|\lspX(x)},
\qquad
(\nu \in \faPqd(\geps), \;x \in \spX).
$$
\end{prop}

\proof
We give the proof  under the assumption that $P$
is of residue type, see Remark \remRT.
By finite dimensionality of $\Vtau$ and $\cAtwoP,$ it follows from Theorem
\refer{t: norm Eis uniform tempered}
and Definition \refer{d: space uniform tempered functions}
that there exists a $p \in \Pi_{\gSr(P), \R}(\faPqd)$ and constants
$\geps >0$ and $s> 0$ such that $\nu \mapsto p(\nu)\nE(P\col \nu)$
is a holomorphic function on $\faPqd(\geps),$ with values in
$\Ci(\spX)\otimes \Hom(\cAtwoP, \Vtau).$ Moreover, it
satisfies the following
estimates. For every $u \in U(\fg)$ there exist constants $n \in \N$ and $C >0$
such that
\begin{equation}
\naam{e: estimate nE with p}
\|p(\nu)\nE(P\col \nu \col u ; x)\| \leq C |(\nu,x)|^n \Theta(x) e^{s|\Re \nu|\lspX(x)},
\qquad(\nu \in \faPqd(\geps), \; x \in \spX).
\end{equation}
If we choose $\geps > 0$ sufficiently small, then by Theorem \refer{t: regularity nE},
the function $\nu \mapsto \nE(P\col \nu)$ is already holomorphic
on $\faPqd(\geps).$ By a straightforward application of Cauchy's integral formula,
involving polydiscs of size $O((1 + \lspX(x))^{-1}),$ it follows
that for $\geps >0$ sufficiently small, the following is true.
For every $u \in U(\fg)$ there exist $n \in \N$ and $C > 0$ such that
the estimate (\refer{e: estimate nE with p}) holds with $p =1.$
\qed

\begin{cor}\restypetwo
\naam{c: uniform tempered estimate dE}
Let $P\in \allparabs.$
Then there
exist constants $\geps > 0$ and $s >0$ and for every $u \in U(\fg)$ constants $n \in \N$
and $C > 0,$ such that the function
$\nu \mapsto \dE(P\col \nu)$
is a holomorphic  $\Ci(\spX, \Hom(\Vtau, \cAtwoP))$-valued function
on $\faPqd(\geps)$ satisfying the estimate
$$
\|\dE(P\col \nu \col u;x)\| \leq C |(\nu, x)|^n  \Theta(x) e^{s |\Re \nu|\lspX(x)},
\qquad
(\nu \in \faPqd(\geps), \;x \in \spX).
$$
\end{cor}

\proof
We give the proof under the assumption
that $P$ is of residue type, see  Remark \remRT.
In
view of Definition \refer{d: dual Eisenstein new}, the result
follows from
Proposition \refer{p: uniform tempered estimate nE},
\qed

\begin{cor}\restypetwo
\naam{c: second estimate dE}
Let $P \in \allparabs.$
Then, for all $U \in S(\faPqd)$ and $u \in U(\fg),$
there exist constants $m \in \N$ and $C > 0$ such
that
$$
\|\dE(P\col \nu ; U\col u ; x)\| \leq
 C |(\nu,x)|^m \Theta(x), \qquad (\nu \in i\faPqd,\; x \in \spX).
$$
\end{cor}

\proof
We  give the proof under the assumption that
$P$ is of residue type, see Remark \remRT.
The result follows from the estimate of the previous corollary,
by a straightforward
application of Cauchy's integral formula involving polydiscs
of size $O(( 1 + \lspX(x))^{-1}).$
\qed

\section{The spherical Fourier transform}\eqnreset
\naam{s: spherical Fou}
We recall from \bib{Bps2}, Cor.~17.6, that there exists a constant
$N \in \N$ such that
\begin{equation}
\naam{e: Theta in Lone}
(1 + \lspX)^{-N}\Theta^2 \in L^1(\spX).
\end{equation}

Combining the estimate (\refer{e: Schwartz estimate})
with (\refer{e: Theta in Lone}) and
the estimate of
Corollary \refer{c: second estimate dE},
we see that the integral in the following definition converges absolutely.

\begin{defi}\restypetwo
\naam{d: Fou P}
Let $P\in \allparabs.$
If $f \in \cC(\spX\col \tau),$
we define its Fourier transform $\cF_P f: i\faPqd \to \cAtwoP$ by
$$
\Fou_P f (\nu) = \int_X \dE(P\col \nu\col x)\, f(x) \; dx, \qquad (\nu \in i\faPqd).
$$
\end{defi}

The above definition depends on the validity of
the estimate of Corollary \refer{c: second estimate dE}.
Thus, within the chain of reasoning leading up to Theorem \thmendRT{},
the use of this definition requires $P$ to be of residue type,
see Remark \remRT.

\begin{rem}
\naam{r: Fou if G compact center mod H}
If $G$ has compact center modulo $H,$ then $\fa_{G\iq} = \{0\}$
and $\cA_{2,G} = \cA_2(G/H\col \tau).$ Moreover, using
Remark \refer{r: normalized Eisenstein for P is G} we infer that $f \mapsto \cF_Gf(0)$
is the restriction to $\cC(\spX\col \tau)$ of the orthogonal
projection $L^2(\spX\col \tau) \to L^2_d(\spX \col \tau) = \cA_{2,G}.$
\end{rem}

\begin{lemma}
\restypetwo
\naam{l: continuity Fou into Ci}
Let $P \in \allparabs.$
Then for every $U \in S(\faPqd)$ there exists a constant
$m \in \N$ and a continuous seminorm $s$ on $\cC(\spX\col \tau)$ such
that
$$
\sup_{\nu \in i\faPqd} ( 1 + |\nu|)^{-m}\| \Fou_P f(\nu ; U) \|\leq s(f),
$$
for all $f \in \cC(\spX\col \tau).$
In particular, the Fourier transform
$\cF_P$ maps $\cC(\spX\col \tau)$ continuous linearly
into $\Ci(i\faPqd) \otimes \cAtwoP.$
\end{lemma}

\proof
We give the proof under the assumption
that $P$ is of residue type, see Remark \remRT.
The result follows from the estimates (\refer{e: Theta in Lone}) and
(\refer{e: Schwartz estimate})
combined with the estimate of Corollary \refer{c: second estimate dE}.
\qed

\begin{lemma}\restypetwo
\naam{l: Fou after D}
Let $P\in \allparabs.$
Then for every $D \in \DX$ and all $f \in \cC(\spX\col\tau),$
\begin{equation}
\naam{e: Fou after D}
\Fou_P (Df)(\nu) = \umuP(D\col \nu) \Fou_P f(\nu),\qquad (\nu \in i\faPqd).
\end{equation}
\end{lemma}

\proof
We give  the proof under the assumption that $P$ is of residue type,
see Remark \remRT.
{}From \bib{Bas}, Lemma 7.2, we recall that every $D \in \DX$ acts by
a continuous linear endomorphism on $\cC(\spX\col \tau).$ Since
$\Fou_P: \cC(X\col \tau) \to \Cci(i\faPqd) \otimes \cAtwoP$ is continuous,
it suffices to prove the identity (\refer{e: Fou after D}) for $f$ in the dense subspace
$\Cci(\spX\col \tau)$ of $\cC(\spX\col \tau).$ For such $f$ the identity
is an immediate consequence of Lemma \refer{l: D on dE}.
\qed

Let $\Omega$ be the image in $\DX$ of the Casimir operator defined by the bilinear form
$B$ on $\fg,$ see Section \refer{s: notation}.

\begin{lemma}\restypetwo
\naam{l: estimate umu inverse}
Let $P \in \allparabs$
and let  $\geps >0.$
\begin{enumerate}
\itema
$\umu_P(\Omega\col \nu) = -|\nu|^2 \, I + O(|\nu|)$
as $\nu \in \faPqd(\geps),$ $|\nu| \to \infty.$
\vspace{0.5mm}
\itemb
There
 exists a constant $R>0$ such that for every
$\nu \in \faPqd(\geps)$ with $|\nu| \geq R$ the endomorphism $\umu_P(\Omega\col \nu)$
is invertible and the operator norm of its inverse
satisfies the  estimate
\begin{equation}
\naam{e: estimate umu Omega}
|\nu|^2\, \|\umu_P(\Omega \col \nu)^{-1}\| \leq 2.
\end{equation}
\end{enumerate}
\end{lemma}

\proof
We give the proof under the assumption that $P$ is of residue type,
see Remark \remRT.
Let $\inp{\dotvar}{\dotvar}$ denote the complex bilinear form on $\faPqdc$
that extends the dual of the given  bilinear form on $\faPq.$
There exists a first order polynomial function $\umu_1: \faPqdc \to \End(\cAtwoP)$
such that
\begin{equation}
\naam{e: formula umu P Omega}
\umu_P(\Omega \col \nu) = \inp{\nu}{\nu} I + \umu_1(\nu), \qquad (\nu \in \faPqdc).
\end{equation}
Indeed, this follows by application of Corollary
\refer{c: deco cAPtwo in eigenspaces}.
It follows by a straightforward estimation
that, for $\nu \in \faPqd(\geps),$
\begin{equation}
\naam{e: form as norm plus cO}
\inp{\nu}{\nu} = - {|\nu|^2} + O(|\nu|)
\text{as}\quad  |\nu| \to \infty.
\end{equation}
Using that $\|\mu_1(\nu)\| = O(|\nu|)$ we obtain (a) from
(\refer{e: formula umu P Omega}) and (\refer{e: form as norm plus cO}).
{}From (a) it follows that $-|\nu|^{-2} \umu_P(\Omega\col \nu) =
I + O(|\nu|^{-1}).$
Hence, (b) follows.
\qed

In the following result, $\cS(i\faPqd)$ denotes the Euclidean Schwartz space
of $i\faPqd.$

\begin{prop}\restypetwo
\naam{p: Fou continuous on Schwartz space}
Let $P\in \allparabs.$
Then
the
Fourier transform $\Fou_P$ maps $\cC(\spX\col \tau)$ continuous linearly
into $\cS(i\faPqd) \otimes \cAtwoP.$
\end{prop}

\proof
We give the proof under the assumption
that $P$ is of residue type, see Remark \remRT.
Moreover, we use the argumentation of \bib{Bps2}, p.~436, completion of
the proof of Theorem 19.1, with $\Fou_P$ in place of $\Fou_\pi.$
Let us
label the first two displayed formulas
in the mentioned text in \bib{Bps2} by (E1) and (E2),
respectively.
The estimate (\refer{e: estimate umu Omega})
generalizes the estimate of \bib{Bps2}, Lemma 19.4.

Let $R>0$ be as in Lemma \refer{l: estimate umu inverse} (b) for some
choice of $\geps >0$  and let
$u \in S(\faPqd)$ and $M \in \N.$ In view of the last assertion
of Lemma \refer{l: continuity Fou into Ci} it suffices
to prove the analogue of (E1), i.e., it suffices to prove the existence
of a continuous seminorm $s$ on $\cC(\spX\col \tau)$ such that
$$
\| \Fou_Pf(\nu; u) \| \leq ( 1 + |\nu|)^{-M} s(f),
$$
for every $f \in \cC(\spX\col \tau)$ and all $\nu \in i\faPqd$ with $|\nu | \geq R.$
As in \bib{Bps2} this is done by induction on the order of $u,$ by using
Lemma \refer{l: continuity Fou into Ci} instead of (E2)
and by using Lemmas
\refer{l: Fou after D} and \refer{l: estimate umu inverse}
instead of
\bib{Bps2}, Lemmas 19.3 and  19.4, respectively.
\qed

We end this section with a result on the Fourier transform of a compactly
supported smooth function.
If $S > 0$ we write
$$
\spX_S: = \{ x \in \spX \mid \lspX(x) \leq S\}.
$$
Then $\spX_S$ is a $K$-invariant compact subset of $\spX.$ We write
$
C^\infty_S(\spX\col \tau)$ for the closed subspace of $C^\infty(\spX\col \tau)$ consisting
of functions with support contained in $\spX_S.$

\begin{prop}\restypetwo
\naam{p: Fou of compactly supported function}
Let $P \in \allparabs$
and let
$\geps>0$ be as in Corollary \refer{c: uniform tempered estimate dE}.
For every $f \in \Cci(\spX\col \tau),$ the Fourier transform $\Fou_P(f)$
extends to a holomorphic function on $\faPqd(\geps)$ with values in $\cAtwoP.$
Moreover, let $S > 0.$ Then for every $m \in \N$
there exists a continuous seminorm $p_m$ on $C^\infty_S(\spX\col \tau)$
such that, for every $f \in C^\infty_S(\spX\col \tau),$
\begin{equation}
\naam{e: estimate Fou of compactly supported f}
\| \Fou_P f(\nu) \| \leq (1 + |\nu|)^{-m} p_m(f),
\qquad (\nu \in \faPqd(\geps)).
\end{equation}
\end{prop}

\proof
We give the proof under the assumption that $P$ is of residue type,
see Remark \remRT.
The assertion about holomorphy is a straightforward consequence of the
holomorphy of the Eisenstein integral as formulated in Corollary
\refer{c: uniform tempered estimate dE}.
Let $n\in \N$ be the constant of the mentioned corollary  associated with
$u=1.$ Let $S > 0.$
Then it straightforwardly follows from the estimate
of Corollary \refer{c: uniform tempered estimate dE}
that there exists
a continuous seminorm $p_0$ on $C^\infty_S(\spX\col\tau)$
such that, for every $f \in C^\infty_S(\spX\col \tau),$
\begin{equation}
\naam{e: second estimate Fou of compactly supported f}
\|\Fou_P f(\nu)\|\leq (1 + |\nu|)^n p_0(f),\qquad
(\nu \in \faPqd(\geps)).
\end{equation}
Let $R > 0$ be associated with $\geps >0$ as in Lemma \refer{l: estimate umu inverse}.
Then it follows from the above estimate that, for every $k \in \N$ and
for $\nu \in \faPqd(\geps)$ with $|\nu |\geq R,$
\begin{eqnarray}
|\nu|^{2k} \|\Fou_Pf(\nu)\| & = &
 |\nu|^{2k} \|\umu_P(\Omega\col\nu)^{-k}\Fou_P(\Omega^k f)(\nu)\|\\
&\leq & ( 1 + |\nu|)^n \,p_0(2^k \Omega^k f).
\end{eqnarray}
Taking $k\in \N$ such that $n - 2k \leq m$ we see that
there exists a continuous seminorm
$p_m'$ on $C_S(\spX\col \tau)$ such that
for every $f \in C^\infty_S(\spX \col \tau)$
the estimate (\refer{e: estimate Fou of compactly supported f})
holds for all $\nu \in \faPqd(\geps)$ with
$|\nu| \geq  R.$
{}From (\refer{e: second estimate Fou of compactly supported f})
it follows that there exists a constant $C>0$ such that the estimate
(\refer{e: estimate Fou of compactly supported f})
holds with $C p_0$ in place of $p_m,$ for all $\nu \in \faPqd(\geps)$ with $|\nu| < R.$
Take for $p_m$ any continuous
seminorm with
$p_m \geq \max(C p_0, p_m');$ then the desired assertion follows.
\qed

We end this section with another useful result.

\begin{lemma}\restypetwo
\naam{l: Fou vanishes on cA two}
Let $P \in \allparabs$
 and assume that $\faPq \neq 0.$
Then $\Fou_P$ vanishes on $\cA_2(\spX\col \tau).$
\end{lemma}

\proof
We give the proof under the assumption
that $P$ is of residue type,
see Remark \remRT.
Fix $f \in \cA_2(\spX\col \tau).$ Then there exists
a non-trivial polynomial $q$ in one variable such that
$q(\Omega) f = 0.$ In view of Lemma \refer{l: Fou after D} it follows that
$q(\umu_P(\Omega\col \nu)) \Fou_P f(\nu) = 0$ for all $\nu \in i\faPqd.$
{}From Lemma \refer{l: estimate umu inverse} it follows that the polynomial
function $\nu \mapsto \det q(\umu_P(\Omega\col \nu))$ is not identically zero.
Hence, $\Fou_P f$ vanishes on an open dense subset of $i \faPqd.$ By smoothness
of $\Fou_P f$ it follows that $\Fou_P f = 0.$
\qed

\section{The wave packet transform}\eqnreset
\naam{s: Wave packet}
It follows from the estimate in Proposition \refer{p: uniform tempered estimate nE}
that the integral in the following definition is absolutely convergent.
We agree to write $d\nu$ for the Lebesgue measure $d\mu_P(\nu)$ on $i\faPqd,$
normalized as in Section \refer{s: normalization}.

\begin{defi}\restypetwo
\naam{d: wave packet}
Let $P \in \allparabs.$
Then for every $\gf \in \cS(i\faPqdc)\otimes \cAtwoP,$
we define the wave packet  $\Wave_P \gf: \spX \to \Vtau$ by
$$
\Wave_P \gf(x) := \int_{i\faPqd} \nE(P\col\nu\col x)\, \gf(\nu) \; d\nu\qquad (x \in \spX).
$$
\end{defi}

This definition  depends on the validity
of the estimate of Proposition \refer{p: uniform tempered estimate nE},
which in the chain of reasoning leading up to
Theorem \thmendRT{} requires $P$ to be of residue type,
see Remark \remRT.

Note that the wave packet $\Wave_P \gf$  is smooth  $\tau$-spherical function.

\begin{rem}
\naam{r: Wave if G compact center mod H}
If $G$ has compact center modulo $H,$ then $A_{G\iq} = \{0\}$ and
$\cA_{2,G} = \cA_2(G/H \col \tau).$ In this case the measure $d\nu = d\mu_G$
has total volume $1$ (end of Section \refer{s: normalization}),
and using Remark \refer{r: normalized Eisenstein for P is G}
we infer that $\Wave_G \gf = \gf(0).$
Accordingly, $\Wave_G$ is naturally identified with the inclusion map
$\cA_2(\spX\col \tau) \to \Ci(\spX\col \tau).$
\end{rem}

\begin{thm}\restypetwo
\naam{t: Wave continuous on Schwartz space}
Let $P \in \allparabs$
Then the wave packet map $\Wave_P$ maps $\cS(i\faPqd) \otimes \cAtwoP$
continuous linearly into the Schwartz space $\cC(\spX\col \tau).$
\end{thm}

\proof
We give the proof under the assumption that $P$ is of residue type,
see Remark \remRT.
Let $\gL \in \rmL_P(\fb, \tau)$ and fix $\psi \in \cAtwoP(\gL).$
We recall from the proof of Theorem \refer{t: regularity nE} that the family $F$ defined
by
$$
F(\nu, x) = \nE(P\col \nu\col x)\psi
$$
has components with respect to a basis of $\Vtau$ that are functions
of type ${\rm II}'(\gL)$ in the sense of \bib{BCDn}. Hence, by \bib{BCDn}, Thm.~1,
the map $\ga \mapsto \cW_{\ga,F},$
where
$$
\cW_{\ga, F}(x) = \int_{i\faPqd} \ga(\nu) F(\nu , x) \; d\nu,
$$
is continuous linear from $\cS(i\faPqd)$ into $\cC(\spX\col \tau).$
We note that $\cW_{\ga, F} = \Wave_P (\ga \otimes \psi).$ Hence, the result
follows by using linearity,  the finite dimensionality of
$\cAtwoP$
(see Corollary \refer{c: deco cAPtwo in eigenspaces})
and the decomposition (\refer{e: deco cA two P in eigenspaces}) of the latter space.
\qed

Let $P \in \allparabs$ and $D \in \DX.$
In the following lemma we write $\umu_P(D)$ for the endomorphism
of $\cS(i\faPqd) \otimes \cAtwoP$ given by
$$
[\umuP(D)\gf](\nu) = \umuP(D\col \nu)(\gf(\nu)),
$$
for $\gf \in \cS(i\faPqd)\otimes \cAtwoP$ and $\nu \in i\faPqd.$

\begin{lemma}\restypetwo
\naam{l: D after Wave}
Let $P \in \allparabs$
and  $D \in \DX.$ Then
\begin{equation}
\naam{e: D after Wave}
 D \after \Wave_P = \Wave_P \after \umuP(D)\quad \text{on}\quad \cS(i\faPqd)\otimes \cAtwoP.
\end{equation}
\end{lemma}

\proof
We give the proof under  the assumption that
$P$ is of residue type, see
Remark \remRT.
The operator $D$ defines a continuous linear endomorphism of $\cC(\spX\col \tau),$
by \bib{Bas}, Lemma 7.2. In view of Theorem \refer{t: Wave continuous on Schwartz space}
it follows that both sides of the equation are continuous linear maps from
$\cS(i\faPqd)\otimes\cAtwoP$ to $\cC(\spX\col \tau).$ Hence, by density it suffices
to prove the equality when
applied to an element $\gf \in \Cci(i\faPqd) \otimes \cAtwoP.$ But then the result
is an immediate consequence of Lemma \refer{l: action of DX on nE}
by differentiation under the integral sign, in view of Definition \refer{d: wave packet}.
\qed

We equip $\cC(\spX\col \tau)$ with the restriction of the
$L^2$-inner product $\hinp{\dotvar}{\dotvar}$ from $L^2(\spX\col \tau).$
Similarly, for $P \in \allparabs,$ we equip $\cS(i\faPqd)\otimes  \cAtwoP$
with the restriction of the  $L^2$-type inner product
$\hinp{\dotvar}{\dotvar}$
from $L^2(i\faPqd)\otimes \cAtwoP.$
With respect to these structures the
Fourier transform $\Fou_P$ and the wave
packet map
$\Wave_P$ are adjoint in the following sense.

\begin{lemma}\restypetwo
\naam{l: Fou and Wave transpose}
Let $P \in \allparabs.$
Then the
continuous linear operators $\Fou_P: \cC(\spX\col \tau) \to \cS(i\faPqd) \otimes \cAtwoP$
and $\Wave_P:  \cS(i\faPqd) \otimes \cAtwoP \to \cC(\spX\col \tau)$
are adjoint in the sense that, for all $f \in  \cC(\spX\col \tau)$ and
$\gf \in \cS(i\faPqd) \otimes \cAtwoP,$
\begin{equation}
\naam{e: Fou and Wave adjoint}
\hinp{\Fou_P f}{ \gf} =
\hinp{f}{\Wave_P \gf}
\end{equation}
\end{lemma}

\proof
We give the proof under the assumption that $P$ is of residue
type,
see Remark \remRT.
By continuity and density it suffices to prove (\refer{e: Fou and Wave adjoint})
for all $f \in \Cci(\spX\col \tau)$ and $\gf \in \Cci(i\faPqd) \otimes \cAtwoP.$
For such $f$ and $\gf,$  the formula follows by an application of
Fubini's theorem.
\qed

\section{Fourier inversion for Schwartz functions}\eqnreset
\naam{s: Fourier inversion for Schwartz functions}
In this section we  show that the Fourier inversion formula
(\refer{e: inversion formula with T t F}),
established in \bib{BSfi}, implies an inversion formula for Schwartz
functions, formulated in terms of the Fourier transforms and
the wave packet maps introduced in the previous sections.

The crucial first step is the following.

\begin{prop}
\restypetwo
\naam{p: T t F written as Wave after Fou}
Let $F \subset \gD.$
Then for every $W$-invariant even residue weight $t \in \WT(\gS),$
\begin{equation}
\naam{e: T t F in terms of Wave and Fou}
T^t_F = [W:W_F]\, t(\faFqp)\, \Wave_F \Fou_F \text{on} \Cci(\spX\col \tau).
\end{equation}
\end{prop}

\proof
We give the proof under the assumption that  $P_F$
is of residue type, see  Remark \remRT.
In
case $F = \gD$ and
$G$ has compact center modulo $H,$ then $\fa_{F\iq} = \fa_{G\iq} = \{0\}$
and the proof below has to be read according to the conventions indicated in
Remarks \refer{r: on T gD}, \refer{r: Eisenstein for P is G},
\refer{r: normalized Eisenstein for P is G},
\refer{r: Fou if G compact center mod H} and
\refer{r: Wave if G compact center mod H}.

Let $f \in \Cci(\spX \col \tau)$ and $x \in \spXp.$
It follows from (\refer{e: defi T F t})
that
$$
T_F^t f (x) = |W|\, t(\faFqp)\,\int_{i\faFqd + \geps_F}
\int_X K^t_F(\nu\col x \col y) f(y) \; dy\,
 d\mu_F(\nu),
$$
for all $\geps_F \in \faFqp$ sufficiently close to zero.
In view of (\refer{e: K F t in terms of Eis})
and Definition \refer{d: Fou P}, this
equality may be
rewritten as
\begin{equation}
\naam{e: T F t in terms of Eis and Fou}
T_F^t f (x) = [W:W_F]\, t(\faFqp)\,
 \int_{i\faFqd + \geps_F} \nE(P_F\col\nu \col x) \Fou_F f(\nu) \; d\mu_F(\nu).
\end{equation}
Since the expressions on both sides of the equation extend smoothly to all
of $\spX$ in the variable $x,$ it follows that (\refer{e: T F t in terms of Eis and Fou})
holds for all $x \in \spX.$
{}From Theorem \refer{t: regularity nE}
and Proposition \refer{p: Fou of compactly supported function}
it follows that $\nu \mapsto \nE(P_F\col\nu \col x) \Fou_F f(\nu) $
is holomorphic on $\faFqd(\geps),$ for some $\geps > 0.$ Moreover, from the
mentioned
results it also follows
that
for every $N >0,$ there exists a constant $C_N >0$ such that
$$
 \| \nE(P_F\col\nu \col x) \Fou_F f(\nu)\| \leq C_N (1 + |\nu|)^{-N}, \qquad
(\nu \in \faFqd(\geps)).
$$
This estimate allows us to take the limit of
(\refer{e: T F t in terms of Eis and Fou}) for $\geps_F \to 0;$ thus, using
Definition \refer{d: wave packet} and observing that
$d\mu_F(\nu) =d\nu$ on $i\faFqd,$ we obtain (\refer{e: T t F in terms of Wave and Fou}).
\qed

The
proof of the following result involves an induction step
using the long chain of results marked (RT),
see Remark \remRT.

\begin{thm}{\ }
\naam{t: inversion for Schwartz}
\begin{enumerate}
\itema
Every  $P \in \allparabs$ is of residue type.
\itemb
If $t$ is any $W$-invariant even residue weight for $\gS,$ then
\begin{equation}
\naam{e: Fourier inversion in terms of Wave Fou}
f = \sum_{F \subset \gD} [W:W_F] \, t(\faFqp) \Wave_F \Fou_F f,
\end{equation}
for every $f \in \cC(\spX\col \tau).$
\itemc
The pair $(G,H)$ is of residue type
if and only if $G$ has a compact center modulo $H.$
\end{enumerate}
\end{thm}

\proof
We first show that (c) and (b) follow from (a). Thus, assume (a).
Then viewed as a parabolic subgroup, $G$ is of residue type. By
Remark \refer{r: G of residue type} it follows
that the pair
$({}^\circ G, {}^\circ G \cap H)$
is of residue type. Moreover, if $G$ has compact center modulo $H,$ then
$(G,H)$ is of residue type. If the center of $G$ is not compact modulo $H,$
then $(G,H)$ is not of residue type, by
Definition \refer{d: good residue type new} (a). This establishes (c).

We now turn to (b). Let $f$ be a $W$-invariant even residue weight for $\gS.$
Then for each $F \subset \gD,$ the parabolic subgroup $P_F$
is of residue type so that  Proposition \refer{p: T t F written as Wave after Fou} applies.
It now follows from (\refer{e: inversion formula with T t F}) combined with
(\refer{e: T t F in terms of Wave and Fou})  that
(\refer{e: Fourier inversion in terms of Wave Fou})
holds for every $f \in \Cci(\spX\col \tau).$ Finally, the validity of (b)
follows by density of $\Cci(\spX\col \tau)$ and continuity of each of
 the $\Wave_F\Fou_F,$ for $F \subset \gD$
(Proposition \refer{p: Fou continuous on Schwartz space} and
Theorem \refer{t: Wave continuous on Schwartz space} apply with $P = P_F,$
since $P_F$ is of residue
type).

It remains to prove (a). We will do this by induction
on $\dim \Aq,$  the $\gs$-split rank
of $G.$ First, assume that $\dim \Aq = 0.$ Then $\spX$ is compact,
hence the pair $(G,H)$ is of residue type,
see Remark \refer{r: compact is residue type}. It follows that
$G,$ viewed as a parabolic subgroup, is of residue type,
see Remark \refer{r: G of residue type}.
Since $\allparabs = \{G\},$ this establishes (a) in case $\dim \Aq = 0.$

Thus, assume that $\dim \Aq \geq 1$ and that (a) has
been established for all pairs $(G',\gs')$ with $G'$ of $\gs$-split rank
smaller than $\dim \Aq.$

If the center of $G$ is not compact modulo $H,$
then ${}^\circ \Aq:= {}^\circ G \cap \Aq \subsetneq \Aq.$ Hence,
it follows by the inductive hypothesis that every $\gs$-parabolic
subgroup of ${}^\circ G$ containing ${}^\circ \Aq $ is of residue type.
Denote the set of these parabolic subgroups by ${}^\circ\allparabs.$ Then
$G = {}^\circ G \times C,$
where $C = \exp ({\rm center}\,(\fg) \cap \fp)$ and $P \mapsto PC$
is a bijection from ${}^\circ\allparabs$ onto $\allparabs.$ Moreover,
$M_{PC} = M_P$ for every $P \in \allparabs({}^\circ G).$
In view of Definition \refer{d: parabolic of residue type} it follows
that every parabolic subgroup from $\allparabs$
is of residue type as well, whence (a).

Thus, we may assume that $G$ has compact center modulo $H;$ then
$\fa_{\gD\iq} = 0.$
By the inductive hypothesis, the symmetric pairs $(M_F, M_F \cap v H v^{-1}),$
for $F \subsetneq \gD$ and $v \in \NKaq,$ satisfy condition (a). In particular,
$M_F,$ viewed as a parabolic subgroup of $M_F,$ is of residue type relative
to $v H v^{-1}.$ Since ${}^\circ M_F = M_F,$ it follows by
Remark \refer{r: G of residue type} that the pairs $(M_F, M_F \cap v H v^{-1}),$
are all of residue type. In view of
Definition \refer{d: parabolic of residue type}
it now follows that the standard parabolic subgroups $P_F,$ for $F \subsetneq \gD,$
are all of residue type.

Let $t$ be a $W$-invariant even residue weight on $\gS.$
Then $T_\gD^t$ is a continuous linear operator from $\Cci(\spX\col \tau)$ into
the finite dimensional subspace $\cA^t_\gD(\spX\col \tau)$ of $\Ci(\spX\col \tau);$
all functions in this subspace are $\DX$-finite, see the text after
(\refer{e: defi cA t as span}).

Proposition \refer{p: T t F written as Wave after Fou}
applies for every subset $F \subsetneq \gD.$ Hence,
from  (\refer{e: inversion formula with T t F}) and
(\refer{e: T t F in terms of Wave and Fou}) it  follows that
\begin{equation}
\naam{e: T gD t in terms of Wave Fou F}
T_\gD^t = I - \sum_{F \subsetneq \gD} [W:W_F]\, t(\faFqp)\, \Wave_F \Fou_F
\end{equation}
as an operator from $\Cci(\spX\col \tau)$ into $\Ci(\spX\col \tau).$
Applying Proposition \refer{p: Fou continuous on Schwartz space}
and Theorem \refer{t: Wave continuous on Schwartz space}
we infer that $T_\gD^t$
extends to a continuous linear map from
$\cC(\spX\col \tau)$ into $\cA^t_\gD(\spX\col \tau) \cap
\cC(\spX\col \tau);$ moreover, the latter intersection
is continuously contained in $\cA_2(\spX\col \tau).$
By density of $\Cci(\spX\col \tau),$ the validity of the identity
(\refer{e: T gD t in terms of Wave Fou F}) extends to the space
$\cC(\spX\col \tau).$

By repeated application of  Lemma \refer{l: Fou vanishes on cA two},
with $P = P_F,$ $F \subsetneq \gD,$  it follows
from (\refer{e: T gD t in terms of Wave Fou F}) that
$T_\gD^t = I$ on $\cA_2(\spX\col \tau).$
Finally, by application of Lemma \refer{l: Fou and Wave transpose}
to (\refer{e: T gD t in terms of Wave Fou F})
it follows that $T_\gD^t$ is symmetric with
respect to the $L^2$-inner product on $\cC(\spX\col \tau).$ We conclude that
$T^t_\gD$ is the orthogonal projection from $\cC(\spX\col \tau)$ onto
$\cA_2(\spX\col \tau);$
in particular, it follows that the latter space is finite dimensional. Moreover, since
$\cA_2(\spX\col \tau)$ is dense in $L^2_d(\spX\col \tau)$ it follows that
$\cA_2(\spX\col \tau) = L^2_d(\spX\col \tau)$ and that $T^t_\gD$ is the restriction of
the orthogonal projection  $L^2(\spX\col \tau) \to L^2_d(\spX\col \tau).$
{}From this we conclude that $(G,H)$ is of residue type, see
Definition \refer{d: good residue type new} and Remark \refer{r: one K suffices for residue type}.
Hence, $P_\gD = G$ is of residue type.
It follows that all standard parabolic subgroups in $\allparabs$ are of residue type.
Since every $P \in \allparabs$ is associated with a standard one, see
Lemma \refer{l: on parabolics},
assertion (a) follows by application of Lemma \refer{l: residue type and sim}.
\qed

We define the equivalence relation $\sim$ on the collection
of subsets of $\gD$ by $F \sim F' \iff P_F \sim P_{F'}.$

\begin{lemma}
\naam{t: sum of residue weights}
Let $t$ be a $W$-invariant even residue weight on $\gS$
and let $F\subset \gD.$ Then
$$
\sum_{F'\subset \gD \atop F' \sim F}  t(\faFprqp) = |W(\faFq)|^{-1}.
$$
\end{lemma}

\proof
Let $\allparabs(\AFq)$ denote the collection of
all $P \in \allparabs$ with $\gs$-split component $\APq$ equal to $\AFq.$
Moreover, let $\cS$ denote the collection of all subsets $F'\subset \gD$ with
$F' \sim F.$
For every $P \in \allparabs(\AFq)$ there exists a unique
$F_P \subset \gD$ such that $sPs^{-1} = P_{F_P},$
for some $s \in W,$ see Lemma \refer{l: on parabolics}.
Clearly, $F_P \sim F.$ Moreover, the map $\rmp: P \mapsto F_P$
is surjective from $\allparabs(\AFq)$ onto $\cS.$ If $F' \in \cS,$
let $W_{F',F}$ denote the collection of $s \in W$ mapping $\faFq$ onto $\faFprq.$
Then the map $s \mapsto s^{-1} P_{F'} s$ from $W_{F',F}$ onto $\rmp^{-1}(F')$
factors to a bijection
from $W(\faFprq\mid \faFq)$ onto $\rmp^{-1}(F').$
Starting from (\refer{e: WT sum one}) with $Q = P_F$ we now obtain that
\begin{eqnarray*}
1 & = & \sum_{P\in \allparabs(\AFq)} t(\faPqp)\\
  & = & \;\;\;\;\,\sum_{F'\subset \gD \atop F' \sim F}\;\;\;
 \sum_{s \in W(\faFprq\mid \faFq)} t(s^{-1}(\faFprqp))\\
&=& \;\;\;\;\,
\sum_{F'\subset \gD \atop F' \sim F}  \;|W(\faFprq\mid \faFq)|\,t(\faFprqp).
\end{eqnarray*}
For every $F'\subset \gD$ with $F'\sim F,$
the group $W(\faFq)$ acts freely and transitively from the right
on  $W(\faFprq\mid \faFq).$
Hence, $|W(\faFprq\mid \faFq)| = |W(\faFq)|$
and the result follows.
\qed

\begin{lemma}
\naam{l: Wave Fou through class}
Let $P,Q\in \allparabs$ be associated parabolic subgroups
and let $s \in W(\faPq\mid \faQq).$ Then
\begin{enumerate}
\itema
$ \cF_P f(s\nu) = \nC_{P|Q}(s\col \nu)\Fou_Qf (\nu),$ for all
$f \in \cC(\spX\col \tau)$ and $\nu \in i\faQqd;$
\itemb
$\Wave_P\Fou_P = \Wave_Q\Fou_Q$ as endomorphisms of $\cC(\spX\col \tau).$
\end{enumerate}
\end{lemma}

\proof
It follows from Corollary \refer{c: funct eqn dEis}
combined with the Maass--Selberg relations (\refer{e: Maass Selberg}),
that
$$
\dE(P\col s\nu\col x) = \nC_{P|Q}(s \col \nu) \dE(Q\col \nu \col x),
$$
for all $x \in \spX$ and all $\nu \in i\faQqd.$ Now (a) follows by Definition
\refer{d: Fou P}.
The linear bijection $s$ from $i\faQqd$ onto $i\faPqd$ maps the chosen Lebesgue measures
onto each other, see Section \refer{s: normalization}.
Hence, for $f \in \cC(\spX\col\tau),$
\begin{equation}
\naam{e: formula for Wave P Fou P}
\Wave_P \Fou_P f(x) = \int_{i \faQqd} \nE(P\col s\nu\col x) \Fou_P f(s\nu) \; d\nu.
\end{equation}
Applying (a) and Proposition \refer{p: funct eqn Eis} we obtain that
\begin{eqnarray*}
\nE(P\col s\nu\col x) \Fou_P f(s\nu) &=& \nE(P\col s\nu\col x)
\nC_{P|Q}(s\col \nu) \Fou_Q f(\nu)\\
&=& \nE(Q\col \nu\col x) \Fou_Q f(\nu).
\end{eqnarray*}
Substituting the obtained identity in the right-hand side of
(\refer{e: formula for Wave P Fou P})
we obtain (b).
\qed

\begin{rem}
\naam{r: Wave P Fou P through class}
Let $P\in \allparabs.$ Then it follows from part (b) of the above lemma
that the continuous linear endomorphism
$\Wave_P \after \Fou_P$ of $\cC(\spX \col \tau)$  depends on $P$
through its equivalence class in $\allparabs/\!\!\sim.$
\end{rem}

If $P \in \allparabs,$ we agree to write $W^*_P$ for the normalizer of $\faPq$ in $W.$
Then  $W(\faPq) \simeq W^*_P / W_P.$

\begin{thm}
\naam{t: final Fourier inversion on Schwartz space}
Let $f \in \cC(\spX\col \tau).$ Then
$$
f =  \sum_{[P]\in \allparabs/\sim}  [W:W^*_P] \;\Wave_P \Fou_P f.
$$
\end{thm}

\begin{rem}
\naam{r: inversion and constants in Delorme}
In view of Remark \refer{r: Eis and Delorme}, this theorem corresponds
to part (iii) of Thm.~2 in Delorme's paper \bib{Dpl}. Note that in the latter theorem constants
$|W(\faPq)|^{-1},$ for $P \in \allparabs,$ appear in place of the constants $[W:W^*_P].$
This is due to a different normalization of measures, as will be explained in the sequel
\bib{BSpl2} to this paper.
 \end{rem}

\proof
We observe that $[W:W_P]|W(\faPq)|^{-1} = [W: W^*_P],$ for $P \in \allparabs.$
Since every $P \in \allparabs$ is associated with a standard parabolic
subgroup, see Lemma \refer{l: on parabolics},
the result now follows from Theorem
\refer{t: inversion for Schwartz},
Lemma \refer{t: sum of residue weights} and
Remark \refer{r: Wave P Fou P through class}.
\qed

\section{Properties of the Fourier transforms}\eqnreset
\naam{s: properties of the Fourier transforms}
The purpose of this section is to establish relations between the different
Fourier and wave packet transforms $\Fou_P$ and $\Wave_Q,$ as $P,Q \in \allparabs.$
We shall also determine the image of $\Fou_P$ and the kernel of $\Wave_Q.$
The relation between the several Fourier transforms is given
by Lemma \refer{l: Wave Fou through class} (a).

\begin{lemma}
\naam{l: transformation props C}
Let $P,Q,R \in \allparabs$ be associated. Then, for all
$s \in W(\faQq|\faPq)$ and $t \in W(\fa_{R\iq}|\faQq),$
$$
\nC_{R|P}(ts \col \nu) = \nC_{R|Q}(t\col s\nu) \after \nC_{Q|P}(s\col \nu),
\qquad (\nu \in i\faPqd).
$$
\end{lemma}

\proof
The above identity is an immediate consequence of the functional equation
for the Eisenstein integral, see Proposition \refer{p: funct eqn Eis}, and the definition
of the $C$-function, see Definition \refer{d: c functions as constant term}.
\qed

In particular, from the above lemma with $P=Q=R$ combined
with the Maass--Selberg relations, see Theorem \refer{t: Maass Selberg},
 we see that we may define
a unitary representation $\gg_P$ of $W(\faPq)$ in $L^2(i\faPqd)\otimes  \cAtwoP$
by
$$
[\gg_P(s) \gf](\nu) = \nC_{P|P}(s\col s^{-1}\nu) \gf(s^{-1} \nu),\qquad (\nu \in i\faPqd),
$$
for $\gf \in L^2(i\faPqd)\otimes  \cAtwoP.$
The associated collection of $W(\faPq)$-invariants in $L^2(i\faPqd)\otimes  \cAtwoP$
is denoted by $(L^2(i\faPqd)\otimes \cAtwoP)^{W(\faPq)}.$ The orthogonal projection
from the first onto the latter space is denoted by
$$
\rmP_{W(\faPq)}: \;\;
L^2(i\faPqd)\otimes \cAtwoP \to ( L^2(i\faPqd)\otimes  \cAtwoP)^{W(\faPq)} .
$$
The intersection
of the latter space with
$\cS(i\faPqd)\otimes \cAtwoP$ consists of the functions
$\gf\in \cS(i\faPqd)\otimes \cAtwoP$ satisfying
\begin{equation}
\naam{e: defi W faPq invariants}
\gf(s\nu) = \nC_{P|P}(s\col \nu)\gf(\nu), \qquad
(s \in W(\faPq),\; \nu \in  i\faPqd),
\end{equation}
and is denoted by
$(\cS(i\faPqd) \otimes \cAtwoP)^{W(\faPq)}.$

\begin{cor}
\naam{c: image Fou W invariant}
Let $P \in \allparabs.$
The image  of $\cC(\spX\col \tau)$ under the Fourier transform $\Fou_P$
is contained in the space $(\cS(i\faPqd)\otimes  \cAtwoP)^{W(\faPq)}.$
\end{cor}

\proof
Let $f \in \cC(\spX\col \tau).$ Then it follows from
Lemma \refer{l: Wave Fou through class}(a) with $P = Q$ that $\gf := \Fou_P f$
satisfies (\refer{e: defi W faPq invariants}).
\qed

We can now state the first main result of this section.

\begin{thm}
\naam{t: Fou after Wave}
Let $P,Q \in \allparabs.$
\begin{enumerate}
\itema
If $P$ and $Q$ are not associated, then $\Fou_Q\after \Wave_P = 0.$
\itemb
If $P$ and $Q$ are associated, then $[W: W^*_P] \,\Fou_Q\after \Wave_P \after \Fou_P = \Fou_Q$
on $\cC(\spX\col \tau).$
\itemc
If $P$ and $Q$ are associated, then, for each $s \in W(\faQq\mid \faPq),$
every $\gf \in \cS(i\faPqd)\otimes \cAtwoP$
and all $\nu \in i\faPqd,$
$$
\Fou_Q \after \Wave_P\, \gf \,(s\nu)  =
[W: W^*_P]^{-1}\, \nCQP(s \col \nu) \rmP_{W(\faPq)}\gf(\nu).
$$
In particular, $\Fou_P\after \Wave_P = [W: W^*_P]^{-1}\, \rmP_{W(\faPq)}.$
\end{enumerate}
\end{thm}

The proof is analogous to the proof of Theorem 16.6 in \bib{BSmc}, with
adaptations to deal with the present more general situation.
In the course of the proof we  need two lemmas.
The first of these is a straightforward generalization of
Lemma 16.11 in \bib{BSmc}.

\begin{lemma}
\naam{l: surjective differential}
Let $P \in \allparabs,$ let $\fb\subset \fq$ be a $\Cartan$-stable
Cartan subspace containing $\faPq$
and let $\gL \in \stfbPdc.$ Then for $\gl$ in the complement
of a finite union of affine hyperplanes in $\faPqdc,$
the map $D \mapsto d[\gg(D\col \gL + \dotvar)](\gl)$ is surjective
from $\DX$ onto $ \fa_{P\iq\iC}.$
\end{lemma}

\proof
The proof is a straightforward modification of the proof
of Lemma 16.12 in \bib{BSmc}. In that proof one should everywhere put
$\faPq$ in place of $\faq$ and $\stfbP$ in place of $\fb_{\rm k}.$ In particular,
$\pi_{\faq}$ should be replaced by the projection $\pi_{\faPq}: \fbdc \to \faPqdc$
along the subspace $\stfbPdc.$
\qed

The next lemma is a consequence of Lemma \refer{l: data associated to gL},
which in turn heavily relies
on the information about the infinitesimal
characters of discrete series stated in Theorem
\refer{t: infinitesimal characters L two d real and regular new}.

\begin{lemma}
\naam{l: conjugacy of gL nu}
Let $P,Q \in \allparabs,$ let $\fb \subset \fq$ be a $\Cartan$-stable Cartan
subspace containing $\faq$ and let $\gL_1 \in \rmL_P(\fb,\tau)$ and
$\gL_2 \in \rmL_Q(\fb, \tau).$
Let $\nu_1 \in i\faPqd$ be such that $\gL_1 + \nu_1$ is regular with
respect to $\gS(\fb),$  let $\nu_2 \in i\faQqd$ and assume that
$\gL_1 + \nu_1$ and $ \gL_2 + \nu_2$ are conjugate under $W(\fb).$
Then $P$ and $Q$ are associated, and $\nu_1$
and $\nu_2$ are conjugate under $W(\faQq\mid \faPq).$
\end{lemma}

\proof
Let $s \in W(\fb)$ be such that $s(\gL_1 + \nu_1) = \gL_2 + \nu_2.$
Select $v_1 \in \PcW$ and $v_2 \in \QcW$ such that
$\gL_1 \in L_{P_,v_1}(\fb, \tau)$
and $\gL_2 \in L_{Q_,v_2}(\fb, \tau).$
Let $(\hat \fb_1, t_1)$ and $(\hat \fb_2, t_2)$ be
associated with $P, v_1$ and $Q, v_2,$ respectively, as in Lemma \refer{l: data associated to gL}.
Then $t_1 \gL_1 + t_1 \nu_1$ and $t_2 \gL_2 + t_2 \nu_2$ are conjugate
under $t_2st_1^{-1} \in W(\hat \fb_2\mid \hat \fb_1).$
It follows by application of
\bib{Cn}, Lemma 2, that $t_1 \faPq = \Ad(v_1)^{-1} \faPq$ and
$t_2 \faQq = \Ad(v_2)^{-1} \faQq$
are conjugate under $t_2 s t_1^{-1}.$ This implies that $s(\faPq) = \faQq.$
It follows that $s|_{\faPq} \in W(\faQq\mid \faPq),$ see
Lemma \refer{l: conjugacy in Cartan};
hence, $P \sim Q.$
  It also follows by \bib{Cn}, Lemma 2, that $t_2st_1^{-1}$ maps $t_1 \nu_1$
onto $t_2 \nu_2;$ hence, $s \nu_1 = \nu_2.$
\qed

The following lemma collects some properties of the composition $\Fou_Q\after \Wave_P$
needed in the proof of Theorem \refer{t: Fou after Wave}.

\begin{lemma}
\naam{l: properties Fou Wave}
Let $P, Q \in \allparabs.$ Then the composition $\cT: = \Fou_Q\after \Wave_P$
is a continuous linear map from $\cS(i\faPqd) \otimes \cAtwoP$ to $\cS(i\faQqd) \otimes \cAtwoQ.$
Moreover, it satisfies the following properties.
\begin{enumerate}
\itema
$\umu_Q(D)\after \cT = \cT\after \umu_P(D)$ for all $D \in \DX.$
\itemb
$\cT$ maps into $(\cS(i \faQqd) \otimes \cAtwoQ)^{W(\faQq)}.$
\end{enumerate}
\end{lemma}

\proof
The continuity of $\cT$ follows from
Theorem \refer{t: Wave continuous on Schwartz space}
combined with Proposition \refer{p: Fou continuous on Schwartz space}.
Property (a) follows from Lemma \refer{l: D after Wave}
combined with Lemma \refer{l: Fou after D}.
Finally, (b) follows from Corollary \refer{c: image Fou W invariant}.
\qed

\begin{prop}
\naam{p: about general cT}
Let $P,Q \in \allparabs.$ There exists an open dense $W(\faQq)$-invariant subset
$\Omega  \subset i\faQqd$
with the following property.
Let $\cT$ be any continuous linear map from  $\cS(i\faPqd) \otimes \cAtwoP$ to
$\cS(i\faQqd) \otimes \cAtwoQ$ satisfying
the properties of Lemma \refer{l: properties Fou Wave}.
\begin{enumerate}
\itema
If $P$ and $Q$ are not associated, then $\cT = 0.$
\itemb
If $P$ and $Q$ are associated and
$s_0 \in W(\faPq\mid \faQq),$ then there exists a unique smooth function
$\gb: \Omega \to \Hom(\cAtwoP,\cAtwoQ)$ such that
\begin{equation}
\cT\gf(\nu) = \rmP_{W(\faQq)}(\gb\, s_0^*\gf)(\nu),
\end{equation}
for all $\gf \in \Cci(s_0\Omega)\otimes \cAtwoP$ and $\nu \in \Omega.$
\end{enumerate}
\end{prop}

\proof
For every $\nu \in i\faQqd,$ we define the distribution
$u_\nu \in \cD'(i\faPqd) \otimes \Hom(\cAtwoP, \cAtwoQ)$ by
$$
u_\nu(\gf) = \cT(\gf)(\nu), \qquad (\gf \in \Cci(i\faPqd) \otimes \cAtwoP).
$$
Then it follows from condition (a) that
$$
u_\nu \after \umu_P(D) = \umu_Q(D\col \nu) u_\nu, \qquad (D \in \DX).
$$
Let now $\gL_1 \in \rmL_P(\fb,\tau)$ and  $\psi_1 \in \cAtwoP(\gL_1).$
Let $\gL_2 \in \rmL_Q(\fb, \tau)$
and $\psi_2 \in \cAtwoQ(\gL_2).$
We define the distribution $v_\nu \in \cD'(i\faPqd)$ by
$v_\nu(f) = \hinp{ u_\nu (f \otimes \psi_1)}{\psi_2},$
for $f \in \Cci(i\faPqd).$
It follows
that
\begin{equation}
\naam{e: equation for v nu}
[\gg_{\fb}(D \col \gL_1 + \dotvar) - \gg_\fb(D \col \gL_2 + \nu)] v_\nu = 0.
\end{equation}
Each element $\gL$ from the finite set $ \rmL_Q(\fb, \tau)$ is regular
with respect to $\gS(\fm_{Q\iC}, \stfb_Q),$ see Corollary \refer{c: gL is real regular}.
Let $\Omega_0$ be the set of $\nu \in i\faQqd$ such that $\gL +\nu$ is regular with respect
to $\gS(\fb),$
for every $\gL \in \rmL_Q(\fb,\tau).$ Then $\Omega_0$ is the complement
of a finite union of hyperplanes in $i\faQqd,$ hence open dense.

Let $\nu \in \Omega_0$ and let $\gL_1, \gL_2$ be as above.
Moreover, let $\nu_1 \in i\faPqd$
and assume that  $\nu \notin W(\faQq\mid \faPq)\nu_1.$
Then by Lemma \refer{l: conjugacy of gL nu}, the elements
$\gL_1 + \nu_1 $ and $\gL_2 + \nu$ are not conjugate under $W(\fb).$
Hence, there exists a $D \in \DX$ such that the polynomial function
in front of $v_\nu$ in (\refer{e: equation for v nu}) does not vanish at $\nu_1.$ This implies
that $v_\nu$ vanishes in a neighborhood of $\nu_1.$ Let $\gf \in \Cci(i\faPqd) \otimes \cAtwoP.$
Then it follows from the above by linearity that $\cT(\gf)(\nu) = u_\nu(\gf) = 0$ for all
$\nu \in \Omega_0$ with $\nu \notin W(\faQq\mid \faPq)\supp \gf.$ By density and
continuity, this implies that the function $\cT(\gf)$ vanishes on
$i\faQqd \setminus W(\faQq\mid \faPq)\supp \gf.$ Hence,
$$
\supp \cT \gf \subset W(\faQq \mid \faPq) \supp \gf.
$$
If $P$ and $Q$ are not associated, then it follows that the latter set
has empty interior in $i\faQqd,$ hence $\cT\gf = 0$ by continuity.
This establishes (a).

{}From now on, we assume that $P\sim Q.$ Then it follows from the above that
\begin{equation}
\naam{e: support u nu}
\supp u_\nu \subset W(\faPq\mid \faQq) \nu,
\end{equation}
for every $\nu \in i\faQqd.$

Let $\Omega_1$ be the set of $\nu \in i \faQqd$ whose stabilizer in $W(\faQq)$
is trivial.
Then $\Omega_1$ contains
the complement
of a finite union of hyperplanes in $i \faQqd$ hence is open dense in $i\faQqd.$
Since $W(\faQq)$ acts simply transitively on $W(\faPq\mid \faQq)$ from the
right, we see that if $\nu \in \Omega_1,$ the points $s\nu,$ for $s \in W(\faPq\mid \faQq),$
are mutually different.

Let $\Omega_2'$ be the set of $\mu \in i\faPqd$ such that
for every $\gL_1 \in \rmL_P(\fb,\tau)$ the map $D \mapsto
d[\gg_\fb(D \col \gL_1 + \dotvar)](\mu)$ is surjective from $\DX$ onto
$\faPqdc.$ Then $\Omega_2'$ is an open subset of $i\faPqd$
containing the complement of a finite union of
hyperplanes, see Lemma \refer{l: surjective differential}.
It follows that $\Omega_2 = \cap_{s \in W(\faPq\mid \faQq)} s^{-1} \Omega_2'$
is a similar subset of $i\faQqd.$ We define $\Omega: =  \Omega_1 \cap \Omega_2.$

Let now $\nu\in \Omega.$
We claim that the distribution $u_\nu$ has order zero.
To prove the claim, fix $\gL_1, \gL_2, \psi_1, \psi_2$ as before,
and define $v = v_\nu$ as above. Then by linearity, it suffices to
show that $v$ has order zero.
Since $\supp (v) \subset W(\faPq\mid \faQq)\nu,$ by (\refer{e: support u nu}),
it follows by our
assumption on $\Omega_1$ that we may express $v$ uniquely as a sum of
distributions $v_s,$ for $s \in W(\faPq\mid \faQq),$ with $\supp v_s \subset \{s\nu\}.$
{}From (\refer{e: equation for v nu}) it follows that each $v_s$ satisfies the equations
$$
\gf_D v_s = 0,\qquad (D \in \DX),
$$
where $\gf_D : i \faPqd \to \C$ is given by
$\gf_D = \gg_{\fb}(D\col \gL_1 + \dotvar ) - \gg_{\fb}(D\col \gL + \nu).$

It follows from our assumption on $\Omega_2$ that the collection of differentials
$d\gf_D(s\nu),$ for $D \in \DX,$ spans $\faPqdc.$ Now apply \bib{BSmc}, Lemma 16.10,
to conclude that $v_s$ has order zero, for each $s \in W(\faPq\mid \faQq).$
This establishes the claim that $v_\nu$ has order zero.

It also follows from the above that
$$
u_\nu = \sum_{s \in W(\faPq\mid\faQq)} \gd_{s\nu} \otimes E_s(\nu),
$$
with $E_s(\nu)$ a unique element of $\Hom(\cAtwoP, \cAtwoQ),$
for $s \in W(\faPq\mid \faQq).$ We conclude that, for every
$\gf \in \Cci(i\faPqd)\otimes \cA_{2,P}$ and all $\nu \in \Omega,$
\begin{equation}
\naam{e: first formula cT gf}
\cT\gf(\nu) = \sum_{s \in W(\faPq\mid \faQq)} E_s(\nu) \gf(s\nu).
\end{equation}
Fix $s_0 \in W(\faPq\mid \faQq).$
Let $\nu_0 \in \Omega.$
By the assumption on $\Omega_1$
there exists an open neighborhood $U$ of $\nu_0$ in $\Omega$ such that
the sets $sU$ are mutually disjoint, for $s \in W(\faPq\mid \faQq).$ For
$\gf \in \Cci(s_0 U) \otimes \cAtwoP$
we have
$$
\cT(\gf)(\nu) = E_{s_0}(\nu) \gf(s_0 \nu).
$$
We conclude that $E_{s_0}$ is smooth on $U.$ It follows that
$E_{s_0} \in \Ci(\Omega) \otimes \Hom(\cAtwoP, \cAtwoQ).$
{}From the above asserted uniqueness of the $E_s$ and the transformation property
of $\cT\gf$ stated in Lemma \refer{l: properties Fou Wave} (b), it follows that
$$
E_{s_0}(t\nu) = \nC_{Q|Q} (t\col \nu) E_{s_0 t}(\nu),\qquad (t \in W(\faQq)).
$$
If we combine this with (\refer{e: first formula cT gf}) we obtain, for all
$\gf \in \Cci(s_0\Omega) \otimes \cAtwoP,$
and all $\nu \in \Omega,$ that
$$
\cT \gf(\nu)  = \sum_{t \in W(\faQq)} \nC_{Q|Q}(t\col \nu)^{-1} E_{s_0}(t\nu) \gf(s_0 t \nu)
= |W(\faQq)| \, \rmP_{W(\faQq)} (E_{s_0} s_0^*\gf)(\nu).
$$
This establishes the result with $\gb = |W(\faQq)| E_{s_0}.$
\qed
\medno
{\bf Proof of Theorem \refer{t: Fou after Wave}:\ }
If $P\not\sim Q,$ then it follows from Lemma \refer{l: properties Fou Wave}
combined with Proposition \refer{p: about general cT}
that $\Fou_Q \after \Wave_P = 0.$

To prove (b) we note that
if $f \in \Cci(\spX\col \tau),$ then
it follows from Theorem \refer{t: final Fourier inversion on Schwartz space} that
$f = \sum_{[R] \in \allparabs/\sim} [W:W^*_R] \,\Wave_R\Fou_R f.$
Applying $\Fou_Q$ to both sides of this equation,
it follows by (a) that
$$
\Fou_Q f = [W:W^*_R]\,\Fou_Q \Wave_R \Fou_R f
$$
for each $R \in \allparabs$ with $R \sim Q;$ in particular,
we may take $R = P.$  Thus, (b) follows by continuity and density.

We shall first prove (c) under the assumption  that $P = Q$
and $s =1.$
Put $\cT = \Fou_P \after \Wave_P.$
{}From Lemma \refer{l: properties Fou Wave}
and Proposition \refer{p: about general cT} we deduce that
$$
\cT\psi(\nu) = \rmP_{W(\faPq)} (\gb \psi) (\nu),
$$
for all $\psi \in \Cci(\Omega)\otimes \cAtwoP$ and
$\nu \in \Omega;$ here $\Omega \subset i\faPqd$
is an open dense subset, and $\gb \in \Ci(\Omega) \otimes \End(\cAtwoP).$
It follows from Lemma \refer{l: Fou and Wave transpose}
that the operator $\cT$ is symmetric with respect to
the $L^2$-inner product on $\cS(i\faPq) \otimes \cAtwoP.$
Let $\gf \in \cS(i\faPq) \otimes \cAtwoP.$ Then for all
$ \psi \in \Cci(\Omega)\otimes \cAtwoP,$
$$
\hinp{\cT\gf}{\psi} = \hinp{\gf}{\cT\psi} = \hinp{\gf}{\rmP_{W(\faPq)}(\gb \psi)}
= \hinp{\gb^* \rmP_{W(\faPq)} \gf} {\psi}.
$$
This implies that, for all $\gf \in \cS(i\faPqd)\otimes \cAtwoP,$
$\cT\gf = \gb^* \rmP_{W(\faPq)} \gf$ on $\Omega.$ We claim that in fact $\gb^* =
[W:W^*_P]^{-1} \, I$ on
$\Omega.$

To prove the claim we note that it follows from
the established part (b) with $P =Q$ that
$\Fou_P = [W:W^*_P]\,\cT \Fou_P.$ We infer that,
for $f \in \Cci(\spX\col \tau)$ and $\nu \in \Omega,$
\begin{eqnarray*}
\Fou_P f (\nu) &= & [W:W^*_P]\,\cT \Fou_P f (\nu)\\
&=&
 [W:W^*_P]\,\gb(\nu)^* \rmP_{W(\faPq)} \Fou_P f(\nu)\\
&=&  [W:W^*_P]\, \gb(\nu)^* \Fou_P f(\nu).
\end{eqnarray*}
In the last equality we have used Corollary \refer{c: image Fou W invariant}.
The claim now follows by application of Lemma \refer{l: surjectivity ev Fou} below.
We infer that for all $\gf \in \cS(i\faPq) \otimes \cAtwoP,$ we have
$\cT \gf = [W\col W_P^*]^{-1}\, \rmP_{W(\faPq)} \gf$
on $\Omega,$ hence on $i\faPqd,$ by continuity and density.
This establishes part (c) of the theorem for $P = Q$ and $s =1.$
For a general pair of associated
parabolic subgroups $P,Q$ and a general $s \in W(\faPq\mid \faQq)$ the assertion follows
by application of Lemma \refer{l: Wave Fou through class} (a).
\qed

\begin{lemma}
\naam{l: surjectivity ev Fou}
Let $P \in \allparabs$ and let $\nu \in i\faPqd$ have
trivial stabilizer in $W(\faPq).$
Then $f \mapsto \Fou_Pf(\nu)$ maps $\Cci(\spX\col \tau)$ onto $\cAtwoP.$
\end{lemma}

\proof
The proof is a reduction to the lemma below, in a way completely analogous to
the proof of \bib{BSmc}, Lemma 16.13.
\qed

\begin{lemma}
Let $P \in \allparabs$ and let $\nu\in i\faPqd$ have trivial stabilizer in $W(\faPq).$
Then for every $\psi \in \cAtwoP\setminus\{0\},$
the Eisenstein integral $\nE(P\col \nu)\psi$ does not
vanish identically on $\spX.$
\end{lemma}

\proof
The proof is completely analogous to the proof of
\bib{BSmc}, Lemma 16.14, involving the constant term of the Eisenstein
integral along $P.$
\qed

\begin{prop}
\naam{p: kernel Wave}
Let $P \in \allparabs.$
\begin{enumerate}
\itema
The map $\rmP_{W(\faPq)}$ restricts to a continuous linear projection from
$\cS(i\faPqd)\otimes  \cAtwoP$ onto $(\cS(i\faPqd) \otimes \cAtwoP)^{W(\faPq)}.$ This projection is orthogonal
with respect to the given $L^2$-inner product.
\itemb
$\Wave_P \after \rmP_{W(\faPq)} = \Wave_P.$
\itemc
The kernel of $\Wave_P$ equals the kernel of $\rmP_{W(\faPq)}.$
\end{enumerate}
\end{prop}

\proof
It follows from Theorem \refer{t: Fou after Wave} (c)  that
$\rmP_{W(\faPq)} = [W\col W_P^*]\,\Fou_P \after \Wave_P.$
It now follows from application of
Lemma \refer{l: properties Fou Wave},
that $\rmP_{W(\faPq)}$ is a continuous linear endomorphism of
$\cS(i \faPqd) \otimes \cAtwoP,$ with
image contained in $(\cS(i \faPqd) \otimes \cAtwoP)^{W(\faPq)}.$
The latter space is contained in $(L^2(i \faPqd) \otimes \cAtwoP)^{W(\faPq)},$
on which space $\rmP_{W(\faPq)}$ equals the identity. Hence,
$\rmP_{W(\faPq)}$ is a projection and (a) readily follows.

Starting from Theorem \refer{t: Fou after Wave} (b) with
$P = Q,$ we obtain by taking adjoints and applying
Lemma \refer{l: Fou and Wave transpose}, that
$$
\Wave_P \after [W\col W_P^*]\Fou_P \after \Wave_P = \Wave_P.
$$
Assertion (b) now follows by application of Theorem \refer{t: Fou after Wave} (c).

{}From (b) it follows that $\ker\,\rmP_{W(\faPq)} \subset \ker \Wave_P.$
The converse inclusion follows by Theorem \refer{t: Fou after Wave} (c).
\qed

\begin{prop}
\naam{p: image Fou}
Let $P \in \allparabs.$ Then the image of $\Fou_P$ equals
$(\cS(i\faPqd) \otimes  \cAtwoP)^{W(\faPq)}.$
\end{prop}

\proof
That $\image(\Fou_P) \subset (\cS(i\faPqd)\otimes  \cAtwoP)^{W(\faPq)}$ was asserted in
Corollary \refer{c: image Fou W invariant}. For the converse inclusion,
let $\gf \in (\cS(i\faPqd)\otimes  \cAtwoP)^{W(\faPq)}.$ Then
$$
\gf = \rmP_{W(\faPq)} \gf
= \Fou_P( [W:W_P^*]\, \Wave_P \gf) \in \image(\Fou_P),
$$
in view of Theorem \refer{t: Fou after Wave} (c).
\qed

\section{The spherical Plancherel theorem}\eqnreset
\naam{s: the spherical Plancherel theorem}
In this final section we establish the Plancherel theorem
for functions from  $\cC(\spX\col \tau)$ and  $L^2(\spX\col \tau).$
If $P \in \allparabs,$ then by $\cC_P(\spX\col \tau)$
we denote the image of the operator $\Wave_P$ in $\cC(\spX\col \tau).$

\begin{thm}{\ }
\naam{t: spherical Planch on Schwartz}
Let $P \in \allparabs.$
\begin{enumerate}
\itema
The space $\cC_P(\spX\col \tau)$ is closed in $\cC(\spX\col \tau)$
and depends on $P \in \allparabs$ through
its class in $\allparabs/\!\!\sim.$
\itemb
The spaces $\ker \Fou_P$ and $\cC_P(\spX\col \tau)$ are each other's orthocomplement
in $\cC(\spX\col \tau).$
\itemc
The space $\cC(\spX \col \tau)$ admits the following finite direct sum
decomposition
\begin{equation}
\naam{e: cC as dir sum cCR}
\cC(\spX\col \tau)= \bigoplus_{[R] \in \allparabs/\sim} \; \cC_R(\spX\col \tau).
\end{equation}
The summands are mutually orthogonal with respect to the $L^2$-inner product
on $\cC(\spX\col \tau).$
\itemd
For each $P \in \allparabs,$
the operator $[W:W^*_P]\, \Wave_P\after \Fou_P$ is the
projection onto $\cC_P(\spX\col \tau)$ along the remaining summands in the above
direct sum.
\end{enumerate}
\end{thm}

\begin{rem}
\naam{r: comparison with HC and Del}
For the case of the group this result, together
with Propositions \refer{p: kernel Wave} and
\refer{p: image Fou}, is a consequence of Harish-Chandra's
Plancherel theorem for $K$-finite Schwartz functions,
see Remark \refer{r: comparison with nE of HC} and \bib{HC3}, \S 27.
See also \bib{Arthur}, Ch.~III, \S 1,
for a review of the Plancherel theorem for spherical Schwartz
functions.

For general symmetric spaces, the result,
together with Propositions \refer{p: kernel Wave} and \refer{p: image Fou},
is equivalent
to Thm.~2 in Delorme's paper \bib{Dpl}, in view of Remark \refer{r: Eis and Delorme}.
See also Remark \refer{r: inversion and constants in Delorme}.
\end{rem}

\proof
If $Q \in \allparabs,$ $Q\not\sim P,$ then, by Theorem \refer{t: Fou after Wave} (a),
\begin{equation}
\naam{e: Fou Q zero on cC P}
\Fou_Q = 0\quad \text{on}\;\;\cC_P(\spX\col \tau).
\end{equation}
In view of Theorem \refer{t: final Fourier inversion on Schwartz space} this implies that
\begin{equation}
\naam{e: Wave after Fou on cC P}
[W:W^*_P]\, \Wave_P \Fou_P = I \quad \text{on} \quad \cC_P(\spX\col \tau).
\end{equation}
We infer that
$
\cC_P(\spX\col \tau) = \image (\Wave_P\Fou_P).
$
By Remark \refer{r: Wave P Fou P through class} it follows that $\cC_P(\spX\col \tau)$
depends on $P$ through
its class in $\allparabs/\!\!\sim.$ This establishes the second assertion of (a).
{}From Theorem \refer{t: final Fourier inversion on Schwartz space}
we see that $\cC(\spX\col \tau)$ is the vector
sum of the spaces $\cC_{R}(\spX\col \tau),$ for $[R] \in \allparabs/\!\!\sim.$
To establish the orthogonality of the summands, let $P,Q \in \allparabs,$
$P \not\sim Q,$ let $f \in \cC_P(\spX\col \tau)$ and
$\psi \in \cS(i\faQqd) \otimes \cA_{2,Q}.$ Then $\hinp{f}{\Wave_Q\psi} = \hinp{\Fou_Q f}{\psi} = 0,$
by Lemma \refer{l: Fou and Wave transpose} and (\refer{e: Fou Q zero on cC P}). This establishes (c).
Combining (c) with (\refer{e: Fou Q zero on cC P}) and (\refer{e: Wave after Fou on cC P}),
we obtain (d). {}From (c) it follows that $\cC_P(\spX\col \tau)$ is closed,
whence (a).

It remains to establish (b).  {}From (\refer{e: Fou Q zero on cC P})
it follows that $\ker \Fou_P$ contains the part $\cC_0$
of the sum (\refer{e: cC as dir sum cCR}) consisting of the summands with $R \not\sim P.$
On the other hand, $\ker \Fou_P \cap \cC_P(\spX\col \tau) = 0$
by (\refer{e: Wave after Fou on cC P})
and we conclude that $\ker \Fou_P = \cC_0.$ The latter space
equals the orthocomplement of $\cC_P(\spX\col \tau),$ by the orthogonality
of the sum (\refer{e: cC as dir sum cCR}).
\qed

We denote by $L^2_P(\spX\col \tau)$ the closure in $L^2(\spX\col \tau)$
of $\cC_P(\spX\col \tau).$ This space depends on $P$ through its class
in $\allparabs/\!\!\sim.$

\begin{cor}
\naam{c: ortho deco of Ltwotau}
The space $L^2(\spX\col \tau)$ admits the following orthogonal direct sum decomposition
into closed subspaces
$$
L^2(\spX\col \tau)= \bigoplus_{[P] \in \allparabs/\sim} \; L^2_P(\spX\col \tau).
$$
\end{cor}

\proof
Since $\cC(\spX\col \tau)$ is dense in $L^2(\spX\col \tau),$
the result is an immediate consequence of Theorem \refer{t: spherical Planch on Schwartz} (c).
\qed

We recall that a continuous linear map between Hilbert spaces
$T: \cH_1 \to \cH_2$ is called a partial isometry if $T$ is an isometry
from $(\ker T)^\perp$ into $\cH_2.$ In particular, this means
that $\image T$ is a closed subspace of $\cH_2.$

\begin{thm}{\ }
\naam{t: spherical Planch on Ltwo}
Let $P \in \allparabs.$
\begin{enumerate}
\itema
The operator $\Fou_P$ extends uniquely to a continuous linear map from
$L^2(\spX\col \tau)$ to $L^2(i\faPqd)\otimes  \cAtwoP,$ denoted $\Fou_P.$
\itemb
The operator $\Wave_P$ extends uniquely to a continuous linear
map from $L^2(i\faPqd)\otimes  \cAtwoP$ to $L^2(\spX\col \tau),$ denoted $\Wave_P.$
This extension
is the adjoint of the extended operator $\Fou_P.$%
\itemc
The extended operator $[W:W^*_P]^{1/2}\Fou_P$
is a partial isometry from $L^2(\spX\col \tau)$
onto  the space $(L^2(i\faPqd)\otimes  \cAtwoP)^{W(\faPq)},$ with kernel equal
to the orthocomplement of $L^2_P(\spX\col \tau).$
\itemd
The extended operator $[W:W^*_P]^{1/2}\Wave_P$ is a partial isometry
from $L^2(i\faPqd)\otimes  \cAtwoP$  onto $L^2_P(\spX\col \tau)$ with kernel equal
to the orthocomplement of $(L^2(i\faPqd) \otimes \cAtwoP)^{W(\faPq)}.$
\iteme
Let $\repparabs \subset \allparabs$ be a choice of representatives
for the classes in $\allparabs/\!\!\sim.$ Then
$$
I = \sum_{R \in \repparabs} [W:W^*_R]\, \Wave_R \Fou_R \text{on}
L^2(\spX\col \tau).
$$
\end{enumerate}
\end{thm}

\begin{rem}
In view of Remark \refer{r: Eis and Delorme},
this result corresponds to part (iv) of Theorem 2
in \bib{Dpl}. See also Remark \refer{r: comparison with HC and Del}.
\end{rem}
\proof
Fix $\bfP_\gs$ as in (e).
Let $f \in \cC(\spX\col \tau).$ Then it follows from
Theorem \refer{t: spherical Planch on Schwartz}
(c) and (d), combined with Lemma \refer{l: Fou and Wave transpose}
that
$$
\|f\|^2 = \hinp{f}{f} = \sum_{R \in \bfP_\gs} [W:W^*_R]\, \hinp{f}{\Wave_R \Fou_R f}
= \sum_{R \in \bfP_\gs} [W:W^*_R]\, \|\Fou_R f\|^2.
$$
In particular, this equation
 holds for a choice of $\bfP_\gs$ with $P \in \bfP_\gs.$
It follows  that $\Fou_P$ is continuous with respect
to the $L^2$-topologies. By density of $\cC(\spX\col \tau)$ in $L^2(\spX\col \tau),$
it follows that $\Fou_P$ has a unique continuous linear extension
$L^2(\spX\col \tau) \to L^2(i\faPqd) \otimes \cAtwoP.$ Hence (a).

Since $\Wave_P$ is the adjoint of $\Fou_P$ with respect to the $L^2$-inner products on the
Schwartz spaces, it follows that the adjoint of the extension of $\Fou_P$ is
a continuous linear extension of $\Wave_P$ to a continuous linear operator from
$L^2(i\faPqd) \otimes \cAtwoP$ to $L^2(\spX\col \tau).$ This extension is unique
by density of $\cS(i\faPqd)$ in $L^2(i\faPqd).$ This proves (b).

By continuity and density, the formula in (e) follows from  Theorem
\refer{t: final Fourier inversion on Schwartz space}.
{}From Theorem \refer{t: spherical Planch on Schwartz}  (b) and (c) it follows that, for
$R \in \bfP_\gs$ with $R \not\sim P,$ $\Fou_P = 0$ on $\cC_R(\spX\col \tau),$
hence on $L^2_R(\spX\col ,\tau),$ by continuity and density.
Put $\tFou_P:= [W: W_P^*]^{1/2} \Fou_P$ and $\tWave_P := [W: W_P^*]^{1/2} \Wave_P.$
Using (e) we infer that $\ker \tFou_P = L^2_P(\spX\col \tau)^\perp$
and that $\tWave_P\after \tFou_P$ is the orthogonal projection from
$L^2(\spX\col \tau)$ onto $L^2_P(\spX\col \tau).$ Since $\tWave_P = \tFou_P^*,$
it follows that $\tFou_P$ is isometric from $L^2_P(\spX\col \tau)$ onto $\image \tFou_P$
and that $\tWave_P$ is isometric from $\image \tFou_P$ onto $L^2_P(\spX\col \tau).$
It follows from Theorem \refer{t: Fou after Wave}
and continuity and density that $\tFou_P\after \tWave_P = \rmP_{W(\faPq)}$
on $L^2(i\faPqd) \otimes \cAtwoP.$ Hence,
$$
\image(\tFou_P) = \tFou_P(L^2_P(\spX\col \tau)) =
 \image(\rmP_{W(\faPq)}) = (L^2(i\faPqd) \otimes \cAtwoP)^{W(\faPq)}
$$
and (c) follows. Finally, (d) follows from (c) by taking adjoints.
\qed

\def\adritem#1{\hbox{\small #1}}
\def\distance{\hbox{\hspace{3.5cm}}}
\def\apetail{@}

\def\adderik{\vbox{
\adritem{E.P. van den Ban}
\adritem{Mathematisch Instituut}
\adritem{Universiteit Utrecht}
\adritem{PO Box 80 010}
\adritem{3508 TA Utrecht}
\adritem{Netherlands}
\adritem{E-mail: ban{\apetail}math.uu.nl}
}
}

\def\addhenrik{\vbox{
\adritem{H. Schlichtkrull}
\adritem{Matematisk Institut}
\adritem{K\o benhavns Universitet}
\adritem{Universitetsparken 5}
\adritem{2100 K\o benhavn \O}
\adritem{Denmark}
\adritem{E-mail: schlicht@math.ku.dk}
}
}
\vfill
\hbox{\vbox{\adderik}\vbox{\distance}\vbox{\addhenrik}}
\end{document}